# FUZZY INTERVAL MATRICES, NEUTROSOPHIC INTERVAL MATRICES AND THEIR APPLICATIONS

**Vasantha Kandasamy**
**Florentin Smarandache**

# Fuzzy interval matrices, Neutrosophic interval matrices and their applications

**VASANTHA KANDASAMY**
vasanthakandasamy@gmail.com
http://www.vasantha.net
http://mat.iitm.ac.in/~wbv

**FLORENTIN SMARANDACHE**
smarand@unm.edu
http://www.gallup.unm.edu/~smarandache

# CONTENTS

Chapter Three

## FUZZY MODELS AND NEUTROSOPHIC MODELS USING FUZZY INTERVAL MATRICES AND NEUTROSOPHIC INTERVAL MATRICES

## PREFACE

The new concept of fuzzy interval matrices has been introduced in this book for the first time. The authors have not only introduced the notion of fuzzy interval matrices, interval neutrosophic matrices and fuzzy neutrosophic interval matrices but have also demonstrated some of its applications when the data under study is an unsupervised one and when several experts analyze the problem.

Further, the authors have introduced in this book multi-expert models using these three new types of interval matrices. The new multi expert models dealt in this book are FCIMs, FRIMs, FCInMs, FRInMs, IBAMs, IBBAMs, nIBAMs, FAIMs, FAnIMS, etc. Illustrative examples are given so that the reader can follow these concepts easily.

This book has three chapters. The first chapter is introductory in nature and makes the book a self-contained one. Chapter two introduces the concept of fuzzy interval matrices. Also the notion of fuzzy interval matrices, neutrosophic interval matrices and fuzzy neutrosophic interval matrices, can find applications to Markov chains and Leontief economic models. Chapter three gives the application of fuzzy interval matrices and neutrosophic interval matrices to real-world problems by constructing the models already mentioned. Further these models are mainly useful when the data is an unsupervised one and when one needs a multi-expert model. The new concept of fuzzy interval matrices and neutrosophic interval matrices will find their applications in engineering, medical, industrial, social and psychological problems. We have given a long list of references to help the interested reader.

W.B.VASANTHA KANDASAMY
FLORENTIN SMARANDACHE

## Chapter One

# BASIC CONCEPTS

This chapter recalls some basic definitions and results to make this book a self contained one. This chapter has eleven sections. In section one we just recall the definition of interval matrices. Fuzzy Cognitive Maps (FCMs) are described in section two. Section three gives a brief introduction to neutrosophy. Some basic neutrosophic algebraic structures are described in section four. A brief introduction to neutrosophic graph is given in section five. Section six describes Neutrosophic Cognitive Maps (NCMs). Fuzzy Relational Maps (FRMs) and Neutrosophic Relational Maps (NRMs) definitions are briefly recalled in section seven. Section eight recollects Fuzzy Associative Memories (FAMs) and Neutrosophic Associative Memories (NAMs). Some basic concepts of Bidirectional Associative Memories (BAMs) is given in section nine. Section ten describes Fuzzy Relational Equations (FREs) and the final section describes NREs and their properties.

This book for the first time defines fuzzy interval matrices, neutrosophic interval matrices and fuzzy neutrosophic interval matrices. Most of the fuzzy models work only with single expert. CFCM, CFRM, CNCM and CNRM models can work as multi expert models but, we face certain limitations like the views of experts canceling, leading to 0 which may not be a feasible solution at all times. Fuzzy interval matrices give way to their applications in the construction of multiexpert models using the FCMs, NCMs, FRMs, NRMs, FAMs, NAMs, FREs and NREs.

## 1.1 Definition of Interval Matrices and Examples

In this section we just recall the definition of interval matrices and illustrate it with examples. As the authors are not interested in developing theories about fuzzy interval matrices and neutrosophic interval matrices but more keen on to illustrate how these matrices can find their applications in the real world model when several experts give their opinion.

**DEFINITION 1.1.1:** *Given matrices $B = (b_{ij})$ and $C = (c_{ij})$ of order n such that $b_{ij} \leq c_{ij}$, i, j = 1, 2, …, n. Then the interval matrix $A_I = [B, C]$ is defined by*

$$A_I = [B, C] = \{A = (a_{ij}) \mid b_{ij} \leq a_{ij} \leq c_{ij};\ i, j = 1, 2, \ldots, n\}.$$

*(interval vectors and matrices are vectors and matrices whose elements are interval numbers. The superscript I being used to indicate such a vector or matrix)*

***Example 1.1.1:*** Let $A_I = [B, C]$ where

$$B = \begin{bmatrix} 2 & 7 \\ 0 & -5 \end{bmatrix}$$

and

$$C = \begin{bmatrix} 12 & 20 \\ 15 & 0 \end{bmatrix}.$$

All matrices $D = (d_{ij})$ with

$2 \leq d_{11} \leq 12$,
$7 \leq d_{12} \leq 20$,
$0 \leq d_{21} \leq 15$ and
$-5 \leq d_{22} \leq 0$

are in the interval matrix $A_I = [B, C]$.

Interested reader can refer [98] for more literature.

## 1.2 Definition of Fuzzy Cognitive Maps

In this section we recall the notion of Fuzzy Cognitive Maps (FCMs), which was introduced by Bart Kosko [92] in the year 1986. We also give several of its interrelated definitions. FCMs have a major role to play mainly when the data concerned is an unsupervised one.

Further this method is most simple and an effective one as it can analyse the data by directed graphs and connection matrices.

**DEFINITION 1.2.1:** *An FCM is a directed graph with concepts like policies, events etc. as nodes and causalities as edges. It represents causal relationship between concepts.*

***Example 1.2.1****:* In Tamil Nadu (a southern state in India) in the last decade several new engineering colleges have been approved and started. The resultant increase in the production of engineering graduates in these years is disproportionate with the need of engineering graduates.

This has resulted in thousands of unemployed and underemployed graduate engineers. Using an expert's opinion we study the effect of such unemployed people on the society. An expert spells out the five major concepts relating to the unemployed graduated engineers as

$E_1$ – Frustration
$E_2$ – Unemployment
$E_3$ – Increase of educated criminals
$E_4$ – Under employment
$E_5$ – Taking up drugs etc.

The directed graph where $E_1, \ldots, E_5$ are taken as the nodes and causalities as edges as given by an expert is given in the following Figure 1.2.1:

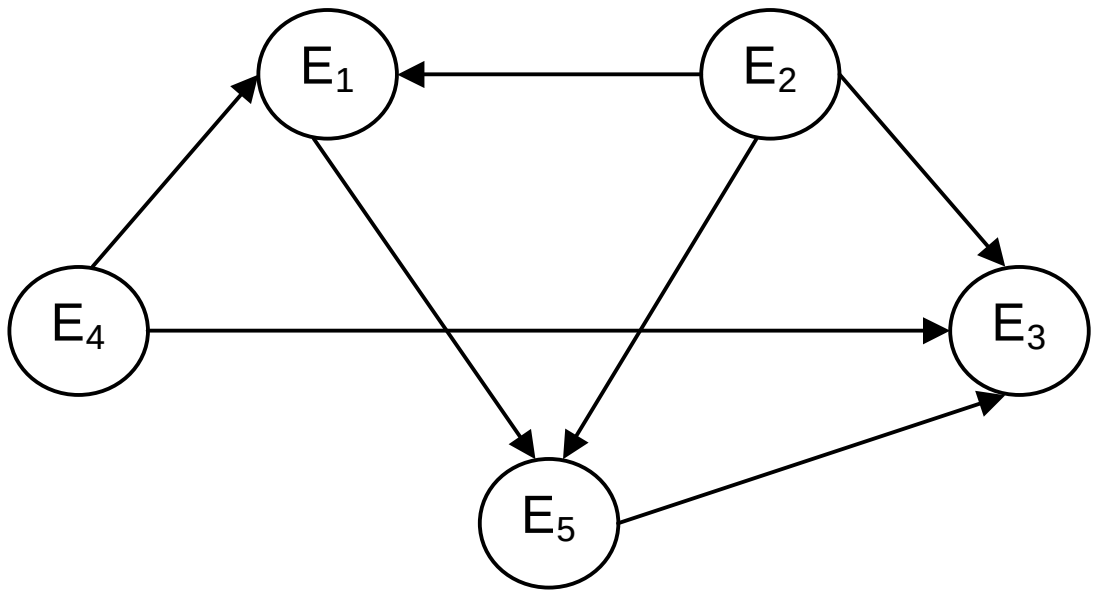

FIGURE: 1.2.1

According to this expert, increase in unemployment increases frustration. Increase in unemployment, increases the educated criminals. Frustration increases the graduates to take up to evils like drugs etc. Unemployment also leads to the increase in number of persons who take up to drugs, drinks etc. to forget their worries and unoccupied time. Under-employment forces then to do criminal acts like theft (leading to murder) for want of more money and so on. Thus one cannot actually get data for this but can use the expert's opinion for this unsupervised data to obtain some idea about the real plight of the situation. This is just an illustration to show how FCM is described by a directed graph.

{If increase (or decrease) in one concept leads to increase (or decrease) in another, then we give the value 1. If there exists no relation between two concepts the value 0 is given. If increase (or decrease) in one concept decreases (or increases) another, then we give the value –1. Thus FCMs are described in this way.}

**DEFINITION 1.2.2:** *When the nodes of the FCM are fuzzy sets then they are called as fuzzy nodes.*

**DEFINITION 1.2.3:** *FCMs with edge weights or causalities from the set {–1, 0, 1} are called simple FCMs.*

**DEFINITION 1.2.4:** *Consider the nodes / concepts $C_1, \ldots, C_n$ of the FCM. Suppose the directed graph is drawn using edge weight $e_{ij} \in \{0, 1, -1\}$. The matrix E be defined by $E = (e_{ij})$ where $e_{ij}$ is the weight of the directed edge $C_i C_j$. E is called the adjacency matrix of the FCM, also known as the connection matrix of the FCM.*

It is important to note that all matrices associated with an FCM are always square matrices with diagonal entries as zero.

**DEFINITION 1.2.5:** *Let $C_1, C_2, \ldots, C_n$ be the nodes of an FCM. $A = (a_1, a_2, \ldots, a_n)$ where $a_i \in \{0, 1\}$. A is called the instantaneous state vector and it denotes the on-off position of the node at an instant i.e.,*

$$a_i = 0 \text{ if } a_i \text{ is off and}$$
$$a_i = 1 \text{ if } a_i \text{ is on}$$

*for $i = 1, 2, \ldots, n$.*

**DEFINITION 1.2.6:** *Let $C_1, C_2, \ldots, C_n$ be the nodes of an FCM. Let $\overrightarrow{C_1C_2}, \overrightarrow{C_2C_3}, \overrightarrow{C_3C_4}, \ldots, \overrightarrow{C_iC_j}$ be the edges of the FCM ($i \neq j$). Then the edges form a directed cycle. An FCM is said to be cyclic if it possesses a directed cycle. An FCM is said to be acyclic if it does not possess any directed cycle.*

**DEFINITION 1.2.7:** *An FCM with cycles is said to have a feedback.*

**DEFINITION 1.2.8:** *When there is a feedback in an FCM, i.e., when the causal relations flow through a cycle in a revolutionary way, the FCM is called a dynamical system.*

**DEFINITION 1.2.9:** *Let $\overline{C_1C_2}, \overline{C_2C_3}, \ldots, \overline{C_{n-1}C_n}$ be a cycle. When $C_i$ is switched on and if the causality flows through the edges of a cycle and if it again causes $C_i$, we say that the dynamical system goes round and round. This is true for any node $C_i$, for $i = 1, 2, \ldots, n$. The equilibrium state for this dynamical system is called the hidden pattern.*

**DEFINITION 1.2.10:** *If the equilibrium state of a dynamical system is a unique state vector, then it is called a fixed point.*

***Example 1.2.2:*** Consider a FCM with $C_1, C_2, \ldots, C_n$ as nodes. For example let us start the dynamical system by switching on $C_1$. Let us assume that the FCM settles down with $C_1$ and $C_n$ on i.e. the state vector remains as (1, 0, 0, …, 0, 1) this state vector (1, 0, 0, …, 0, 1) is called the fixed point.

**DEFINITION 1.2.11:** *If the FCM settles down with a state vector repeating in the form*

$$A_1 \to A_2 \to \ldots \to A_i \to A_1$$

*then this equilibrium is called a limit cycle.*

**DEFINITION 1.2.12:** *Finite number of FCMs can be combined together to produce the joint effect of all the FCMs. Let $E_1, E_2, \ldots, E_p$ be the adjacency matrices of the FCMs with nodes $C_1, C_2, \ldots, C_n$ then the combined FCM is got by adding all the adjacency matrices $E_1, E_2, \ldots, E_p$.*

*We denote the combined FCM adjacency matrix by $E = E_1 + E_2 + \ldots + E_p$.*

**NOTATION:** Suppose $A = (a_1, \ldots, a_n)$ is a vector which is passed into a dynamical system E. Then $AE = (a'_1, \ldots, a'_n)$ after thresholding and updating the vector suppose we get $(b_1, \ldots, b_n)$ we denote that by

$$(a'_1, a'_2, \ldots, a'_n) \hookrightarrow (b_1, b_2, \ldots, b_n).$$

Thus the symbol '$\hookrightarrow$' means the resultant vector has been thresholded and updated.

FCMs have several advantages as well as some disadvantages. The main advantage of this method it is simple. It functions on expert's opinion. When the data happens to be an unsupervised one the FCM comes handy. This is the only known fuzzy technique that gives the hidden pattern of the situation. As we have a very well known theory, which states that the strength of the data depends on, the number of experts' opinion we can use combined FCMs with several experts' opinions.

At the same time the disadvantage of the combined FCM is when the weightages are 1 and –1 for the same $C_i$ $C_j$, we have the sum adding to zero thus at all times the connection matrices $E_1, \ldots, E_k$ may not be conformable for addition.

Combined conflicting opinions tend to cancel out and assisted by the strong law of large numbers, a consensus emerges as the sample opinion approximates the underlying population opinion. This problem will be easily overcome if the FCM entries are only 0 and 1.

We have just briefly recalled the definitions. For more about FCMs please refer Kosko [92, 96, 216].

## 1.3 An Introduction to Neutrosophy

In this section we introduce the notion of neutrosophic logic created by the author Florentine Smarandache [168-172], which is an extension / combination of the fuzzy logic in which indeterminacy is included. It has become very essential that the notion of neutrosophic logic play a vital role in several of the real world problems like law, medicine, industry, finance, IT, stocks and share etc.

Use of neutrosophic notions will be illustrated/ applied in the later sections of this chapter. Fuzzy theory only measures the grade of membership or the non-existence of a membership in the revolutionary way but fuzzy theory has failed to attribute the concept when the relations between notions or nodes or concepts in problems are indeterminate. In fact one can say the inclusion of the concept of indeterminate situation with fuzzy concepts will form the neutrosophic logic.

As in this book the concept of only fuzzy cognitive maps are dealt which mainly deals with the relation / non-relation between two nodes or concepts but it fails to deal the relation between two conceptual nodes when the relation is an indeterminate one. Neutrosophic logic is the only tool known to us, which deals with the notions of indeterminacy, and here we give a brief description of it. For more about Neutrosophic logic please refer [168-172].

**DEFINITION 1.3.1:** *In the neutrosophic logic every logical variable x is described by an ordered triple x = (T, I, F) where T is the degree of truth, F is the degree of false and I the level of indeterminacy.*

(A). To maintain consistency with the classical and fuzzy logics and with probability there is the special case where

$$T + I + F = 1.$$

(B). But to refer to intuitionistic logic, which means incomplete information on a variable proposition or event one has

$$T + I + F < 1.$$

(C). Analogically referring to Paraconsistent logic, which means contradictory sources of information about a same logical variable, proposition or event one has

$$T + I + F > 1.$$

Thus the advantage of using Neutrosophic logic is that this logic distinguishes between relative truth that is a truth is one or a few worlds only noted by 1 and absolute truth denoted by $1^{+}$. Likewise neutrosophic logic distinguishes between relative falsehood, noted by 0 and absolute falsehood noted by ${}^{-}0$.

It has several applications. One such given by [168-172] is as follows:

***Example 1.3.1:*** From a pool of refugees, waiting in a political refugee camp in Turkey to get the American visa, a% have the chance to be accepted – where a varies in the set A, r% to be rejected – where r varies in the set R, and p% to be in pending (not yet decided) – where p varies in P.

Say, for example, that the chance of someone Popescu in the pool to emigrate to USA is (between) 40-60% (considering different criteria of emigration one gets different percentages, we have to take care of all of them), the chance of being rejected is 20-25% or 30-35%, and the chance of being in pending is 10% or 20% or 30%. Then the neutrosophic probability that Popescu emigrates to the Unites States is

NP (Popescu) = ((40-60) (20-25) ∪ (30-35), {10,20,30}), closer to the life.

This is a better approach than the classical probability, where 40 P(Popescu) 60, because from the pending chance – which will be converted to acceptance or rejection – Popescu might get extra percentage in his will to emigrating and also the superior limit of the subsets sum

$$60 + 35 + 30 > 100$$

and in other cases one may have the inferior sum < 0, while in the classical fuzzy set theory the superior sum should be 100 and the inferior sum μ 0. In a similar way, we could say about the element Popescu that Popescu ((40-60), (20-25) ∪ (30-35), {10, 20, 30}) belongs to the set of accepted refugees.

***Example 1.3.2:*** The probability that candidate C will win an election is say 25-30% true (percent of people voting for him), 35% false (percent of people voting against him), and 40% or 41% indeterminate (percent of people not coming to the ballot box, or giving a blank vote – not selecting any one or giving a negative vote cutting all candidate on the list). Dialectic and dualism don't work in this case anymore.

***Example 1.3.3:*** Another example, the probability that tomorrow it will rain is say 50-54% true according to meteorologists who have investigated the past years weather, 30 or 34-35% false according to today's very sunny and droughty summer, and 10 or 20% undecided (indeterminate).

***Example 1.3.4:*** The probability that Yankees will win tomorrow versus Cowboys is 60% true (according to their confrontation's history giving Yankees' satisfaction), 30-32% false (supposing Cowboys are actually up to the mark, while Yankees are declining), and 10 or 11 or 12% indeterminate (left to the hazard: sickness of players, referee's mistakes, atmospheric conditions during the game). These parameters act on players' psychology.

As in this book we use mainly the notion of neutrosophic logic with regard to the indeterminacy of any relation in cognitive maps we are restraining ourselves from dealing with several interesting concepts about neutrosophic logic. As FCMs deals with unsupervised data and the existence or non-existence of cognitive relation, we do not in this book elaborately describe the notion of neutrosophic concepts.

However we just state, suppose in a legal issue the jury or the judge cannot always prove the evidence in a case, in several places we may not be able to derive any conclusions from the existing facts because of which we cannot make a conclusion that no relation exists or otherwise. But existing relation is an indeterminate. So in the case when the concept of indeterminacy exists the judgment ought to be very carefully analyzed be it a civil case or a criminal case. FCMs are deployed only where the existence or non-existence is dealt with but however in our Neutrosophic Cognitive Maps we will deal with the notion of indeterminacy of the evidence also. Thus legal side has lot of Neutrosophic (NCM) applications. Also we will show how NCMs can be used to study factors as varied as stock markets, medical diagnosis, etc.

## 1.4 Some Basic Neutrosophic Structures

In this section we define some new neutrosophic algebraic structures like neutrosophic fields, neutrosophic spaces and neutrosophic matrices and illustrate them with examples. For these notions are used in the definition of neutrosophic cognitive maps which is dealt in the later sections of this chapter.

Throughout this book by 'I' we denote the indeterminacy of any notion/ concept/ relation. That is when we are not in a position to associate a relation between any two concepts then we denote it as indeterminacy.

Further in this book we assume all fields to be real fields of characteristic 0 all vector spaces are taken as real spaces over reals and we denote the indeterminacy by 'I' as i will make a confusion as i denotes the imaginary value viz $i^2 = -1$ that is $\sqrt{-1} = i$.

**DEFINITION 1.4.1:** *Let K be the field of reals. We call the field generated by $K \cup I$ to be the neutrosophic field for it involves the indeterminacy factor in it. We define $I^2 = I$, $I + I = 2I$ i.e., $I + \ldots + I = nI$, and if $k \in K$ then $k.I = kI$, $0I = 0$. We denote the neutrosophic field by K(I) which is generated by $K \cup I$ that is $K(I) = \langle K \cup I \rangle$.*

***Example 1.4.1:*** Let R be the field of reals. The neutrosophic field is generated by $\langle R \cup I \rangle$ i.e. R(I) clearly $R \subset \langle R \cup I \rangle$.

***Example 1.4.2:*** Let Q be the field of rationals. The neutrosophic field is generated by Q and I i.e. $Q \cup I$ denoted by Q(I).

**DEFINITION 1.4.2:** *Let K(I) be a neutrosophic field we say K(I) is a prime neutrosophic field if K(I) has no proper subfield which is a neutrosophic field.*

***Example 1.4.3:*** Q(I) is a prime neutrosophic field where as R(I) is not a prime neutrosophic field for $Q(I) \subset R(I)$.

It is very important to note that all neutrosophic fields are of characteristic zero. Likewise we can define neutrosophic subfield.

**DEFINITION 1.4.3:** *Let K(I) be a neutrosophic field, $P \subset K(I)$ is a neutrosophic subfield of P if P itself is a neutrosophic field. K(I) will also be called as the extension neutrosophic field of the neutrosophic field P.*

Now we proceed on to define neutrosophic vector spaces, which are only defined over neutrosophic fields. We can define two types of neutrosophic vector spaces one when it is a neutrosophic vector space over ordinary field other being neutrosophic vector space over neutrosophic fields. To this end we have to define neutrosophic group under addition.

**DEFINITION 1.4.4:** *We know Z is the abelian group under addition. Z(I) denote the additive abelian group generated by the set Z and I, Z(I) is called the neutrosophic abelian group under '+'.*

Thus to define basically a neutrosophic group under addition we need a group under addition. So we proceed on to define neutrosophic abelian group under addition. *Suppose G is an additive abelian group under '+'. G(I) = $\langle G \cup I\rangle$, additive group generated by G and I, G(I) is called the neutrosophic abelian group under '+'.*

***Example 1.4.4:*** Let Q be the group under '+'; Q (I) = $\langle Q \cup I\rangle$ is the neutrosophic abelian group under addition; '+'.

***Example 1.4.5:*** R be the additive group of reals, R(I) = $\langle R \cup I\rangle$ is the neutrosophic group under addition.

***Example 1.4.6:*** $M_{n\times m}(I) = \{(a_{ij}) \mid a_{ij} \in Z(I)\}$ be the collection of all n × m matrices under '+' $M_{n\times m}(I)$ is a neutrosophic group under '+'.

Now we proceed on to define neutrosophic subgroup.

**DEFINITION 1.4.5:** *Let G(I) be the neutrosophic group under addition. P $\subset$ G(I) be a proper subset of G(I). P is said to be neutrosophic subgroup of G(I) if P itself is a neutrosophic group i.e. P = $\langle P_1 \cup I\rangle$ where $P_1$ is an additive subgroup of G.*

***Example 1.4.7:*** Let Z(I) = $\langle Z \cup I\rangle$ be a neutrosophic group under '+'. $\langle 2Z \cup I\rangle$ = 2Z(I) is the neutrosophic subgroup of Z(I).

In fact Z(I) has several neutrosophic subgroups.

Now we proceed on to define the notion of neutrosophic quotient group.

**DEFINITION 1.4.6:** *Let G (I) = $\langle G \cup I \rangle$ be a neutrosophic group under '+', suppose P (I) be a neutrosophic subgroup of G (I) then the neutrosophic quotient group*

$$\frac{G(I)}{P(I)} = \{a + P(I) \mid a \in G(I)\}.$$

***Example 1.4.8:*** Let Z(I) be a neutrosophic group under addition, Z the group of integers under addition P = 2Z(I) is a neutrosophic subgroup of Z(I) the neutrosophic subgroup of Z(I), the neutrosophic quotient group

$$\frac{Z(I)}{P} = \{a + 2Z(I) \mid a \in Z(I)\} = \{(2n+1) + (2n+1) I \mid n \in Z\}.$$

Clearly $\frac{Z(I)}{P}$ is a group. For P = 2Z (I) serves as the additive identity. Take a, b $\in \frac{Z(I)}{P}$. If a, b $\in$ Z(I) \ P then two possibilities occur.

a + b is odd times I or a + b is odd or a + b is even times I or even if a + b is even or even times I then a + b $\in$ P. if a + b is odd or odd times I a + b $\in \frac{Z(I)}{P = 2Z(I)}$.

It is easily verified that P acts as the identity and every element in

$$a + 2Z\,(I) \in \frac{Z(I)}{2Z(I)}$$

has inverse. Hence the claim.

Now we proceed on to define the notion of neutrosophic vector spaces over fields and then we define neutrosophic vector spaces over neutrosophic fields.

**DEFINITION 1.4.7:** *Let G(I) by an additive abelian neutrosophic group. K any field. If G(I) is a vector space over K then we call G(I) a neutrosophic vector space over K.*

Now we give the notion of strong neutrosophic vector space.

**DEFINITION 1.4.8:** *Let G(I) be a neutrosophic abelian group. K(I) be a neutrosophic field. If G(I) is a vector space over K(I) then we call G(I) the strong neutrosophic vector space.*

**THEOREM 1.4.1:** *All strong neutrosophic vector space over K(I) are a neutrosophic vector space over K; as K $\subset$ K(I).*

*Proof:* Follows directly by the very definitions.

Thus when we speak of neutrosophic spaces we mean either a neutrosophic vector space over K or a strong neutrosophic vector space over the neutrosophic field K(I). By basis we mean a linearly independent set which spans the neutrosophic space.

Now we illustrate with an example.

***Example 1.4.9:*** Let R(I) × R(I) = V be an additive abelian neutrosophic group over the neutrosophic field R(I). Clearly V is a strong neutrosophic vector space over R(I). The basis of V are {(0,1), (1,0)}.

***Example 1.4.10:*** Let V = R(I) × R(I) be a neutrosophic abelian group under addition. V is a neutrosophic vector space over R. The neutrosophic basis of V are {(1,0), (0,1), (I,0), (0,I)}, which is a basis of the vector space V over R.

A study of these basis and its relations happens to be an interesting form of research.

**DEFINITION 1.4.9:** *Let G(I) be a neutrosophic vector space over the field K. The number of elements in the neutrosophic basis is called the neutrosophic dimension of G(I).*

**DEFINITION 1.4.10:** *Let G(I) be a strong neutrosophic vector space over the neutrosophic field K(I). The number of elements in the strong neutrosophic basis is called the strong neutrosophic dimension of G(I).*

We denote the neutrosophic dimension of G(I) over K by $N_k$ (dim) of G (I) and that the strong neutrosophic dimension of G (I) by $SN_{K(I)}$ (dim) of G(I).

Now we define the notion of neutrosophic matrices.

**DEFINITION 1.4.11:** *Let $M_{nxm} = \{(a_{ij}) \;/\; a_{ij} \in K(I)\}$, where K (I), is a neutrosophic field.We call $M_{nxm}$ to be the neutrosophic matrix.*

***Example 1.4.11:*** Let Q(I) = ⟨Q ∪ I⟩ be the neutrosophic field.

$$M_{4\times 3} = \begin{pmatrix} 0 & 1 & I \\ -2 & 4I & 0 \\ 1 & -I & 2 \\ 3I & 1 & 0 \end{pmatrix}$$

is the neutrosophic matrix, with entries from rationals and the indeterminacy I. We define product of two neutrosophic matrices whenever the production is defined as follows:

Let

$$A = \begin{pmatrix} -1 & 2 & -I \\ 3 & I & 0 \end{pmatrix}_{2\times 3}$$

and

$$B = \begin{pmatrix} -I & 1 & 2 & 4 \\ 1 & I & 0 & 2 \\ 5 & -2 & 3I & -I \end{pmatrix}_{3\times 4}$$

$$AB = \begin{bmatrix} -6I+2 & -1+4I & -2-3I & I \\ -2I & 3+I & 6 & 12+2I \end{bmatrix}_{2\times 4}$$

(we use the fact $I^2 = I$).

To define Neutrosophic Cognitive Maps we direly need the notion of neutrosophic matrices. We use square neutrosophic matrices for Neutrosophic Cognitive Maps and use rectangular neutrosophic matrices for Neutrosophic Relational Maps (NRMs).

## 1.5 Some Basic Notions about Neutrosophic Graphs

In this section we for the first time introduce the notion of neutrosophic graphs, illustrate them and give some basic properties. We need the notion of neutrosophic graphs basically to obtain neutrosophic cognitive maps which will be nothing but directed neutrosophic graphs. Similarly neutrosophic relational maps will also be directed neutrosophic graphs.

It is no coincidence that graph theory has been independently discovered many times since it may quite properly be regarded as an area of applied mathematics. The subject finds its place in the work of Euler. Subsequent rediscoveries of graph theory were by Kirchhoff and Cayley. Euler (1707-1782) became the father of graph theory as well as topology when in 1936 he settled a famous unsolved problem in his day called the Konigsberg Bridge Problem.

Psychologist Lewin proposed in 1936 that the life space of an individual be represented by a planar map. His view point led the psychologists at the Research center for Group Dynamics to another psychological interpretation of a graph in which people are represented by points and interpersonal relations by lines. Such relations include love, hate, communication and power. In fact it was precisely this approach which led the author to a personal discovery of graph theory, aided and abetted by psychologists L. Festinger and D. Cartwright.

Here it is pertinent to mention that the directed graphs of an FCMs or FRMs are nothing but the psychological inter-relations or feelings of different nodes; but it is unfortunate that in all these studies the concept of indeterminacy was never given any place, so in this chapter for the first time we will be having graphs in which the notion of indeterminacy i.e. when two vertex should be connected or not is never dealt with. If graphs are to display human feelings then this point is very vital for in several situations certain relations between concepts may certainly remain an indeterminate. So this section will purely cater to the properties of such graphs the edges of certain vertices may not find its connection i.e., they are indeterminates, which we will be defining as neutrosophic graphs.

The world of theoretical physics discovered graph theory for its own purposes. In the study of statistical mechanics by Uhlenbeck the points stands for molecules and two adjacent points indicate nearest neighbor interaction of some physical kind, for example magnetic attraction or repulsion. But it is forgotten in all these situations we may have molecule structures which need not attract or repel but remain without action or not able to predict the action for such analysis we can certainly adopt the concept of neutrosophic graphs.

In a similar interpretation by Lee and Yang the points stand for small cubes in Euclidean space where each cube may or may not be occupied by a molecule. Then two points are adjacent whenever both spaces are occupied. Feynmann proposed the diagram in which the points represent physical particles and the lines represent paths of the particles after collisions. Just at each stage of applying graph theory we may now feel the neutrosophic graph theory may be more suitable for application.

Now we proceed on to define the neutrosophic graph.

**DEFINITION 1.5.1:** *A neutrosophic graph is a graph in which at least one edge is an indeterminacy denoted by dotted lines.*

**NOTATION***:* The indeterminacy of an edge between two vertices will always be denoted by dotted lines.

***Example 1.5.1:*** The following are neutrosophic graphs:

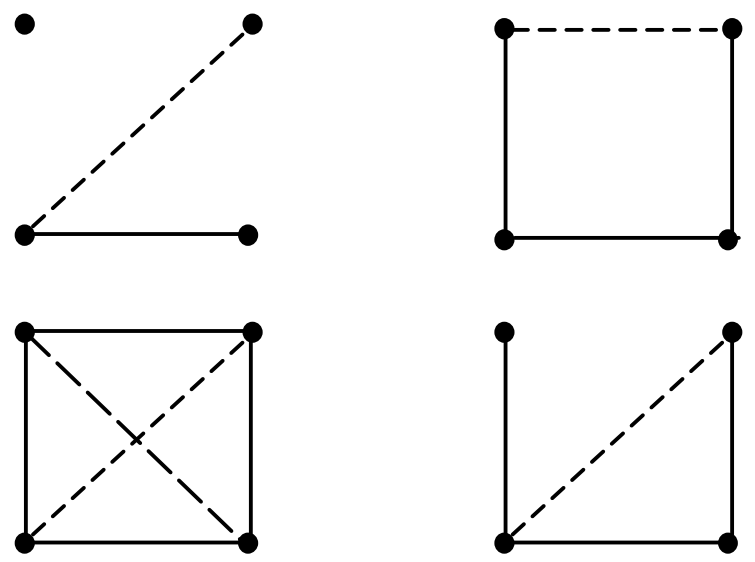

FIGURE: 1.5.1

All graphs in general are not neutrosophic graphs.

***Example 1.5.2:*** The following graphs are not neutrosophic graphs given in Figure 1.5.2.

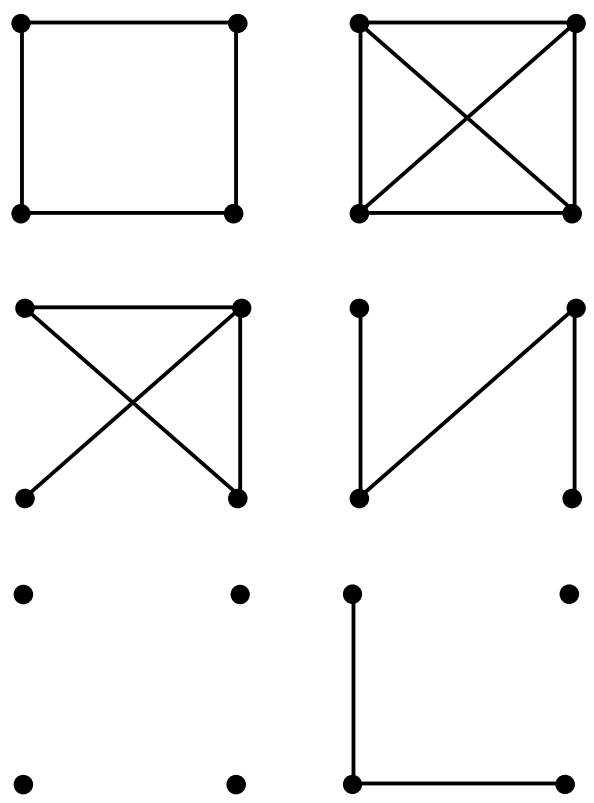

FIGURE: 1.5.2

**DEFINITION 1.5.2:** *A neutrosophic directed graph is a directed graph which has at least one edge to be an indeterminacy.*

**DEFINITION 1.5.3:** *A neutrosophic oriented graph is a neutrosophic directed graph having no symmetric pair of directed indeterminacy lines.*

**DEFINITION 1.5.4:** *A neutrosophic subgraph H of a neutrosophic graph G is a subgraph H which is itself a neutrosophic graph.*

**THEOREM 1.5.1:** *Let G be a neutrosophic graph. All subgraphs of G are not neutrosophic subgraphs of G.*

*Proof:* By an example. Consider the neutrosophic graph given

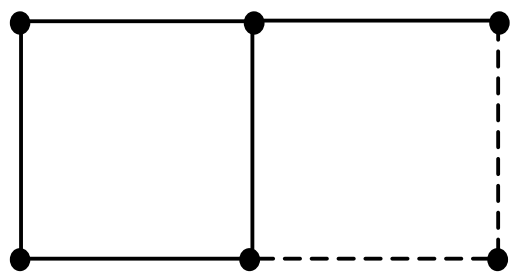

FIGURE: 1.5.3

in Figure 1.5.3.
This has a subgraph given by Figure 1.5.4.

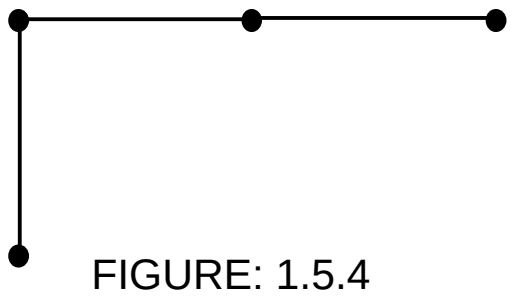

FIGURE: 1.5.4

which is not a neutrosophic subgraph of G.

**THEOREM 1.5.2:** *Let G be a neutrosophic graph. In general the removal of a point from G need not be a neutrosophic subgraph.*

*Proof:* Consider the graph G given in Figure 1.5.5.

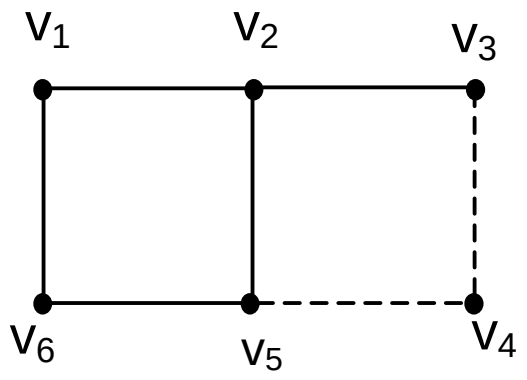

FIGURE: 1.5.5

G \ $v_4$ is only a subgraph of G but is not a neutrosophic subgraph of G.

Thus it is interesting to note that this is a main feature by which a graph differs from a neutrosophic graph.

**DEFINITION 1.5.5:** *Two graphs G and H are neutrosophically isomorphic if*

*i. They are isomorphic.*
*ii. If there exists a one to one correspondence between their point sets which preserve indeterminacy adjacency.*

**DEFINITION 1.5.6:** *A neutrosophic walk of a neutrosophic graph G is a walk of the graph G in which at least one of the lines is an indeterminacy line. The neutrosophic walk is neutrosophic closed if* $v_0 = v_n$ *and is neutrosophic open otherwise.*

*It is a neutrosophic trial if all the lines are distinct and at least one of the lines is a indeterminacy line and a path, if all points are distinct (i.e. this necessarily means all lines are distinct and at least one line is a line of indeterminacy). If the neutrosophic walk is neutrosophic closed then it is a neutrosophic cycle provided its n points are distinct and* $n \geq 3$.

*A neutrosophic graph is neutrosophic connected if it is connected and at least a pair of points are joined by a path. A neutrosophic maximal connected neutrosophic subgraph of G is called a neutrosophic connected component or simple neutrosophic component of G.*

*Thus a neutrosophic graph has at least two neutrosophic components then it is neutrosophic disconnected. Even if one is a component and another is a neutrosophic component still we do not say the graph is neutrosophic disconnected.*

***Example 1.5.3:*** Neutrosophic disconnected graphs are given in Figure 1.5.6.

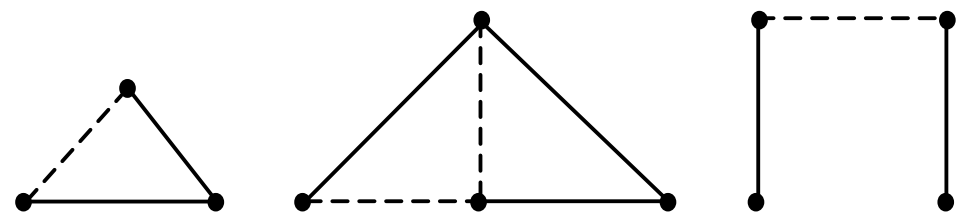

FIGURE: 1.5.6

***Example 1.5.4:*** Graph which is not neutrosophic disconnected is given by Figure 1.5.7.

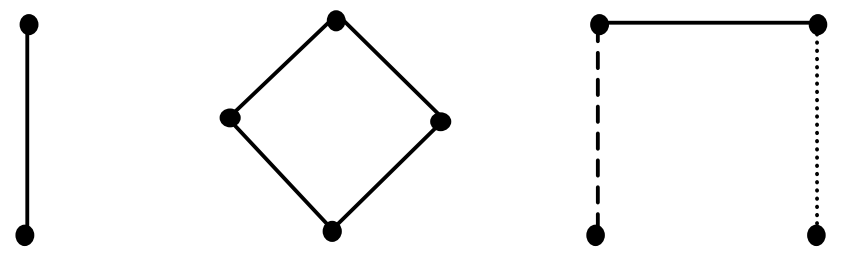

FIGURE: 1.5.7

Several results in this direction can be defined and analyzed.

**DEFINITION 1.5.7:** *A neutrosophic bigraph, G is a bigraph, G whose point set V can be partitioned into two subsets $V_1$ and $V_2$ such that at least a line of G which joins $V_1$ with $V_2$ is a line of indeterminacy.*

This neutrosophic bigraphs will certainly play a role in the study of FRMs and in fact we give a method of conversion of data from FRMs to FCMs. As both the models FRMs and FCMs work on the adjacency or the connection matrix we just define the neutrosophic adjacency matrix related to a neutrosophic graph G given by Figure 1.5.8.

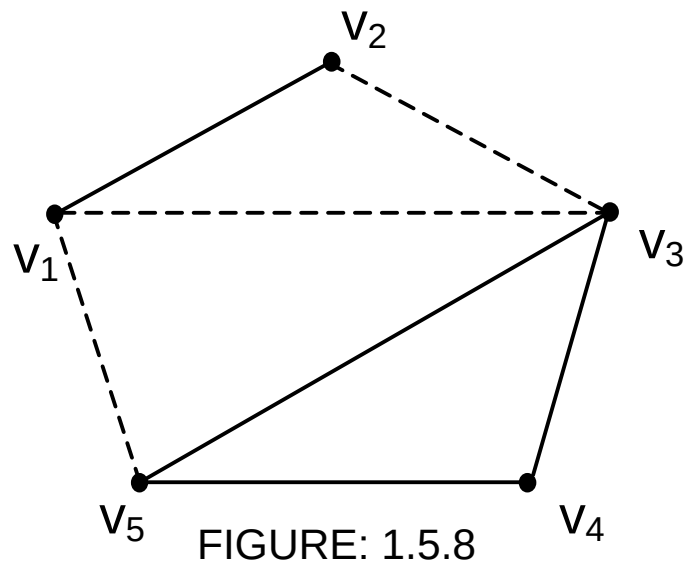

FIGURE: 1.5.8

The neutrosophic adjacency matrix is N(A)

$$N(A) = \begin{bmatrix} 0 & 1 & I & 0 & I \\ 1 & 0 & I & 0 & 0 \\ I & I & 0 & 1 & 1 \\ 0 & 0 & 1 & 0 & 1 \\ I & 0 & 1 & 1 & 0 \end{bmatrix}.$$

Its entries will not only be 0 and 1 but also the indeterminacy *I*.

**DEFINITION 1.5.8:** *Let G be a neutrosophic graph. The adjacency matrix of G with entries from the set (I, 0, 1) is called the neutrosophic adjacency matrix of the graph.*

Thus one finds a very interesting application of neutrosophy graphs in Neutrosophic Cognitive Maps and Neutrosophic Relational Maps. Application of these concepts will be dealt in the last section of this chapter.

Now as our main aim is the study of Neutrosophic Cognitive Maps we do not divulge into a very deep study of neutrosophic graphs or its properties but have given only the basic and the appropriate notions which are essential for studying of this book.

## 1.6 On Neutrosophic Cognitive Maps with Examples

The notion of Fuzzy Cognitive Maps (FCMs) which are fuzzy signed directed graphs with feedback are discussed and described in section 1.5 of this chapter. The directed edge $e_{ij}$ from causal concept $C_i$ to concept $C_j$ measures how much $C_i$ causes $C_j$. The time varying concept function $C_i(t)$ measures the non negative occurrence of some fuzzy event, perhaps the strength of a political sentiment, historical trend or opinion about some topics like child labor or school dropouts etc. FCMs model the world as a collection of classes and causal relations between them.

The edge $e_{ij}$ takes values in the fuzzy causal interval [–1, 1] ($e_{ij}$ = 0 indicates no causality, $e_{ij} > 0$ indicates causal increase; that $C_j$ increases as $C_i$ increases and $C_j$ decreases as $C_i$ decreases, $e_{ij} < 0$ indicates causal decrease or negative causality $C_j$ decreases as $C_i$ increases or $C_j$, increases as $C_i$ decreases. Simple FCMs have edge value in {–1, 0, 1}. Thus if causality occurs it occurs to maximal positive or negative degree.

It is important to note that $e_{ij}$ measures only absence or presence of influence of the node $C_i$ on $C_j$ but till now any researcher has not contemplated the indeterminacy of any relation between two nodes $C_i$ and $C_j$. When we deal with unsupervised data, there are situations when no relation can be determined between some two nodes. So in this section we try to introduce the indeterminacy in FCMs, and we choose to call this generalized structure as Neutrosophic Cognitive Maps (NCMs). In our view this will certainly give a more appropriate result and also caution the user about the risk of indeterminacy.

Now we proceed on to define the concepts about NCMs.

**DEFINITION 1.6.1:** *A Neutrosophic Cognitive Map (NCM) is a neutrosophic directed graph with concepts like policies, events etc. as nodes and causalities or indeterminates as edges. It represents the causal relationship between concepts.*

Let $C_1, C_2, \ldots, C_n$ denote n nodes, further we assume each node is a neutrosophic vector from neutrosophic vector space V. So a node $C_i$ will be represented by $(x_1, \ldots, x_n)$ where $x_k$'s are zero or one or I and $x_k = 1$ means that the node $C_k$ is in the on state and $x_k = 0$ means the node is in the off state and $x_k = I$ means the nodes state is an indeterminate at that time or in that situation.

Let $C_i$ and $C_j$ denote the two nodes of the NCM. The directed edge from $C_i$ to $C_j$ denotes the causality of $C_i$ on $C_j$ called connections. Every edge in the NCM is weighted with a number in the set {-1, 0, 1, I}. Let $e_{ij}$ be the weight of the directed edge $C_iC_j$, $e_{ij} \in \{-1, 0, 1, I\}$. $e_{ij} = 0$ if $C_i$ does not have any effect on $C_j$, $e_{ij} = 1$ if increase (or decrease) in $C_i$ causes increase (or decreases) in $C_j$, $e_{ij} = -1$ if increase (or decrease) in $C_i$ causes decrease (or increase) in $C_j$ . $e_{ij} = I$ if the relation or effect of $C_i$ on $C_j$ is an indeterminate.

**DEFINITION 1.6.2:** *NCMs with edge weight from {-1, 0, 1, I} are called simple NCMs.*

**DEFINITION 1.6.3:** *Let $C_1$, $C_2$, …, $C_n$ be nodes of a NCM. Let the neutrosophic matrix N(E) be defined as N(E) = ($e_{ij}$) where $e_{ij}$ is the weight of the directed edge $C_i$ $C_j$, where $e_{ij} \in \{0, 1, -1, I\}$. N(E) is called the neutrosophic adjacency matrix of the NCM.*

**DEFINITION 1.6.4:** *Let $C_1$, $C_2$, …, $C_n$ be the nodes of the NCM. Let A = ($a_1$, $a_2$,…, $a_n$) where $a_i \in \{0, 1, I\}$. A is called the instantaneous state neutrosophic vector and it denotes the on – off – indeterminate state position of the node at an instant*

*$a_i$ = 0 if $a_i$ is off (no effect)*
*$a_i$ = 1 if $a_i$ is on (has effect)*
*$a_i$ = I if $a_i$ is indeterminate(effect cannot be determined)*

*for i = 1, 2,…, n.*

**DEFINITION 1.6.5:** *Let $C_1$, $C_2$, …, $C_n$ be the nodes of the FCM. Let $\overrightarrow{C_1C_2}$, $\overrightarrow{C_2C_3}$, $\overrightarrow{C_3C_4}$, … , $\overrightarrow{C_iC_j}$ be the edges of the NCM. Then the edges form a directed cycle. An NCM is said to be cyclic if it possesses a directed cyclic. An NCM is said to be acyclic if it does not possess any directed cycle.*

**DEFINITION 1.6.6:** *An NCM with cycles is said to have a feedback. When there is a feedback in the NCM i.e. when the causal relations flow through a cycle in a revolutionary manner the NCM is called a dynamical system.*

**DEFINITION 1.6.7:** *Let $\overrightarrow{C_1C_2}, \overrightarrow{C_2C_3}, \cdots, \overrightarrow{C_{n-1}C_n}$ be cycle, when $C_i$ is switched on and if the causality flow through the edges of a cycle and if it again causes $C_i$, we say that the dynamical system goes round and round. This is true for any node $C_i$, for i = 1, 2,…, n. the equilibrium state for this dynamical system is called the hidden pattern.*

**DEFINITION 1.6.8:** *If the equilibrium state of a dynamical system is a unique state vector, then it is called a fixed point. Consider the NCM with $C_1$, $C_2$,…, $C_n$ as nodes. For example let us start the dynamical system by switching on $C_1$. Let us assume* *that the NCM settles down with $C_1$ and $C_n$ on, i.e. the state vector remain as (1, 0,…, 1) this neutrosophic state vector (1,0,…, 0, 1) is called the fixed point.*

**DEFINITION 1.6.9:** *If the NCM settles with a neutrosophic state vector repeating in the form*

$$A_1 \to A_2 \to \ldots \to A_i \to A_1,$$

*then this equilibrium is called a limit cycle of the NCM.*

Here we give the methods of determining the hidden pattern of NCM.
Let $C_1, C_2,\ldots, C_n$ be the nodes of an NCM, with feedback. Let E be the associated adjacency matrix. Let us find the hidden pattern when $C_1$ is switched on when an input is given as the vector $A_1 = (1, 0, 0,\ldots, 0)$, the data should pass through the neutrosophic matrix N(E), this is done by multiplying $A_1$ by the matrix N(E). Let $A_1N(E) = (a_1, a_2,\ldots, a_n)$ with the threshold operation that is by replacing $a_i$ by 1 if $a_i > k$ and $a_i$ by 0 if $a_i < k$ (k – a suitable positive integer) and $a_i$ by I if $a_i$ is not a integer. We update the resulting concept, the concept $C_1$ is included in the updated vector by making the first coordinate as 1 in the resulting vector. Suppose $A_1N(E) \to A_2$ then consider $A_2N(E)$ and repeat the same procedure. This procedure is repeated till we get a limit cycle or a fixed point.

**DEFINITION 1.6.10:** *Finite number of NCMs can be combined together to produce the joint effect of all NCMs. If $N(E_1)$, $N(E_2),\ldots, N(E_p)$ be the neutrosophic adjacency matrices of a NCM with nodes $C_1, C_2,\ldots, C_n$ then the combined NCM is got by adding all the neutrosophic adjacency matrices $N(E_1),\ldots, N(E_p)$. We denote the combined NCMs adjacency neutrosophic matrix by $N(E) = N(E_1) + N(E_2)+\ldots+ N(E_p)$.*

**NOTATION:** Let $(a_1, a_2, \ldots , a_n)$ and $(a'_1, a'_2, \ldots , a'_n)$ be two neutrosophic vectors. We say $(a_1, a_2, \ldots , a_n)$ is equivalent to $(a'_1, a'_2, \ldots , a'_n)$ denoted by $((a_1, a_2, \ldots , a_n) \sim (a'_1, a'_2, \ldots, a'_n)$ if $(a'_1, a'_2, \ldots , a'_n)$ is got after thresholding and updating the vector $(a_1, a_2, \ldots , a_n)$ after passing through the neutrosophic adjacency matrix N(E).

The following are very important:

***Note 1:*** The nodes $C_1$, $C_2$, …, $C_n$ are nodes are not indeterminate nodes because they indicate the concepts which are well known. But the edges connecting $C_i$ and $C_j$ may be indeterminate i.e. an expert may not be in a position to say that $C_i$ has some causality on $C_j$ either will he be in a position to state that $C_i$ has no relation with $C_j$ in such cases the relation between $C_i$ and $C_j$ which is indeterminate is denoted by I.

***Note 2:*** The nodes when sent will have only ones and zeros i.e. on and off states, but after the state vector passes through the neutrosophic adjacency matrix the resultant vector will have entries from {0, 1, I} i.e. they can be neutrosophic vectors.
The presence of I in any of the coordinate implies the expert cannot say the presence of that node i.e. on state of it after passing through N(E) nor can we say the absence of the node i.e. off state of it the effect on the node after passing through the dynamical system is indeterminate so only it is represented by I. Thus only in case of NCMs we can say the effect of any node on other nodes can also be indeterminates. Such possibilities and analysis is totally absent in the case of FCMs.

***Note 3:*** In the neutrosophic matrix N(E), the presence of I in the $a_{ij}$ the place shows, that the causality between the two nodes i.e. the effect of $C_i$ on $C_j$ is indeterminate. Such chances of being indeterminate is very possible in case of unsupervised data and that too in the study of FCMs which are derived from the directed graphs.
Thus only NCMs helps in such analysis.

Now we shall represent a few examples to show how in this set up NCMs is preferred to FCMs. At the outset before we proceed to give examples we make it clear that all unsupervised data need not have NCMs to be applied to it. Only data which have the relation between two nodes to be an indeterminate need to be modeled with NCMs if the data has no indeterminacy factor between any pair of nodes one need not go for NCMs; FCMs will do the best job.

## 1.7 Definition and Illustration of Fuzzy Relational Maps (FRMs)

In this section, we introduce the notion of Fuzzy relational maps (FRMs); they are constructed analogous to FCMs described and discussed in the earlier sections. In FCMs we promote the correlations between causal associations among concurrently active units. But in FRMs we divide the very causal associations into two disjoint units, for example, the relation between a teacher and a student or relation between an employee or employer or a relation between doctor and patient and so on. Thus for us to define a FRM we need a domain space and a range space which are disjoint in the sense of concepts. We further assume no intermediate relation exists within the domain elements or node and the range spaces elements. The number of elements in the range space need not in general be equal to the number of elements in the domain space.

Thus throughout this section we assume the elements of the domain space are taken from the real vector space of dimension n and that of the range space are real vectors from the vector space of dimension m (m in general need not be equal to n). We denote by R the set of nodes $R_1,\ldots, R_m$ of the range space, where $R = \{(x_1,\ldots, x_m) \mid x_j = 0 \text{ or } 1\}$ for j = 1, 2,…, m. If $x_i = 1$ it means that the node $R_i$ is in the on state and if $x_i = 0$ it means that the node $R_i$ is in the off state. Similarly D denotes the nodes $D_1, D_2,\ldots, D_n$ of the domain space where $D = \{(x_1,\ldots, x_n) \mid x_j = 0 \text{ or } 1\}$ for i = 1, 2,…, n. If $x_i = 1$ it means that the node $D_i$ is in the on state and if $x_i = 0$ it means that the node $D_i$ is in the off state.

Now we proceed on to define a FRM.

**DEFINITION 1.7.1:** *A FRM is a directed graph or a map from D to R with concepts like policies or events etc, as nodes and causalities as edges. It represents causal relations between spaces D and R .*

*Let $D_i$ and $R_j$ denote that the two nodes of an FRM. The directed edge from $D_i$ to $R_j$ denotes the causality of $D_i$ on $R_j$ called relations. Every edge in the FRM is weighted with a number in the set {0, ±1}. Let $e_{ij}$ be the weight of the edge $D_iR_j$, $e_{ij} \in \{0, \pm 1\}$. The weight of the edge $D_i R_j$ is positive if increase in $D_i$ implies increase in $R_j$ or decrease in $D_i$ implies decrease in $R_j$ ie causality of $D_i$ on $R_j$ is 1. If $e_{ij} = 0$, then $D_i$ does not have any effect on $R_j$ . We do not discuss the cases when increase in $D_i$ implies decrease in $R_j$ or decrease in $D_i$ implies increase in $R_j$ .*

**DEFINITION 1.7.2:** *When the nodes of the FRM are fuzzy sets then they are called fuzzy nodes. FRMs with edge weights {0, ±1} are called simple FRMs.*

**DEFINITION 1.7.3:** *Let $D_1, \ldots, D_n$ be the nodes of the domain space D of an FRM and $R_1, \ldots, R_m$ be the nodes of the range space R of an FRM. Let the matrix E be defined as $E = (e_{ij})$ where $e_{ij}$ is the weight of the directed edge $D_iR_j$ (or $R_jD_i$), E is called the relational matrix of the FRM.*

***<u>Note</u>:*** It is pertinent to mention here that unlike the FCMs the FRMs can be a rectangular matrix with rows corresponding to the domain space and columns corresponding to the range space. This is one of the marked difference between FRMs and FCMs.

**DEFINITION 1.7.4:** *Let $D_1, ..., D_n$ and $R_1,\ldots, R_m$ denote the nodes of the FRM. Let $A = (a_1,\ldots,a_n)$, $a_i \in \{0, 1\}$. A is called the instantaneous state vector of the domain space and it denotes the on-off position of the nodes at any instant. Similarly let $B = (b_1,\ldots, b_m)$ $b_i \in \{0, 1\}$. B is called instantaneous state vector of the range space and it denotes the on-off position of the nodes at any instant $a_i = 0$ if $a_i$ is off and $a_i = 1$ if $a_i$ is on for i= 1, 2,…, n Similarly, $b_i = 0$ if $b_i$ is off and $b_i = 1$ if $b_i$ is on, for i= 1, 2,…, m.*

**DEFINITION 1.7.5:** *Let $D_1, \ldots, D_n$ and $R_1,\ldots, R_m$ be the nodes of an FRM. Let $D_iR_j$ (or $R_j D_i$) be the edges of an FRM, j = 1, 2,…,* *m and i= 1, 2,…, n. Let the edges form a directed cycle. An FRM is said to be a cycle if it posses a directed cycle. An FRM is said to be acyclic if it does not posses any directed cycle.*

**DEFINITION 1.7.6:** *An FRM with cycles is said to be an FRM with feedback.*

**DEFINITION 1.7.7:** *When there is a feedback in the FRM, i.e. when the causal relations flow through a cycle in a revolutionary manner, the FRM is called a dynamical system.*

**DEFINITION 1.7.8:** *Let $D_i$ $R_j$ (or $R_j$ $D_i$), $1 \leq j \leq m$, $1 \leq i \leq n$. When $R_i$ (or $D_j$) is switched on and if causality flows through edges of the cycle and if it again causes $R_i$ (or$D_j$), we say that the dynamical system goes round and round. This is true for any node $R_j$ (or $D_i$) for $1 \leq i \leq n$, (or $1 \leq j \leq m$). The equilibrium state of this dynamical system is called the hidden pattern.*

**DEFINITION 1.7.9:** *If the equilibrium state of a dynamical system is a unique state vector, then it is called a fixed point. Consider an FRM with $R_1, R_2,\ldots, R_m$ and $D_1, D_2,\ldots, D_n$ as nodes. For example, let us start the dynamical system by switching on $R_1$ (or $D_1$). Let us assume that the FRM settles down with $R_1$ and $R_m$ (or $D_1$ and $D_n$) on, i.e. the state vector remains as (1, 0, …, 0, 1) in R (or 1, 0, 0, … , 0, 1) in D), This state vector is called the fixed point.*

**DEFINITION 1.7.10:** *If the FRM settles down with a state vector repeating in the form*

$$A_1 \to A_2 \to A_3 \to \ldots \to A_i \to A_1 \text{ (or } B_1 \to B_2 \to \ldots \to B_i \to B_1)$$

*then this equilibrium is called a limit cycle.*

Here we give the methods of determining the hidden pattern of FRM.
Let $R_1, R_2,\ldots, R_m$ and $D_1, D_2,\ldots, D_n$ be the nodes of a FRM with feedback. Let E be the relational matrix. Let us find a hidden pattern when $D_1$ is switched on i.e. when an input is given as vector $A_1 = (1, 0, \ldots, 0)$ in $D_1$, the data should pass through the relational matrix E. This is done by multiplying $A_1$ with the relational matrix E. Let $A_1E = (r_1, r_2, \ldots, r_m)$, after thresholding and updating the resultant vector we get $A_1 E \in R$. Now let $B = A_1E$ we pass on B into $E^T$ and obtain $BE^T$. We update and threshold the vector $BE^T$ so that $BE^T \in D$. This procedure is repeated till we get a limit cycle or a fixed point.

**DEFINITION 1.7.11:** *Finite number of FRMs can be combined together to produce the joint effect of all the FRMs. Let $E_1, \ldots, E_p$ be the relational matrices of the FRMs with nodes $R_1, R_2, \ldots, R_m$ and $D_1, D_2, \ldots, D_n$, then the combined FRM is represented by the relational matrix $E = E_1 + \ldots + E_p$.*

Now we proceed on to describe a NRM.

Neutrosophic Cognitive Maps (NCMs) promote the causal relationships between concurrently active units or decides the absence of any relation between two units or the indeterminance of any relation between any two units. But in Neutrosophic Relational Maps (NRMs) we divide the very causal nodes into two disjoint units. Thus for the modeling of a NRM we need a domain space and a range space which are disjoint in the sense of concepts. We further assume no intermediate relations exist within the domain and the range spaces. The number of elements or nodes in the range space need not be equal to the number of elements or nodes in the domain space.

Throughout this section we assume the elements of a domain space are taken from the neutrosophic vector space of dimension n and that of the range space are neutrosophic vector space of dimension m. (m in general need not be equal to n). We denote by R the set of nodes $R_1, \ldots, R_m$ of the range space, where $R = \{(x_1, \ldots, x_m) \mid x_j = 0 \text{ or } 1 \text{ for } j = 1, 2, \ldots, m\}$.

If $x_i = 1$ it means that node $R_i$ is in the on state and if $x_i = 0$ it means that the node $R_i$ is in the off state and if $x_i = I$ in the resultant vector it means the effect of the node $x_i$ is indeterminate or whether it will be off or on cannot be predicted by the neutrosophic dynamical system.

It is very important to note that when we send the state vectors they are always taken as the real state vectors for we know the node or the concept is in the on state or in the off state but when the state vector passes through the Neutrosophic dynamical system some other node may become indeterminate i.e. due to the presence of a node we may not be able to predict the presence or the absence of the other node i.e., it is indeterminate, denoted by the symbol I, thus the resultant vector can be a neutrosophic vector.

**DEFINITION 1.7.12:** *A Neutrosophic Relational Map (NRM) is a Neutrosophic directed graph or a map from D to R with concepts like policies or events etc. as nodes and causalities as edges. (Here by causalities we mean or include the indeterminate causalities also). It represents Neutrosophic Relations and Causal Relations between spaces D and R .*

*Let $D_i$ and $R_j$ denote the nodes of an NRM. The directed edge from $D_i$ to $R_j$ denotes the causality of $D_i$ on $R_j$ called relations. Every edge in the NRM is weighted with a number in the set {0, +1, –1, I}. Let $e_{ij}$ be the weight of the edge $D_i R_j$, $e_{ij} \in \{0, 1, –1, I\}$. The weight of the edge $D_i R_j$ is positive if increase in $D_i$ implies increase in $R_j$ or decrease in $D_i$ implies decrease in $R_j$ i.e. causality of $D_i$ on $R_j$ is 1. If $e_{ij} = –1$ then increase (or decrease) in $D_i$ implies decrease (or increase) in $R_j$. If $e_{ij} = 0$ then $D_i$ does not have any effect on $R_j$. If $e_{ij} = I$ it implies we are not in a position to determine the effect of $D_i$ on $R_j$ i.e. the effect of $D_i$ on $R_j$ is an indeterminate so we denote it by I.*

**DEFINITION 1.7.13:** *When the nodes of the NRM take edge values from {0, 1, –1, I} we say the NRMs are simple NRMs.*

**DEFINITION 1.7.14:** *Let $D_1, …, D_n$ be the nodes of the domain space D of an NRM and let $R_1, R_2,…, R_m$ be the nodes of the range space R of the same NRM. Let the matrix N(E) be defined as $N(E) = (e_{ij})$ where $e_{ij}$ is the weight of the directed edge $D_i R_j$ (or $R_j D_i$) and $e_{ij} \in \{0, 1, –1, I\}$. N(E) is called the Neutrosophic Relational Matrix of the NRM.*

The following remark is important and interesting to find its mention in this book.

**Remark**: Unlike NCMs, NRMs can also be rectangular matrices with rows corresponding to the domain space and columns corresponding to the range space. This is one of the marked difference between NRMs and NCMs. Further the number of entries for a particular model which can be treated as disjoint sets when dealt as a NRM has very much less entries than when the same model is treated as a NCM.

Thus in many cases when the unsupervised data under study or consideration can be spilt as disjoint sets of nodes or concepts; certainly NRMs are a better tool than the NCMs.

**DEFINITION 1.7.15:** *Let $D_1$, …, $D_n$ and $R_1$,…, $R_m$ denote the nodes of a NRM. Let $A = (a_1,…, a_n)$, $a_i \in \{0, 1, –1\}$ is called the Neutrosophic instantaneous state vector of the domain space and it denotes the on-off position of the nodes at any instant. Similarly let $B = (b_1,…, b_n)$ $b_i \in \{0, 1, –1\}$, B is called instantaneous state vector of the range space and it denotes the on-off position of the nodes at any instant, $a_i = 0$ if $a_i$ is off and $a_i = 1$ if $a_i$ is on for i = 1, 2, …, n. Similarly, $b_i = 0$ if $b_i$ is off and $b_i = 1$ if $b_i$ is on for i = 1, 2,…, m.*

**DEFINITION 1.7.16:** *Let $D_1$,…, $D_n$ and $R_1$, $R_2$,…, $R_m$ be the nodes of a NRM. Let $D_i R_j$ (or $R_j D_i$) be the edges of an NRM, j = 1, 2,…, m and i = 1, 2,…, n. The edges form a directed cycle. An NRM is said to be a cycle if it possess a directed cycle. An NRM is said to be acyclic if it does not possess any directed cycle.*

**DEFINITION 1.7.17:** *A NRM with cycles is said to be a NRM with feedback.*

**DEFINITION 1.7.18:** *When there is a feedback in the NRM i.e. when the causal relations flow through a cycle in a revolutionary manner the NRM is called a Neutrosophic dynamical system.*

**DEFINITION 1.7.19:** *Let $D_i$ $R_j$ (or $R_j$ $D_i$) $1 \leq j \leq m$, $1 \leq i \leq n$, when $R_j$ (or $D_i$) is switched on and if causality flows through edges of a cycle and if it again causes $R_j$ (or $D_i$) we say that the Neutrosophical dynamical system goes round and round. This is true for any node $R_j$ ( or $D_i$ ) for $1 \leq j \leq m$ (or $1 \leq i \leq n$). The equilibrium state of this Neutrosophical dynamical system is called the Neutrosophic hidden pattern.*

**DEFINITION 1.7.20:** *If the equilibrium state of a Neutrosophical dynamical system is a unique Neutrosophic state vector, then it is called the fixed point. Consider an NRM with $R_1$, $R_2$, …, $R_m$ and $D_1$, $D_2$,…, $D_n$ as nodes. For example let us start the dynamical system by switching on $R_1$ (or $D_1$). Let us assume that the NRM settles down with $R_1$ and $R_m$ (or $D_1$ and $D_n$) on, or indeterminate on, i.e. the Neutrosophic state vector remains as (1, 0, 0,…, 1) or (1, 0, 0,…I) (or (1, 0, 0,…1) or (1, 0, 0,…I) in D), this state vector is called the fixed point.*

**DEFINITION 1.7.21:** *If the NRM settles down with a state vector repeating in the form $A_1 \to A_2 \to A_3 \to \ldots \to A_i \to A_1$ (or $B_1 \to B_2 \to \ldots \to B_i \to B_1$) then this equilibrium is called a limit cycle.*

Here we give the methods of determining the hidden pattern in a NRM

Let $R_1$, $R_2$,…, $R_m$ and $D_1$, $D_2$,…, $D_n$ be the nodes of a NRM with feedback. Let N(E) be the Neutrosophic Relational Matrix. Let us find the hidden pattern when $D_1$ is switched on i.e. when an input is given as a vector; $A_1$ = (1, 0, …, 0) in D; the data should pass through the relational matrix N(E). This is done by multiplying $A_1$ with the Neutrosophic relational matrix N(E). Let $A_1N(E) = (r_1, r_2, \ldots, r_m)$ after thresholding and updating the resultant vector we get $A_1E \in R$, Now let $B = A_1E$ we pass on B into the system $(N(E))^T$ and obtain $B(N(E))^T$. We update and threshold the vector $B(N(E))^T$ so that $B(N(E))^T \in D$.

This procedure is repeated till we get a limit cycle or a fixed point.

**DEFINITION 1.7.22:** *Finite number of NRMs can be combined together to produce the joint effect of all NRMs. Let $N(E_1)$, $N(E_2)$,…, $N(E_r)$ be the Neutrosophic relational matrices of the NRMs with nodes $R_1$,…, $R_m$ and $D_1$,…,$D_n$, then the combined NRM is represented by the neutrosophic relational matrix N(E) = $N(E_1) + N(E_2) +\ldots+ N(E_r)$.*

## 1.8 Introduction to Fuzzy Associative Memories

A fuzzy set is a map $\mu : X \to [0, 1]$ where X is any set called the domain and [0, 1] the range i.e., μ is thought of as a membership function i.e., to every element $x \in X$ μ assigns membership value in the interval [0, 1]. But very few try to visualize the geometry of fuzzy sets. It is not only of interest but is meaningful to see the geometry of fuzzy sets when we discuss fuzziness. Till date researchers over looked such visualization [Kosko, 92-96], instead they have interpreted fuzzy sets as generalized indicator or membership functions mappings μ from domain X to range [0, 1]. But functions are hard to visualize. Fuzzy theorist often picture membership functions as two-dimensional graphs with the domain X represented as a one-dimensional axis.

The geometry of fuzzy sets involves both domain X = $(x_1,\ldots, x_n)$ and the range [0, 1] of mappings $\mu : X \to [0, 1]$. The geometry of fuzzy sets aids us when we describe fuzziness, define fuzzy concepts and prove fuzzy theorems. Visualizing this geometry may by itself provide the most powerful argument for fuzziness.

An odd question reveals the geometry of fuzzy sets. What does the fuzzy power set $F(2^X)$, the set of all fuzzy subsets of X, look like? It looks like a cube, What does a fuzzy set look like? A fuzzy subsets equals the unit hyper cube $I^n = [0, 1]^n$. The fuzzy set is a point in the cube $I^n$. Vertices of the cube $I^n$ define a non-fuzzy set. Now with in the unit hyper cube $I^n = [0, 1]^n$ we are interested in a distance between points, which led to measures of size and fuzziness of a fuzzy set and more fundamentally to a measure. Thus within cube theory directly extends to the continuous case when the space X is a subset of $R^n$. The next step is to consider mappings between fuzzy cubes. This level of abstraction provides a surprising and fruitful alternative to the prepositional and predicate calculus reasoning techniques used in artificial intelligence (AI) expert systems. It allows us to reason with sets instead of propositions. The fuzzy set framework is numerical and multidimensional. The AI framework is symbolic and is one dimensional with usually only bivalent expert rules or propositions allowed. Both frameworks can encode structured knowledge in linguistic form. But the fuzzy approach translates the structured knowledge into a flexible numerical framework and processes it in a manner that resembles neural network processing. The numerical framework also allows us to adaptively infer and modify fuzzy systems perhaps with neural or statistical techniques directly from problem domain sample data.

Between cube theory is fuzzy-systems theory. A fuzzy set defines a point in a cube. A fuzzy system defines a mapping between cubes. A fuzzy system S maps fuzzy sets to fuzzy sets. Thus a fuzzy system S is a transformation S: $I^n \to I^P$. The n-dimensional unit hyper cube $I^n$ houses all the fuzzy subsets of the domain space or input universe of discourse X = $\{x_1, \ldots, x_n\}$. $I^p$ houses all the fuzzy subsets of the range space or output universe of discourse, Y = $\{y_1, \ldots, y_p\}$. X and Y can also denote subsets of $R^n$ and $R^p$. Then the fuzzy power sets F $(2^X)$ and F $(2^Y)$ replace $I^n$ and $I^p$.

In general a fuzzy system S maps families of fuzzy sets to families of fuzzy sets thus S: $I^{n_1} \times \ldots \times I^{n_r} \to I^{p_1} \times \ldots \times I^{p_s}$ Here too we can extend the definition of a fuzzy system to allow arbitrary products or arbitrary mathematical spaces to serve as the domain or range spaces of the fuzzy sets. We shall focus on fuzzy systems S: $I^n \to I^P$ that map balls of fuzzy sets in $I^n$ to balls of fuzzy set in $I^p$. These continuous fuzzy systems behave as associative memories. The map close inputs to close outputs. We shall refer to them as Fuzzy Associative Maps or FAMs.

The simplest FAM encodes the FAM rule or association $(A_i, B_i)$, which associates the p-dimensional fuzzy set $B_i$ with the n-dimensional fuzzy set $A_i$. These minimal FAMs essentially map one ball in $I^n$ to one ball in $I^p$. They are comparable to simple neural networks. But we need not adaptively train the minimal

FAMs. As discussed below, we can directly encode structured knowledge of the form, "If traffic is heavy in this direction then keep the stop light green longer" is a Hebbian-style FAM correlation matrix. In practice we sidestep this large numerical matrix with a virtual representation scheme. In the place of the matrix the user encodes the fuzzy set association (Heavy, longer) as a single linguistic entry in a FAM bank linguistic matrix. In general a FAM system F: $I^n \to I^b$ encodes the processes in parallel a FAM bank of m FAM rules ($A_1$, $B_1$), …, ($A_m$ $B_m$). Each input A to the FAM system activates each stored FAM rule to different degree. The minimal FAM that stores ($A_i$, $B_i$) maps input A to $B_i$' a partly activated version of $B_i$. The more A resembles $A_i$, the more $B_i$' resembles $B_i$. The corresponding output fuzzy set B combines these partially activated fuzzy sets $B_1^1, B_2^1, \ldots, B_m^1$. B equals a weighted average of the partially activated sets B = $w_1 B_1^1 + \ldots + w_n B_m^1$ where $w_i$ reflects the credibility frequency or strength of fuzzy association ($A_i$, $B_i$). In practice we usually defuzzify the output waveform B to a single numerical value $y_j$ in Y by computing the fuzzy centroid of B with respect to the output universe of discourse Y.

More generally a FAM system encodes a bank of compound FAM rules that associate multiple output or consequent fuzzy sets $B_i$, …, $B_i^s$ with multiple input or antecedent fuzzy sets $A_i^1$, …, $A_i^r$. We can treat compound FAM rules as compound linguistic conditionals. This allows us to naturally and in many cases easily to obtain structural knowledge. We combine antecedent and consequent sets with logical conjunction, disjunction or negation. For instance, we could interpret the compound association ($A^1$, $A^2$, B), linguistically as the compound conditional "IF $X^1$ is $A^1$ AND $X^2$ is $A^2$, THEN Y is B" if the comma is the fuzzy association ($A^1$, $A^2$, B) denotes conjunction instead of say disjunction.

We specify in advance the numerical universe of discourse for fuzzy variables $X^1$, $X^2$ and Y. For each universe of discourse or fuzzy variable X, we specify an appropriate library of fuzzy set values $A_1^r$, …, $A_k^2$ Contiguous fuzzy sets in a library overlap. In principle a neural network can estimate these libraries of fuzzy sets. In practice this is usually unnecessary.

The library sets represent a weighted though overlapping quantization of the input space X. They represent the fuzzy set values assumed by a fuzzy variable. A different library of fuzzy sets similarly quantizes the output space Y. Once we define the library of fuzzy sets we construct the FAM by choosing appropriate combinations of input and output fuzzy sets Adaptive techniques can make, assist or modify these choices.

An Adaptive FAM (AFAM) is a time varying FAM system. System parameters gradually change as the FAM system samples and processes data. Here we discuss how natural network algorithms can adaptively infer FAM rules from training data. In principle, learning can modify other FAM system components, such as the libraries of fuzzy sets or the FAM-rule weights $w_i$.

In the following section we propose and illustrate an unsupervised adaptive clustering scheme based on competitive learning to blindly generate and refine the bank of FAM rules. In some cases we can use supervised learning techniques if we have additional information to accurately generate error estimates. Thus Fuzzy Associative Memories (FAMs) are transformation. FAMs map fuzzy sets to fuzzy sets. They map unit cubes to unit cubes. In simplest case the FAM system consists of a single association. In general the FAM system consists of a bank of different FAM association. Each association corresponds to a different numerical FAM matrix or a different entry in a linguistic FAM-bank matrix. We do not combine these matrices as we combine or superimpose neural-network associative memory matrices. We store the matrices and access them in parallel. We begin with single association FAMs. We proceed on to adopt this model to the problem.

## 1.9 Some Basic Concepts of BAM

Now we go forth to describe the mathematical structure of the Bidirectional Associative Memories (BAM) model. Neural networks recognize ill defined problems without an explicit set of rules. Neurons behave like functions, neurons transduce an unbounded input activation x(t) at time t into a bounded output signal S(x(t)) i.e. Neuronal activations change with time.

Artificial neural networks consists of numerous simple processing units or neurons which can be trained to estimate sampled functions when we do not know the form of the functions. A group of neurons form a field. Neural networks contain many field of neurons. In our text $F_x$ will denote a neuron field, which contains n neurons, and $F_y$ denotes a neuron field, which contains p neurons. The neuronal dynamical system is described by a system of first order differential equations that govern the time-evolution of the neuronal activations or which can be called also as membrane potentials.

$$\dot{x}_i = g_i(X, Y, ...)$$
$$\dot{y}_j = h_j(X, Y, ...)$$

where $\dot{x}_i$ and $\dot{y}_j$ denote respectively the activation time function of the $i^{th}$ neuron in $F_X$ and the $j^{th}$ neuron in $F_Y$. The over dot denotes time differentiation, $g_i$ and $h_j$ are some functions of X, Y, ... where $X(t) = (x_1(t), ... , x_n(t))$ and $Y(t) = (y_1(t), ... , y_p(t))$ define the state of the neuronal dynamical system at time t.

The passive decay model is the simplest activation model, where in the absence of the external stimuli, the activation decays in its resting value

$$\dot{x}_i = x_i$$
and
$$\dot{y}_j = y_j.$$

The passive decay rate $A_i > 0$ scales the rate of passive decay to the membranes resting potentials $\dot{x}_i = –A_ix_i$. The default rate is $A_i = 1$, i.e. $\dot{x}_i = –A_ix_i$. The membrane time constant $C_i > 0$ scales the time variables of the activation dynamical system. The default time constant is $C_i = 1$. Thus $C_i\dot{x}_i = –A_ix_i$.

The membrane resting potential $P_i$ is defined as the activation value to which the membrane potential equilibrates in the absence of external inputs. The resting potential is an additive constant and its default value is zero. It need not be positive.

$$P_i = C_i \dot{x}_i + A_i x_i$$
$$I_i = \dot{x}_i + x_i$$

is called the external input of the system. Neurons do not compute alone. Neurons modify their state activations with external input and with feed back from one another. Now, how do we transfer all these actions of neurons activated by inputs their resting potential etc. mathematically. We do this using what are called synaptic connection matrices.

Let us suppose that the field $F_X$ with n neurons is synaptically connected to the field $F_Y$ of p neurons. Let $m_{ij}$ be a synapse where the axon from the $i^{th}$ neuron in $F_X$ terminates. $M_{ij}$ can be positive, negative or zero. The synaptic matrix M is a n by p matrix of real numbers whose entries are the synaptic efficacies $m_{ij}$.

The matrix M describes the forward projections from the neuronal field $F_X$ to the neuronal field $F_Y$. Similarly a p by n synaptic matrix N describes the backward projections from $F_Y$ to $F_X$. Unidirectional networks occur when a neuron field synaptically intra connects to itself. The matrix M be a n by n square matrix. A Bidirectional network occur if $M = N^T$ and $N = M^T$. To describe this synaptic connection matrix more simply, suppose the n neurons in the field $F_X$ synaptically connect to the p-neurons in field $F_Y$. Imagine an axon from the $i^{th}$ neuron in $F_X$ that terminates in a synapse $m_{ij}$, that about the $j^{th}$ neuron in $F_Y$. We assume that the real number $m_{ij}$ summarizes the synapse and that $m_{ij}$ changes so slowly relative to activation fluctuations that is constant.

Thus we assume no learning if $m_{ij} = 0$ for all t. The synaptic value $m_{ij}$ might represent the average rate of release of a neuro-transmitter such as norepinephrine. So, as a rate, $m_{ij}$ can be positive, negative or zero.

When the activation dynamics of the neuronal fields $F_X$ and $F_Y$ lead to the overall stable behaviour the bidirectional networks are called as Bidirectional Associative Memories (BAM). As in the analysis of the HIV/AIDS patients relative to the migrancy we state that the BAM model studied presently and predicting the future after a span of 5 or 10 years may not be the same.

For the system would have reached stability and after the lapse of this time period the activation neurons under investigations and which are going to measure the model would be entirely different.

Thus from now onwards more than the uneducated poor the educated rich and the middle class will be the victims of HIV/AIDS. So for this study presently carried out can only give how migration has affected the life style of poor labourer and had led them to be victims of HIV/AIDS.

Further not only a Bidirectional network leads to BAM also a unidirectional network defines a BAM if M is symmetric i.e. $M = M^T$. We in our analysis mainly use BAM which are bidirectional networks. However we may also use unidirectional BAM chiefly depending on the problems under investigations. We briefly describe the BAM model more technically and mathematically.

An additive activation model is defined by a system of n + p coupled first order differential equations that inter connects the fields $F_X$ and $F_Y$ through the constant synaptic matrices M and N.

$$x_i = -A_i x_i + \sum_{j=1}^{p} S_j(y_j) n_{ji} + I_i \qquad (1.9.1)$$

$$y_i = -A_j y_j + \sum_{i=1}^{n} S_i(x_i) m_{ij} + J_j \qquad (1.9.2)$$

$S_i(x_i)$ and $S_j(y_j)$ denote respectively the signal function of the $i^{th}$ neuron in the field $F_X$ and the signal function of the $j^{th}$ neuron in the field $F_Y$.

Discrete additive activation models correspond to neurons with threshold signal functions.

The neurons can assume only two values ON and OFF. ON represents the signal +1, OFF represents 0 or – 1 (– 1 when the representation is bipolar). Additive bivalent models describe asynchronous and stochastic behaviour.

At each moment each neuron can randomly decide whether to change state or whether to emit a new signal given its current activation. The Bidirectional Associative Memory or BAM is a non adaptive additive bivalent neural network. In neural literature the discrete version of the equation (3.1) and (3.2) are often referred to as BAMs.

A discrete additive BAM with threshold signal functions arbitrary thresholds inputs an arbitrary but a constant synaptic connection matrix M and discrete time steps K are defined by the equations

$$x_i^{k+1} = \sum_{j=1}^{p} S_j(y_j^k) m_{ij} + I_i \qquad (1.9.3)$$

$$y_j^{k+1} = \sum_{i=1}^{n} S_i\left(x_i^k\right) m_{ij} + J_i \qquad (1.9.4)$$

where $m_{ij} \in M$ and $S_i$ and $S_j$ are signal functions. They represent binary or bipolar threshold functions. For arbitrary real valued thresholds $U = (U_1, ..., U_n)$ for $F_X$ neurons and $V = (V_1, ..., V_P)$ for $F_Y$ neurons the threshold binary signal functions corresponds to

$$S_i(x_i^k) = \begin{cases} 1 \quad \text{if} \quad x_i^k > U_i \\ S_i(x_i^{k-1}) \text{ if } x_i^k = U_i \\ 0 \quad \text{if} \quad x_i^k < U_i \end{cases} \qquad (1.9.5)$$

and

$$S_j(x_j^k) = \begin{cases} 1 \quad \text{if} \quad y_j^k > V_j \\ S_j(y_j^{k-1}) \text{ if } y_j^k = V_j \\ 0 \quad \text{if} \quad y_j^k < V_j \end{cases} \qquad (1.9.6).$$

The bipolar version of these equations yield the signal value -1 when $x_i < U_i$ or when $y_j < V_j$. The bivalent signal functions allow us to model complex asynchronous state change patterns. At any moment different neurons can decide whether to compare their activation to their threshold. At each moment any of the 2n subsets of $F_X$ neurons or 2p subsets of the $F_Y$ neurons can decide to change state. Each neuron may randomly decide whether to check the threshold conditions in the equations (1.9.5) and (1.9.6). At each moment each neuron defines a random variable that can assume the value ON(+1) or OFF(0 or –1). The network is often assumed to be deterministic and state changes are synchronous i.e. an entire field of neurons is updated at a time. In case of simple asynchrony only one neuron makes a state change decision at a time. When the subsets represent the entire fields $F_X$ and $F_Y$ synchronous state change results.

In a real life problem the entries of the constant synaptic matrix M depends upon the investigator's feelings. The synaptic matrix is given a weightage according to their feelings. If $x \in F_X$ and $y \in F_Y$ the forward projections from $F_X$ to $F_Y$ is defined by the matrix M. $\{F(x_i, y_j)\} = (m_{ij}) = M$, $1 \le i \le n$, $1 \le j \le p$.

The backward projections is defined by the matrix $M^T$. $\{F(y_i, x_i)\} = (m_{ji}) = M^T$, $1 \le i \le n$, $1 \le j \le p$. It is not always true that the backward projections from $F_Y$ to $F_X$ is defined by the matrix $M^T$.

Now we just recollect the notion of bidirectional stability. All BAM state changes lead to fixed point stability. The property holds for synchronous as well as asynchronous state changes.

A BAM system ($F_X$, $F_Y$, M) is bidirectionally stable if all inputs converge to fixed point equilibria. Bidirectional stability is a dynamic equilibrium. The same signal information flows back and forth in a bidirectional fixed point. Let us suppose that A denotes a binary n-vector and B denotes a binary p-vector. Let A be the initial input to the BAM system. Then the BAM equilibrates to a bidirectional fixed point ($A_f$, $B_f$) as

$$\begin{array}{lllll} A & \rightarrow & M & \rightarrow & B \\ A' & \leftarrow & M^T & \leftarrow & B \\ A' & \rightarrow & M & \rightarrow & B' \\ A'' & \leftarrow & M^T & \leftarrow & B' \quad \text{etc.} \\ A_f & \rightarrow & M & \rightarrow & B_f \\ A_f & \leftarrow & M^T & \leftarrow & B_f \quad \text{etc.} \end{array}$$

where A’, A’’, ... and B’, B’’, ... represents intermediate or transient signal state vectors between respectively A and $A_f$ and B and $B_f$. The fixed point of a Bidirectional system is time dependent.

The fixed point for the initial input vectors can be attained at different times. Based on the synaptic matrix M which is developed by the investigators feelings the time at which bidirectional stability is attained also varies accordingly.

## 1.10 Properties of Fuzzy Relations and FREs

In this section we just recollect the properties of fuzzy relations like, fuzzy equivalence relation, fuzzy compatibility relations, fuzzy ordering relations, fuzzy morphisms and sup-i-compositions of fuzzy relation. For more about these concepts please refer [90].

Now we proceed on to define fuzzy equivalence relation. A crisp binary relation R(X, X) that is reflexive, symmetric and transitive is called an equivalence relation. For each element x in X, we can define a crisp set $A_x$, which contains all the elements of X that are related to x, by the equivalence relation.

$$A_x = \{y \mid (x, y) \in R (X, X)\}$$

$A_x$ is clearly a subset of X. The element x is itself contained in $A_x$ due to the reflexivity of R, because R is transitive and symmetric each member of $A_x$, is related to all the other members of $A_x$. Further no member of $A_x$, is related to any element of X not included in $A_x$. This set $A_x$ is referred to an as equivalence class of R (X, X) with respect to x. The members of each equivalence class can be considered equivalent to each other and only to each other under the relation R. The family of all such equivalence classes defined by the relation which is usually denoted by X / R, forms a partition on X.

A fuzzy binary relation that is reflexive, symmetric and transitive is known as a fuzzy equivalence relation or similarity relation. In the rest of this section let us use the latter term. While the max-min form of transitivity is assumed, in the following discussion on concepts; can be generalized to the alternative definition of fuzzy transitivity.

While an equivalence relation clearly groups elements that are equivalent under the relation into disjoint classes, the interpretation of a similarity relation can be approached in two different ways. First it can be considered to effectively group elements into crisp sets whose members are similar to each other to some specified degree. Obviously when this degree is equal to 1, the grouping is an equivalence class. Alternatively however we may wish to consider the degree of similarity that the elements of X have to some specified element $x \in X$. Thus for each $x \in X$, a similarity class can be defined as a fuzzy set in which the membership grade of any particular element represents the similarity of that element to the element x. If all the elements in the class are similar to x to the degree of 1 and similar to all elements outside the set to the degree of 0 then the grouping again becomes an equivalence class. We know every fuzzy relation R can be uniquely represented in terms of its α-cuts by the formula

$$R = \bigcup_{\alpha \in (0,1]} \alpha . {}^{\alpha}R .$$

It is easily verified that if R is a similarity relation then each α-cut, ${}^{\alpha}R$ is a crisp equivalence relation. Thus we may use any similarity relation R and by taking an α - cut ${}^{\alpha}R$ for any value $\alpha \in (0, 1]$, create a crisp equivalence relation that represents the presence of similarity between the elements to the degree α. Each of these equivalence relations form a partition of X. Let $\pi ({}^{\alpha}R)$ denote the partition corresponding to the equivalence relation ${}^{\alpha}R$. Clearly any two elements x and y belong to the same block of this partition if and only if $R (x, y) \geq \alpha$. Each similarity relation is associated with the set $\pi (R) = \{\pi ({}^{\alpha}R) \mid \alpha \in (0,1]\}$ of partition of X. These partitions are nested in the sense that $\pi ({}^{\alpha}R)$ is a refinement of $\pi ({}^{\beta}R)$ if and only $\alpha \geq \beta$.

The equivalence classes formed by the levels of refinement of a similarity relation can be interpreted as grouping elements that are similar to each other and only to each other to a degree not less than α.

Just as equivalences classes are defined by an equivalence relation, similarity classes are defined by a similarity relation. For a given similarity relation R(X, X) the similarity class for each x ∈ X is a fuzzy set in which the membership grade of each element y ∈ X is simply the strength of that elements relation to x or R(x, y). Thus the similarity class for an element x represents the degree to which all the other members of X are similar to x. Expect in the restricted case of equivalence classes themselves, similarity classes are fuzzy and therefore not generally disjoint.

Similarity relations are conveniently represented by membership matrices. Given a similarity relation R, the similarity class for each element is defined by the row of the membership matrix of R that corresponds to that element.

Fuzzy equivalence is a cutworthy property of binary relation R(X, X) since it is preserved in the classical sense in each α-cut of R. This implies that the properties of fuzzy reflexivity, symmetry and max-min transitivity are also cutworthy. Binary relations are symmetric and transitive but not reflexive are usually referred to as quasi equivalence relations.

The notion of fuzzy equations is associated with the concept of compositions of binary relations. The composition of two fuzzy binary relations P (X, Y) and Q (Y, Z) can be defined, in general in terms of an operation on the membership matrices of P and Q that resembles matrix multiplication. This operation involves exactly the same combinations of matrix entries as in the regular matrix multiplication. However the multiplication and addition that are applied to these combinations in the matrix multiplication are replaced with other operations, these alternative operations represent in each given context the appropriate operations of fuzzy set intersections and union respectively. In the max-min composition for example, the multiplication and addition are replaced with the min and max operations respectively.

We shall give the notational conventions. Consider three fuzzy binary relations P (X, Y), Q (Y, Z) and R (X, Z) which are defined on the sets

$X = \{x_i \mid i \in I\}$
$Y = \{y_j \mid j \in J\}$ and
$Z = \{z_k \mid k \in K\}$

where we assume that $I = N_n$ $J = N_m$ and $K = N_5$. Let the membership matrices of P, Q and R be denoted by $P = [p_{ij}]$, $Q = [q_{ij}]$, $R = [r_{ik}]$ respectively, where $p_{ij} = P(x_i, y_j)$, $q_{jk} = Q(y_j, z_k)$ $r_{ij} = R(x_i, z_k)$ for all $i \in I$ $(=N_n)$, $j \in J = (N_m)$ and $k \in K$ $(=N_5)$. This clearly implies that all entries in the matrices P, Q, and R are real numbers from the unit interval [0, 1]. Assume now that the three relations constrain each other in such a way that $P \circ Q = R$ where $\circ$ denotes max-min composition. This means that $\max_{j \in J} \min (p_{ij}, q_{jk}) = r_{ik}$ for all $i \in I$ and $k \in K$. That is the matrix equation $P \circ Q = R$ encompasses $n \times s$ simultaneous equations of the form $\max_{j \in J} \min (p_{ij}, q_{jk}) = r_{ik}$. When two of the components in each of the equations are given and one is unknown these equations are referred to as fuzzy relation equations.

When matrices P and Q are given the matrix R is to determined using $P \circ Q = R$. The problem is trivial. It is solved simply by performing the max-min multiplication – like operation on P and Q as defined by $\max_{j \in J} \min (p_{ij}, q_{jk}) = r_{ik}$. Clearly the solution in this case exists and is unique. The problem becomes far from trivial when one of the two matrices on the left hand side of $P \circ Q = R$ is unknown. In this case the solution is guaranteed neither to exist nor to be unique.

Since R in $P \circ Q = R$ is obtained by composing P and Q it is suggestive to view the problem of determining P (or alternatively Q ) from R to Q (or alternatively R and P) as a decomposition of R with respect to Q (or alternatively with respect to P). Since many problems in various contexts can be formulated as problems of decomposition, the utility of any method for solving $P \circ Q = R$ is quite high. The use of fuzzy relation equations in some applications is illustrated. Assume that we have a method for solving P ° Q = R only for the first decomposition problem (given Q and R).

Then we can directly utilize this method for solving the second decomposition problem as well. We simply write P ° Q = R in the form $Q^{-1}$ o $P^{-1} = R^{-1}$ employing transposed matrices. We can solve $Q^{-1}$ o $P^{-1} = R^{-1}$ for $Q^{-1}$ by our method and then obtain the solution of P ° Q = R by $(Q^{-1})^{-1} = Q$.

We study the problem of partitioning the equations P ° Q = R. We assume that a specific pair of matrices R and Q in the equations P ° Q = R is given. Let each particular matrix P that satisfies P ° Q = R is called its solution and let S (Q, R) = {P | P ° Q = R} denote the set of all solutions (the solution set).

It is easy to see this problem can be partitioned, without loss of generality into a set of simpler problems expressed by the matrix equations $p_i$ o Q = $r_i$ for all $i \in I$ where

$$P_i = [p_{ij} \mid j \in J] \text{ and}$$
$$r_i = [r_{ik} \mid k \in K].$$

Indeed each of the equation in $\max_{j \in J} \min (p_{ij} q_{jk}) = r_{ik}$ contains unknown $p_{ij}$ identified only by one particular value of the index i, that is, the unknown $p_{ij}$ distinguished by different values of i do not appear together in any of the individual equations. Observe that $p_i$, Q, and $r_i$ in $p_i$ ° Q = $r_i$ represent respectively, a fuzzy set on Y, a fuzzy relation on Y × Z and a fuzzy set on Z. Let $S_i$ (Q, $r_i$) = [$p_i$ | $p_i$ o Q = $r_i$] denote, for each $i \in I$, the solution set of one of the simpler problem expressed by $p_i$ ° Q = $r_i$.

Thus the matrices P in S (Q, R) = [P | P ° Q = R ] can be viewed as one column matrix

$$P = \begin{bmatrix} p_1 \\ p_2 \\ \vdots \\ p_n \end{bmatrix}$$

where $p_i \in S_i$ (Q, $r_i$) for all i ∈ I = (=$N_n$). It follows immediately from $\max_{j\in J} \min (p_{ij}\ q_{jk}) = r_{ik}$. That if $\max_{j\in J} q_{jk} < r_{ik}$ for some i ∈ I and some k ∈ K, then no values $p_{ij} \in [0, 1]$ exists (j ∈ J) that satisfy P ∘ Q = R, therefore no matrix P exists that satisfies the matrix equation. This proposition can be stated more concisely as follows if

$$\max_{j\in J} q_{jk} < \max_{j\in J} r_{ik}$$

for some k ∈ K then S (Q, R) = ϕ. This proposition allows us in certain cases to determine quickly that P ∘ Q = R has no solutions its negation however is only a necessary not sufficient condition for the existence of a solution of P ∘ Q = R that is for S (Q, R) ≠ ϕ. Since P ∘ Q = R can be partitioned without loss of generality into a set of equations of the form $p_i \circ Q = r_i$ we need only methods for solving equations of the later form in order to arrive at a solution. We may therefore restrict our further discussion of matrix equations of the form P ∘ Q = R to matrix equation of the simpler form P ∘ Q = r, where p = [$p_j$ | j ∈ J], Q = [$q_{jk}$ | j ∈ J, k ∈ K] and r = {$r_k$ | k ∈ K].

We just recall the solution method as discussed by [43]. For the sake of consistency with our previous discussion, let us again assume that p, Q and r represent respectively a fuzzy set on Y, a fuzzy relation on Y × Z and a fuzzy set on Z. Moreover let J = $N_m$ and K = $N_s$ and let S (Q, r) = {p | p ∘ Q = r} denote the solution set of

$$p \circ Q = r.$$

In order to describe a method of solving p ∘ Q = r we need to introduce some additional concepts and convenient notation. First let $\wp$ denote the set of all possible vectors.

$$p = \{p_j \mid j \in J\}$$

such that $p_j \in [0, 1]$ for all j ∈ J, and let a partial ordering on $\wp$ be defined as follows for any pair $p^1, p^2 \in \wp$ $p^1 \leq p^2$ if and only if $p_i^2 \leq p_j^2$ for all j ∈J. Given an arbitrary pair $p^1, p^2 \in$ $\wp$ such that $p^1 \leq p^2$ let $[p^1, p^2] = \{p \in \wp \mid p^1 \leq p < p^2\}$. For any pair $p^1, p^2 \in \wp$ $(\{p^1, p^2\} \leq \}$ is a lattice.

Now we recall some of the properties of the solution set S (Q, r). Employing the partial ordering on $\wp$, let an element $\hat{p}$ of S (Q, r) be called a maximal solution of p ° Q = r if for all p ∈ S (Q, r), $p \geq \hat{p}$ implies $p = \hat{p}$ if for all p ∈ S (Q, r) $p < \tilde{p}$ then that is the maximum solution. Similar discussion can be made on the minimal solution of p ° Q = r. The minimal solution is unique if $p \geq \hat{p}$ (i.e. $\hat{p}$ is unique).

It is well known when ever the solution set S (Q, r) is not empty it always contains a unique maximum solution $\hat{p}$ and it may contain several minimal solution. Let $\breve{S}$ (Q, r) denote the set of all minimal solutions. It is known that the solution set S (Q, r) is fully characterized by the maximum and minimal solution in the sense that it consists exactly of the maximum solution $\hat{p}$ all the minimal solutions and all elements of $\wp$ that are between $\hat{p}$ and the numeral solution.

Thus S (Q, r) = $\bigcup_{p} [\tilde{p}, \hat{p}]$ where the union is taken for all $\tilde{p} \in \breve{S}$ (Q, r). When S (Q, r) ≠ ϕ, the maximum solution.

$\hat{p} = [\hat{p}_j \mid j \in J]$ of p ° Q = r is determined as follows:

$$\hat{p}_j = \min_{k \in K} \sigma(q_{ik}, r_k) \text{ where } \sigma(q_{jk}, r_k) = \begin{cases} r_k & \text{if } q_{jk} > r_k \\ 1 & \text{otherwise} \end{cases}$$

when $\hat{p}$ determined in this way does not satisfy p ° Q = r then S(Q, r) = ϕ. That is the existence of the maximum solution $\hat{p}$ as determined by $\hat{p}_j = \min_{k \in K} \sigma(q_{ik}, r_k)$ is a necessary and sufficient condition for S (Q, r) ≠ ϕ. Once $\hat{p}$ is determined by $\hat{p}_j = \min_{k \in K} \sigma(q_{ik}, r_k)$, we must check to see if it satisfies the given matrix equations p ° Q = r. If it does not then the equation has no solution (S (Q, r) = ϕ), otherwise $\hat{p}$ in the maximum solution of the equation and we next determine the set $\tilde{S}$ (Q, r) of its minimal solutions.

## 1.11 Binary Neutrosophic Relation and their Properties

In this section we introduce the notion of neutrosophic relational equations and fuzzy neutrosophic relational equations and analyze and apply them to real-world problems, which are abundant with the concept of indeterminacy. We also mention that most of the unsupervised data also involve at least to certain degrees the notion of indeterminacy.

Throughout this section by a neutrosophic matrix we mean a matrix whose entries are from the set N = [0, 1] $\cup$ *I* and by a fuzzy neutrosophic matrix we mean a matrix whose entries are from N' = [0, 1] $\cup$ {n*I* / n $\in$ (0,1]}.

Now we proceed on to define binary neutrosophic relations and binary neutrosophic fuzzy relation.

A binary neutrosophic relation $R_N(x, y)$ may assign to each element of X two or more elements of Y or the indeterminate *I*. Some basic operations on functions such as the inverse and composition are applicable to binary relations as well. Given a neutrosophic relation $R_N(X, Y)$ its domain is a neutrosophic set on X $\cup$ *I* domain R whose membership function is defined by domR(x) = $\max_{y \in X \cup I} R_N(x, y)$ for each x $\in$ X $\cup$ *I*.

That is each element of set X $\cup$ *I* belongs to the domain of R to the degree equal to the strength of its strongest relation to any member of set Y $\cup$ *I*. The degree may be an indeterminate *I* also. Thus this is one of the marked difference between the binary fuzzy relation and the binary neutrosophic relation. The range of $R_N(X,Y)$ is a neutrosophic relation on Y, ran R whose membership is defined by ran R(y) = $\max_{x \in X} R_N(x, y)$ for each y $\in$ Y, that is the strength of the strongest relation that each element of Y has to an element of X is equal to the degree of that element's membership in the range of R or it can be an indeterminate *I*.

The height of a neutrosophic relation $R_N(x, y)$ is a number h(R) or an indeterminate *I* defined by $h_N(R) = \max_{y \in Y \cup I} \max_{x \in X \cup I} R_N(x, y)$. That is $h_N(R)$ is the largest membership grade attained by any pair (x, y) in R or the indeterminate *I*.

A convenient representation of the neutrosophic binary relation $R_N(X, Y)$ are membership matrices $R = [\gamma_{xy}]$ where $\gamma_{xy} \in R_N(x, y)$. Another useful representation of a binary neutrosophic relation is a neutrosophic sagittal diagram. Each of the sets X, Y represented by a set of nodes in the diagram, nodes corresponding to one set are clearly distinguished from nodes representing the other set. Elements of X' × Y' with non-zero membership grades in $R_N(X, Y)$ are represented in the diagram by lines connecting the respective nodes. These lines are labeled

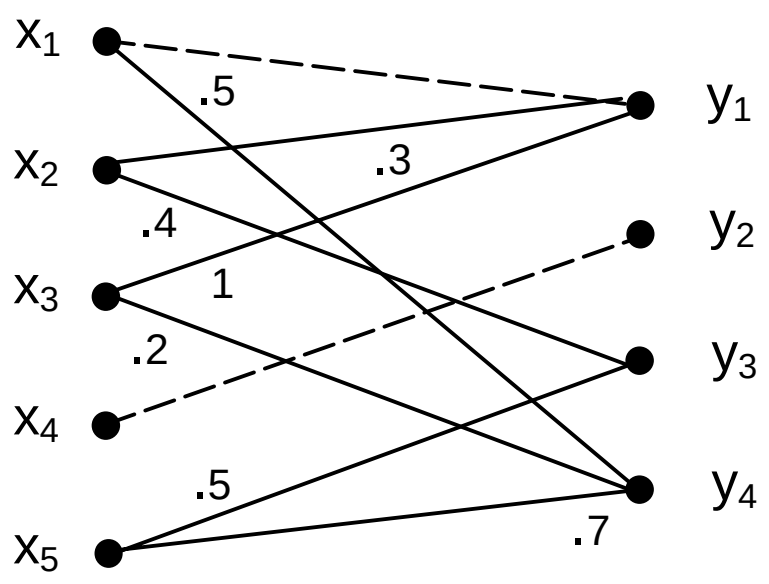

FIGURE: 1.11.1

with the values of the membership grades.

An example of the neutrosophic sagittal diagram is a binary neutrosophic relation $R_N(X, Y)$ together with the membership neutrosophic matrix is given below.

$$\begin{array}{c} \\ x_1 \\ x_2 \\ x_3 \\ x_4 \\ x_5 \end{array}\begin{array}{c} \begin{array}{cccc} y_1 & y_2 & y_3 & y_4 \end{array} \\ \begin{bmatrix} I & 0 & 0 & 0.5 \\ 0.3 & 0 & 0.4 & 0 \\ 1 & 0 & 0 & 0.2 \\ 0 & I & 0 & 0 \\ 0 & 0 & 0.5 & 0.7 \end{bmatrix} \end{array}.$$

The inverse of a neutrosophic relation $R_N(X, Y) = R(x, y)$ for all $x \in X$ and all $y \in Y$. A neutrosophic membership matrix $R^{-1} = [r^{-1}_{yx}]$ representing $R^{-1}_N(Y, X)$ is the transpose of the matrix R for $R_N(X, Y)$ which means that the rows of $R^{-1}$ equal the columns of R and the columns of $R^{-1}$ equal rows of R. Clearly $(R^{-1})^{-1} = R$ for any binary neutrosophic relation.

Consider any two binary neutrosophic relation $P_N(X, Y)$ and $Q_N(Y, Z)$ with a common set Y. The standard composition of these relations which is denoted by $P_N(X, Y) \bullet Q_N(Y, Z)$ produces a binary neutrosophic relation $R_N(X, Z)$ on $X \times Z$ defined by $R_N(x, z) = [P \bullet Q]_N(x, z) = \max_{y \in Y} \min[P_N(x, y), Q_N(x, y)]$ for all $x \in X$ and all $z \in Z$.

This composition which is based on the standard $t_N$-norm and $t_N$-co-norm, is often referred to as the max-min composition. It can be easily verified that even in the case of binary neutrosophic relations $[P_N(X, Y) \bullet Q_N(Y, Z)]^{-1} = Q_N^{-1}(Z, Y) \bullet P_N^{-1}(Y, X)$. $[P_N(X, Y) \bullet Q_N(Y, Z)] \bullet R_N(Z, W) = P_N(X, Y) \bullet [Q_N(Y, Z) \bullet R_N(Z, W)]$, that is, the standard (or max-min) composition is associative and its inverse is equal to the reverse composition of the inverse relation. However, the standard composition is not commutative, because $Q_N(Y, Z) \bullet P_N(X, Y)$ is not well defined when $X \neq Z$. Even if $X = Z$ and $Q_N (Y, Z) \circ P_N (X, Y)$ are well defined still we can have $P_N (X, Y) \circ Q (Y, Z) \neq Q (Y, Z) \circ P (X, Y)$.

Compositions of binary neutrosophic relation can the performed conveniently in terms of membership matrices of the relations. Let $P = [p_{ik}]$, $Q = [q_{kj}]$ and $R = [r_{ij}]$ be membership matrices of binary relations such that $R = P \circ Q$. We write this using matrix notation

$$[r_{ij}] = [p_{ik}] \text{ o } [q_{kj}]$$

where $r_{ij} = \max_k \min (p_{ik}, q_{kj})$.

A similar operation on two binary relations, which differs from the composition in that it yields triples instead of pairs, is known as the relational join. For neutrosophic relation $P_N (X, Y)$ and $Q_N (Y, Z)$ the relational join P * Q corresponding to the neutrosophic standard max-min composition is a ternary relation $R_N (X, Y, Z)$ defined by $R_N (x, y, z) = [P * Q]_N (x, y, z) = \min [P_N (x, y), Q_N (y, z)]$ for each $x \in X$, $y \in Y$ and $z \in Z$.

This is illustrated by the following Figure 1.11.2.

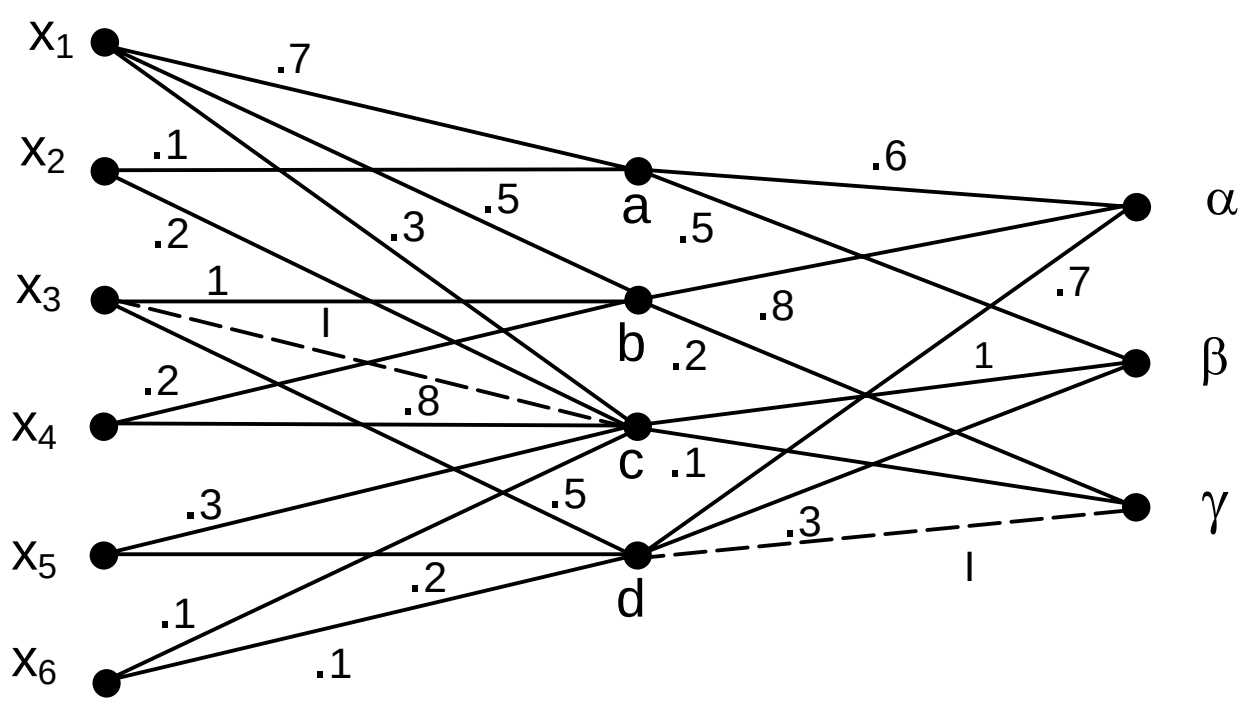

FIGURE: 1.11.2

In addition to defining a neutrosophic binary relation there exists between two different sets, it is also possible to define neutrosophic binary relation among the elements of a single set X.

A neutrosophic binary relation of this type is denoted by $R_N(X, X)$ or $R_N(X^2)$ and is a subset of $X \times X = X^2$.

These relations are often referred to as neutrosophic directed graphs or neutrosophic digraphs. [221-222]

Neutrosophic binary relations $R_N(X, X)$ can be expressed by the same forms as general neutrosophic binary relations. However they can be conveniently expressed in terms of simple diagrams with the following properties:

I. Each element of the set X is represented by a single node in the diagram.
II. Directed connections between nodes indicate pairs of elements of X for which the grade of membership in R is non zero or indeterminate.
III. Each connection in the diagram is labeled by the actual membership grade of the corresponding pair in R or in indeterminacy of the relationship between those pairs.

The neutrosophic membership matrix and the neutrosophic sagittal diagram is as follows for any set X = {a, b, c, d, e}.

$$\begin{bmatrix} 0 & I & 0.3 & 0.2 & 0 \\ 1 & 0 & I & 0 & 0.3 \\ I & 0.2 & 0 & 0 & 0 \\ 0 & 0.6 & 0 & 0.3 & I \\ 0 & 0 & 0 & I & 0.2 \end{bmatrix}.$$

Neutrosophic membership matrix for x is given above and the neutrosophic sagittal diagram is given below.

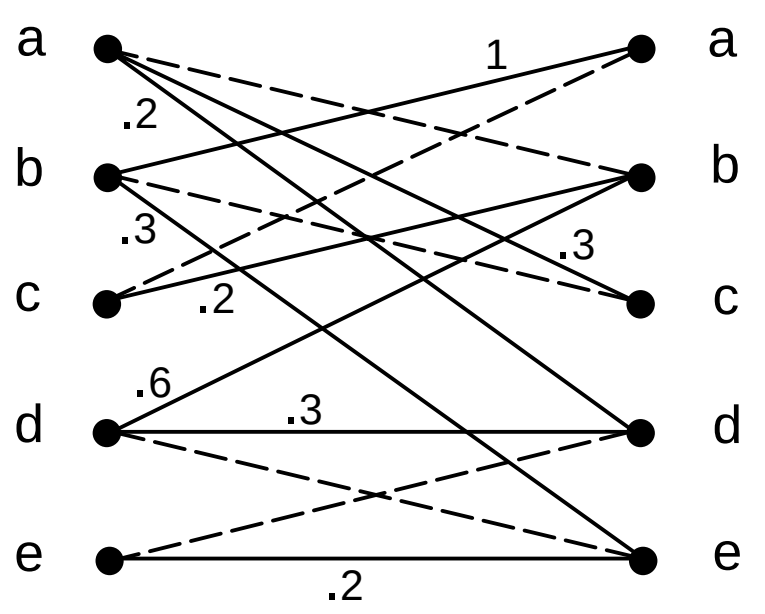

FIGURE 1.11.3

Neutrosophic diagram or graph is left for the reader as an exercise.

The notion of reflexivity, symmetry and transitivity can be extended for neutrosophic relations $R_N$ (X, Y) by defining them in terms of the membership functions or indeterminacy relation.

*Thus $R_N$ (X, X) is reflexive if and only if $R_N$ (x, x) = 1 for all x $\in$ X.* If *this is not the case for some x $\in$ X the relation is irreflexive.*

*A weaker form of reflexivity, if for no x in X; $R_N$(x, x) = 1 then we call the relation to be anti-reflexive referred to as $\in$-reflexivity, is sometimes defined by requiring that*

$$R_N (x, x) \geq \in \text{ where } 0 < \in < 1.$$

*A fuzzy relation is symmetric if and only if*

$$R_N (x, y) = R_N (y, x) \text{ for all } x, y, \in X.$$

*Whenever this relation is not true for some x, y* $\in$ *X the relation is called asymmetric. Furthermore when* $R_N(x, y) > 0$ *and* $R_N(y, x) > 0$ *implies that x = y for all x, y* $\in$ *X the relation R is called anti-symmetric.*

*A fuzzy relation* $R_N(X, X)$ *is transitive (or more specifically max-min transitive) if*

$$R_N(x, z) \geq \max_{y \in Y} \min [R_N(x, y), R_N(y, z)]$$

*is satisfied for each pair* $(x, z) \in X^2$*. A relation failing to satisfy the above inequality for some members of X is called non-transitive and if* $R_N(x, x) < \max_{y \in Y} \min [RN(x, y), RN(y, z)]$ *for all* $(x, x) \in X^2$*, then the relation is called anti-transitive.*

*Given a relation* $R_N(X, X)$ *its transitive closure* $\overline{R}_{NT}(x, X)$ *can be analyzed in the following way.*

The transitive closure on a crisp relation $R_N(X, X)$ is defined as the relation that is transitive, contains

$$R_N(X, X) < \max_{y \in Y} \min [R_N(x, y) R_N(y, z)]$$

for all $(x, x) \in X^2$, then the relation is called anti-transitive. Given a relation $R_N(x, x)$ its transitive closure $\overline{R}_{NT}(X, X)$ can be analyzed in the following way.

The transitive closure on a crisp relation $R_N(X, X)$ is defined as the relation that is transitive, contains $R_N$ and has the fewest possible members. For neutrosophic relations the last requirement is generalized such that the elements of transitive closure have the smallest possible membership grades, that still allow the first two requirements to be met.

Given a relation $R_N(X, X)$ its transitive closure $\overline{R}_{NT}(X, X)$ can be determined by a simple algorithm.

Now we proceed on to define the notion of neutrosophic equivalence relation.

**DEFINITION 1.11.1:** *A crisp neutrosophic relation* $R_N(X, X)$ *that is reflexive, symmetric and transitive is called an neutrosophic*

*equivalence relation. For each element x in X, we can define a crisp neutrosophic set $A_x$ which contains all the elements of X that are related to x by the neutrosophic equivalence relation.*

*Formally $A_x$ = [ y | (x, y) $\in R_N$ (X, X)}. $A_x$ is clearly a subset of X. The element x is itself contained in $A_x$, due to the reflexivity of R because R is transitive and symmetric each member of $A_x$ is related to all other members of $A_x$. Further no member of $A_x$ is related to any element of X not included in $A_x$. This set $A_x$ is clearly referred to as an neutrosophic equivalence class of $R_N$ (X, x) with respect to x. The members of each neutrosophic equivalence class can be considered neutrosophic equivalent to each other and only to each other under the relation R.*

But here it is pertinent to mention that in some X; (a, b) may not be related at all to be more precise there may be an element a $\in$ X which is such that its relation with several or some elements in X \ {a} is an indeterminate. The elements which cannot determine its relation with other elements will be put in as separate set.

A neutrosophic binary relation that is reflexive, symmetric and transitive is known as a neutrosophic equivalence relation.

Now we proceed on to define Neutrosophic intersections neutrosophic t – norms ($t_N$ – norms)

Let A and B be any two neutrosophic sets, the intersection of A and B is specified in general by a neutrosophic binary operation on the set N = [0, 1] $\cup$ *I*, that is a function of the form

$$i_N: N \times N \to N.$$

For each element x of the universal set, this function takes as its argument the pair consisting of the elements membership grades in set A and in set B, and yield the membership grade of the element in the set constituting the intersection of A and B. Thus,

$$(A \cap B)(x) = i_N [A(x), B(x)] \text{ for all } x \in X.$$

In order for the function $i_N$ of this form to qualify as a fuzzy intersection, it must possess appropriate properties, which ensure that neutrosophic sets produced by $i_N$ are intuitively acceptable as meaningful fuzzy intersections of any given pair of neutrosophic sets. It turns out that functions known as $t_N$-norms, will be introduced and analyzed in this section. In fact the class of $t_N$- norms is now accepted as equivalent to the class of neutrosophic fuzzy intersection. We will use the terms $t_N$ – norms and neutrosophic intersections inter changeably.

Given a $t_N$ – norm, $i_N$ and neutrosophic sets A and B we have to apply:

$$(A \cap B)(x) = i_N [A(x), B(x)]$$

for each $x \in X$, to determine the intersection of A and B based upon $i_N$.

However the function $i_N$ is totally independent of x, it depends only on the values A (x) and B(x). Thus we may ignore x and assume that the arguments of $i_N$ are arbitrary numbers $a, b \in [0, 1] \cup I = N$ in the following examination of formal properties of $t_N$-norm.

A neutrosophic intersection/ $t_N$-norm $i_N$ is a binary operation on the unit interval that satisfies at least the following axioms for all $a, b, c, d \in N = [0, 1] \cup I$.

| | |
|---|---|
| $1_N$ | $i_N(a, 1) = a$ |
| $2_N$ | $i_N(a, I) = I$ |
| $3_N$ | $b \leq d$ implies |
| | $i_N(a, b) \leq i_N(a, d)$ |
| $4_N$ | $i_N(a, b) = i_N(b, a)$ |
| $5_N$ | $i_N(a, i_N(b, d)) = i_N(a, b), d)$. |

We call the conditions $1_N$ to $5_N$ as the axiomatic skeleton for neutrosophic intersections / $t_N$ – norms. Clearly $i_N$ is a continuous function on $N \setminus I$ and $i_N(a, a) < a \ \forall a \in N \setminus I$

$$i_N(I\ I) = I.$$

If $a_1 < a_2$ and $b_1 < b_2$ implies $i_N(a_1, b_1) < i_N(a_2, b_2)$. Several properties in this direction can be derived as in case of t-norms.

The following are some examples of $t_N$ –norms

1. $i_N(a, b) = \min(a, b)$

$i_N$ (a, *I*) = min (a, *I*) = *I* will be called as standard neutrosophic intersection.

2. $i_N$ (a, b) = ab for a, b $\in$ N \ *I* and $i_N$ (a, b) = *I* for a, b $\in$ N where a = *I* or b = *I* will be called as the neutrosophic algebraic product.

3. Bounded neutrosophic difference.

   $i_N$ (a, b) = max (0, a + b – 1) for a, b $\in$ *I*

   $i_N$ (a, *I*) = *I* is yet another example of $t_N$ – norm.

1. Drastic neutrosophic intersection
2.

$$i_N (a, b) = \begin{cases} a \text{ when } b = 1 \\ b \text{ when } a = 1 \\ I \text{ when } a = I \\ \quad\quad \text{or } b = I \\ \quad\quad \text{or } a = b = I \\ 0 \quad \text{otherwise.} \end{cases}$$

As *I* is an indeterminate adjoined in $t_N$ – norms. It is not easy to give then the graphs of neutrosophic intersections. Here also we leave the analysis and study of these $t_N$ – norms (i.e. neutrosophic intersections) to the reader.

The notion of neutrosophic unions closely parallels that of neutrosophic intersections. Like neutrosophic intersection the general neutrosophic union of two neutrosophic sets A and B is specified by a function

$$\mu_N: N \times N \to N \text{ where } N = [0, 1] \cup I.$$

The argument of this function is the pair consisting of the membership grade of some element x in the neutrosophic set A and the membership grade of that some element in the neutrosophic set B, (here by membership grade we mean not only the membership grade in the unit interval [0, 1] but also the indeterminacy of the membership). The function returns the membership grade of the element in the set A $\cup$ B.

Thus (A $\cup$ B) (x) = $\mu_N$ [A (x), B(x)] for all x $\in$ X. Properties that a function $\mu_N$ must satisfy to be initiatively acceptable as neutrosophic union are exactly the same as properties of functions that are known. Thus neutrosophic union will be called as neutrosophic t-co-norm; denoted by $t_N$ – co-norm.

A neutrosophic union / $t_N$ – co-norm $\mu_N$ is a binary operation on the unit interval that satisfies at least the following conditions for all a, b, c, d $\in$ N = $[0, 1] \cup I$

$C_1$ $\mu_N (a, 0) = a$

$C_2$ $\mu_N (a, I) = I$

$C_3$ $b \leq d$ implies $\mu_N (a, b) \leq \mu_N (a, d)$

$C_4$ $\mu_N (a, b) = \mu_N (b, a)$

$C_5$ $\mu_N (a, \mu_N (b, d)) = \mu_N (\mu_N (a, b), d)$.

Since the above set of conditions are essentially neutrosophic unions we call it the axiomatic skeleton for neutrosophic unions / $t_N$-co-norms. The addition requirements for neutrosophic unions are

i. $\mu_N$ is a continuous functions on $N \setminus \{I\}$

ii. $\mu_N (a, a) > a$.

iii. $a_1 < a_2$ and $b_1 < b_2$ implies $\mu_N (a_1. b_1) < \mu_N (a_2, b_2)$; $a_1, a_2, b_1, b_2 \in N \setminus \{I\}$.

We give some basic neutrosophic unions.
Let $\mu_N : [0, 1] \times [0, 1] \to [0, 1]$

$\mu_N (a, b) = \max (a, b)$
$\mu_N (a, I) = I$ is called as the standard neutrosophic union.
$\mu_N (a, b) = a + b – ab$ and
$\mu_N (a, I) = I$ .

This function will be called as the neutrosophic algebraic sum.

$$\mu_N (a, b) = \min (1, a + b) \text{ and } \mu_N (a, I) = I$$

will be called as the neutrosophic bounded sum. We define the notion of neutrosophic drastic unions

$$\mu_N (a, b) = \begin{cases} a \text{ when } b = 0 \\ b \text{ when } a = 0 \\ I \text{ when } a = I \\ \quad \text{or } b = I \\ 1 \ \text{ otherwise.} \end{cases}$$

Now we proceed on to define the notion of neutrosophic Aggregation operators. Neutrosophic aggregation operators on neutrosophic sets are operations by which several neutrosophic sets are combined in a desirable way to produce a single neutrosophic set.

Any neutrosophic aggregation operation on n neutrosophic sets ($n \geq 2$) is defined by a function $h_N$: $N^n \to N$ where N = [0, 1] $\cup$ *I* and $N^n = \underbrace{N \times ... \times N}_{n-\text{times}}$ when applied to neutrosophic sets $A_1$, $A_2$,…, $A_n$ defined on X the function $h_N$ produces an aggregate neutrosophic set A by operating on the membership grades of these sets for each x $\in$ X (Here also by the term membership grades we mean not only the membership grades from the unit interval [0, 1] but also the indeterminacy *I* for some x $\in$ X are included). Thus

$$A_N (x) = h_N (A_1 (x), A_2 (x),…, A_n(x))$$

for each x $\in$ X.

**Chapter Two**

# INTRODUCTION TO FUZZY INTERVAL MATRICES AND NEUTROSOPHIC INTERVAL MATRICES AND THEIR GENERALIZATIONS

Zadeh has introduced to the world the notion of fuzzy theory in the year (1965). Though four decades has passed lots of development have been made in fuzzy theory or more so in their applications.

Till date the concept of fuzzy interval matrices has not been introduced. This book ventures to introduce for the first time the notion of fuzzy interval matrices and neutrosophic interval matrices and give a few of their application to fuzzy and neutrosophic models.

The notion of fuzzy matrices and neutrosophic matrices are used in models like FCM, NCM, FRM, NRM, BAM, FAM and so on. Now we venture to define fuzzy interval matrices, fuzzy interval bimatrices, neutrosophic interval matrices, neutrosophic interval bimatrices, neutrosophic interval n-matrices fuzzy neutrosophic interval matrices and fuzzy neutrosophic interval n-matrices.

This chapter has three sections. In section one the authors introduce the new notion of fuzzy interval matrices and describe them with examples. The notion of fuzzy interval bimatrices and their generalizations are introduced in section two. The neutrosophic interval matrices and the concept of fuzzy neutrosophic interval matrices are introduced in section three. Here the generalization of these concepts are also carried out.

## 2.1 Fuzzy Interval Matrices

In the section we introduce the notion of fuzzy interval matrices [A, B] where A and B are fuzzy m × n matrices i.e., A = $(a_{ij})$, B = $(b_{ij})$ and $a_{ij}$, $b_{ij} \in [0, 1]$ further $a_{ij} \leq b_{ij}$ for $1 \leq i \leq m$ and $1 \leq j \leq n$. Thus the fuzzy interval matrices or interval fuzzy matrices contain all matrices C = $(c_{ij})$ such that $(c_{ij})$ takes its value in between $a_{ij}$ and $b_{ij}$ i.e., C ∈ [A, B] implies $a_{ij} \leq c_{ij} \leq b_{ij}$.

It we take m = n then we say the fuzzy interval matrices are square matrices i.e., [A, B] is the fuzzy interval square matrices. Here it is important to mention that in the definition of the fuzzy interval matrices be it square or rectangle one, we can take any interval [a, b] such that $0 \leq a < b \leq 1$; i.e., $[a, b] \subseteq [0, 1]$.

So it is up to the convenience of the expert to choose the whole interval or a proper subinterval. Further it is important to mention here that we in our models will take also the interval [–1, 1] and call it as fuzzy interval matrices defined on the fuzzy interval [–1, 1] this will be finding its application when we work with the FCM and FRM models.

Just we will illustrate these concepts by examples.

***Example 2.1.1:*** Let [A, B] be a fuzzy interval 2 × 3 matrix, where

$$A = \begin{bmatrix} 0 & 0.2 & 0.6 \\ 0.1 & 0.4 & 0.21 \end{bmatrix}$$

and

$$B = \begin{bmatrix} 0.6 & 0.6 & 1 \\ 1 & 0.8 & 0.9 \end{bmatrix}$$

We see [A, B] is a fuzzy interval 2 × 3 matrices.

$$C = \begin{bmatrix} 0.2 & 0.4 & 0.7 \\ 0.6 & 0.5 & 0.8 \end{bmatrix}$$

is a matrix belonging to the fuzzy interval matrix [A, B] .

$$D = \begin{bmatrix} 0.8 & 1 & 0 \\ 0.1 & 0.2 & 0.8 \end{bmatrix}.$$

Clearly D $\notin$ [A, B].

But still D is a fuzzy matrix but does not belong to this interval of fuzzy matrices [A, B].

Now we proceed on to give an example of a fuzzy interval square matrices.

***Example 2.1.2:*** Now consider the 3 × 3 square fuzzy interval matrix given by [A, B] where

$$A = \begin{bmatrix} 0 & 0.3 & 0.2 \\ 0.1 & 0 & 0.3 \\ 0.4 & 0.5 & 0 \end{bmatrix}$$

and

$$B = \begin{bmatrix} 0.6 & 0.7 & 0.5 \\ 1 & 0 & 0.8 \\ 0.9 & 0.6 & 1 \end{bmatrix}$$

defined in the interval [0, 1] we say C = ($C_{ij}$) $\in$ [A, B] if $0 \le C_{11} \le 0.6$, $0.3 \le C_{12} < 0.7$, $0.2 \le C_{13} \le 0.5$, $0.1 \le C_{21} \le 1$, $0 \le C_{22} \le 0$, $0.3 \le C_{23} \le 0.8$, $0.4 \le C_{31} \le 0.9$, $0.5 \le C_{32} \le 0.6$ and $0 \le C_{33} \le 1$. All 3 × 3 fuzzy matrices need not belong to this fuzzy interval 3 × 3 matrix for take

$$E = \begin{bmatrix} 1 & 0 & 1 \\ 0 & 1 & 0 \\ 0.8 & 0.5 & 0.3 \end{bmatrix}.$$

Clearly E is a 3 × 3 fuzzy matrix but E does not belong to this fuzzy interval matrix [A, B]. Let [A, B] be a fuzzy interval matrix defined on the interval [0, 1] where A = ($a_{ij}$) and B = ($b_{ij}$). We define the fuzzy matrix

$$M = \left( \frac{a_{ij} + b_{ij}}{2} \right)$$

to be the medianal fuzzy matrix of the fuzzy interval matrix [A, B]. Every fuzzy interval matrix [A, B] defined on the interval [0, 1] (or [a, b] ⊆ [0 1], it 0 ≤ a < b ≤ 1) has a unique medianal fuzzy matrix of the fuzzy interval matrix [a, b].

Every fuzzy interval matrix [A, B] defined on the interval [0, 1] (or [a, b] ⊆ [0 1] i.e. 0 ≤ a < b ≤ 1) has a unique medianal matrix. We call the matrix A = ($a_{ij}$) to be the minimal fuzzy matrix of the fuzzy interval matrix [A, B] and B = ($b_{ij}$) is defined to be the maximal fuzzy matrix of the fuzzy interval matrix [A, B]. Clearly the minimal and maximal matrices are unique we further say [A, B] is a fuzzy closed interval matrix.

(A, B) will be called as the fuzzy open interval matrices and [A, B) and (A, B] will be called as the half closed open and a half open closed fuzzy interval matrix. When we have open interval fuzzy matrix or half open-closed interval fuzzy matrix or closed open interval fuzzy matrix we have only quasi medianal matrix. We have only quasi minimal matrix and quasi maximal matrix for the fuzzy interval matrix which is a fuzzy open interval matrix.

When we have a fuzzy open-closed interval matrix (A, B] then we have a quasi minimal fuzzy matrix a quasi medianal fuzzy matrix but a unique maximal fuzzy matrix. When the fuzzy interval matrix is a fuzzy closed open interval matrix [A, B) then we have a unique minimal fuzzy matrix given by A = ($a_{ij}$) but only a quasi medianal and quasi maximal fuzzy matrices.

These concepts will be helpful to us when we apply these matrices in real fuzzy models. How ever we have not defined explicitly, but illustrate with example to show how these closed fuzzy interval matrix or open fuzzy interval matrix or open-closed fuzzy interval matrix or closed-open fuzzy interval matrix by the following examples:

When we say [A, B] is a fuzzy closed interval matrix defined on the interval [0, 1] (or [a, b]) we mean if $(a_{ij}) = A$ then atleast one of the $a_{ij}$ is 0 in case it is defined on [0, 1] (or atleast one of $a_{ij} = a$ if defined on [a, b]), and for $B = (b_{ij})$ atleast one of $b_{ij} =1$ if defined on [0, 1] (or atleast) one of $b_{ij} = b$ if defined on [a, b]).

Now we illustrate this by the following example:

***Example 2.1.3:*** Let [A, B] be a fuzzy interval matrix defined on the interval say [0, 0.7], if

$$A = \begin{bmatrix} 0 & 0 & 0.2 & 0.3 \\ 0.1 & 0 & 0.3 & 0.4 \\ 0 & 0.1 & 0.6 & 0.2 \end{bmatrix}$$

and

$$B = \begin{bmatrix} 0.7 & 0.6 & 0.6 & 0.5 \\ 0.3 & 0 & 0.7 & 0.6 \\ 0 & 0.7 & 0.6 & 0.4 \end{bmatrix}.$$

Clearly [A, B] is a closed interval fuzzy matrix.
If

$$A = \begin{bmatrix} 1 & 0.1 & 0.3 & 0.4 \\ 0.2 & 0.2 & 0.6 & 0.1 \\ 0.1 & 0.2 & 0.4 & 0.1 \end{bmatrix}$$

then clearly [A, B] is not a closed interval fuzzy matrix as none of the $a_{ij}$ is 0, where $A = (a_{ij})$.
If

$$B = \begin{bmatrix} 0.2 & 0.6 & 0.3 & 0.6 \\ 0.4 & 0.5 & 0.5 & 0.5 \\ 0.4 & 0.6 & 0.6 & 0.3 \end{bmatrix}$$

then the fuzzy interval matrix [A, B] defined on the interval [0, 0.7] is not a closed fuzzy interval matrix as none of the $b_{ij}$ is 0.7 where B = ($b_{ij}$). Thus we now proceed on to define fuzzy open interval matrix [A, B] defined on the interval [0, 1] (or on the interval [a, b] $\subseteq$ [0, 1]).

Let (A B) be a fuzzy open interval matrix defined on the interval [0, 1] (or on the interval [a, b], [a, b] $\subseteq$ [0, 1] i.e., $0 \leq a < b \leq 1$). We say (A, B) is a open interval fuzzy matrix if A= ($a_{ij}$) and B = ($b_{ij}$) then $a_{ij} \neq 0$ for any i and j and $b_{ij} \neq 1$ for any i, j [i.e., in case it is defined on the interval [a, b] none of the $a_{ij}$ = a and none of the $b_{ij}$ = b].

Now we illustrate this by the following example:

***Example 2.1.4:*** Let (A, B) be a open fuzzy interval matrix defined on the interval (0.1, 0.8). Given A = ($a_{ij}$) and B = ($b_{ij}$) are $4 \times 4$ matrices with entries from the open interval (0.1, 0.8) i.e.

$$A = \begin{bmatrix} 0.2 & 0.3 & 0.2 & 0.2 \\ 0.4 & 0.2 & 0.3 & 0.2 \\ 0.3 & 0.4 & 0.2 & 0.3 \\ 0.2 & 0.3 & 0.4 & 0.2 \end{bmatrix} = (a_{ij}).$$

None of the elements in A i.e. ($a_{ij}$) is 0.1 and

$$B = \begin{bmatrix} 0.3 & 0.7 & 0.6 & 0.6 \\ 0.6 & 0.3 & 0.5 & 0.5 \\ 0.6 & 0.5 & 0.3 & 0.6 \\ 0.6 & 0.5 & 0.7 & 0.3 \end{bmatrix} = (b_{ij}),$$

we see none of the $b_{ij}$ is.8. Thus the fuzzy interval matrix (A, B) is a open interval matrix defined on (0.1, 0.8).

Clearly [A, B] is a closed interval matrix on the interval [0.2, 0.7].

Thus a open interval fuzzy matrix is open only for some intervals and for some intervals it can be closed. For if we take in A = ($a_{ij}$), a fuzzy matrix the min ($a_{ij}$) = a say, and in B = ($b_{ij}$), a fuzzy matrix if max ($b_{ij}$) = b then the fuzzy interval matrix [A, B] is a closed interval matrix on the interval [a, b].

On similar lines we can define the notion of half open-closed fuzzy interval matrix and half closed open fuzzy interval matrix defined on (0, 1] (or on the interval (a, b] where $0 \leq a < b \leq 1$). What are the possible ways of defining operations on the fuzzy interval matrices defined on the fuzzy interval [0, 1] (or [a, b]). Let [A B] be the fuzzy interval matrix defined on the interval [0, 1] (or on [a, b], $0 \leq a < b \leq 1$].

Let C = ($c_{ij}$) and D = ($d_{ij}$) be two fuzzy matrices in the fuzzy interval matrices [A, B] i.e., $a_{ij} \leq c_{ij} \leq b_{ij}$ where A = ($a_{ij}$) and B = ($b_{ij}$) and $a_{ij} \leq d_{ij} \leq b_{ij}$. Now C + D, the sum of the matrices defined as

$$\left(\frac{c_{ij} + d_{ij}}{2}\right).$$

Clearly C + D $\in$ [A, B] i.e., to the fuzzy interval matrix. We do not say this is the only way by which addition or any operation is to be defined on the fuzzy matrices of the fuzzy interval matrices [A, B].

Likewise, for product we can define the max min rule i.e., D.C = R = ($r_{ij}$) only when the fuzzy interval matrix is a square n×n fuzzy matrices Clearly DC = R is in the interval of fuzzy interval matrix [A, B]. R is defined by $r_{ij}$ = max (min ($d_{ik}$, $c_{kj}$)]. We will have compatibility problem, if the matrices are not square matrices. For if [A, B] be a fuzzy interval matrix defined on the interval [0, 1].

Let

$$C = \begin{bmatrix} 0 & 0.1 & 0.3 & 0.1 \\ 0.2 & 0.3 & 0 & 0.2 \\ 0 & 0.4 & 0.3 & 0.5 \\ 0.2 & 0.4 & 0.6 & 1 \end{bmatrix}$$

and

$$D = \begin{bmatrix} 0.3 & 0.2 & 0.5 & 0.1 \\ 0 & 0.5 & 0.2 & 0.3 \\ 1 & 0 & 0.7 & 1 \\ 0 & 0.6 & 0.9 & 0.2 \end{bmatrix}$$

C.D = R = $(r_{ij})$ where $r_{ij}$ = max [min $(c_{ik}, d_{kj})$].

$r_{11}$ = max [min (0, 0.3), min (0.1, 0), min (0.3, 1) min (0.1, 0)]
= max [0, 0, 0.3 0] = 0.3.

| | | | | | | | | |
|---|---|---|---|---|---|---|---|---|
| $r_{12}$ | = | 0.1, | $r_{13}$ | = | 0.3, | $r_{14}$ | = | 0.3, |
| $r_{21}$ | = | 0.2, | $r_{22}$ | = | 0.3, | $r_{23}$ | = | 0.2, |
| $r_{24}$ | = | 0.3, | $r_{31}$ | = | 0.3, | $r_{32}$ | = | 0.5, |
| $r_{33}$ | = | 0.5, | $r_{34}$ | = | 0.3, | $r_{41}$ | = | 0.6, |
| $r_{42}$ | = | 0.6, | $r_{43}$ | = | 0.4, | $r_{44}$ | = | 0.6. |

$$\text{Thus } R = (r_{ij}) = \begin{bmatrix} 0.3 & 0.1 & 0.3 & 0.3 \\ 0.2 & 0.3 & 0.2 & 0.3 \\ 0.3 & 0.5 & 0.5 & 0.3 \\ 0.6 & 0.6 & 0.9 & 0.6 \end{bmatrix},$$

and R $\in$ [A, B]. Here also it has become pertinent to mention that this operation product defined by us is not the only operation, and interested researcher can define the product in the interval matrices in a different way. Here we see under the operations which we have defined the interval fuzzy matrices need to be only fuzzy square matrices. Suppose we have two fuzzy interval matrices say [A, B] and [$A_1$, $B_1$] defined on two intervals or on the same fuzzy interval [0, 1]. Suppose we take a fuzzy matrix C $\in$ [A, B] and $C_1$ in [$A_1$, $B_1$] and if the fuzzy interval matrices [A, B] are $m \times n$ matrices and that of the fuzzy matrices in the fuzzy interval matrix [$A_1$, $B_1$] are $p \times q$ matrices. If $n \neq p$ we cannot define any operation relating C and $C_1$. The only known operation relating C and $C_1$ is max min operation provided n = p. Otherwise to the best of our knowledge we do not know any other compatible operations.

For if C = $(C_{ij})$ is a m× n fuzzy matrix in [A, B] and $C_1 = (C^1 ij)$ is n × q fuzzy matrix in the interval of fuzzy matrices $[A_1, B_1]$,

$$C \bullet C_1 = \max_K [\min (C_{1k}, C^1_{kj}) = S = (s_{ij}).$$

Clearly S is not in the fuzzy interval matrix [A, B] or in the fuzzy interval matrix $[A_1, B_1]$.

Before we proceed onto define other concepts we just illustrate this product by a concrete example and derive the consequences due to the newly defined product.

***Example 2.1.5:*** Let [A, B] be a 5 × 3 fuzzy interval matrix defined on the interval [0, 1] and $[A_1, B_1]$ be a 3 × 4 fuzzy interval matrix defined on the interval [0, 1]. Let C ∈ [A, B] and $C_1 \in [A_1, B_1]$ where C is a 5 × 3 fuzzy matrix and $C_1$ is a 3 × 4 fuzzy matrix with entries from the fuzzy interval [0, 1].

$$C = \begin{bmatrix} 0 & 0.1 & 0.2 \\ 0.5 & 0 & 0.6 \\ 0.2 & 0.5 & 0.4 \\ 0.3 & 0.1 & 0 \\ 1 & 0 & 0.2 \end{bmatrix}_{5\times 3} \in [A, B]$$

$$C_1 = \begin{bmatrix} 0.2 & 0 & 1 & 0.3 \\ 0 & 0.1 & 0.4 & 1 \\ 0.5 & 0 & 0.3 & 0.2 \end{bmatrix}_{3\times 4} \in [A_1, B_1]$$

$$C \bullet C_1 = (S_{ij}) = \max_K [\min (C_{1k}, C'_{kj}) = S$$

$$S = \begin{bmatrix} 0.2 & 0.1 & 0.2 & 0.2 \\ 0.5 & 0 & 0.5 & 0.3 \\ 0.4 & 0.1 & 0.4 & 0.5 \\ 0.2 & 0.1 & 0.3 & 0.3 \\ 0.2 & 0 & 1 & 0.3 \end{bmatrix}_{5\times 4}.$$

Clearly S $\not\subseteq$ [A, B] and S $\not\subseteq$ [$A_1$, $B_1$]. Thus we see the compatibility cannot be obtained hence we are forced to get a new fuzzy interval matrix which are 5 × 4 fuzzy matrices. Thus we can say this method of product on different fuzzy interval matrices can pave way for another new fuzzy interval matrices. Thus if $[A, B]_{m\times n} \bullet [A_1, B_1]_{n\times p}$ then $[A_2, B_2]_{m\times p}$, of course defined on the same fuzzy interval [0, 1] (or [a, b]; $0 \leq a < b \leq 1$ i.e. [a, b] $\subseteq$ [0, 1]). This type of operations on interval of matrices will certainly find very nice applications. One such is in the interval fuzzy relational equations dealt in chapter 3 of this book.

Having defined fuzzy interval matrices now we proceed on to define the notion of fuzzy interval bimatrices. Before we proceed to define fuzzy interval bimatrices we are forced to define the notion of just interval bimatrices and interval n-matrices, n a positive integer greater than or equal to 2.

## 2.2 Interval Bimatrices and their Generalizations

We have just recalled the definition of interval matrices in the earlier chapter. Now in this section we proceed onto define the notion of interval bimatrices and interval m-matrices m $\geq$ 2, m an integer when m = 2 we get the interval bimatrices. We illustrate these concepts with examples.

**DEFINITION 2.2.1:** *Let [A, B] = [$A_1$, $B_1$] $\cup$ [$A_2$, $B_2$] where [$A_1$ $B_1$] and [$A_2$ $B_2$] interval matrices defined on the same interval [a, b] or on different intervals [$a_1$, $b_1$] and [$a_2$, $b_2$]. We call [A, B] = [$A_1$, $B_1$] $\cup$ [$A_2$, $B_2$] to be the interval bimatrix defined on the same interval [a, b] or on different intervals [$a_1$, $b_1$] $\cup$ [$a_2$, $b_2$]; with no mathematical meaning attached to the symbol '$\cup$'; we say [A, B] is the interval bimatrix defined on the bi-interval [$a_1$, $b_1$] $\cup$[$a_2$, $b_2$] or on the bi-interval [a, b] $\cup$ [a, b].*

*Now we have several other factors to be observed. First if both the interval matrices are m × n matrices then we say [A, B] = [$A_1$, $B_1$] $\cup$ [$A_2$, $B_2$] is a m × n interval bimatrix; if m $\neq$ n we say the interval bimatrix is a rectangular interval bimatrix. If*

*$m = n$, we call $[A, B] = [A_1, B_1] \cup [A_2, B_2]$ to be a square interval bimatrix defined on the bi-interval $[a, b] \cup [a, b]$ or $[a_1, b_1] \cup [a_2, b_2]$.*

*If $[A_1, B_1]$ is a $m \times n$ rectangular interval matrix and $[A_2, B_2]$ is a $p \times q$ rectangular interval matrix then we call $[A, B] = [A_1, B_1] \cup [A_2, B_2]$ to be a mixed rectangular interval bimatrix defined on the bi-interval $[a_1, b_1] \cup [a_2, b_2]$.*

*If the interval bimatrix $[A, B] = [A_1, B_1] \cup [A_2, B_2]$ is such that $[A_1, B_1]$ is a $m \times m$ square interval matrix and $[A_2, B_2]$ is a $p \times p$ ($p \neq m$) square interval matrix then we call $[A, B]$ a mixed interval square bimatrix defined on the bi-interval $[a_1, b_1] \cup [a_2, b_2]$. If one of the interval matrix in a interval bimatrix is a square interval matrix and the other interval matrix $[A_2, B_2]$ is a rectangular interval matrix then we just call $[A, B]$ the mixed interval bimatrix defined on the bi-interval $[a_1, b_1] \cup [a_2, b_2]$.*

If $[A, B] = [A_1, B_1] \cup [A_2, B_2]$ is an interval bimatrix defined on the bi-interval, $[a_1, b_1] \cup [a_2, b_2]$ and if $M = M_1 \cup M_2$ is a bimatrix of the interval bimatrix where $M_1 \in [A_1, B_1]$ and $M_2 \in [A_2, B_2]$ and $N = N_1 \cup N_2$ is a bimatrix of the interval bimatrix $[A, B]$ then the operations $M + N = (M_1 + N_1) \cup (M_2 + N_2)$ and $MN = M_1N_1 \cup M_2N_2$ are performed as in case of interval matrices.

If compatibility of operations is assured in the interval matrices $[A_1, B_1]$ and $[A_2, B_2]$ then certainly the compatibility of operations are guaranteed in case of the interval bimatrices $[A, B] = [A_1, B_1] \cup [A_2, B_2]$.

All properties which exist in case of interval matrices can be easily extended in case of interval bimatrices with appropriate modifications. Further it has become pertinent to mention here that the notion of closed interval bimatrix, half open closed interval bimatrix, half closed open bimatrix or open interval bimatrix or mixed interval bimatrix can be defined as in case of interval matrix.

For proper understanding, we just mention when both the interval matrices of the interval bimatrix $[A, B] = [A_1, B_1] \cup [A_2, B_2]$ where $[A, B]$ is defined on the closed bi-interval $[a_1, b_1] \cup [a_2, b_2]$ then we say the interval bimatrix is defined on the closed bi-interval $[a_1, b_1] \cup [a_2, b_2]$ where if $A_1 = (a^1_{ij})$, $B_1 = (b^1_{ij})$, $A_2 = (a^2_{ij})$ and $B_2 = (b^2_{ij})$ then $a_1 \leq a^1_{ij} \leq b_1$; $a_1 \leq b^1_{ij} \leq b_1$, $a_2 \leq a^2_{ij} \leq b_2$ and $a_2 \leq b^2_{ij} \leq b_2$. Like wise the other concepts are defined.

Having defined these concepts now we proceed on to give examples of them.

***Example 2.2.1:*** Let $[A, B] = [A_1, B_1] \cup [A_2, B_2]$ be a interval rectangular $3 \times 4$ bimatrix defined on the bi-interval $[a_1, b_1] \cup [a_2, b_2]$ where $[a_1, b_1] = [a_2, b_2] = [–R, R]$, R the reals. Let $M \in [A, B]$ then $M = M_1 \cup M_2$, where $M_1 \in [A_1, B_1]$ and $M_2 \in [A_2, B_2]$ with

$$M_1 = \begin{bmatrix} 3 & 2 & 0 & -1 \\ 0 & 1 & 0 & 2 \\ 1 & 0 & +1 & 0 \end{bmatrix} \in [A_1, B_1]$$

and

$$M_2 = \begin{bmatrix} 0 & 1 & 2 & -1 \\ 2 & 0 & 0 & 1 \\ 0 & -1 & 2 & 0 \end{bmatrix} \in [A_2, B_2]$$

$$M = M_1 \cup M_2 \in [A, B].$$

$N = N_1 \cup N_2$ where

$$N_1 = \begin{bmatrix} 0 & 1 & 0 & 1 \\ 1 & 0 & 0 & 1 \\ 0 & 0 & 1 & 1 \end{bmatrix}$$

and

$$N_2 = \begin{bmatrix} 3 & 3 & 2 & -2 \\ 2 & 1 & 0 & 3 \\ 1 & -1 & 3 & 0 \end{bmatrix}$$

$$\begin{aligned} M + N &= (M_1 \cup M_2) + (N_1 \cup N_2) \\ &= (M_1 + N_1) \cup (M_2 + N_2) \end{aligned}$$

$$= \begin{bmatrix} 3 & 3 & 0 & 0 \\ 1 & 1 & 0 & 3 \\ 1 & 0 & 2 & 1 \end{bmatrix} \cup \begin{bmatrix} 2 & 1 & 2 & -1 \\ 2 & -1 & 0 & 1 \\ 0 & -1 & 2 & -3 \end{bmatrix} \in [A, B]$$

'+' is the closed binary operation on the 3 × 4 rectangular interval bimatrix. How ever the usual matrix product is not defined on this particular interval bimatrix.

***Example 2.2.2:*** Let $[A, B] = [A_1, B_1] \cup [A_2, B_2]$ be a square 4 × 4 interval bimatrix defined on the bi-interval $[0, \infty] \cup [0, \infty]$. Let $M = M_1 \cup M_2 \in [A, B]$ where $M_1 \in [A_1, B_1]$ and $M_2 \in [A_2, B_2]$.

$$M_1 = \begin{bmatrix} 2 & 0 & 1 & 5 \\ 0 & 1 & 1 & 2 \\ 3 & 0 & 0 & 4 \\ 1 & 0 & 5 & 0 \end{bmatrix} \in [A_1, B_1]$$

and

$$M_2 = \begin{bmatrix} 0 & 1 & 2 & 0 \\ 1 & 2 & 3 & 4 \\ 0 & 1 & 2 & 0 \\ 1 & 0 & 0 & 2 \end{bmatrix} \in [A_2, B_2].$$

Take $N_1 \cup N_2 \in [A, B]$ then $M_1N_1 \cup M_2N_2 \in [A, B]$. Thus this interval bimatrix is closed with respect to both matrix addition '+' and matrix multiplication ×. But clearly [A, B] under + is not a bigroup only a bimonoid.

Thus the algebraic structure which the interval bimatrix has depends both on the operations and the interval on which it is defined apart from the compatibility of the operations.

Now we give yet another example of a mixed interval bimatrix.

***Example 2.2.3:*** Let $[A, B] = [A_1, B_1] \cup [A_2, B_2]$ be defined on the bi-interval $[-\infty, \infty] \cup [0, \infty]$.

Let $M = M_1 \cup M_2 \in [A, B]$ where

$$M_1 = \begin{bmatrix} -3 & 0 & 1 \\ 0 & -2 & 2 \\ 4 & 1 & 0 \end{bmatrix} \in [A_1, B_1]$$

is the collection of 3 × 3 matrices and $M_2$ is a 4 × 3 matrix where

$$M_2 = \begin{bmatrix} 0 & 1 & 2 \\ 0 & 9 & 4 \\ 0 & 0 & 6 \\ 8 & 0 & 0 \end{bmatrix} \in [A_2, B_2].$$

Thus we see [A, B] is a mixed interval bimatrix defined on the bi-interval $[-\infty\ \infty] \cup [0, \infty] = Z^+ \cup \{0\}$. One can define compatible binary operations on the mixed interval bimatrix [A, B]. We see the interval matrix $[A_1, B_1]$ is a ring on the interval $[-\infty,\infty]$ under matrix addition and matrix multiplication where as $[A_2, B_2]$ is only a monoid under '+' and not compatible with respect to '×'.

Having seen and defined the notion of interval bimatrix now we proceed on to define the generalization of interval bimatrices. First we define a interval trimatrix or a interval 3-matrix and illustrate it with an example.

**DEFINITION 2.2.2:** *Let $[A, B] = [A_1, B_1] \cup [A_2, B_2] \cup [A_3, B_3]$ where $[A_i, B_i]$ is an interval matrix defined on the interval $[a_i, b_i]$, i = 1, 2, 3. The symbol '$\cup$' is just only a notational convenience. Thus if $C^i \in [A_i, B_i]$ then $C^i = (c^i_{pj})$ where i = 1, 2, 3 and $a^i_{pj} \le c^i_{pj} \le b^i_{pj}$, i = 1, 2, 3 and $A_i = (a^i_{pj})$ and $B_i = (b^i_{pj})$, i = 1, 2, 3, [A, B] is defined as the interval trimatrix on the tri-interval $[a_1, b_1] \cup [a_2, b_2] \cup [a_3, b_3]$.*

Any element in the interval trimatrix [A, B] is a trimatrix of the form $C = C^1 \cup C^2 \cup C^3$ where $C^i \in [A_i, B_i]$, i = 1, 2, 3.

For more about trimatrices please refer chapter 1 and [].

Now we can as in case of interval bimatrices which are square interval bimatrix, rectangular interval bimatrix, mixed interval bimatrix and so on we in case of interval trimatrix define the notion of square interval trimatrix by which we have all trimatrices $M = M_1 \cup M_2 \cup M_3$ in [A, B], to be a square $n \times n$ matrix where each $M_i$ is an element of the interval matrix $[A_i, B_i]$, i = 1, 2, 3. We say the interval trimatrix $[A, B] = [A_1, B_1] \cup [A_2, B_2] \cup [A_3, B_3]$ to be a rectangular $m \times n$ matrix if each matrix $M_i$ in the interval matrix $[A_i, B_i]$ is a $m \times n$ matrix for i = 1, 2, 3. The interval of the trimatrix can be defined on the tri-interval $[a_1, b_1] \cup [a_2, b_2] \cup [a_3, b_3]$. Like wise we say a trimatrix is a mixed square trimatrix if each matrix in the interval matrix $[A_i, B_i]$ is a $n_i \times n_i$ square matrix i = 1, 2, 3. Clearly $n_1 \neq n_2$ or $n_1 \neq n_3$ Like wise a mixed rectangular interval of trimatrix and just mixed interval of trimatrix are defined.

Now we proceed on to give some examples of trimatrices.

***Example 2.2.4:*** Let $[A, B] = [A_1, B_1] \cup [A_2, B_2] \cup [A_3, B_3]$ where each $[A_i, B_i]$ is a $3 \times 3$ square interval matrix defined on $[a_i, b_i]$, i = 1, 2, 3. Let the tri-interval be $[a_1, b_1] \cup [a_2, b_2] \cup [a_3, b_3] = [0, \infty] \cup [-R, R] \cup [-Z, Z]$. If $M \in [A, B]$ then $M = M_1 \cup M_2 \cup M_3$ where

$$M_1 = \begin{bmatrix} 0 & 1 & 2 \\ 3 & 4 & 5 \\ 6 & 0 & 1 \end{bmatrix} \in [A_1, B_1]$$

$$M_2 = \begin{bmatrix} -3/7 & 0 & 7 \\ 0 & 0.8 & -1 \\ 1/5 & -1 & 0 \end{bmatrix} \in [A_2, B_2]$$

and

$$M_3 = \begin{bmatrix} -4 & 6 & 23 \\ 1 & 0 & 5 \\ 3 & -7 & 0 \end{bmatrix} \in [A_3, B_3].$$

$M = M_1 \cup M_2 \cup M_3$ is a trimatrix of the interval of trimatrix [A, B] defined on the tri-interval. This interval of trimatrix is called as the square $3 \times 3$ interval of trimatrix.

We give yet another example of mixed rectangular trimatrix.

***Example 2.2.5:*** Let $[A, B] = [A_1, B_1] \cup [A_2, B_2] \cup [A_3, B_3]$ be a mixed rectangular interval trimatrix where the interval matrix $[A_1, B_1]$ contains all $3 \times 2$ matrices with entries on the interval $[0, \propto] = R^+$, $[A_2, B_2]$ is the $2 \times 5$ rectangular interval matrix with entries from $Z = [-\propto, \propto]$, the set of integers. $[A_3, B_3]$ is a $6 \times 1$ rectangular interval matrix with entries from the negative rationals including 0 denoted by $Q^- \cup \{0\}$.

Thus if $M = M_1 \cup M_2 \cup M_3 \in [A, B]$ the mixed rectangular interval trimatrix then

$$M = \begin{bmatrix} \sqrt{3} & 9 \\ 5 & 0 \\ 3 & \sqrt{19} \end{bmatrix} \cup \begin{bmatrix} 0 & 2 & 5 & -3 & 1 \\ 4 & -1 & 0 & -2 & 0 \end{bmatrix} \cup \begin{bmatrix} -3 \\ -7/2 \\ 0 \\ -5 \\ -6/7 \\ 5/3 \end{bmatrix}.$$

Thus [A, B] is a mixed rectangular interval trimatrix defined on the tri-interval $R^+ \cup Z \cup Q^- \cup \{0\}$.

Now we proceed on to generalize this notion to the case of interval n-matrix $n \geq 2$, for when $n = 2$, we get the interval bimatrix and when $n = 3$ we get the interval trimatrix.

**DEFINITION 2.2.3:** *Let [A, B] = $[A_1, B_1] \cup [A_2, B_2] \cup \ldots \cup [A_n, B_n]$ where each $[A_i, B_i]$ is a interval matrix defined on the interval $[a_i, b_i]$, i = 1, 2, …, n. Thus any element in the set [A, B] is a n- matrix given by $M = M_1 \cup M_2 \cup \ldots \cup M_n$ where each $M_i$ is the matrix of the interval matrix $[A_i, B_i]$, i = 1, 2, …, n.[A, B] is defined to be the interval n matrix defined on the n-interval $[a_1, b_1] \cup [a_2, b_2] \cup \ldots \cup [a_n, b_n]$.*

All properties defined in the case of interval bimatrix, and or interval trimatrix can be easily defined and extended in the case of interval n matrix ($n \geq 2$).

Now we just illustrate by an example the interval 6-matrix, which is a mixed matrix.

***Example 2.2.6:*** Let [A, B] = $[A_1, B_1] \cup [A_2, B_2] \cup \ldots \cup [A_6, B_6]$ where each $[A_i, B_i]$ is a interval matrix for i = 1, 2, …, 6. $[A_1, B_1]$ is a $3 \times 3$ interval matrix defined on the set of positive integers; $[A_2, B_2]$ is a $4 \times 2$ interval matrix defined on the set of rationals; $[A_3, B_3]$ is a $5 \times 3$ interval matrix with entries from positive reals, $[A_4, B_4]$ is a $2 \times 5$ matrix with entries from $Z_{10}$ (set of integers modulo 10), $[A_5, B_5]$ is a interval $7 \times 1$ column matrix with entries from R (set of reals) and $[A_6, B_6]$ is a interval $2 \times 2$ square matrix with entries from Q.

Any M in [A, B] = $[A_1, B_1] \cup \ldots \cup [A_6, B_6]$ will be of the form; $M = M_1 \cup M_2 \cup \ldots \cup M_6$ where each $M_i$ is from the interval matrix $[A_i, B_i]$, i = 1, 2, …, 6.

$$M = \begin{bmatrix} 0 & 2 & 3 \\ 1 & 0 & 1 \\ 5 & 4 & 0 \end{bmatrix} \cup \begin{bmatrix} 3 & 0 \\ -7/2 & 6 \\ 2/9 & 1 \\ 5 & -6 \end{bmatrix} \cup$$

$$\begin{bmatrix} 6 & 0.3 & 0.4 \\ 1 & 0 & 2 \\ 0.1 & 5 & \sqrt{7} \\ \sqrt{6} & 2 & 0 \\ 1 & 0 & \sqrt{6} \\ \sqrt{5} & \sqrt{2} & \sqrt{2} \end{bmatrix} \cup \begin{bmatrix} 3 & 2 & 0 & 1 & 5 \\ 6 & 8 & 9 & 0 & 7 \end{bmatrix} \cup$$

$$\begin{bmatrix} 8 \\ \sqrt{3} \\ 5 \\ 0 \\ 1 \\ 2 \\ \sqrt{7} \end{bmatrix} \cup \begin{bmatrix} 7/2 & 0 \\ 9/5 & 5 \end{bmatrix}.$$

Thus [A, B] is a interval 6-matrix defined on the six intervals $(Z^+ \cup \{0\}) \cup Q \cup (R^+ \cup \{0\}) \cup Z_{10} \cup R \cup Q$. All operations which are compatible on [A, B] as a whole can be defined.

Now we just say $[A, B] = [A_1, B_1] \cup [A_2, B_2] \cup \ldots \cup [A_n, B_n]$, a mixed square interval n-matrix, if every n-matrix in the interval matrices $[A_i, B_i]$ are square interval matrices, for i = 1, 2, …, n. Here $[A, B] = [A_1, B_1] \cup \ldots \cup [A_n, B_n]$, where each of the interval matrices are square matrices of order $i \times i$. $i \in Z^+$.

Similarly a mixed rectangular interval n-matrix $[A, B] = [A_1, B_1] \cup [A_2, B_2] \cup \ldots \cup [A_n, B_n]$ where each $[A_i, B_i]$ is a rectangular interval matrix; i = 1, 2, …, n. $[A, B] = [A_1, B_1] \cup \ldots \cup [A_n, B_n]$ is called rectangular $p \times m$ interval n-matrix, if each interval matrix $[A_i, B_i]$ is a $p \times m$ interval matrix; i = 1, 2, …, n. $[A, B] = [A_1, B_1] \cup \ldots \cup [A_n, B_n]$ is called the rectangular $p \times m$ interval n-matrix if each interval matrix $[A_i, B_i]$ is a $p \times m$ interval matrix, i = 1, 2, …, n.

Now having defined the notion of interval bimatrices, interval trimatrices and interval n-matrices in general; we now proceed onto define fuzzy interval bimatrices and fuzzy interval n-matrices. As even the very notion of fuzzy interval matrices have not been defined and defined only now in this book for the first time we here define fuzzy interval n-matrices; $n \geq 2$ and illustrate them with examples.

**DEFINITION 2.2.4:** *Let [A, B] be an interval matrix with entries from [0, 1] or [a, b] with $0 \leq a < b \leq 1$. [A, B] is called the fuzzy interval matrix defined on the interval [a, b] or [0, 1]. If [A, B] contains only $n \times n$ fuzzy matrices then we call [A, B] the fuzzy interval of fuzzy rectangular matrices.*

We just illustrate by an example.

***Example 2.2.7:*** Let [A, B] be a interval of 2 × 2 fuzzy square matrices where

$$A = \begin{bmatrix} 0 & 0.2 \\ 0.3 & 0.5 \end{bmatrix}$$

and

$$B = \begin{bmatrix} 1 & 0.8 \\ 0.9 & 1 \end{bmatrix},$$

[A, B] is the fuzzy interval 2 × 2 square matrix.

***Example 2.2.8:*** Let [A, B] be an interval of 3 × 1 fuzzy rectangular matrix defined on the fuzzy interval [0, 0.6] where

$$A = \begin{pmatrix} 0 \\ 0.1 \\ 0.2 \end{pmatrix}$$

and

$$B = \begin{pmatrix} 0.5 \\ 0.6 \\ 0.4 \end{pmatrix}.$$

.

**DEFINITION 2.2.5:** *Let $[A_1, B_1]$ and $[A_2, B_2]$ be any two interval fuzzy matrices defined on the fuzzy intervals say $[a_1, b_1]$ and $[a_2, b_2]$ respectively. Let $[A, B] = [A_1, B_1] \cup [A_2, B_2]$ where '$\cup$' is just a symbol used for notational convenience, then [A, B] contains all fuzzy bimatrices of the form $M_1 \cup M_2$ where $M_1$ is the fuzzy matrix from the fuzzy interval matrix $[A_1, B_1]$ defined on the interval $[a_1, b_1]$ and $M_2$ is the fuzzy matrix from the fuzzy interval matrix $[A_2, B_2]$ defined on the interval $[a_2, b_2]$, $M = M_1 \cup M_2$ is a fuzzy bimatrix defined on the fuzzy bi-interval $[a_1, b_1] \cup [a_2, b_2]$ where M is called the fuzzy bimatrix of the fuzzy interval bimatrix [A, B]. ($[a_1, b_1] \subseteq [0, 1]$ and $[a_2, b_2] \subseteq [0, 1]$)*

Now we illustrate a fuzzy interval bimatrix by some examples before we proceed onto define some more properties about interval fuzzy bimatrices on the fuzzy bi-interval $[a_1, b_1] \cup [a_2, b_2]$.

***Example 2.2.9:*** Let $[A, B] = [A_1, B_1] \cup [A_2, B_2]$ be a interval fuzzy bimatrices defined on the fuzzy bi-interval $[0.2, 1] \cup [0, 0.7]$ where $[A_1, B_1]$ contains all $3 \times 3$ interval fuzzy matrices defined on the fuzzy intervals $[0.2, 1] \subset [0, 1]$ and $[A_2, B_2]$ contains all $3 \times 3$ interval fuzzy matrices defined on the fuzzy interval $[0, 0.7] \subseteq [0, 1]$.

We call [A, B] the interval fuzzy square $3 \times 3$ bimatrix defined on the bi-interval $[0.2, 1] \cup [0, 0.7]$. Any element M in [A, B] will be of the form $M = M_1 \cup M_2$ where

$$M_1 = \begin{bmatrix} 0.2 & 0.7 & 1 \\ 0.5 & 0.2 & 0.7 \\ 1 & 0.3 & 0.4 \end{bmatrix} \text{ and } \begin{bmatrix} 0 & 0.7 & 0.2 \\ 0.6 & 0 & 0 \\ 0.5 & 0.6 & 0 \end{bmatrix} = M_2.$$

Now on this fuzzy interval of bimatrices [A, B] we can perform both matrix addition defined by max function and matrix multiplication defined by max min function.

Now we give yet another example.

***Example 2.2.10:*** Let [A, B] = $[A_1, B_1] \cup [A_2, B_2]$ be a interval fuzzy bimatrix defined on the fuzzy bi-interval $[0, 0.5] \cup [0.3, 1]$. $[A_1, B_1]$ is the set of all m × n interval fuzzy matrices defined on the fuzzy interval [0, 0.5] and $[A_2, B_2]$ is the set of all m × n fuzzy interval matrices defined on the fuzzy [0.3, 1]. We call [A, B] the m × n rectangular fuzzy interval bimatrix defined on the bi-interval $[0, 0.5] \cup [0.3, 1]$. If we take m = 2 and n = 4, we would get the set of all 2 × 4 rectangular fuzzy interval bimatrix defined on the bi-interval $[0, 0.5] \cup [0.3, 1]$.

Just any element $M = M_1 \cup M_2$ in [A, B] will be of the form

$$M = \begin{bmatrix} 0 & 0.4 & 0.3 & 0 \\ 0.2 & 0 & 0.1 & 0.5 \end{bmatrix} \cup \begin{bmatrix} 0.3 & 1 & 1 & 0.8 \\ 1 & 0.8 & 0.7 & 0.5 \end{bmatrix}.$$

We give yet another example of a fuzzy interval bimatrix defined on the fuzzy interval.

***Example 2.2.11:*** Let [A, B] = $[A_1, B_1] \cup [A_2, B_2]$ be any fuzzy interval bimatrix defined on the fuzzy bi-interval $[a_1, b_1] \cup [a_2, b_2]$. To be more specific let $[A_1, B_1]$ contain all fuzzy interval 2 × 2 square fuzzy matrices with entries from the fuzzy interval [0, 0.7] and $[A_2, B_2]$ contains all fuzzy interval 5 × 5 square fuzzy matrices with entries from the fuzzy interval [0.5, 1]. Thus any element in the fuzzy interval bimatrix [A, B] will be of the form $M = M_1 \cup M_2$ where

$$M_1 = \begin{bmatrix} 0 & 0.2 \\ 0.7 & 0.5 \end{bmatrix} \text{ and } \begin{bmatrix} 0.5 & 0.6 & 1 & 1 & 0.7 \\ 0.6 & 0.6 & 0.7 & 1 & 1 \\ 1 & 1 & 1 & 0.6 & 0.8 \\ 0.9 & 1 & 1 & 0.8 & 1 \\ 0.8 & 0.6 & 0.7 & 1 & 0.9 \end{bmatrix} = M_2.$$

We call this situation in which the interval matrices are fuzzy square matrices, as mixed square fuzzy interval bimatrices defined on the fuzzy bi-interval $[0, 7] \cup [0.5, 1]$.

Next we proceed onto give another example of a fuzzy interval bimatrix.

***Example 2.2.12:*** Let [A, B] = $[A_1, B_1] \cup [A_2, B_2]$ be a fuzzy interval bimatrix defined on the fuzzy bi-interval $[a_1, b_1] \cup [a_2, b_2]$, where $[A_1, B_1]$ is the fuzzy interval rectangular $3 \times 4$ matrix defined on the fuzzy interval [0, 0.6] and $[A_2, B_2]$ is the fuzzy interval rectangular $6 \times 2$ matrix defined on the fuzzy interval [0.1, 1].

We call the fuzzy interval bimatrix [A, B] to be the fuzzy interval mixed rectangular bimatrix defined on the bi-interval $[0, 0.6] \cup [0.1, 1]$. Any element M in [A, B] will be of the form

$$M = \begin{bmatrix} 0 & 0.1 & 0.2 & 0 \\ 0.1 & 0 & 0.4 & 0.52 \\ 0.6 & 0.3 & 0 & 0.51 \end{bmatrix} \text{ and } \begin{bmatrix} 1 & 0.2 \\ 0.3 & 1 \\ 1 & 1 \\ 0.3 & 0.9 \\ 0.6 & 0.2 \\ 0.8 & 1 \end{bmatrix} = M_1 \cup M_2,$$

both the fuzzy interval matrices are rectangular fuzzy matrices.

Next we illustrate the final example of the fuzzy interval mixed bimatrix.

***Example 2.2.13:*** Let [A, B] = $[A_1, B_1] \cup [A_2, B_2]$ be a fuzzy interval mixed bimatrix defined on the fuzzy bi-interval $[a_1, b_1] \cup [a_2, b_2]$. To be more specific where $[A_1, B_1]$ contains the collection of $4 \times 4$ fuzzy interval matrix with entries from [0.6, 1] and $[A_2, B_2]$ contains all $1 \times 5$ fuzzy interval matrix with entries from [0, 0.7]. Any element M in the fuzzy interval bimatrix which we defined as the mixed fuzzy interval bimatrix will be of the form $M = M_1 \cup M_2$ where

$$M_1 \cup M_2 = \begin{bmatrix} 0.6 & 1 & 1 & 1 \\ 1 & 0.5 & 0.6 & 0.7 \\ 0.8 & 0.8 & 1 & 0.6 \\ 0.9 & 1 & 0.9 & 1 \end{bmatrix} \cup \begin{bmatrix} 0 & 0.1 & 0.7 & 0 & 0.5 \end{bmatrix}.$$

Thus by these examples we have illustrated and defined 5 types of fuzzy interval bimatrices, namely.

1. Fuzzy interval square bimatrices where both the fuzzy interval matrices will contain only $m \times m$ fuzzy square matrices.
2. Fuzzy interval rectangular bimatrices where both the fuzzy interval matrices will contain $p \times q$ fuzzy rectangular matrices.
3. Fuzzy interval mixed square bimatrices where one of the fuzzy interval matrix will be a $m \times m$ fuzzy square matrix where as the other will be a $t \times t$ fuzzy square matrix $t \neq m$.
4. Fuzzy interval mixed rectangular bimatrix will contain fuzzy interval $m \times n$ rectangular matrices and $p \times q$ fuzzy interval rectangular matrices $p \neq m$ and or $q \neq n$.
5. The fuzzy interval mixed bimatrices will contain both fuzzy interval square matrices and fuzzy interval rectangular matrices.

Thus when we speak of bimatrix they will fall under only these five categories. Further these will find their applications in Fuzzy Cognitive Bimaps (FCBMs) and Fuzzy Relational Bimaps (FRBMs).

Now we proceed onto define the notion of fuzzy interval trimatrices on fuzzy tri-intervals.

**DEFINITION 2.2.6:** *$[A, B] = [A_1, B_1] \cup [A_2, B_2] \cup [A_3, B_3]$ denotes the fuzzy interval trimatrices defined on the fuzzy tri-interval $[a_1, b_1] \cup [a_2, b_2] \cup [a_3, b_3] \subseteq [0, 1] \cup [0, 1] \cup [0, 1]$. Here each of the $[A_i, B_i]$ are fuzzy interval matrices defined on the fuzzy interval $[a_i, b_i]$; $i = 1, 2, 3$.*

*Thus any element M in [A, B] will be of the form $M = M_1 \cup M_2 \cup M_3$ where M is a fuzzy trimatrix, each $M_i$ is a fuzzy matrix from the fuzzy interval matrix $[A_i, B_i]$; (i =1, 2, 3) defined on the fuzzy intervals $[a_i, b_i]$.*

We now illustrate this by the following example:

***Example 2.2.14:*** Let $[A, B] = [A_1, B_1] \cup [A_2, B_2] \cup [A_3, B_3]$ be a fuzzy interval trimatrix defined on the tri-intervals $[a_1, b_1] \cup [a_2, b_2] \cup [a_3, b_3]$; where each of $[A_i, B_i]$ is a interval fuzzy square matrix or interval fuzzy rectangular matrix.

We denote any element M of [A, B] by $M = M_1 \cup M_2 \cup M_3$ where $M_i$ is a fuzzy matrix of the fuzzy interval matrix $[A_i, B_i]$, i = 1, 2, 3.

***Example 2.2.15:*** Let $[A, B] = [A_1, B_1] \cup [A_2, B_2] \cup [A_3, B_3]$ be a fuzzy interval trimatrix defined on the fuzzy intervals [0, 0.8], [0.2, 0.9] and [0.4, 1] respectively.

Let the interval fuzzy matrix $[A_1, B_1]$ contain all $2 \times 2$ fuzzy matrix taking entries from [0, 0.8]. The fuzzy interval matrix $[A_2, B_2]$ contains all $1 \times 4$ rectangular fuzzy matrix with entries from [0.2, 1] and the fuzzy interval matrix $[A_3, B_3]$ contains all $5 \times 3$ fuzzy matrices taking entries from the interval [0.4, 1]. If M is any element of [A, B] then $M = M_1 \cup M_2 \cup M_2$ where

$$M_1 = \begin{bmatrix} 0 & 0.6 \\ 0.7 & 0.2 \end{bmatrix} \in [A_1, B_1]$$

$$M_2 = (0.2, 0.7, 0.9, 0.6) \in [A_2, B_2]$$

and

$$M_3 = \begin{bmatrix} 0.4 & 1 & 1 \\ 1 & 0.6 & 0.7 \\ 1 & 0.6 & 0.9 \\ 0.6 & 0.8 & 0.7 \\ 0.8 & 0.9 & 1 \end{bmatrix} \in [A_3, B_3].$$

Now we proceed onto define and generalize this concept to fuzzy interval n-matrix.

**DEFINITION 2.2.7:** *Let [A, B] = [$A_1$, $B_1$] ∪ [$A_2$, $B_2$] ∪ … ∪ [$A_n$, $B_n$]; (where n >3 where each [$A_i$, $B_i$] is a fuzzy interval square matrix or a fuzzy interval rectangular matrix defined on the fuzzy interval [$a_i$, $b_i$] ⊆ [0, 1]; i =1, 2, …, n. We call [A, B] the fuzzy interval n-matrix.*

*Any element of the fuzzy interval n-matrix is a fuzzy n-matrix given by M = $M_1$ ∪ $M_2$ ∪ … ∪ $M_n$ where each $M_i$ is a fuzzy matrix of the fuzzy interval matrix [$A_i$, $B_i$] defined on the fuzzy interval [$a_i$, $b_i$]; i = 1, 2, …, n.*

We just illustrate this by the following example:

***Example 2.2.16:*** Let [A, B] = [$A_1$, $B_1$] ∪ … ∪ [$A_8$, $B_8$] be a fuzzy interval 8-matrix defined on the 8-fuzzy interval [$a_1$, $b_1$] ∪ … ∪ [$a_8$, $b_8$], where [$A_i$, $B_i$] is a fuzzy interval matrix defined on the interval [$a_i$, $b_i$]; i = 1, 2, …, 8.

Now if every fuzzy interval matrix [$A_i$, $B_i$] of the fuzzy interval n-matrix [A, B] = [$A_1$, $B_1$] ∪ … ∪ [$A_n$ $B_n$], is a m × m square fuzzy interval matrix then we call [A, B] to be a square fuzzy interval n-matrix.

If every fuzzy interval matrix [$A_i$, $B_i$] of the fuzzy interval n-matrix is a m × n rectangular fuzzy interval matrix then we call [A, B] the rectangular fuzzy interval n-matrix. If every fuzzy interval matrix [$A_i$, $B_i$] of the fuzzy interval n-matrix is either a square matrix or a rectangular fuzzy interval matrix then we call [A, B] a mixed fuzzy interval n-matrix. If each of fuzzy interval matrix [$A_i$, $B_i$] is a $m_i \times m_i$ square fuzzy matrix where $m_i \neq m_j$ if i ≠ j then we call [A, B] as the mixed square fuzzy interval n-matrix.

If each of the fuzzy interval matrix [$A_i$, $B_i$] is a $m_i \times n_i$, rectangular fuzzy matrix $m_i \neq m_j$ if (i ≠ j), then we call the fuzzy interval n-matrix the mixed rectangular fuzzy interval n-matrix.

These fuzzy interval n-matrices will find its applications when we need to analyze multi expert opinions.

## 2.3 Neutrosophic Interval Matrices and their Generalizations

Now we proceed onto define the notion of neutrosophic interval matrices for the first time. Only the notion of interval matrices have been defined, even the notion of fuzzy interval matrices have been defined only in this book.

First we define neutrosophic interval matrices. We call a m × n matrix with entries from a neutrosophic field to be a neutrosophic matrix. For more about neutrosophic matrices, neutrosophic bimatrices and their generalizations please refer [].

**DEFINITION 2.3.1:** *Let A = ($a_{ij}$) and B = ($b_{ij}$) be two m × n neutrosophic matrices with entries from the neutrosophic rational field ⟨Q ∪ I⟩. [A, B] will be called the neutrosophic interval matrices, where if C = ($c_{ij}$) is any other m × n neutrosophic matrix with ($a_{ij}$) ≤ ($c_{ij}$) ≤ ($b_{ij}$) then we say C = ($c_{ij}$) ∈ [A, B] i.e. the matrix C = ($c_{ij}$) is a matrix of the neutrosophic interval matrix [A, B].*

Now we illustrate this by the following example:

***Example 2.3.1:*** Let [A, B] be a neutrosophic interval 3 × 2 matrices, where

$$A = \begin{bmatrix} 2+5I & 0 \\ 5 & 7-I \\ 8 & 8I-9 \end{bmatrix}$$

and

$$B = \begin{bmatrix} 10I+4 & 12 \\ 10 & 9I+20 \\ 18 & 10I+7 \end{bmatrix}.$$

Now

$$C = \begin{bmatrix} 9I+5 & 8 \\ 4+2I & 6-I \\ 12 & 3I+1 \end{bmatrix}$$

belongs to the neutrosophic interval of matrices [A, B] for we see $(a_{ij}) \leq (c_{ij}) \leq (b_{ij})$; suppose

$$D = \begin{bmatrix} 21I + 40 & -45 \\ 0 & 8I - 25 \\ 70 & 24I + 48 \end{bmatrix}$$

be any $3 \times 2$ neutrosophic matrix. Clearly D does not belong to the interval of neutrosophic matrices, [A, B].

We also define neutrosophic interval of matrices in a different way however both the definitions are equivalent.

**DEFINITION 2.3.2:** *[A, B] is called the neutrosophic interval of matrices if A and B are $m \times n$ matrices with entries from the neutrosophic interval $\{[a + bI, c + dI] / a \in \{-n, p\}, b \in \{-m, t\}, c \in \{-r, s\}$ and $d \in \{-u, v\}$ where n, m, p, r, t, u, v and s are reals and the minimal neutrosophic elements got from these intervals [a + bI, c + dI] form the entries of the neutrosophic matrix A and the maximal entries of these elements got from the interval [a + bI, c + dI] form the entries of the neutrosophic matrix B. Thus if $C = (c_{ij})$ is a any matrix the elements $(c_{ij}) \in [a + bI, c + dI]$ and $(a_{ij}) \leq (c_{ij}) \leq (b_{ij})$.*

We illustrate this by the following simple example.

***Example 2.3.2:*** Let [A, B] be a neutrosophic interval matrix where

$$A = \begin{bmatrix} 2 + 5I & 7 + 10I \\ 2 + 6I & 5 + 11I \end{bmatrix}$$

and

$$B = \begin{bmatrix} 18 + 28I & 40 + 32I \\ 16 + 12I & 32 + 14I \end{bmatrix}.$$

Now we can say the neutrosophic number 2 + 5I is the minimum number and 40 + 32I is the maximum number. Thus all matrices which takes value in the neutrosophic interval [2 + 5I, 40 + 32I] will be in the neutrosophic interval matrix [A, B].

For instance take

$$C = \begin{bmatrix} 3+7I & 40+30I \\ 2+7I & 20+12I \end{bmatrix}.$$

Clearly C $\in$ [A, B].
Suppose take

$$D = \begin{bmatrix} 50+3I & 18+42I \\ 2+7I & 5+11I \end{bmatrix}.$$

Clearly D $\notin$ [A, B] for 50 + 3I $\notin$ [2 + 5I, 40 + 32I] also 18 + 42I $\notin$ [2 + 5I, 40 + 32I] so D is not an element of the neutrosophic interval matrix [A, B].

Now we say the neutrosophic interval matrix is a rectangular neutrosophic interval matrix if A and B are m $\times$ n neutrosophic matrices. We say the neutrosophic interval matrix is a square neutrosophic interval matrix if both A and B are n $\times$ n neutrosophic matrices. The interval of definition is determined by the minimum element in A = ($a_{ij}$) and the maximum element in B = ($b_{ij}$).

Now we proceed to define the operations on the neutrosophic interval of matrices. Usual compatible matrix operations on neutrosophic matrices can be defined as the binary operations provided the resultant neutrosophic matrix falls with in the neutrosophic interval of matrices. In case they do not fall in the neutrosophic interval of matrices then we say the closure property is not satisfied. We are not always guaranteed of the closure axiom in case of all neutrosophic interval of matrices under all matrix operations.

Now we proceed onto define the notion of fuzzy neutrosophic interval of matrices.

**DEFINITION 2.3.3:** *We say a neutrosophic matrix A to be a fuzzy neutrosophic matrix if the entries of A are from [a + bI, c + dI] where a, b, c, d $\in$ [0, 1].*

We just illustrate it by a simple example.

***Example 2.3.3**:* Let

$$A = \begin{bmatrix} 0.2+I & 0.7+1I \\ 1+0.I & 0.2+0.3I \\ 1+I & 0 \end{bmatrix}.$$

Clearly A is a fuzzy neutrosophic matrix defined on the interval [0, 1+I].

We can define square fuzzy neutrosophic matrix, rectangular fuzzy neutrosophic matrix, row fuzzy neutrosophic vector / matrix and column fuzzy neutrosophic vector / matrix.

**DEFINITION 2.3.4:** *[A, B] will be called as the fuzzy neutrosophic interval matrices if each of A and B are m $\times$ n fuzzy neutrosophic matrix and $(a_{ij}) \leq (b_{ij})$ where A = $(a_{ij})$ and B = $(b_{ij})$, $a_{ij}$, $b_{ij}$ $\in$ [a + bI, c + dI] with a, b, c, d $\in$ [0, 1], $1 \leq i \leq m$ and $1 \leq j \leq n$. If m = n then we see the fuzzy neutrosophic interval matrix is a fuzzy neutrosophic interval square matrix. If m $\neq$ n, then we call [A, B] the fuzzy neutrosophic interval rectangular matrix.*

*We say a fuzzy neutrosophic matrix C = $(c_{ij})$ $\in$ [A, B] where A = $(a_{ij})$ and B = $(b_{ij})$ if and only if $a_{ij} \leq c_{ij} \leq b_{ij}$, $1 \leq {}_{i}, {}_{j} \leq n$. If [A, B] is a fuzzy neutrosophic interval of square matrix if A = $(a_{ij})$ and B = $(b_{ij})$ are p $\times$ p square fuzzy neutrosophic matrices.*

We illustrate this by the following example:

***Example 2.3.4:*** Let [A, B] be a fuzzy neutrosophic 3 $\times$ 2 interval matrix where

$$A = \begin{bmatrix} 0.5+I & 0.2I \\ 0.3I & 0.4+0.3I \\ I & 0 \end{bmatrix}$$

and

$$B = \begin{bmatrix} 0.9+I & 0.8I \\ I & 0.7+I \\ 1+I & 0 \end{bmatrix}.$$

The minimal element is 0 and maximal element is 1 + I:

These concepts will find their applications in neutrosophic Cognitive Maps (NCMs) and Neutrosophic Relational Maps (NRMs) models.

Now we proceed on to define neutrosophic interval bimatrix.

**DEFINITION 2.3.5*:*** *Let $[A_1, B_1]$ and $[A_2, B_2]$ be two neutrosophic interval matrices. The neutrosophic interval bimatrix $[A, B] = [A_1, B_1] \cup [A_2, B_2]$ is defined to be the collection of all neutrosophic bimatrices $M = M_1 \cup M_2$ where $M_1$ is the neutrosophic matrix from the neutrosophic interval matrix $[A_1, B_1]$ and $M_2$ is the neutrosophic matrix from the neutrosophic interval matrix $[A_2, B_2]$. If both the matrices $M_1$ and $M_2$ are $n \times n$ square matrices then we call [A, B] to be the neutrosophic interval square bimatrix.*

*If both the neutrosophic interval matrices $[A_1, B_1]$ and $[A_1, B_2]$ are rectangular $n \times m$ neutrosophic interval matrices then we call [A, B] the neutrosophic interval rectangular bimatrix. If one of $[A_1, B_1]$ is a neutrosophic interval square matrix and $[A_2, B_2]$ is a neutrosophic interval rectangular matrix, then $[A, B] = [A_1, B_1] \cup [A_2, B_2]$ will be defined as the neutrosophic interval mixed bimatrix.*

*If one of $[A_1, B_1]$ is a neutrosophic interval of $n \times n$ square matrix and $[A_2, B_2]$ is a neutrosophic interval of $p \times p$ square matrix $p \neq n$ then we call $[A, B] = [A_1, B_1] \cup [A_2, B_2]$ to be a neutrosophic interval of mixed square bimatrices.*

*If one of $[A_1, B_1]$ is a neutrosophic interval of rectangular matrix say $m \times n$ $(m \neq n)$ and $[A_2, B_2]$ is a neutrosophic interval $p \times q$ $(p \neq q, p \neq m)$ rectangular matrix. Then we call $[A, B] = [A_1, B_1] \cup [A_2, B_2]$ to be the neutrosophic interval of mixed rectangular bimatrices.*

Thus we can categorize the neutrosophic interval of bimatrices into these 5 classes. We just give a few examples of neutrosophic interval of bimatrices.

***Example 2.3.5:*** Let $[A, B] = [A_1, B_1] \cup [A_2, B_2]$ where $[A_1, B_1]$ is a $2 \times 2$ neutrosophic interval of square matrices and $[A_2, B_2]$ be the neutrosophic interval of $4 \times 1$ rectangular matrices. Thus $[A, B]$ is a neutrosophic interval of mixed bimatrices. Any element M in $[A, B] = [A_1, B_1] \cup [A_2, B_2]$ will be of the form $M = M_1 \cup M_2$ where

$$M_1 = \begin{bmatrix} 2I & 0 \\ 2+4I & 5-I \end{bmatrix}$$

and

$$M_2 = [20, 25 + 14I, 0, 5 + I].$$

Next we give an example of a neutrosophic interval of mixed rectangular matrices.

***Example 2.3.6:*** Let $[A, B] = [A_1, B_1] \cup [A_2, B_2]$ where $[A_1, B_1]$ is the neutrosophic interval of $3 \times 5$ rectangular matrices and $[A_2, B_2]$ is the neutrosophic interval of $4 \times 2$ rectangular matrices. Then $[A, B]$ is the neutrosophic interval of mixed rectangular bimatrices. Any element M in $[A, B]$ will be a mixed rectangular neutrosophic bimatrix where $M = M_1 \cup M_2$, with

$$M_1 = \begin{bmatrix} 0 & 1+I & 7+I & 2I & 5 \\ 3I & 1-I & 8-5I & 4 & 3-I \\ 2 & 0 & 7I & 3 & 4-I \end{bmatrix} \in [A_1, B_1]$$

and

$$M_2 = \begin{bmatrix} 2I & 0 \\ 3+I & 4 \\ 5I-2 & 7I \\ 20 & 3I+5 \end{bmatrix} \in [A_2, B_2].$$

We give yet another example of a neutrosophic interval of rectangular bimatrix.

***Example 2.3.7:*** Let $[A, B] = [A_1, B_1] \cup [A_2, B_2]$ where $[A_1, B_1]$ is the neutrosophic interval of $2 \times 5$ rectangular matrix and $[A_2, B_2]$ is the neutrosophic interval of $2 \times 5$ matrix, where

$$A_1 = \begin{bmatrix} 3+4I & 0 & 2+I & 2 & 7+4I \\ 5 & 3 & 0 & 3 & 5+2I \end{bmatrix}$$

$$B_1 = \begin{bmatrix} 20+4I & 17 & 20+I & 3+20I & 7+8I \\ 12+7I & 45+24I & 14 & 5 & 8+6I \end{bmatrix}.$$

The real minimum is 0, real maximum is 17, neutrosophic minimum is 2 + I and neutrosophic maximum is 45 + 24I.
Thus

$$C = \begin{bmatrix} 0 & 4+5I & 4+I & 3 & 7+5I \\ 14 & 10+5I & 10 & 4 & 8+9I \end{bmatrix}$$

is an element of the neutrosophic interval matrix $[A_1, B_1]$.
Let

$$A_2 = \begin{bmatrix} 0 & 0.3I & 0.5 & 0.2+0.3I & 0.10I \\ 0.4+I & 0.2I & 0 & 0.2 & 0.3I \end{bmatrix}$$

and

$$B_2 = \begin{bmatrix} 1 & 1+I & 0.8 & 0.5+0.8I & I \\ 1+I & 0.9+I & 1+0.8 & 0.6 & 0.9I \end{bmatrix}.$$

Clearly real minimal element or entry is 0, real maximal entry is 1; neutrosophic maximal entries are 1 + I, I; neutrosophic minimal entries are 0.2I, 0.2 + 0.3I.
Let

$$D = \begin{bmatrix} 0 & 0.9I & I+0.6 & 1+I & 1 \\ 0.8+I & 0.8I & 0.9I & 1+0.5I & 0.6I \end{bmatrix},$$

D is an element of the neutrosophic interval matrix $[A_2, B_2]$ Thus E = C $\cup$ D is an element of the neutrosophic interval bimatrix [A, B].
Take

$$C' = \begin{bmatrix} I & 90 & 45I & 0 & 3+I \\ 29I+70 & 1 & 5I & 100 & 2+90I \end{bmatrix}.$$

The neutrosophic matrix C' does not belong to the neutrosophic interval matrix $[A_1, B_1]$.
Take the neutrosophic matrix

$$D' = \begin{bmatrix} 5I & 9 & 2 & 0 & I \\ 0.7I & 0.6+I & 0.2+3I & 7I & 0 \end{bmatrix}$$

D' is not an element of the neutrosophic interval matrix $[A_2, B_2]$. Thus D' $\cup$ D is not an element of the neutrosophic interval bimatrix $[A, B] = [A_2, B_2] \cup [A_1, B_1]$.

Note in the bimatrix $M = M_1 \cup M_2$ even if one of the neutrosophic matrices $M_1$ or $M_2$ is not in the neutrosophic interval matrix $[A_1, B_1]$ or $[A_2, B_2]$ respectively then M does not belong to the neutrosophic interval bimatrix [A, B].

Now having given the definition of a neutrosophic interval bimatrix now we proceed on to define the notion of neutrosophic interval trimatrix.

The new notion of neutrosophic interval bimatrix will find its application in neutrosophic Relational Bimaps (NRBMs) and Neutrosophic Cognitive Bimaps (NCBMs), which are described in chapter 3 of this book.

**DEFINITION 2.3.6:** *Let [A, B] = [$A_1$, $B_1$] $\cup$ [$A_2$, $B_2$] $\cup$ [$A_3$, $B_3$] be such that each [$A_i$, $B_i$] is a neutrosophic interval matrix, for i = 1, 2, 3. Any neutrosophic trimatrix M = $M_1 \cup M_2 \cup M_3$ where $M_i \in [A_i, B_i]$, i = 1, 2, 3 belongs to [A, B]. Thus [A, B] denotes the collection of all neutrosophic trimatrices M = $M_1 \cup M_2 \cup M_3$, where $M_i \in [A_i, B_i]$, i = 1, 2, 3.*

*This collection [A, B] is called the neutrosophic interval of trimatrices. The neutrosophic interval of trimatrices can be square trimatrices, rectangular trimatrices or mixed square trimatrices or mixed rectangular trimatrices or just mixed trimatrices.*

Thus these also fall under only five categories.

We illustrate a few of them with examples.

***Example 2.3.8:*** Let [A, B] = [$A_1$, $B_1$] $\cup$ [$A_2$, $B_2$] $\cup$ [$A_3$, $B_3$] where [$A_i$, $B_i$] are neutrosophic interval square 2 × 2 matrices for i = 1, 2, 3. Here we give $A_i$ and $B_i$ for i = 1, 2, 3 so that one can easily find out the neutrosophic 2 × 2 square trimatrices belonging to the neutrosophic interval of square matrices.

$$A_1 = \begin{bmatrix} 0.2I & 0 \\ 0.7 & 0.1+I, \end{bmatrix}, B_1 = \begin{bmatrix} I & 1 \\ 1+I & 0.8+I \end{bmatrix}$$

$$A_2 = \begin{bmatrix} 5 & 7I \\ 12I & 0 \end{bmatrix}, B_2 = \begin{bmatrix} 25 & 35I \\ 7+15I & 45 \end{bmatrix}$$

$$A_3 = \begin{bmatrix} 3+I & 2-I \\ I+1 & 5-3I \end{bmatrix} \text{ and } B_3 = \begin{bmatrix} 20+20I & 2+I \\ 15+5I & 10I+10 \end{bmatrix}.$$

Now we know given any matrix $M_i$ in the neutrosophic interval of matrices [$A_i$, $B_i$], i = 1, 2, 3 which of the $M_i$ are in the interval of matrices [$A_i$, $B_i$] and which of the $M_i$ do not belong to the neutrosophic interval of matrices.

Let $M = M_1 \cup M_2 \cup M_3$ if all the $M_i \in [A_i, B_i]$, i = 1, 2, 3 then we declare M is an element of the neutrosophic trimatrix of the neutrosophic interval of trimatrices [A, B].

For take $M = M_1 \cup M_2 \cup M_3$ where

$$M_1 = \begin{bmatrix} 0.8I & I \\ 0.8 & 0.2+0.5I \end{bmatrix},$$

$$M_2 = \begin{bmatrix} 20 & 20I \\ 15I & 10 \end{bmatrix}$$

and

$$M_3 = \begin{bmatrix} 1+I & 2 \\ 5+5I & 2+3I \end{bmatrix}.$$

Clearly $M = M_1 \cup M_2 \cup M_2$ is a neutrosophic trimatrix belonging to the neutrosophic interval of trimatrices.

Take $N = N_1 \cup N_2 \cup N_3$ where

$$N_1 = \begin{bmatrix} 8I & 0 \\ 9 & 0.2I \end{bmatrix},$$

$$N_2 = \begin{bmatrix} -I & 0.5I \\ 0 & 90 \end{bmatrix}$$

and

$$N_3 = \begin{bmatrix} 3+I & 2-I \\ 3I-2 & 2+5I \end{bmatrix}.$$

Clearly N is a neutrosophic trimatrix but N is not an element of the neutrosophic interval of trimatrix [A, B].

Thus we have just seen an example of a neutrosophic interval of square trimatrices.

Now we give yet another example of a neutrosophic interval of mixed trimatrices.

***Example 2.3.9:*** Let [A, B] = $[A_1, B_1] \cup [A_2, B_2] \cup [A_3, B_3]$ be a neutrosophic interval of trimatrices where $[A_1, B_1]$ is a neutrosophic interval of $2 \times 3$ matrices where

$$A_1 = \begin{bmatrix} 2 & 3I & 5 \\ 0 & 2I & I \end{bmatrix}$$

and

$$B_1 = \begin{bmatrix} 9 & 20I & 15 \\ 5I & 12I & 7I \end{bmatrix}.$$

Thus this neutrosophic interval of $2 \times 3$ matrices has the minimal real to be 0, maximal real to be 15, minimal neutrosophic element to be I and maximal neutrosophic element to be 20I so if we take any neutrosophic $2 \times 3$ matrices with any of the entries as reals greater than 15 or negative reals or neutrosophic elements less than I or greater than 20I then the $2 \times 3$ neutrosophic matrix will not belong to the neutrosophic interval matrix $[A_1, B_1]$.

Now $[A_2, B_2]$ is the neutrosophic interval of $5 \times 1$ neutrosophic matrices where

$$A_2 = \begin{bmatrix} 2+I \\ 0 \\ 7-I \\ 5+2I \\ 1+I \end{bmatrix}$$

and

$$B_2 = \begin{bmatrix} 20+4I \\ 11+9I \\ 40+21I \\ 5+2I \\ 15-12I \end{bmatrix}.$$

Now the real minimum is 0, neutrosophic positive minimum is 1 + I, neutrosophic negative minimum is 7 – I. No real entries can find its place in the neutrosophic interval of matrices. The neutrosophic positive maximum is 40 + 21I and neutrosophic negative maximum is 15 – 12I.

Thus if 60 + 17I is an entry in any $5 \times 1$ neutrosophic matrix than that matrix does not belong to the neutrosophic interval of $5 \times 1$ matrices $[A_2, B_2]$. $[A_3, B_3]$ is a neutrosophic interval of $3 \times 3$ square matrices where

$$A_2 = \begin{bmatrix} 0.2I & 0 & 0.2 \\ 0.3I & 0.I & 0 \\ 0 & 0 & 0.I \end{bmatrix}$$

and

$$B_3 = \begin{bmatrix} 0.9I & 0.8 & 0.6 \\ 0.8I & 0.8I & 0.7 \\ 0.4 & 0 & 0.6I \end{bmatrix}.$$

Thus 0 is the real minimum 0.8 is the real maximum 0.I is the neutrosophic minimum and 0.9I is the neutrosophic maximum. If a $3 \times 3$ neutrosophic matrix M has an entry to be 1 + I then M $\not\in [A_3, B_3]$. Thus $[A, B] = [A_1, B_1] \cup [A_2, B_2] \cup [A_3, B_3]$ is the neutrosophic interval of mixed trimatrices.

Now we proceed onto extend this notion to the definition neutrosophic interval of n-matrices n > 3. For when n = 2 we get the class of neutrosophic interval of bimatrices and when n = 3 we get the class of neutrosophic interval of trimatrices.

**DEFINITION 2.3.7:** *Let $[A, B] = [A_1, B_1] \cup [A_2, B_2] \cup \ldots \cup [A_n, B_n]$ where each $[A_i, B_i]$ is a neutrosophic interval matrix, i = 1, 2, …, n. Thus any element M in [A, B] will be a neutrosophic n-matrix, $M = M_1 \cup M_2 \cup \ldots \cup M_n$, where $M_i \in [A_i, B_i]$; i = 1, 2, …, n. Thus $[A, B] = [A_1, B_1] \cup [A_2, B_2] \cup \ldots \cup [A_n, B_n]$ denotes the collection of all neutrosophic n-matrices satisfying the condition each $M_i$ is an element of the neutrosophic interval* *matrix $[A_i, B_i]$, i = 1, 2, …, n and [A, B] is defined to be the neutrosophic interval n-matrix.*

*When all the matrices $M_i$ in $M = M_1 \cup \ldots \cup M_n$, $M \in [A, B]$ is a $m \times m$ neutrosophic square matrix then we call [A, B] to be the neutrosophic interval $m \times m$ square n-matrices. If in the neutrosophic interval of n-matrices $[A, B] = [A_1, B_1] \cup [A_2, B_2] \cup \ldots \cup [A_n, B_n]$ if each of a neutrosophic matrices $M_i$, in each of the neutrosophic interval of matrices $[A_i, B_i]$ is a $m \times p$ ($m \neq p$) rectangular matrices then we call [A, B] to be the neutrosophic interval of $m \times p$ rectangular n-matrices. If in the neutrosophic interval of n-matrices $[A, B] = [A_1, B_1] \cup [A_2, B_2] \cup \ldots \cup [A_n, B_n]$ each of the neutrosophic matrices $M_i$, in each of the neutrosophic interval of matrices $[A_i, B_i]$ is a square $p_i \times p_i$ matrices; i = 1, 2, …, n, $p_i \neq p_j$ if $i \neq j$; $i \leq i, j \leq n$, then we call [A B] to be the neutrosophic interval of mixed square n-matrices. If in the neutrosophic interval of matrices $[A, B] = [A_1, B_1] \cup [A_2, B_2] \cup \ldots \cup [A_n, B_n]$, each of the neutrosophic matrices $M_i$, in each of the neutrosophic interval of matrices $[A_i, B_i]$, i = 1, 2, …, n are $t_i \times r_i$ rectangular neutrosophic matrices $t_i \neq r_i$ (i = 1, 2, …, n} then we call [A, B] the neutrosophic interval of mixed rectangular n-matrices.*

*Finally in the neutrosophic interval of n-matrices $[A, B] = [A_1, B_1] \cup [A_2, B_2] \cup \ldots \cup [A_n, B_n]$, each of the neutrosophic matrices $M_i$ in each of the neutrosophic interval of matrices $[A_i, B_i]$ is either a square matrix or a rectangular matrix for i = 1, 2, …, n; then we call [A, B] to be the neutrosophic interval of mixed n-matrices.*

Since the very notion of neutrosophic matrices is new and still new is the notion of neutrosophic interval of matrices and still abstract is the notion of neutrosophic interval of n-matrices we have tried to explain the 5 types of neutrosophic interval n-matrices.

Now we proceed onto illustrate them with examples.

***Example 2.3.10:*** Let $[A, B] = [A_1, B_1] \cup [A_2, B_2] \cup \ldots \cup [A_6, B_6]$ where [A, B] is a neutrosophic interval of mixed 6-matrices; $[A_1, B_1]$ is the neutrosophic interval of $2 \times 3$ matrices where

$$A_1 = \begin{bmatrix} I & 0 & 2I \\ 1 & 2I & 0 \end{bmatrix}$$

and

$$B_1 = \begin{bmatrix} 2I & 5I & 10I \\ 3 & 4I & 9 \end{bmatrix}$$

the maximal real value is 9, the minimal real value is 0, the maximal neutrosophic value is 10I and the minimal neutrosophic value is I.

$[A_2\ B_2]$ is the neutrosophic interval of 4 × 1 neutrosophic matrices where

$$A_2 = \begin{bmatrix} 10+I \\ 1+I \\ 2+5I \\ 1+2I \end{bmatrix}$$

and

$$B_2 = \begin{bmatrix} 11+I \\ 13+I \\ 19+20I \\ 4+5I \end{bmatrix},$$

where the minimal neutrosophic element is 1 + I and the maximal neutrosophic element is 19 + 20I.

$[A_3, B_3]$ is the neutrosophic interval of 3 × 1 matrices where

$$A_3 = [0.I, 0, 0.2] \text{ and } B_3 = [I, 1, 0.9],$$

the minimal real element is 0, the maximal real element is 1, the minimal neutrosophic element is 0.I and the maximal neutrosophic element is I.

$[A_4, B_4]$ is the neutrosophic interval of 3 × 3 neutrosophic matrices where

$$A_4 = \begin{bmatrix} 0 & 2 & 5I \\ -5 & 0 & 4 \\ -2 & -I & 1 \end{bmatrix}$$

and

$$B_4 = \begin{bmatrix} 2 & 20 & 15I \\ 0 & 5I & 45 \\ 2+5I & 10I & 2 \end{bmatrix}$$

The minimal neutrosophic element is –I and the maximal neutrosophic element is 15I. The minimal real element is –5 and the maximal real entry is 45.

[$A_5$, $B_5$] is the neutrosophic interval of (4 × 2) neutrosophic matrices where

$$A_5 = \begin{bmatrix} 2+I & 0 \\ 0 & 3+2I \\ 5+5I & I \\ 2 & 11+I \end{bmatrix}$$

and

$$B_5 = \begin{bmatrix} 4+5I & 5+I \\ 2+I & 4+7I \\ 25+16I & 40I \\ 45 & 22+4I \end{bmatrix}$$

The neutrosophic maximal element is 25 + 16I. The minimal neutrosophic element is 2+I (or I). The minimal real is 0. The maximal real is 45.

[$A_6$, $B_6$] is neutrosophic interval of 2 × 2 square matrices, where

$$A_6 = \begin{bmatrix} 0.6 & 0.2+0.I \\ 0 & 0.2I \end{bmatrix}.$$

and

$$B_6 = \begin{bmatrix} 0.6+I & 0.9+I \\ 1+I & 0.9+0.2I \end{bmatrix}$$

The real minimum is 0.
The real maximum is 0.6
The neutrosophic maximum is 1 + I
The neutrosophic minimum is 0.2I and 0.2 + 0.I

Thus $[A, B] = [A_1, B_1] \cup [A_2, B_2] \cup \ldots \cup [A_6, B_6]$ is a neutrosophic interval of mixed 6-matrices.

We have given an example for neutrosophic interval of mixed n-matrices. We will give more example of neutrosophic interval of mixed square matrices.

***Example 2.3.11:*** Let $[A, B] = [A_1, B_1] \cup [A_2, B_2] \cup [A_4, B_4]$ where [A, B] is the neutrosophic interval of mixed square 4-matrices.

$[A_1, B_1]$ is the neutrosophic interval of $2 \times 2$ square matrices where

$$A_1 = \begin{bmatrix} I & 5 \\ 0 & 3I \end{bmatrix}$$

and

$$B_1 = \begin{bmatrix} 12I & 16 \\ 20 & 9I \end{bmatrix}.$$

The minimal real is 0.
The maximal real is 20.
The minimal neutrosophic element is I and
the maximal neutrosophic element is 12I. Clearly

$$C = \begin{bmatrix} 20I & -I \\ 41 & 0 \end{bmatrix}$$

is a $2 \times 2$ neutrosophic square matrix but $C \notin [A_1, B_1]$ for $20\,I \notin [I, 12I]$ and $41 \notin [0, 20]$.

Now $[A_2, B_2]$ is the neutrosophic interval of $4 \times 4$ neutrosophic square matrices where

$$A_2 = \begin{bmatrix} 0 & 0.I & 0.1 & 0.2 \\ 0.2I & 0 & 0 & 0.1 \\ 0.3I & 0.1 & 0 & 0.4 \\ 0 & 0 & 0 & 0.3I \end{bmatrix}$$

and

$$B_2 = \begin{bmatrix} 1 & I & 0.6 & 1 \\ I & 0.9I & 0.9 & 0.5 \\ I & 0.8 & 0.4I & 1 \\ I & 0 & 0 & 0.8I \end{bmatrix},$$

we see the minimal real is 0, the minimal neutrosophic element is 0.I. The maximal real is 1 and the maximal neutrosophic entry is I.

$[A_3, B_3]$ is the neutrosophic interval matrix of $6 \times 6$ matrices where

$$A_3 = \begin{bmatrix} 0 & 1 & 3I & 4 & 0 & 8I \\ 9I & 6 & 12I & 0 & 9 & 12 \\ 31 & 8 & 10 & 9I & 4 & 0 \\ 0 & 1 & 2 & 6I & 0 & 9I \\ 1 & 2 & 6 & 8I & 0 & 0 \\ 4I & 0 & 0 & 9I & 8 & 12 \end{bmatrix}$$

$$B_3 = \begin{bmatrix} 20 & 1 & 3I & 4 & 0 & 18I \\ 10I & 6 & 21I & 0 & 9 & 12 \\ 85 & 8 & 10 & 9I & 4+8I & 1 \\ 90I & 12 & 21 & 6I & 8 & 9I \\ 120 & 20I & 6+6I & 18I & 10 & 6 \\ 14I & I & 6 & 9I & 8 & 12 \end{bmatrix}.$$

We see the real minimum is 0 the real maximum is 120 the neutrosophic minimum is 3I and the neutrosophic maximum is 90I.

The neutrosophic interval matrix $[A_4, B_4]$ is such that it is a square $2 \times 2$ neutrosophic matrix.

Now let

$$A_4 = \begin{bmatrix} 0.3 & 0 \\ 0.2I & 0.31I \end{bmatrix}$$

and

$$B_4 = \begin{bmatrix} 1+I & 1 \\ 0.8I & 0.32I \end{bmatrix}.$$

The minimum real is 0. The maximum real is 1. The minimum neutrosophic value is 0.2I and the maximum neutrosophic value is I + 1. Thus we see $[A, B] = [A_1, B_1] \cup [A_2, B_2] \cup [A_3, B_3] \cup [A_4, B_4]$ is a neutrosophic interval mixed square 4-matrix.

Now we proceed on to define the notion of fuzzy neutrosophic interval bimatrix, fuzzy neutrosophic interval trimatrix and fuzzy neutrosophic interval n-matrix ($n > 3$). These concepts will be very much helpful for these are used in the fuzzy neutrosophic models, Fuzzy neutrosophic relational Bimaps (FNRBM) and Fuzzy neutrosophic relational n-maps (FNRnMs).

We have just introduced the notion of fuzzy neutrosophic interval matrices.

Now we will first define the notion of fuzzy neutrosophic interval bimatrices as we have already introduced the notion of fuzzy neutrosophic interval matrices.

**DEFINITION 2.3.8:** *Let $[A, B] = [A_1, B_1] \cup [A_2, B_2]$ where $[A_1, B_1]$ and $[A_2, B_2]$ are fuzzy neutrosophic interval matrices, where $A_1 = \left(a_{ij}^1\right)$, $B_1 = \left(b_{ij}^1\right)$, $A_2 = \left(a_{ij}^2\right)$ and $B_2 = \left(b_{ij}^2\right)$, the minimal elements in $A_i = \left(a_{ij}^i\right)$ will be the least element of the entries in* *the fuzzy neutrosophic interval matrices [$A_i$, $B_i$] and the maximal elements of $B_i$ = $\left(b_{ij}^{i}\right)$ will be the greatest element of fuzzy neutrosophic interval matrices [$A_i$, $B_i$], i = 1, 2.*

*Thus [A, B] will contain elements M = $M_1 \cup M_2$ which are fuzzy neutrosophic bimatrices with $M_1 \in$ [$A_1$, $B_1$] and $M_2 \in$ [$A_2$, $B_2$]. [A, B] is called the fuzzy neutrosophic interval bimatrix.*

We first illustrate this by an example.

***Example 2.3.12:*** Let [A, B] = [$A_1$, $B_1$] $\cup$ [$A_2$, $B_2$] be a fuzzy neutrosophic interval bimatrix where [$A_1$, $B_1$] is the fuzzy neutrosophic interval of 2 × 3 matrices and [$A_2$, $B_2$] is a fuzzy neutrosophic interval of 2 × 2 matrices

$$A_1 = \begin{bmatrix} 0.I & 0.2 & 0 \\ 0 & 0.3I & 0.1 \end{bmatrix}$$

and

$$B_1 = \begin{bmatrix} I & 1 & 0.9 \\ 0.8 & 0.9I & 1 \end{bmatrix},$$

where 0 is the real minimum and 0.I is the neutrosophic minimum and 1 is the real maximum and I is the neutrosophic maximum.

Here

$$A_2 = \begin{bmatrix} 0 & 0.2I \\ 0.3I & 0.1 \end{bmatrix}$$

and

$$B_2 = \begin{bmatrix} 0.9I & 0.8I \\ 0.3I & 0.7 \end{bmatrix}.$$

Here 0 is the real minimum, 0.7 is the real maximum and 0.2I is the neutrosophic minimum and 0.9I is the neutrosophic maximum. [A, B] = [$A_1$, $B_1$] $\cup$ [$A_2$, $B_2$] is the fuzzy neutrosophic interval bimatrix we call this fuzzy neutrosophic interval bimatrix to be a mixed fuzzy neutrosophic interval bimatrix or fuzzy neutrosophic interval mixed bimatrix.

Now we proceed on to define the 5 district types of fuzzy neutrosophic interval matrices.

A fuzzy neutrosophic interval bimatrix $[A, B] = [A_1, B_1] \cup [A_2, B_2]$ is said to be a fuzzy neutrosophic interval $m \times m$, square bimatrix if both the fuzzy neutrosophic interval matrices $[A_1, B_1]$ and $[A_2, B_2]$ are $m \times m$ fuzzy neutrosophic square matrices.

We call the fuzzy neutrosophic interval bimatrix $[A, B] = [A_1, B_1] \cup [A_2, B_2]$ to be $m \times p$ rectangular bimatrix if both $[A_1, B_1]$ and $[A_2, B_2]$ are $m \times p$ ($m \neq p$) rectangular fuzzy neutrosophic interval matrices.

The fuzzy neutrosophic interval bimatrix $[A, B] = [A_1, B_1] \cup [A_2, B_2]$ is said to be a mixed rectangular fuzzy neutrosophic bimatrix if $[A_1, B_1]$ is a fuzzy neutrosophic interval $m \times p$ ($m \neq p$) rectangular matrix and $[A_2, B_2]$ is a fuzzy neutrosophic interval $t \times q$ ( $t \neq q$) rectangular matrix ($t \neq m$).

The fuzzy neutrosophic interval bimatrix $[A, B] = [A_1, B_1] \cup [A_2, B_2]$ is said to be a fuzzy neutrosophic interval mixed square bimatrix if $[A_1, B_1]$ is a fuzzy neutrosophic interval $m \times m$ square matrix and $[A_2, B_2]$ is a fuzzy neutrosophic interval of $n \times n$ square matrix $m \neq n$.

A fuzzy neutrosophic interval bimatrix $[A, B] = [A_1, B_1] \cup [A_2, B_2]$ is said to be a fuzzy neutrosophic interval mixed bimatrix if $[A_1, B_1]$ is a fuzzy neutrosophic interval square matrix and $[A_2, B_2]$ is a fuzzy neutrosophic interval rectangular matrix.

Now we proceed on to define the new notion of fuzzy neutrosophic interval trimatrices and illustrate them with nice examples.

**DEFINITION 2.3.9:** *Let [A, B] = $[A_1, B_1] \cup [A_2, B_2] \cup [A_3, B_3]$ where each of the $[A_i, B_i]$ are fuzzy neutrosophic interval matrices for i = 1, 2, 3. Let [A, B] denote the collection of all trimatrices $M = M_1 \cup M_2 \cup M_2$ where $M_i$ is a fuzzy neutrosophic matrix from the fuzzy neutrosophic interval matrix $[A_i, B_i]$; i = 1, 2, 3.*

*We call [A, B] to be the fuzzy neutrosophic interval trimatrix and any element in [A, B] will be a trimatrix $M = M_1 \cup M_2 \cup M_3$.*

Now we illustrate this situation by a simple example.

***Example 2.3.13:*** Let [A, B] = $[A_1, B_1] \cup [A_2, B_2] \cup [A_3, B_3]$, where [A, B] is a fuzzy neutrosophic interval trimatrix with $[A_1, B_1]$ a 3 × 1 fuzzy neutrosophic rectangular matrix where

$$A_1 = \begin{bmatrix} 0.3I \\ 0 \\ 0.2I \end{bmatrix}$$

and

$$B_1 = \begin{bmatrix} I \\ 1 \\ 0.7I \end{bmatrix}.$$

The minimum real is 0, maximum real is 1 minimum neutrosophic is 0.2I and maximum neutrosophic value is I. $[A_2, B_2]$ is a 2 × 2 fuzzy neutrosophic interval matrix with

$$A_2 = \begin{bmatrix} 0.6 & 0.3I \\ 0.2 & 0.4I \end{bmatrix}$$

and

$$B_2 = \begin{bmatrix} 1 & I \\ 0.7 & 0.4I \end{bmatrix}$$

the maximum real is 1, the minimum real is 0.2, the maximum neutrosophic value is I and the minimum neutrosophic value is

0.3I. $[A_3, B_3]$ is a $2 \times 4$ fuzzy neutrosophic interval rectangular matrix where

$$A_3 = \begin{bmatrix} 0.01I & 0.1 & 0.08 & 0.03 \\ 0.5 & 0.2I & 0.1 & 0.08I \end{bmatrix}$$

and

$$B_3 = \begin{bmatrix} I & 1 & 0.3 & 0.8 \\ 0.9 & I & 1 & 0.8I \end{bmatrix}.$$

The maximum real value is 1, the maximum neutrosophic value is I. The minimum real value is 0.03 and the minimum neutrosophic value is 0.01 I. Thus $[A, B] = [A_1, B_1] \cup [A_2, B_2] \cup [A_3, B_3]$ is a fuzzy neutrosophic interval mixed trimatrices. Interested reader can construct more examples of such fuzzy neutrosophic interval trimatrices.

Now we proceed on to define the notion of fuzzy neutrosophic interval n-matrices (n > 3).

**DEFINITION 2.3.10:** *Let $[A, B] = [A_1, B_1] \cup [A_2, B_2] \cup \ldots \cup [A_n, B_n]$, where $[A_i, B_i]$ is a fuzzy neutrosophic interval matrix for i = 1, 2, …, n. [A, B] be the collection of all fuzzy neutrosophic n-matrices of the form $M = M_1 \cup M_2 \cup \ldots \cup M_n$ where each $M_i \in [A_i, B_i]$ i.e., $M_i$ is the fuzzy neutrosophic matrix from the fuzzy neutrosophic interval of matrices, true for i = 1, 2, …, n. We call [A, B] the fuzzy neutrosophic interval n matrices.*

First we illustrate this situation for n = 5, before we define finer properties about fuzzy neutrosophic interval n-matrices.

***Example 2.3.14:*** Let $[A, B] = [A_1, B_1] \cup [A_2, B_2] \cup \ldots \cup [A_5, B_5]$ be the fuzzy neutrosophic interval 5-matrices where $[A_1, B_1]$ is the fuzzy neutrosophic interval of $2 \times 2$ matrices with

$$A_1 = \begin{bmatrix} 0 & 0.01I \\ 0.04 & 0.09I \end{bmatrix}$$

and

$$B_1 = \begin{bmatrix} 1 & I \\ 0.9 & 0.9I \end{bmatrix}.$$

0 is the real minimum, 1 is the real maximum. I the neutrosophic maximum and 0.01I the neutrosophic minimum. $[A_2, B_2]$ is a fuzzy neutrosophic interval 1 × 6 matrix, where

$$A_2 = [(0, 01I, 0.1, 0.04I, 0.003I, 0.6)]$$

and

$$B_2 = [(0.01, I, 1, I, 0.06I, 0.9)]$$

0 is the real minimum, 1 is the real maximum. 0.003I is the neutrosophic minimum and I is the neutrosophic maximum.

The fuzzy neutrosophic interval 3 × 2 matrix $[A_3, B_3]$ with

$$A_3 = \begin{bmatrix} 0.0004I & 0.02 \\ 0.02I & 0.04 \\ 0.003 & 0.1 \end{bmatrix}$$

and

$$B_3 = \begin{bmatrix} 0.9I & 0.8 \\ 0.8I & 0.6 \\ 0.5 & 0.6 \end{bmatrix}.$$

The real minimum is 0.003, the real maximum is 0.8, the neutrosophic minimum being 0.0004I and the neutrosophic maximum being 0.9I.

$[A_4, B_4]$ is a fuzzy neutrosophic interval 4 × 1 column matrix with

$$A_4 = \begin{bmatrix} 0.003I \\ 0.004 \\ 0.001 \\ 0.013I \end{bmatrix}$$

and

$$B_4 = \begin{bmatrix} 0.04I \\ 0.8 \\ 0.7 \\ 0.15I \end{bmatrix}.$$

The real minimum is 0.004 and the real maximum is 0.8. The neutrosophic minimum is 0.003I and the neutrosophic maximum being 0.15I.

$[A_5, B_5]$ the fuzzy neutrosophic interval $2 \times 4$ matrices with

$$A_5 = \begin{bmatrix} 0.00I & 0.004 & 0.007 & 0.0014 \\ 0.07I & 0.0006I & 0.001 & 0.0009 \end{bmatrix}$$

and

$$B_5 = \begin{bmatrix} 0.09I & 0.7 & 0.63 & 0.14 \\ 0.07I & 0.61I & 0.07 & 0.114 \end{bmatrix}.$$

The real minimum is 0.0009 and the real maximum being 0.63. The neutrosophic minimum 0.0006I and the neutrosophic maximum being 0.61I.

Thus $[A, B] = [A_1, B_1] \cup [A_2, B_2] \cup \ldots \cup [A_5, B_5]$ is a fuzzy neutrosophic interval mixed 5-matrix.

Now having used the term 'mixed' in the definition of fuzzy neutrosophic interval n-matrix we proceed to explain the 5 categories of fuzzy neutrosophic interval n-matrices ($n > 3$).

Let $[A, B] = [A_1, B_1] \cup [A_2, B_2] \cup \ldots \cup [A_n, B_n]$ be a fuzzy neutrosophic interval n-matrices. We say [A, B] is a fuzzy neutrosophic interval square n-matrices, if each of the fuzzy neutrosophic interval matrices $[A_i, B_i]$ is a $m \times m$ square matrix for i = 1, 2, 3, …, n, thus any element M in [A, B] is a $m \times m$ square n-matrix i.e. $M = M_1 \cup M_2 \cup \ldots \cup M_n$ and each $M_i \in$ $[A_i, B_i]$ and $M_i$ is a $m \times m$ fuzzy neutrosophic matrix for i = 1, 2, …, n.

We call [A, B] a fuzzy neutrosophic interval rectangular n-matrix if each of the fuzzy neutrosophic interval matrix $[A_i, B_i]$ is a $m \times p$ ($m \neq p$) matrix for i = 1, 2, 3, …, n, where any element M in [A, B] is expressed in the form $M = M_1 \cup M_2 \cup \ldots \cup M_n$ where $M_i \in [A_i, B_i]$ is a fuzzy neutrosophic $m \times p$ matrix for i = 1, 2, …, n and M is the $m \times p$ rectangular fuzzy neutrosophic n-matrix.

Now we proceed onto define fuzzy neutrosophic interval mixed n-matrices.

Let $[A, B] = [A_1, B_1] \cup [A_2, B_2] \cup \ldots \cup [A_n, B_n]$ be a fuzzy neutrosophic interval n-matrices. Here each $[A_i, B_i]$ is a fuzzy neutrosophic interval matrix, for i = 1, 2, …, n. If some of $[A_i, B_i]$ are fuzzy neutrosophic interval square matrices and some other (or rest of) $[A_j\ B_j]$ are fuzzy neutrosophic interval rectangular matrices $1 \leq i, j \leq n$. Then we call [A, B] to be a fuzzy neutrosophic interval mixed n-matrices.

Let $[A, B] = [A_1, B_1] \cup [A_2, B_2] \cup \ldots \cup [A_n, B_n]$ be a fuzzy neutrosophic interval n-matrices if each of $[A_i, B_i]$ is a fuzzy neutrosophic interval $p_i \times p_i$ square matrices, with $p_i \neq p_j$ if $i \neq j$, i = 1, 2, 3, …, n then we call [A, B] the fuzzy neutrosophic interval mixed square n-matrices.

$[A, B] = [A_1, B_1] \cup [A_2, B_2] \cup \ldots \cup [A_n, B_n]$ be a fuzzy neutrosophic interval n-matrices, if each of $[A_i, B_i]$ is a fuzzy neutrosophic $p_i \times q_i$ rectangular matrices ($p_i \neq q_i$) and $p_i \neq p_j$ if $i \neq j$ for all i = 1, 2, 3, …, n, then we call [A, B] to be a fuzzy neutrosophic interval mixed rectangular n-matrices.

Now we have defined and seen examples of several types of fuzzy neutrosophic interval n-matrices built using both the concept of n-matrices and fuzzy neutrosophic interval matrices. These will find their applications in neutrosophic fuzzy models.

**Chapter Three**

# FUZZY MODELS AND NEUTROSOPHIC MODELS USING FUZZY INTERVAL MATRICES AND NEUTROSOPHIC INTERVAL MATRICES

In the earlier chapter we have defined several types of interval n-matrices (n = 1, 2, …; n < ∞). When we use fuzzy interval n-matrices and fuzzy neutrosophic interval n-matrices they will find their applications in Fuzzy Cognitive Maps (FCMs), which are not simple FCMs, Fuzzy Relational Maps (FRMs) which are not simple FRMs, FCBMs and FRBMs. To this end we make some modifications. Further we define some more concepts on these fuzzy interval matrices and fuzzy interval bimatrices.

Like wise the concept of fuzzy neutrosophic interval matrices and fuzzy neutrosophic interval bimatrices will find their applications in the fuzzy neutrosophic models like NCMs NCBMs, NRMs and NRBMs. In fact their generalizations also find their applications. In order to apply them we first give some justification for defining some minimal, maximal and optimal elements in the fuzzy interval of matrices.

We without loss of generality assume in this chapter by a fuzzy interval matrix we mean those interval matrices [A, B] where the entries are taken from the fuzzy interval [–1, 1].

This chapter has nine sections; each of these sections are described by a brief introduction.

## 3.1 Description of FCIMs Model

In this section for the first time we introduce the notion of a new model, called Fuzzy Cognitive Interval maps model (FCIMs-model) when multi-experts opinion are available. This model is generalized to FCInMs model if n sets of experts give their opinion on the same problem. We have described these models with real world problems.

Suppose we are working with the fuzzy cognitive maps defined in chapter I of this book, we further make an assumption that all the p number of experts spells out their opinion on n number of nodes. Using the directed graph suppose we get a fuzzy matrix $P_i$ given by the $i^{th}$ expert and further the FCM is not a simple FCM, $1 \leq i \leq p$.

The entries in the $n \times n$ matrices will be from the fuzzy interval $[-1, 1]$. Now in general if we take the collection of all $n \times n$ fuzzy matrices with entries from the fuzzy interval $[-1, 1]$ we have the interval of fuzzy $n \times n$ matrix. [A, B] would be an infinite collection in general. Further this fuzzy interval matrices associated with any FCM model will satisfy the following condition.

1) The fuzzy interval matrix will always be a square matrix.
2) The fuzzy interval matrix will always have the main diagonal entries to be zero.
3) The number of fuzzy matrices in the fuzzy interval $n \times n$ matrices is though infinite by all means for us in our system we can have only a finite number of fuzzy matrices associated with a FCM and with its associated fuzzy interval $n \times n$ square matrices.

This fuzzy interval of $n \times n$ square matrix associated with the FCMs of p-experts will be called as the Fuzzy Cognitive

Interval Maps (FCIMs) of the multi experts dynamical system, as the related connection matrices forms an interval of $n \times n$ square matrices.

This fuzzy interval of matrices satisfying the conditions 1, 2 and 3 has lots of advantage over the infinite collection. Suppose we have some p experts who have given their opinion on n concepts. Then we will have only p, $n \times n$ square fuzzy matrices with entries from the fuzzy interval [–1, 1]. Now we using these p, $n \times n$ square fuzzy matrices form an associated interval of square $n \times n$ matrices by the following method.

Now we known if they are the fuzzy connection matrices all the main diagonal terms are zero. Now in order to obtain the fuzzy interval $n \times n$ matrix [A, B] using the p-matrices we have to construct A and B for A and B may not exist in general. Now we call A the minimal matrix (element) of [A, B] and B the maximal (element) of [A, B]. We give the method by which A is built using the p matrices. Suppose $A = (a_{ij})$, $a_{ii} = 0$ for $1 \leq i \leq n$; $1 \leq i, j \leq n$. So all the diagonal elements are zero. We make the observation of the element $a_{12}$ in all the p matrices $P_1, \ldots, P_p$ where $P^t = (a^t_{ij})$, $1 \leq t \leq p$; we take the minimum value from the p entries $a^1_{12}, a^2_{12}, \ldots, a^p_{12}$ and put in the new matrix as the value of $a_{12}$ likewise for every $a_{ij}$, $1 \leq i, j \leq n$.

This newly formed matrix may not in general be any of the matrices given by the p experts. We call this matrix the minimal element A of the fuzzy interval of the $n \times n$ matrices got from the p experts.

Like wise we form the maximal matrix $B = (b_{ij})$ by taking the maximal element. If B is the maximal matrix, the elements of the fuzzy interval of matrices, in general B need not be a connection matrix given by any of the p experts. Now having obtained the minimal and maximal fuzzy matrices as A and B; we form the fuzzy interval matrix [A, B]. Clearly by the very construction of A and B, all the p connection fuzzy matrices given by the p experts will lie in the fuzzy interval of matrices [A, B].

Now having constructed the minimal and maximal element on the fuzzy interval of matrices we construct the optimal fuzzy matrix O as follows;

$$O = \frac{A+B}{2} = \frac{(a_{ij} + b_{ij})}{2} = (o_{ij})$$

may be in [A, B] if O is not in [A, B] we include O in the fuzzy interval of matrices [A, B] and call it as the optimal fuzzy matrix of the interval of fuzzy matrices and the associated weighted directed graph will be called as the optimal weighted directed graph.

Now we can work with the minimal matrix A, maximal matrix B and optimal matrix O and compare our results. We would also adjoin the matrix $\overline{A}$ which will be the average matrix of the combined matrices of the p experts excluding the minimal, optimal and the maximal matrices provided they are not the opinion given by any of the p experts. We can find also the resultant of state vectors using the average connection matrix $\overline{A}$ if $\overline{A} \in$ [A, B], well other wise we will adjoin $\overline{A}$ also to the interval of fuzzy matrices as our interval contains only the matrices related with the p experts opinion.

These four matrices A, B, O and $\overline{A}$ may or may not in general be some of the p experts opinion, The fuzzy interval matrix which is formed will have A to be minimal and B to be the maximal element [A, B], the fuzzy interval matrix i.e., if M $\in$ [A, B] and if M = $(m_{ij})$ then $(a_{ij}) \leq (m_{ij}) \leq (b_{ij})$, $1 \leq i, j \leq n$ also for O =$(o_{ij})$ then, $a_{ij} \leq o_{ij} \leq b_{ij}$, $1 \leq i, j \leq n$. Further for $\overline{A} = \left(\overline{a_{ij}}\right)$. We have $a_{ij} \leq \overline{a_{ij}} \leq b_{ij}$ for $1 \leq i, j \leq n$.

Now we have roughly sketched how the FCIM model and the associated fuzzy interval matrix [A, B] is constructed which we choose to call as the fuzzy interval of FCIM matrices or matrix from the FCIM model or dynamical system.

Only while working with a real model, we can see how the fuzzy interval matrix of the FCIM functions and its influence on the very system.

***Note:*** If we use only simple FCMs we will not have much to say for the entries in the fuzzy interval matrix would be [-1, 0, 1] so no meaning can be attributed to the optimal and average

fuzzy matrices. Further the minimal and maximal matrix will have only entries from [-1, 0, 1] so nothing more can be achieved, due to these limitations only we have choosen to work with FCMs which are not simple and these will certainly give more accurate prediction/solution than the simple FCMs.

Now we proceed on to illustrate a FCIM model. We will use only three experts opinion, we would find the fuzzy interval matrix [A, B] and find O and $\overline{A}$ and work with a state vector and find its effect on A, B, O and $\overline{A}$, apart from the three experts opinion connection matrices. First we describe the problem and give the three experts directed graph and the related or associated connection matrices.

Here we give a model called the symptom disease model in children taking 5 major attributes using three experts. [157, 198, 201, 216].

Let us take the 5 attributes which the child shows as symptoms of the disease as $C_1, C_2, \ldots, C_5$:

$C_1$ – Fever with cold / cough
$C_2$ – Fever with vomiting/ loose motion / loss of appetite
$C_3$ – Respiratory diseases
$C_4$ – Gastroenteritis
$C_5$ – Tuberculosis.

The directed graph given by the doctor who is used as an expert. The directed graph is a weighted one; given by the figure 3.1.1.

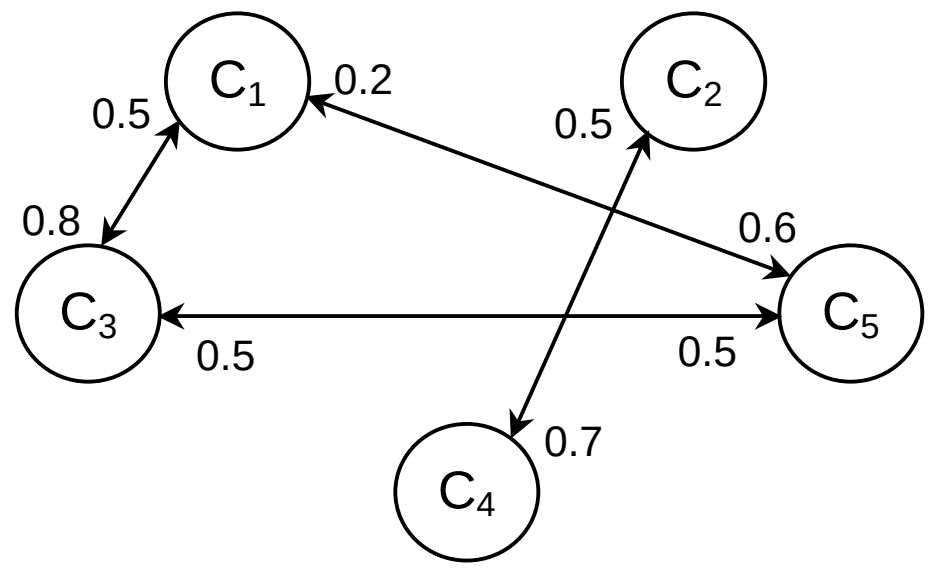

FIGURE 3.1.1

The fuzzy matrix $M_1$ (which is not a simple FCM) associated with the directed weighted graph is given below:

$$M_1 = \begin{array}{c} \\ C_1 \\ C_2 \\ C_3 \\ C_4 \\ C_5 \end{array}\begin{array}{c} \begin{array}{ccccc} C_1 & C_2 & C_3 & C_4 & C_5 \end{array} \\ \left[\begin{array}{ccccc} 0 & 0 & 0.8 & 0 & 0.6 \\ 0 & 0 & 0 & 0.7 & 0 \\ 0.5 & 0 & 0 & 0 & 0.5 \\ 0 & 0.5 & 0 & 0 & 0 \\ 0.2 & 0 & 0.5 & 0 & 0 \end{array}\right] \end{array}.$$

The opinion given by the second expert who is also a doctor and the weighted directed graph given by him is as follows:

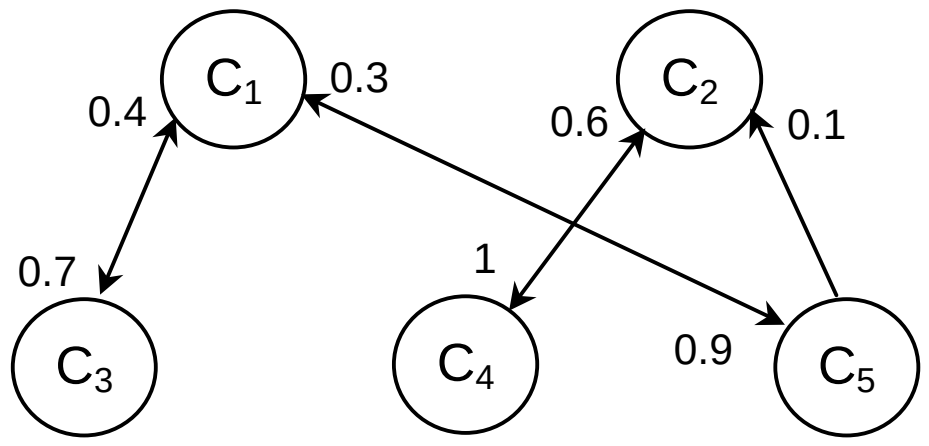

FIGURE 3.1.2

The fuzzy matrix $M_2$ from the directed graph given by the second expert is as follows:

$$M_2 = \begin{array}{c} \\ C_1 \\ C_2 \\ C_3 \\ C_4 \\ C_5 \end{array}\begin{array}{c} \begin{array}{ccccc} C_1 & C_2 & C_3 & C_4 & C_5 \end{array} \\ \left[\begin{array}{ccccc} 0 & 0 & 0.4 & 0 & 0.9 \\ 0 & 0 & 0 & 1 & 0 \\ 0.4 & 0 & 0 & 0 & 0 \\ 0 & 0.6 & 0 & 0 & 0 \\ 0.3 & 0.1 & 0 & 0 & 0 \end{array}\right] \end{array}.$$

The weighted directed graph given by the 3rd expert who is also a doctor is given by the figure 3.1.3:

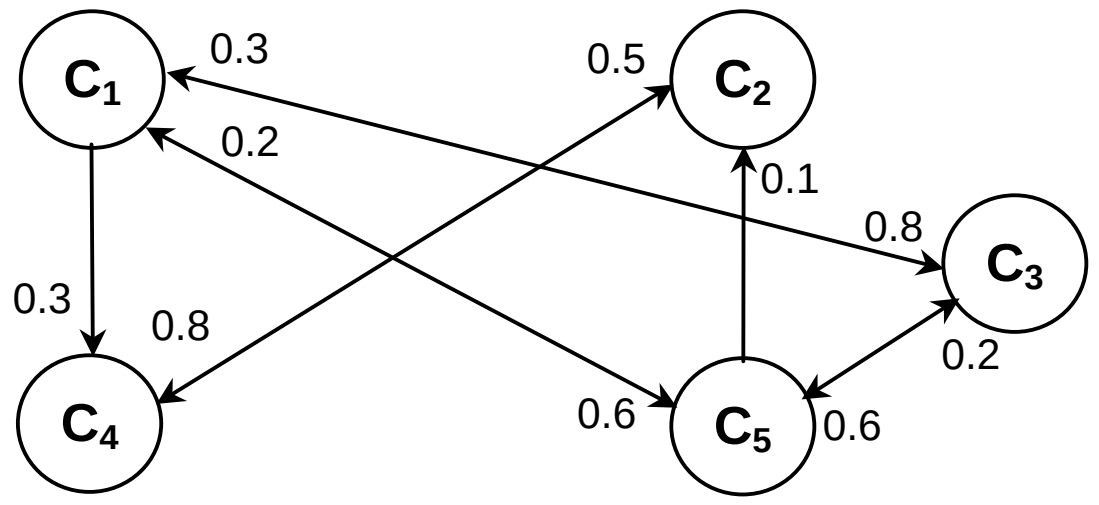

FIGURE 3.1.3

The related fuzzy matrix $M_3$ given by the 3rd expert is as follows:

$$M_3 = \begin{array}{c} \\ C_1 \\ C_2 \\ C_3 \\ C_4 \\ C_5 \end{array} \begin{array}{c} \begin{array}{ccccc} C_1 & C_2 & C_3 & C_4 & C_5 \end{array} \\ \begin{bmatrix} 0 & 0 & 0.8 & 0.3 & 0.6 \\ 0 & 0 & 0 & 0.8 & 0 \\ 0.3 & 0 & 0 & 0 & 0.6 \\ 0 & 0.5 & 0 & 0 & 0 \\ 0.2 & 0.1 & 0.2 & 0 & 0 \end{bmatrix} \end{array}$$

Using the three matrices $M_1$, $M_2$ and $M_3$ we obtain the minimal matrix A and the maximal matrix B so that we can form the fuzzy interval matrix of this FCIM model.

$$A = \begin{array}{c} \\ C_1 \\ C_2 \\ C_3 \\ C_4 \\ C_5 \end{array} \begin{array}{c} \begin{array}{ccccc} C_1 & C_2 & C_3 & C_4 & C_5 \end{array} \\ \begin{bmatrix} 0 & 0 & 0.4 & 0 & 0.6 \\ 0 & 0 & 0 & 0.7 & 0 \\ 0.3 & 0 & 0 & 0 & 0.5 \\ 0 & 0.5 & 0 & 0 & 0 \\ 0 & 0 & 0.2 & 0 & 0 \end{bmatrix} \end{array}$$

The maximal matrix B of the fuzzy interval matrix is given below:

$$B = \begin{array}{c} \\ C_1 \\ C_2 \\ C_3 \\ C_4 \\ C_5 \end{array} \begin{array}{c} \begin{array}{ccccc} C_1 & C_2 & C_3 & C_4 & C_5 \end{array} \\ \begin{bmatrix} 0 & 0 & 0.8 & 0.3 & 0.9 \\ 0 & 0 & 0 & 1 & 0 \\ 0.5 & 0 & 0 & 0 & 0.6 \\ 0 & 0.6 & 0 & 0 & 0 \\ 0.2 & 0.1 & 0.5 & 0 & 0 \end{bmatrix} \end{array}$$

Thus we have the fuzzy interval matrix associated with the FCIM given by [A, B], clearly $M_1, M_2, M_3 \in [A, B]$. Now the optimal matrix of the FCIM is given by O where

$$O = \begin{array}{c} \\ C_1 \\ C_2 \\ C_3 \\ C_4 \\ C_5 \end{array} \begin{array}{c} \begin{array}{ccccc} C_1 & C_2 & C_3 & C_4 & C_5 \end{array} \\ \begin{bmatrix} 0 & 0 & 0.6 & 0.15 & 0.75 \\ 0 & 0 & 0 & 0.85 & 0 \\ 0.4 & 0 & 0 & 0 & 0.55 \\ 0 & 0.55 & 0 & 0 & 0 \\ 0.1 & 0.05 & 0.35 & 0 & 0 \end{bmatrix} \end{array}$$

$O \in [A, B]$.

Now we find the CFCM using. $M_1$, $M_2$ and $M_3$ and the average fuzzy matrix $\overline{A}$; $\overline{A} \in [A, B]$.

$$\overline{A} = \frac{M_1 + M_2 + M_3}{3}$$

$$\overline{A} = \begin{bmatrix} 0 & 0 & 0.66 & 0.1 & 0.7 \\ 0 & 0 & 0 & 0.83 & 0 \\ 0.4 & 0 & 0 & 0 & 0.53 \\ 0 & 0.53 & 0 & 0 & 0 \\ 0.13 & 0.06 & 0.33 & 0 & 0 \end{bmatrix}.$$

Now we will study the effect of the state vector, symptom suffered is $C_1$ = fever with cold and cough alone in the on state and all other nodes are in the off state i.e., X = (1, 0, 0, 0, 0). We shall study the hidden pattern of X on the dynamical system FCIM given by the 3 expert doctors, the minimal A, maximal B, optimal O and the CFCM average matrix $\overline{A}$ using the max-min operation.

Hidden pattern of X given by the dynamical system $M_1$;

$XM_1 = (0, 0, 0.8, 0, 0.6)$

after updating and thresholding we get

$XM_1 \hookrightarrow X_1 = (1, 0, 0.8, 0, 0.6)$,

('$\hookrightarrow$' this symbol denotes the vector has been updated and thresholded).

Now the effect of $X_1$ on $M_1$ is given by

$X_1M_1 \hookrightarrow (1, 0, 0.8, 0, 0.6)$
$= X_2$ (say) ( = $X_1$).

Thus we arrive at the fixed point. According to the first expert who is doctor the symptom fever together with the symptom of cold and cough implies the child will have the maximum probability of suffering with some respiratory disease i.e., 0.8 degree it suffers from some respiratory disease and 0.6 degree it suffers from TB. However the difference is only 0.2.

Now we find the hidden pattern for X = (1, 0, 0, 0, 0) the same state vector using the second experts opinion.

$XM_2 = (0, 0, 0.4, 0, 0.9)$
$\hookrightarrow (1, 0, 0.4, 0, 0.9)$
$= X_1$ (say).

$X_1 M_2 \hookrightarrow (1, 0.1, 0.4, 0, 0.9)$
$= X_2$ (say)

$X_2 M_2 \hookrightarrow (1, 0.1, 0.4, 0.1, 0.9)$
$= X_3$ (say)

$X_3 M_2 \hookrightarrow (1, 0.1, 0.4, 0.1, 0.9)$
$= X_4 ( = X_3)$.

The hidden pattern of X is a fixed point given by the resultant vector $X_4 = (1, 0.1, 0.4, 0.1, 0.9)$. According to this expert who is also a doctor, we see that if a child suffers the symptom of fever with cold /cough it has the maximum possibility it suffers from tuberculosis / primary complexes and gives only 0.4 less than half of the possibility that it suffers from respiratory diseases. However it does not rule out that it may also suffer the symptom of fever with vomiting and gastroenteritis.

His argument that most of the infants born in India suffer from malnutrition and many suffer from the primary complexes (if the tuberculosis systems suffered by the child) he further adds when the child has cold/cough with fever many a time it will vomit due to cough and cold / cough with fever also may make the child suffer indigestion due to fever, he says.

Now we proceed on to work with the 3rd experts opinion on the same state vector

$X = (1\ 0\ 0\ 0\ 0)$ .

The effect of X on $M_3$ is given by

$X M_3 \hookrightarrow (1, 0, 0.8, 0.3, 0.6)$
$= X_1$ (say)
$X_1 M_3 \hookrightarrow (1, 0.3, 0.8, 0.3, 0.6)$
$= X_2$ (say)
$X_2 M_3 \hookrightarrow (1, 0.3, 0.8, 0.3, 0.6)$
$= X_3$ ( = $X_2$ say).

We see the hidden pattern of the resultant state vector is a fixed point given by

$X_3 = (1, 0.3, 0.8, 0.3, 0.6)$.

This doctor also feels like the first doctor the symptoms suffered by the child may be mainly due to respiratory disease. However he does not rule out the factor that the child may suffer from 0.6 degree T.B (Tuberculosis or primary complexes). He gives 0.3 degree to the factor that the child may have vomiting with fever and also 0.3 degree it may be due to gastroenteritis.

Now we study the effect of the same state vector X = (1, 0, 0, 0, 0) on the minimal fuzzy matrix A of the fuzzy interval matrix.

| | | |
|---|---|---|
| $XA$ | $\hookrightarrow$ | $(1, 0, 0.4, 0, 0.6)$ |
| | $=$ | $X_1$ (say) |
| $X_1A$ | $\hookrightarrow$ | $(1, 0, 0.4, 0, 0.6)$ |
| | $=$ | $X_2$ ( $= X_1$ say). |

Thus we see the resultant is a fixed point. The minimal matrix expresses; 0.6 degree, the child's symptom (say) may be due to Tuberculosis and 0.4 degree the symptom the child show is due to it may suffer due to respiratory diseases and totally rules out the possibility the child may suffer the symptom from vomiting or gastroenteritis.

Now we proceed on to study the effect of the state vector X = (1, 0, 0, 0, 0) on the maximal fuzzy matrix B of the fuzzy interval matrix [A, B].

| | | |
|---|---|---|
| $XB$ | $\hookrightarrow$ | $(1, 0, 0.8, 0.3, 0.9)$ |
| | $=$ | $X_1$ (say) |
| $X_1 B$ | $\hookrightarrow$ | $(1, 0.3, 0.8, 0.3, 0.9)$ |
| | $=$ | $X_2$ (say) |
| $X_2 B$ | $=$ | $(1, 0.3, 0.8, 0.3, 0.9)$ |
| | $=$ | $X_3$ ( $= X_2$). |

Thus we see the resultant is more or less the same as the 3rd expert but the maximal matrix gives maximal degree to T. B or primary complex disease i.e., according to this maximal system the symptom with fever and cold and cough is the main cause for T.B followed by respiratory diseases, however does not rule out the possibility that the child may have 0.3 degree of suffering from the symptom of vomiting and gastroenteritis.

Now we see the predictions of the optimal matrix O from the fuzzy interval of matrices for the same state vector X = (1, 0, 0, 0, 0)

| | | |
|---|---|---|
| XO | = | (1, 0, 0.6, 0.15, 0.75) |
| | = | $X_1$ (say) |
| $X_1O$ | = | (1, 0.15, 0.6, 0.15, 0.75) |
| | = | $X_2$ (say) |
| $X_2O$ | = | (1, 0.15, 0.6, 0.15, 0.75) |
| | = | $X_3$ (say) = $X_2$. |

Thus we see from the optimal matrix, if a child suffers with a symptom of fever with cold or cough alone there is 0.75 degree probability it suffers from TB (Tuberculosis) and also it may suffer from the respiratory disease upto 0.6 degree but the optimal matrix gives a very moderate degree (say) 0.15 to the symptom it may be have vomiting and also symptom of gastroenteritis.

Now finally we see the resultant of the state vector X = (1, 0, 0, 0, 0) on the CFCM average $\overline{A}$;

$$\text{i.e., } \overline{A} = \frac{M_1 + M_2 + M_3}{3};$$

| | | |
|---|---|---|
| $X\overline{A}$ | $\hookrightarrow$ | (1, 0, 0.66, 0.1, 0.7) |
| | = | $X_1$ (say) |
| $X_1\overline{A}$ | $\hookrightarrow$ | (1, 0.06, 0.66, 0.1, 0.7) |
| | = | $X_2$ (say). |

Now

| | | |
|---|---|---|
| $X_2\overline{A}$ | $\hookrightarrow$ | (1, 0.6, 0.66, 0.1, 0.7) |
| | = | $X_3$ ( = $X_2$). |

The hidden pattern of the dynamical system is a fixed point given by $X_3$ = (1, 0.06, 0.66, 0.1, 0.7). Thus we can say the average value of the CFCM happens to be the best prediction for the symptom fever with cold and cough for the doubt whether the child suffers from Tuberculosis or respiratory disease is very small amounting only to 0.04 degree. Mostly they feel due to cold; fever and cough the child may also suffer a very mild symptom of vomiting followed by a mild symptom of gastroenteritis. Thus the adaptation of fuzzy interval model can give an optimal solution. The above model is only an illustration to make the reader follow how the model works. That is why for compatibility we have used only very limited number of experts and also very limited number of nodes or concepts associated with symptom disease problem.

Here a very moderate rationalist solution is given by the fuzzy matrix $\overline{A}$ belong to the fuzzy interval matrix [A, B]. By this one cannot always conclude only the average of the combined FCMs will give the best or optimal result. It can be anything in general, the minimal matrix A or the maximal matrix B or the optimal matrix O or the matrix associated with any experts opinion. So we cannot make a generalization by this result.

## 3.2 Description and Illustration of FRIM Model

Now we proceed onto define how the fuzzy interval matrix related with a fuzzy relational map (FRM) is defined and illustrate it with some examples. Let the FRIM be formed using m nodes related with the domain space and n nodes related with the range space. We obtain some t experts opinion using the (m, n) nodes. We take only the FRMs, which are not simple FRMs. When we use t experts we use the concept of interval matrices and so choose to call the models as fuzzy relational Interval maps models i.e., (FRIMs models). Thus the FRIMs are weighed directed graphs. All the connection matrices associated with each of the FRMs is a fuzzy $m \times n$ matrix. Thus we have t number of $m \times n$ fuzzy matrices. How to make these fuzzy matrices into an fuzzy interval of m × n matrices. We obtain the fuzzy interval of m ×n matrices [A, B] as follows.

We formulate or build the minimal matrix A in the following way, if A = $(a_{ij})_{m\times n}$ what are the values of $(a_{ij})_{m \times n}$, 1 ≤ i ≤ m, 1 ≤ j ≤ n. $a_{11}$ is given the least value by observing all the t matrices; suppose $T^k = (t_{ij}^k)$, k = 1, 2, …, t choose the least value from $t_{11}^1, t_{11}^2, ..., t_{11}^t$ and mark it as $a_{11}$. For the value $a_{12}$ observe $t_{12}^1, t_{12}^2, ..., t_{12}^t$ and choose the least of them and mark it as $a_{12}$. Like wise fill all the m × n entries of the matrix A = $(a_{ij})$, A will be the minimal matrix of the interval of matrices using the t experts connection matrices $T^1, …, T^t$. Now how to find the values of B = $(b_{ij})$, the maximal element of the interval of matrices. The $b_{ij}$'s; 1 ≤ i, ≤ m, 1 ≤ j ≤ n are filled in the following way.

Choose the largest values from the t matrices $T^1, T^2, … T^t$. $T^k = (t_{ij}^k)$, k = 1, 2, …, t. Suppose we want to fill $b_{11}$ choose the largest (maximal) value from $t_{11}^1, t_{11}^2, ..., t_{11}^t$ and make it as $b_{11}$.

Now for the value $b_{12}$ choose the maximal value from $t_{12}^1, t_{12}^2, ..., t_{12}^t$.

Several of the values will be maximal or the largest values. Like wise fill all the m×n entries of the matrix B = $(b_{ij})$.

Now all the matrices $T^1, T^2, …, T^t \in [A, B]$. A will be called as the minimal m × n fuzzy matrix of the fuzzy interval matrix associated with the FRIM, similarly B will be called as the maximal m × n fuzzy matrix of the fuzzy interval matrix associated with the FRIM. Now the optimal matrix O is formed by

$$O = \frac{(a_{ij} + b_{ij})}{2};$$

O is a fuzzy m×n matrix and clearly O is an element of the fuzzy interval of m × n matrices. Now take

$$\overline{T} = \frac{T^1 + T^2 + ... + T^t}{t}$$

to be the average of the combined FRMS. $\overline{T} \in [A, B]$. So the associated fuzzy interval of m × n matrices with the FRIMs includes O and $\overline{T}$. We can for all our problems work with the t dynamical system $T^1, T^2, \ldots, T^t$ as well as A, B, O and $\overline{T}$ and arrive at the conclusions. It is still important to mention, this will enable us to compare any experts view with the minimal matrix A, maximal matrix B, optimal matrix O and the average matrix $\overline{T}$. Thus the FRIM model with t experts opinion can atmost have t + 4 number of m × n fuzzy matrices. We shall illustrate this also by a very simple model. These models are constructed mainly for illustrative purposes.

Suppose we are interested in studying a model relating teachers and students and we wish to use the FRIMs model. We just give only a simple small model so that the reader can easily understand the way in which it functions so that the reader can adopt it in any problem which he/she desires to study.

Let us consider the relationship between the teacher and the student suppose we take the domain space as the concepts belonging to the teacher, say $D_1, D_2, \ldots, D_5$ and the range space denote the concepts belonging to the student say $R_1$, $R_2$ and $R_3$.

We describe the nodes $D_1, D_2, \ldots, D_5$ of the domain space and $R_1$, $R_2$ and $R_3$ of the range space as follows;

- $D_1$ – teaching is good
- $D_2$ – teaching is poor
- $D_3$ – teaching is mediocre
- $D_4$ – teacher is kind (moulds the character of students in a right way)
- $D_5$ – teacher is harsh (or rude).

We can have more concepts associated with the teacher like non-reactive, as concerned, indifferent and so on. As our only motivation is to give an illustrative model we have restrained ourselves to work only with 5 attributes with teacher and only 3 concepts / nodes related with the student.

The nodes of the range space associated with the student are:

$R_1$ – Good student
$R_2$ – Bad student
$R_3$ – Average student.

The relational directed graph given by the first expert who is a headmaster of the school.

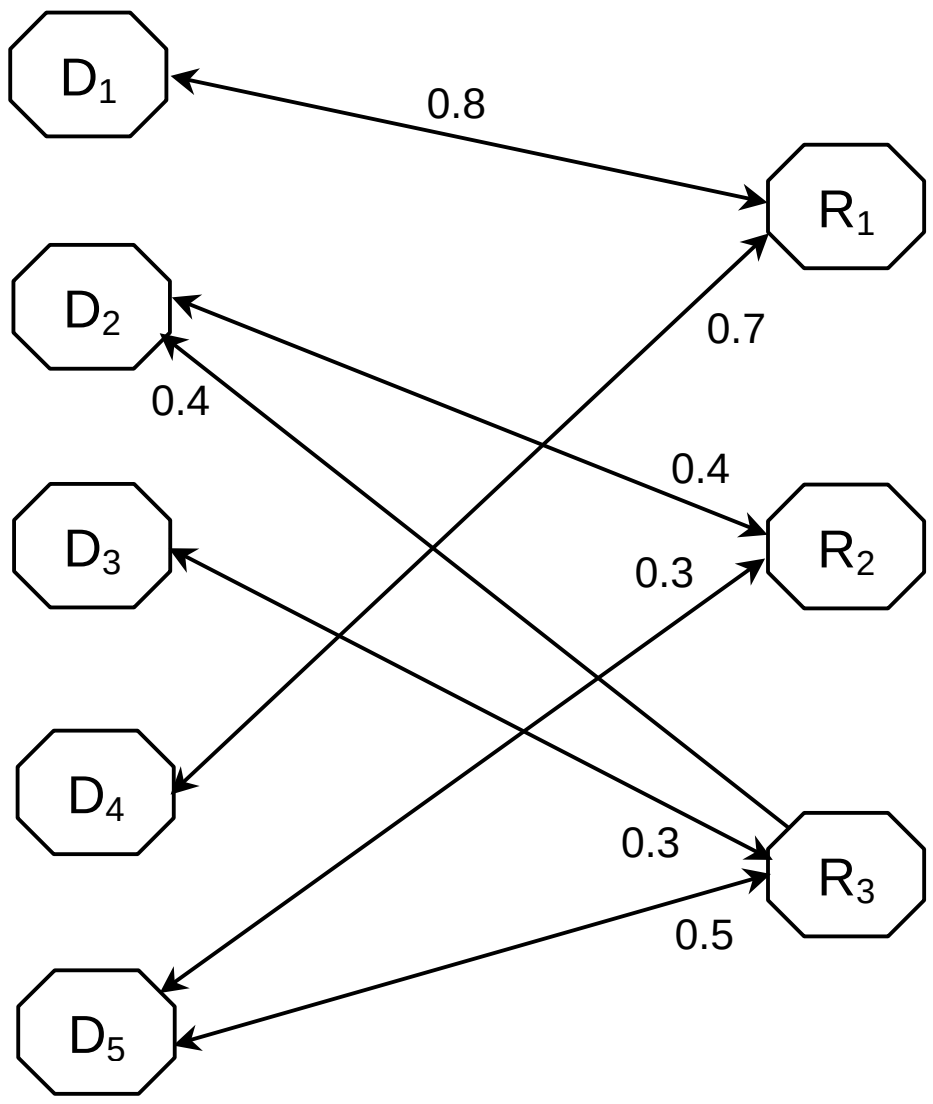

FIGURE: 3.2.1

Using this weighed directed graph we obtain the relational fuzzy $5 \times 3$ matrix $P_1$.

$$P_1 = \begin{array}{c} \\ D_1 \\ D_2 \\ D_3 \\ D_4 \\ D_5 \end{array} \begin{array}{c} \begin{array}{ccc} R_1 & R_2 & R_3 \end{array} \\ \begin{bmatrix} 0.8 & 0 & 0 \\ 0 & 0.4 & 0.4 \\ 0 & 0 & 0.3 \\ 0.7 & 0 & 0 \\ 0 & 0.3 & 0.5 \end{bmatrix} \end{array}.$$

Next we use the opinion of a retired teacher as the second expert. The directed weighed graph given by him is as follows.

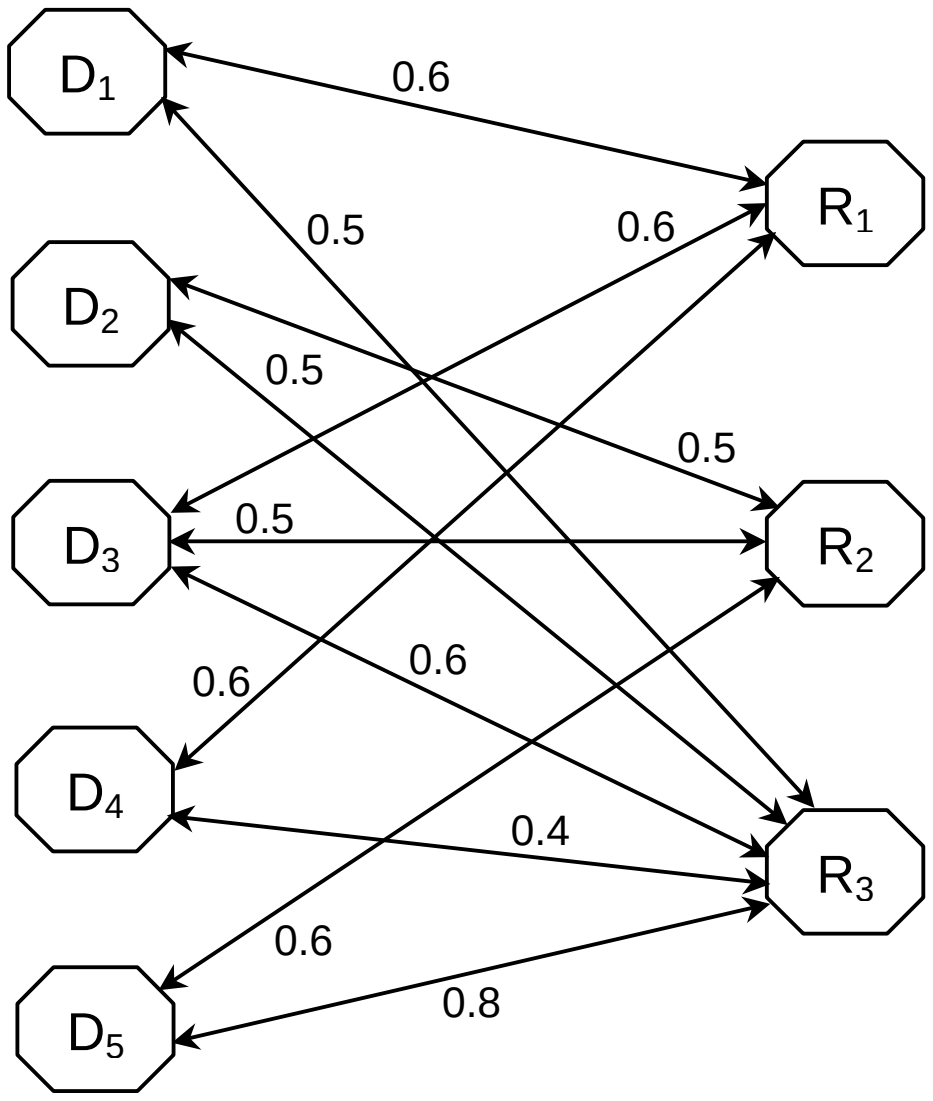

FIGURE: 3.2.2

Using this directed weighed graph we obtain the following fuzzy connection matrix.

$$P_2 = \begin{array}{c} \\ D_1 \\ D_2 \\ D_3 \\ D_4 \\ D_5 \end{array} \begin{array}{c} \begin{array}{ccc} R_1 & R_2 & R_3 \end{array} \\ \left[\begin{array}{ccc} 0.6 & 0 & 0.5 \\ 0 & 0.5 & 0.5 \\ 0.6 & 0.5 & 0.6 \\ 0.6 & 0 & 0.4 \\ 0 & 0.6 & 0.8 \end{array}\right] \end{array}.$$

Now we give the weighed directed graph of the 3rd expert who is an educationalist.

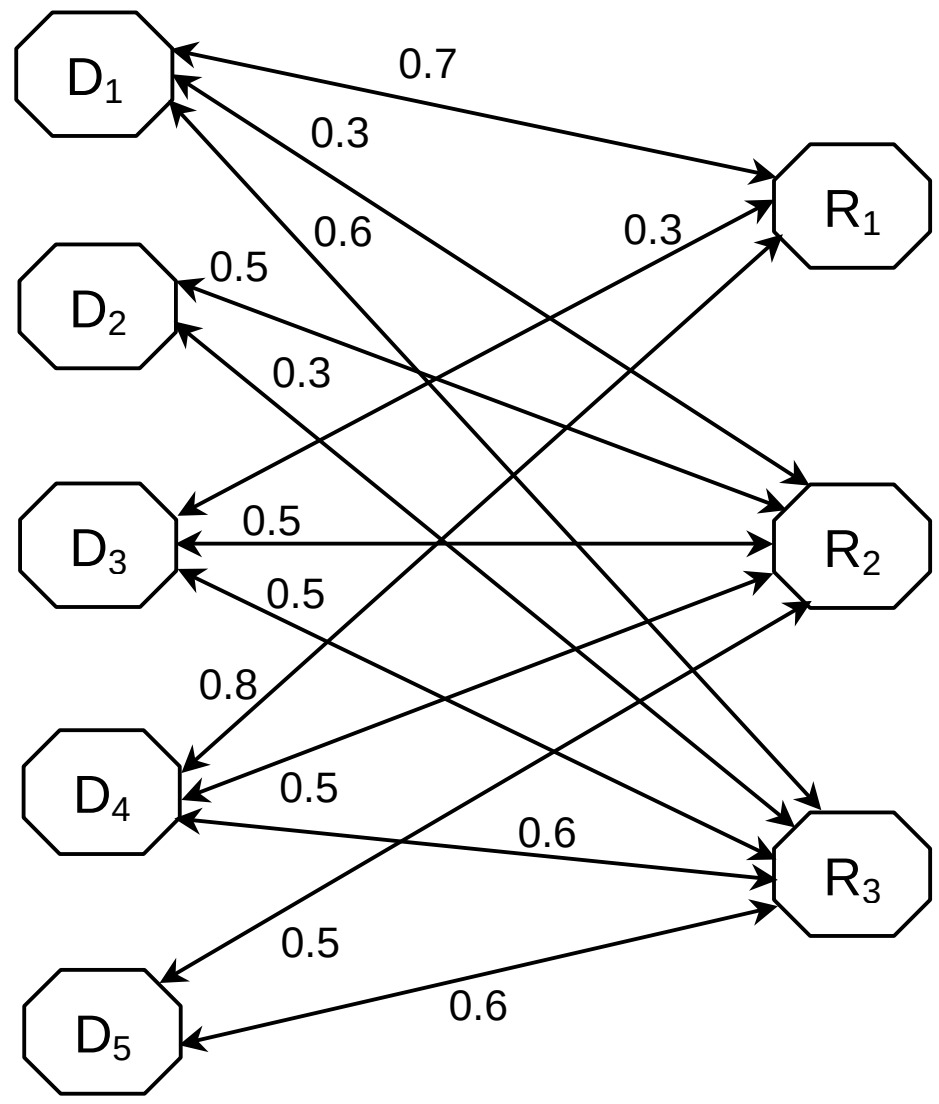

FIGURE: 3.2.3

**The related fuzzy relational matrix using the directed graph is as follows**

$$P_3 = \begin{array}{c} \\ D_1 \\ D_2 \\ D_3 \\ D_4 \\ D_5 \end{array} \begin{array}{c} \begin{array}{ccc} R_1 & R_2 & R_3 \end{array} \\ \begin{bmatrix} 0.7 & 0.3 & 0.6 \\ 0 & 0.5 & 0.3 \\ 0.3 & 0.5 & 0.5 \\ 0.8 & 0.5 & 0.6 \\ 0 & 0.5 & 0.6 \end{bmatrix} \end{array}.$$

Now we want to build the fuzzy interval of matrices using the three matrices $P_1$, $P_2$ and $P_3$. The minimal fuzzy matrix A of the fuzzy interval of matrices [A, B] containing the relational matrices $P_1$, $P_2$ and $P_3$ given by the 3 experts is

$$A = \begin{array}{c} \\ D_1 \\ D_2 \\ D_3 \\ D_4 \\ D_5 \end{array} \begin{array}{c} \begin{array}{ccc} R_1 & R_2 & R_3 \end{array} \\ \begin{bmatrix} 0.6 & 0 & 0 \\ 0 & 0.4 & 0.3 \\ 0 & 0 & 0.3 \\ 0.6 & 0 & 0 \\ 0 & 0.3 & 0.5 \end{bmatrix} \end{array}.$$

The maximal matrix B of the fuzzy interval of matrices associated with the FRIM is given by

$$B = \begin{array}{c} \\ D_1 \\ D_2 \\ D_3 \\ D_4 \\ D_5 \end{array} \begin{array}{c} \begin{array}{ccc} R_1 & R_2 & R_3 \end{array} \\ \begin{bmatrix} 0.8 & 0.3 & 0.6 \\ 0 & 0.5 & 0.5 \\ 0.63 & 0.5 & 0.6 \\ 0.8 & 0.5 & 0.6 \\ 0 & 0.6 & 0.8 \end{bmatrix} \end{array}.$$

Thus we have [A, B] to be fuzzy interval of matrices associated with the FRIM. Clearly, $P_1$, $P_2$, $P_3$ $\in$ [A, B]. Now the optimal matrix O is given by

$$O = \begin{bmatrix} 0.7 & 0.15 & 0.3 \\ 0 & 0.45 & 0.4 \\ 0.32 & 0.25 & 0.45 \\ 0.7 & 0.25 & 0.3 \\ 0 & 0.45 & 0.65 \end{bmatrix}.$$

We see O, the optimal matrix is also an element of the fuzzy interval of matrices [A, B]. Let $\overline{P}$ be the average of the 3 FRM or $\overline{P}$ is the average of the combined FRMs.

$$\text{i.e., } \overline{P} = \frac{P_1 + P_2 + P_3}{3},$$

$$P = \begin{bmatrix} 0.7 & 0.1 & 0.37 \\ 0 & 0.47 & 0.4 \\ 0.31 & 0.33 & 0.47 \\ 0.7 & 0 & 0.33 \\ 0 & 0.47 & 0.63 \end{bmatrix};$$

a fuzzy matrix from the fuzzy interval of matrices associated with the FRIM. We will see using the 7 matrices to obtain the resultant, given by the state vector; with the only node teacher is harsh in the on state i.e., the node $D_5$ is in the on state and all other nodes are in the off state is obtained.

$$\begin{aligned} Y &= (0\ 0\ 0\ 0\ 1) \in D \\ YP_1 &= (0, 0.3, 0.5) \\ &= X \text{ (say)} \\ XP_1^T &\hookrightarrow (0, 0.4, 0.3, 0, 1) \\ &= Y_1 \text{ (say)} \\ Y_1 P_1 &= (0, 0.3, 0.5) \\ &= X_1 \ (= X). \end{aligned}$$

Thus the hidden pattern is a fixed point. According to this, harsh teacher cannot make a good student. Harsh teacher can 0.5 degree produce an average student and 0.3 degree produce a bad student. The teacher being harsh among school going children mars the production of any good children, according to this expert who is a headmaster.

Now we work with the same state vector but use the dynamical system formed by the second expert. Let Y = (0 0 0 0 1),

$$\begin{aligned} YP_2 &= (0, 0.6, 0.8) \\ &= X \text{ (say)} \\ XP_2^T &= (0.5, 0.5, 0.6, 0.4, 1) \\ &= Y_1 \text{ (say)} \\ Y_1 P_2 &= (0.6, 0.6, 0.8) \\ &= X_1 \text{ (say)} \\ X_1 P_2^T &= (0.6, 0.5, 0.6, 0.6, 1) \\ &= Y_2 \text{ (say)} \end{aligned}$$

$$\begin{aligned} Y_2 P_2 &= (0.6, 0.6, 0.8) \\ &= X_2 \, (=X_1). \end{aligned}$$

Thus a harsh teacher according to this expert can produce an average student but have equal degree of producing bad and good students.

Now we study the effect of the state vector Y = (0 0 0 0 1) on the dynamical system $P_3$ given by the 3rd expert.

$$\begin{aligned} YP_3 &= (0, 0.5, 0.6) \\ &= X \text{ (say)} \\ XP_3^T &\hookrightarrow (0.6, 0.5, 0.5, 0.6, 1) \\ &= Y_1 \text{ (say)} \\ Y_1 P_3 &= (0.6, 0.5, 0.6) \\ &= X_1 \text{ (say)} \\ X_1P_3^T &\hookrightarrow (0.6, 0.5, 0.5, 0.6, 0.1) \\ &= Y_2 \text{ (say)} \\ Y_2 P_3 &= (0.6, 0.5, 0.6) \\ &= X_2 \, (=X_1). \end{aligned}$$

Thus the hidden pattern is a fixed point. According to this expert, a harsh teacher can produce both good and average student to the same degree viz. 0.6 and also because of his harshness he produces 50% of bad students.

Now we work with the minimal fuzzy matrix A of the fuzzy interval matrix [A, B] which was constructed using the matrices $P_1$, $P_2$ and $P_3$ using the state vector,

$$\begin{aligned} Y &= (0\ 0\ 0\ 0\ 1) \\ YA &= (0, 0.3, 0.5) \\ &= X \text{ (say)} \\ XA^T &\hookrightarrow (0, 0.3, 0.3, 0, 1) \\ &= Y_1 \text{ (say)} \\ YA &= (0, 0.3, 0.5) \\ &= X_1 \, (= X). \end{aligned}$$

Thus the hidden pattern is a fixed point. According to the minimal matrix A, we get a harsh teacher produces a average student and no good student but has a chance of producing 0.3 bad students.

Now we work with the same state vector Y = (0 0 0 0 1) on the maximal fuzzy matrix B.

$$
\begin{aligned}
YB &= (0, 0.6, 0.8) \\
&= X \text{ (say)} \\
XB^T &\hookrightarrow (0.6, 0.5, 0.6, 0.6, 1) \\
&= Y_1 \text{ (say)}.
\end{aligned}
$$

$Y_1B$ = (0.6, 0.6, 0.8), according to the maximal fuzzy matrix of the fuzzy interval of matrices. We see a harsh teacher is certain to produce the average student and equally a harsh teacher produces a bad student and a good student.

Now we proceed on to find the effect of the state vector Y = (0 0 0 0 1) on the optimal matrix O ∈ [A B].

$$
\begin{aligned}
YO &= (0, 0.45, 0.65) \\
&= X \\
XO^T &\hookrightarrow (0.3, 0.45, 0.45, 0.3, 1) \\
&= Y_1 \text{ (say)} \\
Y_1O &= (0.32, 0.45, 0.65) \\
&= X_1 \\
X_1O^T &\hookrightarrow (0.32, 0.45, 0.45, 0.32, 1) \\
&= Y_2 \text{ (say)} \\
Y_2O &= (0.32, 0.45, 0.65) \\
&= Y_2 \ (= X_1).
\end{aligned}
$$

Thus the hidden pattern of the state vector given by the optimal system is a harsh teacher produces always an average student, and 0.45 degree of making a bad student and 0.32 degree only in making a good student.

Now we find the effect of the state vector Y by the average dynamical system $\overline{P}$.

$Y\overline{P}$ = (0, 0.47, 0.63)
= X

$X\overline{P}^T$ ↪ (0.37, 0.47, 0.47, 0.33, 1)
= $Y_1$ (say)

$Y_1\overline{P}$ = (0.37, 0.47, 0.63)
= $X_1$ (say)

$X_1\overline{P}^T$ ↪ (0.37, 0.47, 0.47, 0.37, 1)
= $Y_2$ (say)

$Y_2\overline{P}$ = (0.37, 0.47, 0.63)
= $X_2$ (= $X_1$).

Thus the hidden pattern is a fixed point. The harsh teacher invariably produces an average student but can produce a bad student up to 0.47 degree and a good student only upto degree 0.37. Thus we see a harsh teacher by their harshness, rudeness and making the student fear (always) and are certain to produce an average student and the possibility of making bad students. Just we have shown how a FRIM related fuzzy interval matrix functions.

## 3.3 Description of FCIBM Model and its Generalization

Next we proceed on to show how the Fuzzy Cognitive Interval bimatrix maps (FCIBM) model works with a fuzzy interval square bimatrices, also the new concept of fuzzy interval of square bimatrices is introduced. Now we show how they are applied in FCIBM's model. Suppose we have several experts analyzing the problem and if each one of them accept to work with two sets of concepts with different numbers say m and n, n ≠ m and some opt to give opinion on m concepts and others opt to give opinion on n concepts; how to find a model, which has the capacity to work on two sets of experts opinions taking the 2 sets of concepts. The FCIBMs constructed in [221] can cater only to two experts at a time, so we now construct a new model called the fuzzy interval bimatrix model associated with a FCIBM.

Let us consider a problem which has m attributes/concepts chosen by a set of t experts and n attributes / concepts chosen by another set of p experts ($n \neq m$), we may have overlaps of concepts, both the sets of workers work on the same problem. We first take the 't' experts opinion on the m concepts and form the $m \times m$ connection FCM matrices. We take only non simple FCM models. Using these t number of $m \times m$ matrices we construct the minimal $m \times m$ fuzzy matrix $A_1$, maximal $m \times m$ fuzzy matrix $B_1$, the optimal matrix $O_1$, where

$$O_1 = \frac{A_1 + B_1}{2} = \left( \frac{\left(a_{ij}^1\right) + \left(b_{ij}^1\right)}{2} \right).$$

The average of all the CFCM matrices $M_1^1 + M_2^1 + \ldots + M_t^1 = \overline{M}^1$ and

$$\overline{M}^1 = \frac{M_1^1 + M_2^1 + \ldots + M_t^1}{t}.$$

Now the fuzzy interval $m \times m$ matrix $[A_1, B_1]$ contains all the t number of $m \times m$ matrices together with $O_1$ and $\overline{M}^1$.

In a similar way now using the p experts opinion and using the n concepts we work and obtain the $[A_2, B_2]$ fuzzy interval $n \times n$ matrices. Thus $[A_2, B_2]$ will be a fuzzy interval $n \times n$ matrices having $A_2 = (a^2_{ij})$ to be the minimal fuzzy matrix of the fuzzy interval matrix $[A_2, B_2]$ and $B_2 = (b^2_{ij})$ to be the maximal fuzzy matrix, constructed using the method given in the earlier section.

$$O^2 = \frac{A_2 + B_2}{2} = \left( \frac{\left(a_{ij}^2\right) + \left(b_{ij}^2\right)}{2} \right) \text{ and } \overline{M^2} = \frac{M_1^2 + \ldots + M_p^2}{p}$$

the fuzzy interval $n \times n$ matrix. $[A_2, B_2]$ will contain $O_2$ and $\overline{M}_2$ together with the p number of $n \times n$ matrices $M_1^2, M_2^2, \ldots, M_p^2$.

Now set [A, B] = [$A_1$, $B_1$] ∪ [$A_2$, $B_2$]; [A, B] will be a fuzzy interval mixed square bimatrix. [A, B] is defined as the fuzzy interval mixed square bimatrix of the FCIBM or associated with it. We will just show how the model functions, suppose [A, B] = [$A_1$, $B_1$] ∪ [$A_2$, $B_2$] be a fuzzy interval square bimatrix associated with the FCIBM, suppose [$A_1$, $B_1$] is a fuzzy interval 4 × 4 square matrix associated with an FCIM on some problem and [$A_2$, $B_2$] is the fuzzy interval 5 × 5 square matrix associated with an FCIM for the same problem. Thus [A, B] = [$A_1$, $B_1$] ∪ [$A_2$, $B_2$] is a fuzzy interval mixed square bimatrix associated with the FCIBM model.

Let $M = M_1 \cup M_2$

$$= \begin{bmatrix} 0 & 0.6 & 0.4 & 0 \\ 0.7 & 0 & 0.2 & 0.1 \\ 0 & 0.6 & 0 & 0 \\ 0.2 & 0.1 & 0 & 0 \end{bmatrix} \cup \begin{bmatrix} 0 & 0.7 & 0.2 & 0.1 & 0 \\ 0.5 & 0 & 0 & 0 & 0.3 \\ 0.7 & 0 & 0 & 0 & 0.6 \\ 0.4 & 0.5 & 0 & 0 & 0.5 \\ 0 & 0 & 0.7 & 0.5 & 0 \end{bmatrix} \in [A, B].$$

Let $M_1$ and $M_2$ be the opinion of two experts one who works with 4 concepts and other with 5 concepts on the same problem.

One wishes to study the problem when the node $C_1$ is in the on state in the first fuzzy dynamical system $M_1$ and the node $D_2$ is in the on state in the second fuzzy dynamical system $M_2$.

Let X = (1 0 0 0) ∪ (0 1 0 0 0) be the state bivector where all other nodes are in the off state except the nodes $C_1$ and $D_2$ in the bivector $X = X_1 \cup X_2$.

Now to find the hidden bipattern of the state bivector $X = X_1 \cup X_2$ on the dynamical bisystem $M = M_1 \cup M_2$.

$$\begin{aligned} XM &= X_1 M_1 \cup X_2 M_2 \\ &= (1, 0.6, 0.4, 0) \cup (0.5, 1, 0, 0, 0.3) \\ &= Y_1 \cup Y_2 \text{ (say).} \end{aligned}$$

$$\text{Let } Y = Y_1 \cup Y_2.$$

$$
\begin{aligned}
YM &= (Y_1 \cup Y_2)\,(M_1 \cup M_2) \\
&= Y_1\,M_1 \cup Y_2\,M_2 \\
&= (1, 0.6, 0.4, 0.1) \cup (0.5, 1, 0.3, 0.3, 0.3) \\
&= Z_1 \cup Z_2 \\
&= Z \text{ (say).} \\
ZM &= Z_1\,M_1 \cup Z_2\,M_2 \\
&= (1, 0.6, 0.4, 0.1) \cup (0.5, 1, 0.3, 0.3, 0.3).
\end{aligned}
$$

Thus the hidden bipattern is a fixed bipoint. Thus we have now shown how the FCIBM system associated with a fuzzy interval mixed square bimatrix model works.

Now we have to cater to a very natural question, which is as follows: Suppose there are several experts working on the same problem but using different number of concepts in the FCIM model, then certainly the fuzzy interval mixed square bimatrix model will not serve anymore so we should try to build a general new model.

So we describe the construction of the most generalized FCIM model which we denote by FCItM, $t \geq 3$.

Suppose we have a set of experts working on the same problem P using the FCIBM model.

We have $p_1$ experts working with $n_1$ concepts on the problem P and $p_2$ experts working with $n_2$ concepts on the same problem P. ($n_1 \neq n_2$) and so on and $p_t$ experts working with $n_t$ concepts on the same problem P, using the FCIBM model with $p_i \neq p_j$ (if $i \neq j$). Then a single dynamical system associated with the FCIM will not serve the purpose, so we have to construct a special dynamical system to work with the problem simultaneously. To this end we do the following construction which we term as Fuzzy Cognitive Interval n-matrix model (FCInM).

Using the $p_1$ experts and $n_1$ concepts on the problem first we form the fuzzy interval $n_1 \times n_1$ square matrix associated with the FCIM. Let us denote this fuzzy interval $n_1 \times n_1$ square matrix associated with the FCIM on the problem P by $[A_1, B_1]$ i.e., we form the FCIM described earlier.

Likewise we form the $[A_i, B_i]$ fuzzy interval $n_i \times n_i$ square matrices associated with FCIM using the $p_i$ experts. We do this for the collection of all the set of $p_1, p_2 \ldots, p_t$ experts i.e., i = 1, 2, …, t.

Let $[A, B] = [A_1, B_1] \cup [A_2, B_2] \cup \ldots \cup [A_t, B_t]$ where [A, B] according to the definition, is a fuzzy interval mixed square matrices and each of the fuzzy interval square matrices $[A_i, B_i]$ contains almost $p_i + 4$ number of $n_i \times n_i$ square fuzzy matrices and $[A_i, B_i]$ is the associated fuzzy interval $n_i \times n_i$ square matrix of the FCIM for the problem P i.e., the FCIM model associated with P. This is true for i = 1, 2, …, t. Thus we see [A, B] is the collection of all fuzzy mixed square t matrices.

Any element M in the set [A, B] will be of the form $M = M_1 \cup M_2 \cup M_3 \cup \ldots \cup M_t$ where $M_i \in [A_i, B_i]$; i = 1, 2, …, t . We call [A, B] the fuzzy interval mixed square t matrices associated with the FCItMs of the problem P. When t = 1 we get the fuzzy interval matrix associated with the FCIM of the problem P. When t = 2 we get the fuzzy interval bimatrix associated with the FCIBM of the problem P and so on.

Now we shall sketch the working when t = 5. Let $[A, B] = [A_1\ B_1] \cup [A_2\ B_2] \cup \ldots \cup [A_5\ B_5]$ be a fuzzy interval mixed square 5 matrix of the FCI5M model on some problem P.

Let $M = M_1 \cup M_2 \cup \ldots \cup M_5$, if X is a state 5 vector whose resultant we are interested in finding out, we set $X = X_1 \cup X_2 \cup \ldots \cup X_5$.

$$\begin{aligned} XM &= (X_1 \cup \ldots \cup X_5)(M_1 \cup M_2 \cup \ldots \cup M_5) \\ &= X_1 M_1 \cup X_2 M_2 \cup \ldots \cup X_5 M_5 \\ &= Y_1 \cup Y_2 \cup \ldots \cup Y_5 \\ &= Y \text{ (say).} \end{aligned}$$

$$\begin{aligned} YM &= Y_1 M_1 \cup \ldots \cup Y_5 M_5 \\ &= Z_1 \cup Z_2 \cup \ldots \cup Z_5 \\ &= Z \text{ (say).} \end{aligned}$$

Then we find ZM; we proceed with the same process until we get a fixed point or a limit cycle of the FCI5M.

Now having explained how fuzzy interval mixed square m matrices are used in the FCItM model now we proceed onto study the application of fuzzy interval mixed rectangular t matrices.

## 3.4 FRIBM model and its Application in the Fuzzy Interval Mixed Rectangular Bimatrices

We have seen in this chapter how the FRIM model functions using fuzzy interval of rectangular matrices. We have also illustrated the working of the model using problem of student teacher relations.

It may so happen that suppose some experts study the problem, with one set of experts who choose $m_1$ concepts for the domain space and $n_1$ concepts for the range space and for the same problem another set of experts choose $m_2$ concepts for the domain space and $n_2$ concepts for the range space ($m_1 \neq m_2$), then how to make a model for it. When the number of concepts given by all experts; both in the range space and the domain space is the same we saw we could use the fuzzy interval rectangular matrix associated with the FRM model called the FRIM model which we have discussed earlier.

So when we have two sets of experts say $p_1$ and $p_2$ and the first set of experts work with $m_1$ concepts / attributes in the domain space and $n_1$ concepts / attributes in the range space and use the FRIM model for the problem. Using the other set of $p_2$ experts for the same problem, who use $m_2$ concepts or attributes in the domain space ($m_2 \neq m_1$ and $n_1 \neq n_2$) and $n_2$ concepts in the range space and adopt the FRIM model, we construct the fuzzy interval mixed rectangular bimatrix for the FRIBM-model by the following way.

First we form the fuzzy interval of $m_1 \times n_1$ rectangular matrices using the $p_1$ experts associated FRIM model and denote it by $[A_1, B_1]$. $[A_1 B_1]$ contains atmost $p_i + 4$ number fuzzy $m_1 \times n_1$ rectangular matrices.

We form the fuzzy interval $m_2 \times n_2$ matrices $[A_2, B_2]$ of the $p_2$ experts associated with the FRIM model for the same problem. Clearly the fuzzy interval matrix $[A_2, B_2]$ contains atmost $p_2$ + 4 numbers of, $m_2 \times n_2$ fuzzy matrices. Set [A, B] = $[A_1, B_1] \cup [A_2, B_2]$, we know [A, B] denotes the fuzzy interval of mixed rectangular bimatrices. Any element M in [A, B] will be of the from $M = M_1 \cup M_2$ where M is a fuzzy mixed rectangular bimatrix.

We call [A, B] to be the fuzzy interval of mixed rectangular bimatrices associated with the FRBM and denote it by FRIBM or is known as the fuzzy interval of mixed bimatrices associated with the FRBM model.

We just illustrate how this model works. If $X = X_1 \cup X_2$ is a state bivector whose resultant we are interested in finding, where $X_1 \in$ to the domain space of the FRIM of the first set of $p_1$ experts and $X_2$ belongs to the domain space of the FRIM of the second set of $p_2$ experts. The hidden pattern of the state bivector $X = X_1 \cup X_2$ on the dynamical bisystem $M = M_1 \cup M_2$ where $M \in [A, B]$, $M_1 \in [A_1, B_1]$ and $M_2 \in [A_2, B_2]$ of the FRIBM is determined as follows.

$$
\begin{aligned}
XM &= X_1M_1 \cup X_2M_2 \\
&= Y_1 \cup Y_2 \\
&= Y \text{ (say)}. \\
YM^T &= (Y_1 \cup Y_2)\ (M_1^T \cup M_2^T) \\
&= Y_1M_1^T \cup Y_2M_2^T \\
&= Z_1 \cup Z_2 \\
&= Z \text{ (say)}. \\
ZM &= (Z_1 \cup Z_2)\ (M_1 \cup M_2) \\
&= Z_1M_1 \cup Z_2M_2 \\
&= T_1 \cup T_2 \\
&= T \text{ (say)}. \\
TM^T &= (T_1 \cup T_2)\ (M_1^T \cup M_2^T) \\
&= T_1M_1^T \cup T_2M_2^T \\
&= S_1 \cup S_2 \\
&= S \text{ (say)}.
\end{aligned}
$$

We proceed on till we arrive at a fixed bipoint or a fixed bicycle.

Having seen the working of the FRIBM model using the associated fuzzy interval of mixed rectangular bimatrices of an FRIBM we now try to find out a general solution for the following natural question, to this effect we describe how a FRItM model is constructed.

Suppose we have a problem P, we have at hand t sets of $p_i$ experts; (i = 1, 2, …, t) and each of the $p_i$ experts work with $m_i$ concepts in the domain space and $n_i$ concepts in the range space i = 1, 2, …, t, $m_i \neq m_j$ (if $i \neq j$) if $t > 2$ certainly the fuzzy interval mixed rectangular bimatrix model associated with a FRIBM cannot be adopted.

So now we build a new model, which we choose to call as the fuzzy interval of mixed rectangular t matrix, associated with the FRItM model.

Let $[A_i, B_i]$ denote the fuzzy interval of $m_i \times n_i$ rectangular matrices given by $p_i$ the experts associated with the FRIM model. The number of fuzzy $m_i \times n_i$ matrices in $[A_i, B_i]$ will be atmost $p_i + 4$. Let us find such collection of fuzzy interval of $m_i \times n_i$ matrices for the t sets of experts, $p_1, \ldots, p_t$. i.e., for i = 1, 2, …, t. Let $[A, B] = [A_1, B_1] \cup [A_2 , B_2] \cup [A_3, B_3] \cup \ldots \cup [A_t, B_t]$.

Any element M in [A, B] is a fuzzy mixed rectangular t matrices i.e., $M = M_1 \cup M_2 \cup \ldots \cup M_t$. where $M_i \in [A_i, B_i]$; i = 1, 2, …, t further $A = A_1 \cup A_2 \cup \ldots \cup A_t$ and $B = B_1 \cup B_2 \cup \ldots \cup B_t$ will serve as the minimal and maximal fuzzy mixed rectangular t-matrices respectively;

Like wise $O = O_1 \cup O_2 \cup \ldots \cup O_t$ will be the optimal fuzzy mixed rectangular t matrix and $\overline{M} = M_1 \cup M_2 \cup \ldots \cup M_t$ will denote the average of the combined FRIM fuzzy interval of mixed rectangular t matrix.

Clearly A, B, O, $\overline{M} \in [A, B]$; we define [A, B] to be the fuzzy interval of rectangular t matrices associated with the FRIMs and denoted by FRItMs.

So if we have any state t-vector, $X = X_1 \cup X_2 \cup \ldots \cup X_t$ associated with the domain spaces of the FRItMs; the hidden pattern of X on any fuzzy rectangular t matrix $M = M_1 \cup M_2 \cup \ldots \cup M_t$ is given by

$$
\begin{aligned}
XM &= (X_1 \cup X_2 \cup \ldots \cup X_t)\,(M_1 \cup M_2 \cup \ldots \cup M_t) \\
&= X_1 M_1 \cup X_2 M_2 \cup \ldots \cup X_t M_t \\
&= Y_1 \cup Y_2 \cup \ldots \cup Y_t \\
&= Y \text{ (say).}
\end{aligned}
$$

$$
\begin{aligned}
YM^T &= (Y_1 \cup Y_2 \cup \ldots \cup Y_t)\;(M_1^T \cup M_2^T \cup \ldots \cup M_t^T) \\
&= Y_1 M_1^T \cup Y_2 M_2^T \cup \ldots \cup Y_t M_t^T \\
&= Z_1 \cup Z_2 \cup \ldots \cup Z_t \\
&= Z \text{ (say).} \\
ZM &= Z_1 M_1 \cup Z_2 M_2 \cup \ldots \cup Z_t M_t \text{ and so on.}
\end{aligned}
$$

We proceed on until we arrive at a fixed t point or a limit t-cycle.

## 3.5 Description of FAIM model and its Generalization

Next we proceed on to study how these concepts of fuzzy interval n-matrices, n = 1, 2, … ($n < \infty$) can be applied to the case of FAM models (Fuzzy Associative Maps model). We have just given a brief description of how the FAM model functions in Chapter 1.

Just we recall if A and B are fit vectors then

$$A \text{ o } M = B$$

where M is a fuzzy $n \times p$ matrix known as the FAM matrix. We know in most of the cases the fuzzy $n \times p$ matrix works as a bidirection system so we can have

$$B \text{ o } M^T = A'.$$

Now we can use both bidirectional FAMs or just FAMs. We can also have n = p. Now we would say how one can use fuzzy interval matrices in FAMs and the need for it come when several experts give their view on the same problem.

FAMs give only one way of working i.e., FAMs function at a time on a single experts opinion. Now suppose we have at hand say some t experts who give their opinion on the same problem using only n and p attributes i.e., using the same sets of attributes. How to find a means of comparison find the minimal or the least element, the maximal element, the optimal element and above all the common opinion derived from all the experts giving every one the same degree of participation; which is very important in case when a decision is to be taken by a set of t people who hold the same amount of share of or equal status or equal responsibility in analyzing the same problem.

Now we show how the fuzzy interval matrices can be used when we work with t experts but all of them work with the same number of attributes along the column p and along the rows n.

Let us assume in our problem we have just t distinct experts, p number of attributes along the column and n number of attributes along the row. So that the associated fuzzy matrix of the FAM for every expert will be only a $n \times p$ fuzzy rectangular matrix.

Now we take all the t fuzzy matrices of the FAM given by the t experts. Now we use the following method: Let $T_1, \ldots, T_t$ be the t number of $n \times p$ fuzzy matrices given by the t experts.

The matrices A and B are constructed as follows using the matrices $T_1, T_2, \ldots, T_t$. Let $T_r = \left(t^r_{ij}\right)$, $1 \le i \le n$ and $1 \le j \le p$, r = 1, 2, …, t.

First fill in the matrix $A = (a_{ij})$ as mentioned below,

$$a_{11} = \min \{t^1{}_{11}, t^2{}_{11}, \ldots , t^t{}_{11}\}$$
$$a_{12} = \min \{t^1{}_{12}, t^2{}_{12}, \ldots , t^t{}_{12}\}$$

thus $a_{ij} = \min \{t^1{}_{ij}, t^2{}_{ij}, \ldots , t^t{}_{ij}\}$;

$$t^r_{ij} \in \left(t^r_{ij}\right) = T_r ;$$

$1 \le r \le t$, $1 \le i \le n$ and $1 \le j \le p$.

Now $a = (a_{ij})$ will be called as the minimal element for the set of t matrices, $T_1, \ldots, T_r$. Similarly define $B = (b_{ij})$ as follows:

$$b_{ij} = \max \{t^1{}_{ij}, t^2{}_{ij}, \ldots, t^t{}_{ij}\}; t^r_{ij} \in T_r ; 1 \le r \le t.$$
$$b_{11} = \max \{t^1{}_{11}, t^2{}_{11}, \ldots, t^t{}_{11}\},$$

$$b_{12} = \max\{t^1_{12}, t^2_{12}, \ldots, t^t_{12}\},$$

and so on.

B = ($b_{ij}$) will be defined as the maximal of all the t-fuzzy n × p matrices. i.e., we will have for any ($t^r_{ij}$) = $T_r$, 1 ≤ r ≤ t, $a_{ij} \leq t^r_{ij} \leq b_{ij}$, 1 ≤ i ≤ n and 1 ≤ j ≤ p. Thus we will have an interval [A, B] to contain atmost t + 2 number fuzzy n × p matrices i.e., [A, B] = {A, $T_1$, …, $T_t$, B}.

Now we define the fuzzy n × p optimal matrix as

$$O = \frac{A+B}{2} = \left(o_{ij}\right) = \left(\frac{a_{ij}+b_{ij}}{2}\right),$$

1 ≤ i ≤ n and 1 ≤ j ≤ p. Similarly we form the average fuzzy n × p matrices

$$\overline{M} = \left(\overline{m}_{ij}\right) = \frac{T_1 + \ldots + T_t}{t} = \left(\frac{t^1_{ij} + \ldots + t^t_{ij}}{t}\right).$$

$$\text{i.e., each } \left(\overline{m}_{ij}\right) = \left(\frac{t^1_{ij} + \ldots + t^t_{ij}}{t}\right),$$

1 ≤ i ≤ n and i ≤ j ≤ p.

Thus now the fuzzy interval of matrices [A, B] contains both O and M. Thus [A, B] contains atmost t + 4 number of n × p fuzzy matrices. We say atmost for it may so happen that one or all of the matrices A, B, O or $\overline{M}$ may coincide with one of the $T_1$, $T_2$, …, $T_t$. Now we call this fuzzy interval of matrices [A, B] to be the fuzzy interval of matrices associated with the FAM model and denote it by FAIM. Now we will illustrate the uses of this fuzzy interval of n × p matrices associated with the FAIM model.

Suppose we take the expert r, we wish to find the effect of any fit vector on the system $T_r$. If $A_r$ is the fit vector we can calculate $B_r$ using the equation $A_r$ o $T_r$ = $B_r$, this is true for r = 1, 2, …, t. Now if we take the fuzzy optimal matrix O we can find the effect of the fit vector $A_r$ on O given by $A_r$ o O = $B^o_r$. Likewise we can find the effect or the resultant of the state vector on the fuzzy minimal matrix A and the maximal matrix B. i.e., $A_r$ o A = $B^A_r$ and $A_r$ o B = $B^B_r$. Now we can compare the four resultant fit vector $B_r$, $B_r^o$, $B_r^A$ and $B_r^B$ and arrive at the best solution / optimal solution.

In the same manner if $B_r$ is the fit vector under consideration we can calculate the resultant fit vector using the equation;

$$B_r \text{ o } T_i^T = A_r^i$$
$$B_r \text{ o } O^T = A_r^o$$
$$B_r \text{ o } A^T = A_r^A$$
and $$B_r \text{ o } B^T = A_r^B$$

for i = 1, 2, …, t; and choose the best solution out of it.

Now a very natural question may arise, suppose the t experts do not wish to take the same set of nodes / concepts for the rows and columns and they are divided into two groups say the first group containing $t_1$ experts and the second group containing $t_2$ experts and the first group works with $n_1$ row of nodes and $p_1$ column of nodes where as the second group of $t_2$ experts work with $n_2$ row nodes and $p_2$ column nodes; then certainly the fuzzy interval $n \times p$ matrices [A, B] associated with the FAIM using t experts opinion cannot be used. So when such a situation arises we give a new model to tackle the problem.

By using the method just explained above we for the $t_1$ experts find its associated fuzzy interval of $n_1 \times p_1$ matrices. Let $[A_1, B_1]$ denote the dynamical system of the FAIM model; and at the same time we find for the $t_2$ experts its associated $n_2 \times p_2$ matrices $[A_2, B_2]$, which gives the FAIM model associated with the $t_2$ experts.

Let us denote by $[A, B] = [A_1, B_1] \cup [A_2, B_2]$, the collection of all mixed rectangular fuzzy bimatrices i.e., any element M in [A, B] will be of the form $M = M_1 \cup M_2$ where $M_1 \in [A_1, B_1]$ and $M_2 \in [A_2, B_2]$. Further $A = A_1 \cup A_2$ and $B = B_1 \cup B_2$ which will be the minimal and maximal fuzzy bimatrices respectively of the fuzzy interval of bimatrices. The optimal fuzzy bimatrices in [A, B] is $O = O_1 \cup O_2$ where $O_1$ and $O_2$ are the optimal fuzzy matrices of the fuzzy interval of matrices $[A_1, B_1]$ and $[A_2, B_2]$ respectively.

The average bimatrix T in [A, B] will be the average of the fuzzy matrices of the $t_1$ experts say $T_1^1, T_2^1, ..., T_{t_1}^1$ denoted by $\overline{T}_1$ of the fuzzy interval of matrices [$A_1$, $B_1$] and the average of the fuzzy matrices of the $t_2$ experts say $T_1^2, T_2^2, ..., T_{t_2}^2$ denoted by $\overline{T}_2$ of the fuzzy interval of matrices [$A_2$, $B_2$]. Thus T = $\overline{T}_1 \cup \overline{T}_2$. We call the fuzzy interval of bimatrix to be the associated model given by two sets of experts $t_1$ and $t_2$ by FAIBM.

Now before we proceed on to define the notion of the fuzzy interval of n-matrices given by n sets of experts we will illustrate how a fuzzy interval of rectangular matrices [A, B] associated with a FAIM model given by t-experts functions.

We give the functioning of the FAIM model related to the psychological and social problems faced by rural women affected with HIV/AIDS. Here we have taken 7 attributes related to the rural women affected with HIV/AIDS, $W_1$, $W_2$, …, $W_7$ and the 10 attributes $R_1$, $R_2$, …, $R_{10}$; which has made a women HIV/AIDS infected. We have taken 3 experts opinion and let $M_1$, $M_2$ and $M_3$ denote the fuzzy vector matrix given by them. We have taken along the rows the problems related to the HIV/AIDS infected rural women and along the columns the cause of rural women becoming HIV/AIDS infected.

The fuzzy vector matrix $M_1$ given by the first expert is as follows:

$$M_1 = \begin{array}{c} \\ w_1 \\ w_2 \\ w_3 \\ w_4 \\ w_5 \\ w_6 \\ w_7 \end{array} \begin{array}{c} \begin{array}{cccccccccc} R_1 & R_2 & R_3 & R_4 & R_5 & R_6 & R_7 & R_8 & R_9 & R_{10} \end{array} \\ \begin{bmatrix} 0.9 & 0.8 & 0.7 & 0 & 0 & 0 & 0 & 0 & 0 & 0.7 \\ 0.5 & 0.8 & 0.6 & 0 & 0 & 0 & 0 & 0 & 0 & 0 \\ 0 & 0.3 & 0.6 & 0 & 0 & 0 & 0 & 0 & 0 & 0 \\ 0 & 0 & 0 & 0.6 & 0 & 0 & 0 & 0 & 0 & 0 \\ 0 & 0 & 0 & 0 & 0.9 & 0.6 & 0.7 & 0 & 0 & 0 \\ 0 & 0 & 0 & 0 & 0 & 0.7 & 0.5 & 0 & 0 & 0 \\ 0 & 0 & 0 & 0 & 0.6 & 0 & 0 & 0 & 0 & 0 \end{bmatrix} \end{array}.$$

For the same set of nodes the fuzzy vector matrix $M_2$ given by the second expert is follows:

$$M_2 = \begin{array}{c}
 \\ w_1 \\ w_2 \\ w_3 \\ w_4 \\ w_5 \\ w_6 \\ w_7
\end{array}
\begin{array}{c}
\begin{array}{cccccccccc} R_1 & R_2 & R_3 & R_4 & R_5 & R_6 & R_7 & R_8 & R_9 & R_{10} \end{array} \\
\begin{bmatrix}
0.8 & 0.7 & 0.6 & 0 & 0 & 0 & 0 & 0.2 & 0 & 0.4 \\
0.6 & 0.8 & 0.7 & 0 & 0 & 0 & 0 & 0 & 0 & 0 \\
0 & 0.3 & 0.6 & 0.1 & 0 & 0 & 0 & 0 & 0 & 0 \\
0 & 0 & 0 & 0.7 & 0 & 0 & 0 & 0 & 0.1 & 0 \\
0 & 0 & 0 & 0 & 0.8 & 0.7 & 0.7 & 0 & 0 & 0 \\
0 & 0 & 0 & 0 & 0 & 0.7 & 0.6 & 0 & 0 & 0 \\
0 & 0 & 0 & 0 & 0.6 & 0.3 & 0 & 0 & 0 & 0
\end{bmatrix}
\end{array}$$

Now we proceed on to give the fuzzy vector matrix $M_3$ related with the third expert.

$$M_3 = \begin{array}{c}
 \\ w_1 \\ w_2 \\ w_3 \\ w_4 \\ w_5 \\ w_6 \\ w_7
\end{array}
\begin{array}{c}
\begin{array}{cccccccccc} R_1 & R_2 & R_3 & R_4 & R_5 & R_6 & R_7 & R_8 & R_9 & R_{10} \end{array} \\
\begin{bmatrix}
0.6 & 0.7 & 0.8 & 0.1 & 0 & 0 & 0 & 0 & 0.1 & 0 \\
0.7 & 0.8 & 0.9 & 0 & 0 & 0 & 0 & 0 & 0 & 0 \\
0 & 0.3 & 0.6 & 0 & 0.1 & 0 & 0 & 0 & 0 & 0.1 \\
0 & 0 & 0 & 0.6 & 0 & 0 & 0.1 & 0 & 0 & 0 \\
0 & 0 & 0 & 0 & 0.7 & 0.6 & 0.8 & 0 & 0 & 0 \\
0 & 0 & 0 & 0 & 0 & 0.7 & 0.8 & 0 & 0 & 0 \\
0 & 0 & 0 & 0 & 0.5 & 0.2 & 0 & 0 & 0.1 & 0
\end{bmatrix}
\end{array}$$

Now using these three fuzzy matrices we have to form the interval of fuzzy matrices containing $M_1$, $M_2$ and $M_3$ for the FAIM model.

To form any fuzzy interval of matrices we need the minimal matrix A and the maximal matrix B. So that they are built in such a way that $M_1, M_2, M_3 \in [A, B]$.

To attain it in a unique way we take in $A = (a_{ij})$; we construct elements $a_{ij}$ by the following method; $a_{ij}$ to be the least of $m_{ij}^1, m_{ij}^2, m_{ij}^3$ where $(m_{ij}^k) \in M_k$, k = 1, 2, 3; i.e.,

$$M_1 = (m_{ij}^1),\ M_2 = (m_{ij}^2) \text{ and } M_3 = (m_{ij}^3),$$

true for $1 \leq i \leq 7$ and $1 \leq j \leq 10$. Thus A will serve as the minimal matrix and all $m_{ij}^{k}$; k = 1, 2, 3 will be such that $a_{ij} \leq m_{ij}^{k}$, k = 1, 2, 3; $1 \leq i < 7$ and $1 \leq j \leq 10$.

Now to construct the maximal matrix B, we follow the same procedure; if B = $(b_{ij})$ then choose

$$b_{ij} = \max\left\{m_{ij}^{1}, m_{ij}^{2}, m_{ij}^{3}\right\};$$

$1 \leq i \leq 7$ and $1 \leq j \leq 10$.

Thus B will be the maximal fuzzy matrix and $m_{ij}^{k} \leq b_{ij}$ for k = 1, 2, 3. Thus we see

$$a_{ij} \leq m^{k}_{ij} \leq b_{ij},$$

k = 1, 2, 3. [A, B] is the fuzzy interval of 7 × 10 matrix, such that $M_1$, $M_2$ and $M_3$ belong to the interval [A, B].

Now we form the optimal matrix O as

$$O = \left(\frac{a_{ij} + b_{ij}}{2}\right) = (o_{ij}),$$

$1 \leq i \leq 7$ and $1 \leq j \leq 10$.

The minimal fuzzy matrix A of the interval of matrices is given below,

$$A = \begin{array}{c} \\ w_1 \\ w_2 \\ w_3 \\ w_4 \\ w_5 \\ w_6 \\ w_7 \end{array}\begin{array}{c} \begin{array}{cccccccccc} R_1 & R_2 & R_3 & R_4 & R_5 & R_6 & R_7 & R_8 & R_9 & R_{10} \end{array} \\ \left[\begin{array}{cccccccccc} 0.6 & 0.7 & 0.6 & 0.1 & 0 & 0 & 0 & 0 & 0 & 0 \\ 0.5 & 0.8 & 0.6 & 0 & 0 & 0 & 0 & 0 & 0 & 0 \\ 0 & 0.3 & 0.6 & 0 & 0 & 0 & 0 & 0 & 0 & 0 \\ 0 & 0 & 0 & 0.6 & 0 & 0 & 0 & 0 & 0 & 0 \\ 0 & 0 & 0 & 0 & 0.6 & 0 & 0 & 0 & 0 & 0 \\ 0 & 0 & 0 & 0 & 0 & 0.7 & 0.5 & 0 & 0 & 0 \\ 0 & 0 & 0 & 0 & 0 & 0.5 & 0 & 0 & 0 & 0 \end{array}\right] \end{array}.$$

The maximal fuzzy matrix B of the fuzzy interval of matrix which denotes the maximum of all the entries $(m^{1}_{ij})$, $(m^{2}_{ij})$ and $(m^{3}_{ij})$ of the fuzzy matrices $M_1$, $M_2$ and $M_3$ is as follows:

$$B = \begin{array}{c} \\ w_1 \\ w_2 \\ w_3 \\ w_4 \\ w_5 \\ w_6 \\ w_7 \end{array} \begin{array}{c} \begin{array}{cccccccccc} R_1 & R_2 & R_3 & R_4 & R_5 & R_6 & R_7 & R_8 & R_9 & R_{10} \end{array} \\ \begin{bmatrix} 0.9 & 0.8 & 0.8 & 0 & 0 & 0 & 0 & 0.2 & 0.1 & 0.7 \\ 0.7 & 0.8 & 0.9 & 0 & 0 & 0 & 0 & 0 & 0 & 0 \\ 0 & 0.3 & 0.6 & 0.1 & 0.1 & 0 & 0 & 0 & 0 & 0.1 \\ 0 & 0 & 0 & 0.7 & 0 & 0 & 0.1 & 0 & 0.1 & 0 \\ 0 & 0 & 0 & 0 & 0.9 & 0.7 & 0.8 & 0 & 0 & 0 \\ 0 & 0 & 0 & 0 & 0 & 0.7 & 0.8 & 0 & 0 & 0 \\ 0 & 0 & 0 & 0 & 0.6 & 0.3 & 0 & 0 & 0.1 & 0 \end{bmatrix} \end{array}.$$

Now the optimal fuzzy matrix O in defined by

$$O = (o_{ij}) = \left(\frac{a_{ij} + b_{ij}}{2}\right).$$

Clearly $O \in [A, B]$.

The optimal fuzzy matrix O of the interval matrix [A, B] is as follows

$$\begin{array}{c} \\ w_1 \\ w_2 \\ w_3 \\ w_4 \\ w_5 \\ w_6 \\ w_7 \end{array} \begin{array}{c} \begin{array}{cccccccccc} R_1 & R_2 & R_3 & R_4 & R_5 & R_6 & R_7 & R_8 & R_9 & R_{10} \end{array} \\ \begin{bmatrix} 0.75 & 0.75 & 0.7 & 0 & 0 & 0 & 0 & 0.1 & 0.05 & 0.35 \\ 0.6 & 0.8 & 0.75 & 0 & 0 & 0 & 0 & 0 & 0 & 0 \\ 0 & 0.3 & 0.6 & 0.05 & 0.05 & 0 & 0 & 0 & 0 & 0.05 \\ 0 & 0 & 0 & 0.65 & 0 & 0 & 0.05 & 0 & 0.05 & 0 \\ 0 & 0 & 0 & 0 & 0.8 & 0.65 & 0.75 & 0 & 0 & 0 \\ 0 & 0 & 0 & 0 & 0 & 0.7 & 0.65 & 0 & 0 & 0 \\ 0 & 0 & 0 & 0 & 0.3 & 0.4 & 0 & 0 & 0.05 & 0 \end{bmatrix} \end{array}$$

Now

$$\overline{M} = \frac{M_1 + M_2 + M_3}{3}$$

$$= \frac{\left(m_{ij}^1 + m_{ij}^2 + m_{ij}^3\right)}{3}$$

$$= (m_{ij}),\ 1 \le i \le 7 \text{ and } 1 \le j \le 10.$$

The matrix $\overline{M}$ is as follows

| | $R_1$ | $R_2$ | $R_3$ | $R_4$ | $R_5$ | $R_6$ | $R_7$ | $R_8$ | $R_9$ | $R_{10}$ |
|---|---|---|---|---|---|---|---|---|---|---|
| $w_1$ | 0.77 | 0.71 | 0.7 | 0.03 | 0 | 0 | 0 | 0.07 | 0.03 | 0.33 |
| $w_2$ | 0.6 | 0.8 | 0.71 | 0 | 0 | 0 | 0 | 0 | 0 | 0 |
| $w_3$ | 0 | 0.3 | 0.6 | 0.03 | 0.03 | 0 | 0 | 0 | 0 | 0.03 |
| $w_4$ | 0 | 0 | 0 | 0.83 | 0 | 0 | 0.03 | 0 | 0.03 | 0 |
| $w_5$ | 0 | 0 | 0 | 0 | 0.8 | 0.63 | 0.77 | 0 | 0 | 0 |
| $w_6$ | 0 | 0 | 0 | 0 | 0 | 0.7 | 0.63 | 0 | 0 | 0 |
| $w_7$ | 0 | 0 | 0 | 0 | 0.57 | 0.37 | 0 | 0 | 0.03 | 0 |

Suppose we take the vector B as given by the experts as B = (0 1 1 0 0 0 0 0 1 0).

Using max. min, backward direction method of calculation in the FAIM fuzzy vector matrix given by the first expert, we get

$$B \text{ o } M_1^T = A. \text{ i.e.,}$$
$$= \text{max. min,. } (m_{ij}. b_j)$$
$$1 \leq i \leq 10$$
$$= (a_i).$$

Thus $A_1$ = (0.8, 0.8, 0.6, 0, 0, 0, 0), since 0.8 is the largest value of the fit vector in A and it is associated with the two nodes vide $W_1$ and $W_2$; child / widower child marriage etc. find its first place also the vulnerability of rural uneducated women find the same state as that of $W_1$. Further the second place is given to $W_3$ disease untreated till it is chronic or they are in last stages, all other states are in the off state.

Suppose we consider the resultant fit vector

$$A_1 = (0.8, 0.8, 0.6, 0, 0, 0, 0);$$

$$A \text{ o } M_1 = B \text{ where}$$
$$= \max (a_i, m_{ij});$$
$$1 \leq j \leq 7$$
$$= B_1$$
$$= (0.9, 0.8, 0.8, 0.6, 0.9, 0.7, 0.7, 0, 0, 0.7).$$

Thus we see the major cause for being affected by HIV/AIDS is $R_1$ and $R_5$ having the maximal value 0.9 closely followed by $R_2$ and $R_3$ whose value from the fit vector $B_1$ is 0.8.

Thus we see the maximum value corresponds to people in rural areas thinking female children to be a burden; so the sooner they get married off the better relief economically, $R_5$ states women in general do not suffer any guilt and fear for life.

Now the value 0.8 corresponding to $R_2$ and $R_3$ read as follows.

Poverty and bad habits of men are the major cause of women being HIV/AIDS victims.

The next largest value being 0.7 taken by the three attributes $R_6$, $R_7$ $R_{10}$. $R_6$ – majority of the women have not changed religion and developed faith in god after the disease; $R_7$ – no moral responsibility on the part of husbands and they infect their wives willfully; and $R_{10}$ – husband hide their disease from their family so the wife becomes HIV/AIDS affected.

The next largest value being 0.6 taken by the attribute $R_4$ – infected women are left uncared by relatives even their husbands. However we see only the nodes $R_8$ and $R_9$ according to this first experts system $M_1$ remains off. Frequent natural abortion / death of born infants and STD/VD infected husbands.

Now we study the second experts opinion using the same fit vector B = (0 1 1 0 0 0 0 0 1 0) using the system $M_2$.

$$
\begin{aligned}
B \circ M_2^T &= A_2 \\
&= \max_{1 \le i \le 10} \min (m_{ij}\, b_j) \\
&= (a_i).
\end{aligned}
$$

Thus $A_2$ = (0.7, 0.8, 0.6, 0.1, 0, 0, 0).

According to this expert $W_2$ takes the largest value followed by $W_1$ and this is closely followed by $W_3$. However the states $W_5$, $W_6$ and $W_7$ are in the off state and $W_4$ takes a very small value viz. 0.1.

We see both the experts agree on the factor that HIV/AIDS infection is due to lack of awareness as they are uneducated and are from rural areas.

Now

$$A_2 \text{ o } M_2 = B_2 = \max.\min (b_j, m_{ij}); \quad 1 \leq i \leq 10$$
$$= (0.7, 0.8, 0.7, 0.1, 0, 0, 0, 0.2, 0.1, 0.4).$$

Thus from this resultant fit vector we see $R_2$ take the maximum value viz. 0.8 closely followed by $R_1$ and $R_3$ which take the value 0.7. The next largest value being 0.4 however $R_5$, $R_6$ and $R_7$ remains in the off state and $R_4$ and $R_9$ take very small value viz. 0.1 and $R_8$ the value 0.2.

Now using the same fit vector B we work with the 3rd expert's fuzzy matrix $M_3$.

$$B \text{ o } M_3^T = A_3 = \max.\min (m_{ij}\ b_j); \quad 1 \leq i \leq 10$$
$$= (a_i).$$

Thus $A_3$ = (0.8, 0.9, 0.6, 0, 0, 0, 0.1) is the resultant fit vector. According to the 3rd expert $W_2$ takes the maximum value 0.9 i.e., unawareness among the rural uneducated women about HIV/AIDS. This is closely followed by the value 0.8 which is taken by $W_1$. The node $W_3$ takes the value 0.6.

However $W_4$, $W_5$ and $W_6$ remain in the off state and the nodes $W_7$ takes a very small value namely 0.1.

Now we study the effect of the fit vector $A_3$ on $M_3$.

$$A_3 \text{ o } M_3 = B_3 = \max.\min (a_i\ m_{ij}); \quad 1 \leq j \leq 10$$
$$= (0.7, 0.8, 0.9, 0.1, 0.1, 0.2, 0, 0, 0.1, 0.1).$$

Now the resultant state vector according to the 3rd expert $M_3$ takes the largest value 0.7 for $R_3$ closely followed by a value 0.8 taken by $R_2$ and $R_2$ is closely followed by $R_1$ which takes the value 0.7, $R_7$, $R_8$ takes 0 value i.e., they remain in the off state. $R_4$, $R_5$, $R_9$ and $R_{10}$ takes the value 0.1 and $R_6$ takes a slightly bigger value 0.2. Thus almost all the three experts model gives large values only for the 3 nodes $R_1$, $R_2$ and $R_3$.

Now we study the effect of the same fit vector B = (0 1 1 0 0 0 0 0 1 0) on the minimal matrix A.

$$
\begin{aligned}
B \text{ o } A^T &= A_m \\
&= \text{max. min } (m_{ij}\ b_j);\ 1 \le i \le 10 \\
&= (a_i).
\end{aligned}
$$

$$\text{i.e., } A = (0.7, 0.8, 0.6, 0, 0, 0, 0).$$

Now the resultant fit vector $A_m$ shows that the maximum value is taken by the node $W_2$ closely followed by 0.7 which is taken by $W_1$.

$W_1$ is closely followed by 0.6 by the node $W_3$. All the other four nodes $W_4$, $W_5$, $W_6$ and $W_7$ take the value 0.

Now we see the resultant of the fit vector $A_m$ on A.

$$
\begin{aligned}
A_m \text{ o } A &= B_m \\
&= \text{max. min } (a_j\ m_{ij}), \\
&\quad\ 1 \le j \le 10 \\
&= (0.6, 0.8, 0.6, 0, 0, 0, 0, 0, 0, 0).
\end{aligned}
$$

Thus we see the maximum value is taken by $R_2$, 0.8, and the $R_1$ and $R_3$ take the value 0.6. Thus we see if the minimal matrix A is used then we get all the nodes $R_4$, $R_5$, $R_6$, $R_7$, $R_8$, $R_9$ and $R_{10}$ remain in the off sate. Thus poverty and not the owners of property is the root cause for the spread of HIV/AIDS. However all other factors seems to be playing no significant role. The other contributing factor being $R_1$, female children is a burden to them and $R_3$ bad habits formed by the men is one of the major reason which is a root cause for women becoming HIV/AIDS infected.

Now for the same fit vector B = (0, 1, 1, 0, 0, 0, 0, 0, 1, 0) we study its effect on the maximal fuzzy matrix B.

$$
\begin{aligned}
B \text{ o } B^T &= A^{max} \\
&= \max \min (b_j, m_{ij}), \\
&\quad 1 \le j \le 10 \\
&= (0.8, 0.8, 0.6, 0.1, 0, 0, 0.1).
\end{aligned}
$$

Now we find that $W_1$ and $W_2$ takes the maximal value 0.8.

This is followed by 0.6 taken by $W_3$ and $W_5$ and $W_6$ remains in the off state. However $W_4$ and $W_7$ takes a very little value 0.1.

Now we find the effect of the fit vector $A^{max}$.

$$
\begin{aligned}
A^{max} \text{ o } B &= B^{max} \\
&= \max \min (b_j, m_{ij}); \\
&\quad 1 \le j \le 10 \\
&= (0.7, 0.8, 0.8, 0.1, 0.1, 0, 0, 0.2, 0.1, 0.7).
\end{aligned}
$$

The maximum value is taken by the nodes $R_2$ and $R_3$, viz. 0.8. Closely followed by 0.7, which is obtained by $R_1$ and $R_{10}$. However $R_6$ and $R_7$ remain in the off state and the least value is taken by $R_3$, $R_5$ and $R_9$, to be 0.1, followed closely by the next least value 0.2 taken by $R_8$.

Now we wish to study the effect of the same fit vector B = (0 1 1 0 0 0 0 0 1 0) on the optimal fuzzy matrix O,

$$
\begin{aligned}
B \text{ o } O^T &= A^o \\
&= \max \min (o_{ij}, b_j); \\
&\quad 1 \le j \le m \\
&= (0.75, 0.8, 0.6, 0.05, 0, 0, 0.05) \\
&= A^o.
\end{aligned}
$$

The maximum value is got by $W_2$, very closely followed by 0.75 taken by $W_1$. Now $W_3$ takes the value of 0.6 closely followed by $W_4$ and $W_7$ as 0.05.

Now we want to find the effect of the fit vector $A^o$.

$$A^o \text{ o } O = \max_{1 \le i \le 10} \min (a_i, o_{ij});$$
$$= (0.75, 0.8, 0.75, 0.05, 0.05, 0.04, 0.05, 0.1, 0.05, 0.35).$$

Thus form the resultant fit vector we see $R_2$ gets the maximum value 0.8 very closely followed by $R_1$ and $R_3$ taking the value 0.75. Now $R_{10}$ takes a value 0.35 none of the state vector is off so according to the optimal fuzzy matrix all the attributes get affected to a very negligible degree i.e., 0.05.

Now we find the effect of the fit vector B = (0 1 1 0 0 0 0 0 1 0) on the average fuzzy matrix of the fuzzy interval of matrices [A, B].

$$B \text{ o } \overline{M}^T = \max_{1 \le j \le 10} \min (\overline{M}_{ij}, b_j);$$
$$= \overline{A}$$
$$= (0.71, 0.8, 0.6, 0.03, 0, 0, 0.03).$$

Thus according to the average of the 3 experts opinion we see the largest value goes to $W_2$, the value being 0.8 the next value is take by the attribute is 0.71 taken by $W_1$, followed by $W_3$ which takes the value 0.6. However $W_5$ and $W_6$ remains in the off state, and very negligible value is taken by $W_4$ and $W_7$ viz. 0.03.

Now we will study the effect of fit vector $\overline{A}$ on $\overline{M}$.

$$\overline{A} \text{ o } \overline{M} = \overline{B}$$
$$= \max_{1 \le j \le 10} \min \left(\overline{a}_i, \overline{m}_{ij}\right);$$
$$= (b_j)$$
$$= (0.77, 0.8, 0.71, 0.03, 0.03, 0.03, 0.03, 0.07, 0.03, 0.03)$$
$$= \overline{B}.$$

Thus from the resultant fit vector $\overline{B}$ we see the maximum value is taken by $R_2$ viz. 0.8 very closely followed by $R_{10}$, 0.77 and still closely followed by 0.71. All nodes comes to on state by a very negligible value less than 0.03. However $R_{10}$ takes some value namely 0.33.

Now having worked with the same state vector B = (0, 1, 1, 0, 0, 0, 0, 0, 1, 0) we see how on these seven fuzzy matrices in the interval of fuzzy matrices [A, B] gives the fit resultant vectors.

Now for this fit vector B = (0 1 1 0 0 0 0 0 1 0) and the fuzzy vector matrix $M_1$, $M_2$, $M_3$. A, B, O and $\overline{M}$ we just write the resultant fit vectors and leave it as the work of the reader to make a comparison.

| | | |
|---|---|---|
| $A_1$ | = | (0.8, 0.8, 0.6, 0, 0, 0, 0); |
| $B_1$ | = | (0.9, 0.8, 0.8, 0.6, 0.9, 0.7, 0.7, 0, 0, 0.7) |
| $A_2$ | = | (0.7, 0.8, 0.6, 0.1, 0, 0, 0) |
| $B_2$ | = | (0.7, 0.8, 0.7, 0.1, 0, 0, 0, 0.2, 0.1, 0.4) |
| $A_3$ | = | (0.8, 0.9, 0.6, 0, 0, 0, 0.1); |
| $B_3$ | = | (0.7, 0.8, 0.9, 0.1, 0.1, 0.2, 0, 0, 0.1, 0.1) |
| $A^m$ | = | (0.7, 0.8, 0.6, 0, 0, 0, 0); |
| $B^m$ | = | (0.6, 0.8, 0.6, 0, 0, 0, 0, 0, 0, 0) |
| $A^{max}$ | = | (0.8, 0.8, 0.6, 0.1, 0, 0, 0.1); |
| $B^{max}$ | = | (0.7, 0.8, 0.8, 0.1, 0.1, 0, 0, 0.2, 0.1, 0.7) |
| $A^o$ | = | (0.75, 0.8, 0.6, 0.5, 0, 0, 0); |
| $B^o$ | = | (0.75, 0.8, 0.75, 0.5, 0.4, 0.5, 0.5, 0.1, 0.5, 0.35) |
| and | | |
| $\overline{A}$ | = | (0.71, 0.8, 0.6, 0.03, 0, 0, 0.03); |
| $\overline{B}$ | = | (0.77, 0.8, 0.71, 0.03, 0.03, 0.03, 0.03, 0.07, 0.03, 0.33). |

We have mainly given this model for illustrating how a fuzzy interval matrix associated with a FAIM system functions.

Now we proceed on to work with a problem P in which we have t set of experts giving opinion in the form of fuzzy vector matrices each of which is a $m_i \times n_i$ matrix, i = 1, 2, …, t and $m_i \neq m_j$ (if $i \neq j$). Further these $t_1$ experts give opinion as a fuzzy vector matrix model with $m_1 \times n_1$ associated fuzzy vector matrix, $t_2$ experts give a FAIM model with $m_2 \times n_2$ associated fuzzy vector matrix and $t_t$ experts give a FAIM model with a $m_t \times n_t$ associated fuzzy vector matrix.

Now if we want to make the FAIM multi expert opinion model as a single system what is to be done, so that this system functions as a whole.

We see FAIM's given by the experts are of varying order also number of experts in each set may or may not be the same. Under these circumstances how to form a physical model which will be unbiased and give the same degree of importance to one and all.

We use the fuzzy interval mixed rectangular t matrix to cater to this need. We shall briefly explain how the model functions. Let us assume we have N experts divided into t sets who want to analyze the problem P using the FAIM model. Each of the t sets contains $t_1, t_2, \ldots, t_p$ experts i.e., $N = t_1 + t_2 + \ldots + t_p$. $t_i$ of them work with a FAIM whose associated fuzzy vector matrix is a $m_i \times n_i$ matrix i.e., $i = 1, 2, \ldots, p$. Now we have to construct a model which will function as a common tool for all the p sets of experts i.e., the N experts who analyze the problem P using FAIM.

Now the $t_1$ set of experts work with a $m_1 \times n_1$ fuzzy vector matrix using a FAIM, $t_2$ set of experts work with a $m_2 \times n_2$ fuzzy vector matrix and so on.

Thus the $t_p$ experts work with the $m_p \times n_p$ fuzzy vector matrix. Clearly $m_i \neq m_j$ $(i \neq j)$ or $n_i \neq n_j$ (if $i \neq j$) i.e., none of the $m_i \times n_i$ matrix is the same as the $m_j \times n_j$ matrix if $i \neq j$.

Here onwards for simple working we assume the FAIM associated with the $t_j$ experts has $t_j$ number of $m_j \times n_j$ fuzzy vector matrices. This is true for $j = 1, 2, \ldots, p$. We want to form the fuzzy interval $m_j \times n_j$ matrix using these $t_j$ number of fuzzy vector matrices. Let the $m_j \times n_j$ fuzzy matrices associated with the $t_j$ experts be $T_1^j, T_2^j, \ldots, T_j^j$. Now the fuzzy interval matrix formed for these $t_j$ experts must contain the $t_j$ fuzzy matrices $T_1^j$, …, $T_{t_j}^j$. Let $T_r^j = \left(m_{qs}^{jr}\right)$; $1 \leq r \leq t_j$; $1 \leq q \leq m_j$ and $1 \leq s \leq n_j$. Now form the minimal matrix $A_j = (a^j_{qs})$; $1 \leq q \leq m_j$;

As $\left(a^{j}_{qs}\right) = \min\left\{m^{j1}_{qs}, m^{j2}_{qs}, ..., m^{jj}_{qs}\right\}$

for $1 \le q \le m_j$ and $1 \le s \le n_j$;

i.e., $\left(a^{j}_{11}\right) = \min\left\{m^{j1}_{11}, m^{j2}_{11}, ..., m^{jj}_{11}\right\}$

$$\left(a^{j}_{12}\right) = \min\left\{m^{j1}_{12}, m^{j2}_{12}, ..., m^{jj}_{12}\right\}$$

and so on.

Now $A_j = \left(a^{j}_{qs}\right)$ is such that $a^{j}_{qs} \le m^{jr}_{qs}$; $1 \le r \le t_j$.

Likewise we construct $B_j = \left(b^{j}_{qs}\right)$, the maximal fuzzy matrix for the interval of fuzzy matrices as

$$\left(b^{j}_{qs}\right) = \max\left\{m^{j1}_{qs}, m^{j2}_{qs}, ..., m^{jt_j}_{qs}\right\}$$

for $1 \le q \le m_j$ and $1 \le s \le n_j$. Thus we have

$$\left(b^{j}_{11}\right) = \max\left\{m^{j1}_{11}, m^{j2}_{11}, ..., m^{jt_j}_{11}\right\}$$

$$\left(b^{j}_{12}\right) = \max\left\{m^{j1}_{12}, m^{j2}_{12}, ..., m^{jt_j}_{12}\right\}$$

and so on. Now

$$B_j = \left(b^{j}_{qs}\right)$$

is the maximal fuzzy matrix of the interval and all $m^{jr}_{qs} \le b^{j}_{qs}$ for $1 \le q \le m_i$ and $1 \le s \le n_i$.

Now $[A_j, B_j]$ is the fuzzy interval of $m_j \times n_j$ fuzzy matrices containing the fuzzy matrices $T^{j}_{1}, T^{j}_{2}, \ldots, T^{j}_{t_j}$ and has $A_j$ to be the minimal element i.e., the minimal fuzzy matrix and $B_j$ to be the maximal fuzzy matrix of the fuzzy interval matrix $[A_j, B_j]$.

Let

$$O_j = \frac{\left(A_j + B_j\right)}{2}$$

be the optimal fuzzy matrix in the fuzzy interval of matrices. We form

$$\overline{T}^{j} = \frac{T^{j}_{1} + ... + T^{j}_{t_j}}{t_j}$$

which is called the average or mean fuzzy matrix of all the $t_j$ experts associated fuzzy matrices. Clearly $O_j$ and $\overline{T}^{j} \in [A_j, B_j]$.

The same procedure is repeated for all the sets of p experts; j = 1, 2, …, p when j = 1, $[A_1, B_1]$ will be the fuzzy interval of $m_1 \times n_1$ matrices associated with the FAIM of the $t_1$ experts, $[A_2, B_2]$ will be the fuzzy interval of $m_2 \times n_2$ matrices associated with the FAIM of the $t_2$ experts and so on. Thus we will have $[A_1, B_1]$, $[A_2, B_2]$, …, $[A_p, B_p]$.

Now set $[A, B] = [A_1, B_1] \cup [A_2, B_2] \cup ... \cup [A_p, B_p]$, clearly [A, B] contains the set of all fuzzy mixed rectangular p-matrices i.e., any element M in [A, B] will be of the form $M = M_1 \cup M_2 \cup \ldots \cup M_p$ where $M_i \in [A_i, B_i]$, i = 1, 2, …, p with

$$A = A_1 \cup A_2 \cup \ldots \cup A_p,$$
$$B = B_1 \cup B_2 \cup \ldots \cup B_p,$$
$$O = O_1 \cup O_2 \cup \ldots \cup O_p \text{ and}$$
$$\overline{T} = \overline{T}^1 \cup \overline{T}^2 \cup ... \cup \overline{T}^p.$$

Now A will be called as the minimal element or minimal fuzzy p-matrix of the fuzzy interval of mixed rectangular p-matrices. B the maximal fuzzy p-matrix, O the optimal fuzzy p-matrix and $\overline{T}$ the average fuzzy p-matrix of the fuzzy interval of p-matrices [A, B].

Now when we use this model we give equal importance to every expert and also to the minimal or maximal or optimal or average value which is highly dependent on each and every fuzzy matrix of the FAIM associated with the expert.

Thus $[A, B] = [A_1, B_1] \cup [A_2, B_2] \cup \ldots \cup [A_p, B_p]$ will be called as the associated fuzzy interval of mixed rectangular matrices of the $FAIN_pM$ or FAIpM model or dynamical system given by p-sets of experts. Now we will show the system when p = 2 i.e., only two sets of experts exists before we show this for any general p, p > 2.

The working of the fuzzy interval of bimatrices associated with the FAIBM of two sets of experts. Let P be the problem under investigation by two sets of experts say $p_1$ and $p_2$, where the set $p_1$ contains $m_1$ experts and the set $p_2$ contains $m_2$ experts. The $m_1$ experts agree to work with $r_1 \times s_1$ fuzzy matrices i.e., they have $r_1$ attributes along the rows and $s_1$ attributes along the columns and the $m_2$ experts work with $r_2$ attributes along the rows and $s_2$ attributes along the columns ($r_1 \neq r_2$) or ($s_1 \neq s_2$). Now using the method described for the fuzzy interval bimatrix associated with the FAIBM, let us assume $[A, B] = [A_1, B_1] \cup [A_2, B_2]$; where to show the working we assume $[A_1, B_1]$ the fuzzy interval matrix contains $3 \times 4$ fuzzy matrices with $A_1$ the fuzzy matrix with minimal entries and $B_1$ the fuzzy matrix with the maximal entries. Let $[A_2, B_2]$ be the fuzzy interval of matrices which contains $5 \times 3$ fuzzy matrices, both these matrices are associated with the same problem using the FAIM models i.e., one set of experts work with (3, 4) attributes and the other set of experts work with (5, 3) attributes. [A, B] is the fuzzy interval bimatrix model associated with the FAIM i.e., the FAIBM model.

Let $M = M_1 \cup M_2 \in [A, B]$ where $M_1 \in [A_1, B_1]$ and $M_2 \in [A_2, B_2]$ with the fuzzy matrix $M_1$ given by an expert in the first set and the fuzzy vector matrix $M_2$ corresponds to the opinion of the expert from the second set of experts.

$$M = \begin{bmatrix} 0.3 & 0 & 0.1 & 0 \\ 0 & 0.8 & 0.4 & 1 \\ 0.9 & 0 & 0 & 0.2 \end{bmatrix} \cup \begin{bmatrix} 0.8 & 0 & 0.6 \\ 0 & 0.7 & 0 \\ 0.2 & 0.5 & 0 \\ 0 & 0 & 0.7 \\ 0 & 0.3 & 0.9 \end{bmatrix}.$$

Suppose the expert wishes to work with the fit bivector B

$$\begin{aligned} B &= (1\ 0\ 1\ 0) \cup (0\ 0\ 1) \\ &= B_1 \cup B_2. \end{aligned}$$

To find the effect of the fit bivector B on the fuzzy vector bimatrix $M = M_1 \cup M_2$

$$\begin{aligned} B \text{ o } M^T &= (B_1 \cup B_2)\,(M_1 \cup M_2)^T \\ &= (B_1 \text{ o } M^T{}_1) \cup (B_2 \text{ o } M^T{}_2) \\ &= \max\min(m^1_{ij}, b^1_j) \cup \max\min(m^2_{ij}, b^2_j) \\ &= (a^1_j) \cup (a^2_j). \end{aligned}$$

$= (0.3, 0.4, 0.9) \cup (0.6, 0, 0, 0.7, 0.9)$
$= A_1 \cup A_2$
$= A.$

$$AoM = [A_1 \cup A_2] \text{ o } [M_1 \cup M_2]$$
$$= (A_1 \text{ o } M_1) \cup (A_2 \text{ o } M_2)$$
$$= \max\min(a_j^1, m_{ij}) \cup \max\min(a_j^2, m_{ij})$$
$$= (b_j^1) \cup (b_j^2)$$
$$= (0.9, 0.4, 0.4, 0.4) \cup (0.6, 0.3, 0.9).$$

Thus we get the resultant bivector to be a binary bipair bivector given by {(0.3, 0.4, 0.9) ∪ (0.6, 0, 0, 0.7, 0.9) = A and B = (0.9, 0.4, 0.4, 0.4, 0.4) ∪ (0.6, 0.3, 0.9)}. Thus we see one can arrive at a resultant for any fit bivector.

Now we will show how a fuzzy interval n matrices associated with a FAInM works.

Let $[A, B] = [A_1, B_1] \cup [A_2, B_2] \cup \ldots \cup [A_n, B_n]$ be a fuzzy interval of n matrices associated with a FAInM where we have n-sets of experts working on the same problem P. Let M be a fuzzy vector n-matrix from the fuzzy interval of n-matrices i.e., $M \in [A, B]$ and $M = M_1 \cup M_2 \cup \ldots \cup M_n$ where $M_i \in [A_i, B_i]$ is a fuzzy vector matrix of the FAIM model given by an expert; this is true for i = 1, 2, …, n. Suppose each $M_i$ in $[A_i, B_i]$ is a $m_i \times s_i$ matrix; i = 1, 2, …, n. Suppose we are given a fit n vector $X = X_1 \cup X_2 \cup \ldots \cup X_n$ where each $X_i$ is a $s_i \times 1$ row matrix, i = 1, 2, …, n. i.e.,

$$X = (x_{i_1}^1) \cup (x_{i_2}^2) \cup \ldots \cup (x_{i_n}^n);$$

$1 \le i_1 \le s_1, 1 \le i_2 \le s_2, \ldots, 1 \le i_n \le s_n$.

$$X \text{ o } M^T = (X_1 \cup X_2 \cup \ldots \cup X_n) \text{ o } (M_1 \cup M_2 \cup \ldots \cup M_n)^T$$
$$= (X_1 \text{ o } M^T{}_1) \cup (X_2 \text{ o } M^T{}_2) \cup \ldots \cup (X_n \text{ o } M^T{}_n)$$
$$= \max\{\min(m_{ij}^1, x_j^1)\} \cup \{\max\{\min(m_{ij}^2, x_j^2)\} \cup \ldots \cup \{\max\{\min(m_{ij}^n, x_j^n)\}$$
$$= (y_i^1) \cup (y_i^2) \cup \ldots \cup (y_i^n)$$
$$= Y_1 \cup Y_2 \cup \ldots \cup Y_n$$
$$= Y \text{ (say)};$$

where each $Y_i$ is a fit vector and Y is a fit n vector.

$$
\begin{aligned}
Y \text{ o } M &= (Y_1 \cup Y_2 \cup \ldots \cup Y_n) \text{ o } (M_1 \cup M_2 \cup \ldots \cup M_n) \\
&= (Y_1 \text{ o } M_1) \cup (Y_2 \text{ o } M_2) \cup \ldots \cup (Y_n \text{ o } M_n) \\
&= \cup \max \{\min \{y_i, m_{ij})\} \\
&= \max \{\min (y_j^1, m_{ij}^1)\} \cup \{\max (\min (y_j^2, m_{ij}^2)\} \\
&\quad \cup \ldots \cup \{\max \min (y_j^n, m_{ij}^n)\} \\
&= Z_1 \cup Z_2 \cup \ldots \cup Z_n \\
&= (Z_i^1) \cup (Z_i^2) \cup \ldots \cup (Z_i^n) .
\end{aligned}
$$

Thus for a given fit n vector X we get the pair of resultant fit n vectors given by

$$\{Y, Z\} = \{(Y_1 \cup Y_2 \cup \ldots \cup Y_n), (Z_1 \cup Z_2 \cup \ldots \cup Z_n)\}.$$

Having seen the use of fuzzy interval n matrices in the FAInM model now we show how fuzzy interval matrices can be used in the case of fuzzy relational equations. A brief description of the notions about fuzzy relational equations have been recalled in chapter one.

## 3.6 Use of Fuzzy Interval Matrices in Fuzzy Relational Equations Model

Now we show how the fuzzy interval matrices can be applied to fuzzy relational equations model. Now we can have several experts opinion expressed in the form of fuzzy relational equation matrices on a particular problem. We can work with each expert and find the solution. However so far we do not have a tool which can comprehend or consolidate the opinion together or make means for comparison. Now using fuzzy interval matrices we are able to achieve the following:

Suppose we have some n experts, giving their opinion about a problem P for these n experts opinion we can give the minimal expected resultant vector matrix, maximal resultant, the optimal solution and the average solution from this consolidated model.

For more literature about fuzzy relational equation models please refer [219].

Let us first briefly recall how the fuzzy relational equations system functions.

Suppose we have a situation in which given the matrices P and Q to determine the matrix R. R is determined by the equation P o Q = R.

$$\text{i.e., } P \circ Q = \max \min (p_{ij}, q_{jk}) = (r_{ik}) = R.$$

Clearly in this case the solution exists and is unique. The problem becomes more complicated when one of the two matrices on the left hand side of P o Q = R is unknown. In this case the solution is neither guranteed nor unique. Now we wish to state in these cases also the model which we are going to form will serve a good purpose.

For this we first define for fuzzy interval of matrices a composition rule; all the while when studying the fuzzy interval matrices associated with FCM or FRM or FAM models we only made use of the fuzzy interval of matrices. Now here to adopt the fuzzy relational equations define compositions when ever compatible among the fuzzy interval of matrices.

Suppose we have 3 fuzzy interval of matrices say [A, B], [X, Y] and [R, S] where for every matrix M in [A B] and N in [X, Y] we have

$$M \circ N = P = \max \min (m_{ij}, n_{jk}) = (p_{ik}) \in [R, S],$$

then we say the fuzzy intervals of matrices [A, B] when composed with the fuzzy interval of matrices [X, Y] give a fuzzy interval of resultant matrix [R, S] denoted by

$$[A, B] \circ [X, Y] = [R, S].$$

Now [A, B] can be a row matrix or a column matrix. Likewise [X, Y] and [R, S], so we can also speak of systems of linear equations and their solutions in terms of fuzzy interval of matrices. We may have several solutions.

This method will help us to find the optimal solution or so. Thus when we say a set of 3 fuzzy interval matrices are compatible under the composition, we only mean the following if [A, B], [X, Y] and [S, T] are 3 fuzzy interval matrices. We have for every $M \in [A, B]$ and $N \in [X, Y]$ such that

$$\begin{aligned} M \text{ o } N &= \max \min. \{m_{ij}, n_{jk}\} \\ &= (r_{ik}) = R \text{ is in } [S,T]. \end{aligned}$$

We denote this symbolically by [A, B] o [X, Y] = [S, T].

Now this rule will be used by us while defining fuzzy relational equations by several experts.

To this end we define the problem and the functioning of the model. Let P be the problem under study X and Y, be two sets such that to each element of X two or more elements are associated in Y with degree of membership from the interval [0, 1]. Thus R[X, Y] denotes this relation, R[X, Y] will be a fuzzy matrix which we choose to call as the fuzzy membership matrix; if X has m elements and Y has n elements the fuzzy membership matrix R [X, Y] will be a $m \times n$ fuzzy matrix for, its entries are form [0, 1].

Now we may have several experts to give their membership fuzzy matrices R [X, Y] and this R [X, Y] which is a $m \times n$ fuzzy matrix will vary from person to person or from expert to expert. Now how to transform these membership of fuzzy matrices into an interval of fuzzy matrices. Let some p experts give the fuzzy membership matrix relating the concepts X and Y given by $P_1, P_2, \ldots, P_p$; we have to form the interval of fuzzy membership matrices [A, B] in which $P_1, P_2, \ldots, P_p$ are elements i.e., $P_1, P_2, \ldots, P_p \in [A, B]$.

Let $P_i = \{(p^i_{rk})\}$, i = 1, 2, …, p, $1 \le r \le m$, $1 \le k \le n$. Choose elements for the fuzzy matrix $A = (a_{rk})$ as follows $1 \le r \le m$, $1 \le k \le n$.

$$a_{11} = \min \{ p^1_{11}, p^2_{11}, \ldots, p^p_{11} \},$$

clearly $a_{11} \le p^i_{11}$; $1 \le i \le p$,

$$a_{12} = \min \{ p^1_{12}, p^2_{12}, \ldots, p^p_{12} \},$$

clearly $a_{12} \le p_{12}^{i}$ ; $1 \le i \le p$, so on.

$a_{ij} = \min \{ p_{ij}^{1}, p_{ij}^{2}, ..., p_{ij}^{p} \}$,

clearly $a_{ij} \le p_{ij}^{t}$ ; $1 \le t \le p$, $1 \le i \le m$; $1 \le j \le n$.

A is called the minimal fuzzy matrix of the fuzzy membership matrices $P_1, P_2, \ldots, P_p$.

Now we define the maximal fuzzy matrix $B = (b_{ij})$ as follows:

$b_{11} = \max \{ p_{11}^{1}, p_{11}^{2}, ..., p_{11}^{p} \}$,

so that $b_{11} \ge p_{11}^{i}$ ; $1 \le i \le p$,

$b_{12} = \max \{ p_{12}^{1}, p_{12}^{2}, ..., p_{12}^{p} \}$,

so that $b_{12} \ge p_{12}^{i}$ , $1 \le i \le p$ and so on.

$b_{ij} = \max \{ p_{ij}^{1}, p_{ij}^{2}, ..., p_{ij}^{p} \}$,

so that $b_{ij} \ge p_{ij}^{t}$ , $1 \le t \le p$.

Thus $B = (b_{ij})$ serves as the fuzzy maximal matrices for the fuzzy membership matrices $P_1, P_2, \ldots, P_p$. Hence [A, B] is a fuzzy interval of membership of matrices containing $P_1, P_2, \ldots, P_p$. We can call A to give the least of all membership and B to give the greatest (maximal) of all memberships. Now we define the optimal membership fuzzy matrix of this fuzzy membership interval of matrix as

$$O = \frac{A + B}{2}.$$

That is

$$O = (o_{ij}) = \left( \frac{(a_{ij}) + (b_{ij})}{2} \right),$$

$1 \le i \le m$, $1 \le j \le n$. Clearly $O \in [A, B]$. We define the average fuzzy memberships matrix of the fuzzy interval of matrices [A, B] to be $\overline{P}$ where

$$\overline{P} = \frac{P_1 + ... + P_p}{p} = \frac{\left(p_{ij}^{1}\right) + ... + \left(p_{ij}^{p}\right)}{p} ;$$

$1 \le i \le m$ and $1 \le j \le n$.

Clearly $\overline{P} \in [A, B]$.

Thus the fuzzy interval of membership matrices with p experts has atmost p + 4 matrices.

Now for this model to function as a fuzzy relational equation model we have to define only the rule of composition of fuzzy interval of matrices.

[A, B] is a fuzzy interval of m × n membership matrices of p experts. Let R be the expected values of the resultant and Q the probable values with which the expert works.

Suppose the p experts choose to give some totally q number of probable values. Clearly the q probable values which will be fuzzy matrices, we will be forming a fuzzy interval of probable matrix solutions. For these we as before find the fuzzy interval of matrices denoted by [X, Y]; thus [X, Y] will be called as the probable fuzzy interval of matrices with X, the calculated minimal probable value fuzzy matrix and Y the maximal probable value fuzzy matrix. The optimal probable value of fuzzy matrix

$$O_p = \frac{X+Y}{2}.$$

The average fuzzy matrix will be sum of the q matrices divided by q which we will denote by Q. Thus [X, Y] the fuzzy probable interval of matrices has atmost q + 4 fuzzy probable matrices.

Now the fuzzy interval of resultant matrices [R, S] is formed as follows for the fuzzy relational equation model. Now take any fuzzy membership matrix P from [A, B] and take the fuzzy probable matrix Q from [X, Y]

$$P \text{ o } Q = \max \min_{j} (p_{ij}\ q_{jk})$$

$$= (r_{ik}) = R.$$

**How to form the fuzzy interval of resultant matrices. Here it is pertinent to mention that even if some expected resultant matrices i.e., they may be searching for such solution. But, how is the interval of matrices formed.**

For the minimal fuzzy matrix of the interval of resultant fuzzy matrices we form R = A o X and for the maximal fuzzy matrix of the interval of resultant fuzzy matrices we form S = B o Y. The optimal resultant solution is given by O o $O_p = O_R$, $O_R$ will be known as the optimal fuzzy resultant matrix of the interval [R, S]. Now for this we cannot have the average. For [R, S] may be an empty fuzzy interval only we have to find the resultants and fill it, or at times it may contain several of the expected fuzzy matrices by the experts or the problem posers. Another situation may be some times solutions of these equations may not be very acceptable. Thus [A, B] o [X, Y] = [R, S] represents the fuzzy dynamical system of fuzzy interval of matrices of the fuzzy relational equations $o_I$ denotes symbolically that we have composition of interval of matrices. i.e., [A, B] $o_I$ [X, Y] = [P o Q] P ∈ [A, B] and Q [X, Y] with A o X = R and B o Y = S]; we denote this symobollically by FRIE (Fuzzy Relational Interval Equations).

What are the advantages of formulating such a model of fuzzy relational interval equations. Here we have made an assumption P o Q = R' where P and Q are given (P ∈ [A, B], Q ∈ [X, Y]), so P o Q = R' always exists for any given P which is the fuzzy membership matrix given by an expert we can always find various elements in [R, S] say like the optimal value probable fuzzy matrix $O_p$, or X the minimal probable fuzzy matrix or Y the maximal probable fuzzy matrix or any other probable fuzzy matrix which may give us a solution; at least close to the desired solution from the interval of resultant fuzzy matrices [R, S].

On the other hand one can be given the probable fuzzy matrix Q from [X, Y] and made to choose the membership fuzzy matrix P from [A, B] so that Q o P gives the best possible solution. It can be even the minimal membership fuzzy matrix A or the maximal fuzzy membership matrix B or the opitimal fuzzy membership matrix O or the average experts opinion fuzzy membership matrix $\overline{P}$ and find Q o $\overline{P}$.

Now we have worked out or discussed with a very nice case given P and Q finding the resultant which always exists and is unique; where P is the membership fuzzy matrix and Q is the probable fuzzy matrix.

It may so happen that we may have the fuzzy interval matrix composition equation.

$$P \text{ o } Q = R' \; ; P \in [A, B], Q \in [X, Y] \text{ and } R' \in [R, S].$$

Suppose any two equations is known i.e., the expected value and the probable value how to find P the fuzzy matrix which is essential to build the structure. The solutions at times may not exists at times even if the solutions exists we may not have a unique solution.

In such precarious situations or which we may choose to call as critical situation we can make use of the optimal or average or minimal or maximal fuzzy matrices of the interval to get a solution. This is one of the major advantages in working with fuzzy intervals (FRIM models) than merely with just fuzzy matrices.

Now we will illustrate the problem by an example so that the reader becomes familiar with working of the fuzzy interval model.

Suppose we want to study the effect of globalization on the silk weavers of Kancheepuram who are bonded labourers using the fuzzy interval of membership, i.e., using the fuzzy interval i.e., the FRIM model.

For more about the description of the problem please refer [202, 206, 219]

Suppose we take four attributes related with the owners of bonded labourers as;

$O_1$ – Globalization / introduction of modern methods using machines in textile industries
$O_2$ – Profit or no loss to the industries
$O_3$ – Availability of raw goods or materials
$O_4$ – Demand for finished goods / products.

The problems related with the bonded labourers are taken under the six major heads;

$B_1$ – No basic education so has no knowledge of any other trade which has made them bonded and lead a life of poverty and penury.

$B_2$ – The advent of power looms and globalization has made them still poorer with no good pay

$B_3$ – Salary they earn per month is insufficient to maintain the family.

$B_4$ – No possibility of any savings so they become more and more bonded by borrowing from the owners they live in permanent depth.

$B_5$ – Even if govt. interferes and frees they do not find and cannot find any proper work for the govt. only frees them but never gives them any alternative employment.

$B_6$ – The hours they work in a day is more than 8 hours or so.

The expert is free to set up any limit for both these attributes related with the bonded labourers and their owners.

Now we give just 3 experts opinion on this problem which are transformed into the fuzzy relational equation. The fuzzy relational equation matrix $P_1$ related with the first expert.

$$P_1 = \begin{array}{c} \\ B_1 \\ B_2 \\ B_3 \\ B_4 \\ B_5 \\ B_6 \end{array} \begin{array}{c} \begin{array}{cccc} O_1 & O_2 & O_3 & O_4 \end{array} \\ \begin{bmatrix} 0.8 & 0 & 0 & 0 \\ 0.8 & 0.3 & 0.3 & 0 \\ 0.1 & 0.2 & 0.3 & 0.4 \\ 0 & 0.1 & 0.1 & 0.1 \\ 0.8 & 0.1 & 0.2 & 0.4 \\ 0.2 & 0.4 & 0.4 & 0.9 \end{bmatrix} \end{array}.$$

The fuzzy relational equation matrix $P_2$ given by the second expert is as follows:

$$P_2 = \begin{array}{c} \\ B_1 \\ B_2 \\ B_3 \\ B_4 \\ B_5 \\ B_6 \end{array} \begin{array}{c} \begin{array}{cccc} O_1 & O_2 & O_3 & O_4 \end{array} \\ \begin{bmatrix} 0.7 & 0.1 & 0 & 0 \\ 0.9 & 0.2 & 0.3 & 0 \\ 0 & 0.1 & 0.2 & 0.3 \\ 0 & 0 & 0.1 & 0.1 \\ 0.9 & 0 & 0.1 & 0.4 \\ 0.1 & 0.2 & 0.4 & 0.7 \end{bmatrix} \end{array}.$$

The fuzzy relational equation matrix given by the 3rd expert be denoted by $P_3$.

$$P_3 = \begin{array}{c} \\ B_1 \\ B_2 \\ B_3 \\ B_4 \\ B_5 \\ B_6 \end{array} \begin{array}{c} \begin{array}{cccc} O_1 & O_2 & O_3 & O_4 \end{array} \\ \begin{bmatrix} 0.9 & 0 & 0 & 0 \\ 0.5 & 0.3 & 0.4 & 0.1 \\ 0.2 & 0.2 & 0.2 & 0.3 \\ 0 & 0 & 0.1 & 0.2 \\ 0.7 & 0.2 & 0.2 & 0.4 \\ 0.2 & 0.3 & 0.3 & 0.8 \end{bmatrix} \end{array}.$$

These three fuzzy relational matrix can be realized as the fuzzy membership matrix given by the 3 experts.

Now we will find the fuzzy interval of fuzzy membership matrices using the three fuzzy membership matrices $P_1$, $P_2$ and $P_3$. Let

$$P_1 = (p^1_{ij}),$$

$$P_2 = (p^2_{ij})$$

and $$P_3 = (p^3_{ij}); \quad 1 \leq i \leq 6 \text{ and } 1 \leq j \leq 4.$$

Let $A = (a_{ij})$ be the minimal fuzzy membership matrix formed using the fuzzy membership matrices $P_1$, $P_2$ and $P_3$.

Define

$$a_{ij} = \min \{ p^1_{ij}, p^2_{ij}, p^3_{ij} \};$$

$1 \leq i \leq 6$ and $1 \leq j \leq 4$.
i.e.,

$$a_{11} = \min\{0.8, 0.7, 0.9\} = 0.7$$
$$a_{12} = \min\{0, 0.1, 0\} = 0$$
$$a_{13} = \min\{0, 0, 0\} = 0$$
$$a_{14} = \min\{0, 0, 0\} = 0$$

so on

$$a_{64} = \min\{0.9, 0.7, 0.8\} = 0.7.$$

Thus we see $A = (a_{ij})$; $1 \leq i \leq 6$ and $1 \leq j \leq 4$ is such that $a_{ij} \leq p_{ij}^{t}$; $t = 1, 2, 3$, with $1 \leq i \leq 6$ and $1 \leq j \leq 3$. So that

$$A = \begin{array}{c} \\ B_1 \\ B_2 \\ B_3 \\ B_4 \\ B_5 \\ B_6 \end{array} \begin{array}{c} \begin{array}{cccc} O_1 & O_2 & O_3 & O_4 \end{array} \\ \begin{bmatrix} 0.7 & 0 & 0 & 0 \\ 0.5 & 0.2 & 0.3 & 0 \\ 0 & 0.1 & 0.2 & 0.3 \\ 0 & 0 & 0.1 & 0.1 \\ 0.7 & 0 & 0.1 & 0.4 \\ 0.1 & 0.2 & 0.3 & 0.7 \end{bmatrix} \end{array}.$$

A is the minimal fuzzy membership matrix.

Now we calculate $B = (b_{ij})$, the maximal fuzzy membership matrix in the following way

$$b_{ij} = \max\{p_{ij}^{1}, p_{ij}^{2}, p_{ij}^{3}\}$$

for $1 \leq i \leq 6$ and $1 \leq j \leq 3$.

$$b_{11} = \max\{p_{11}^{1}, p_{11}^{2}, p_{11}^{3}\}$$
$$= \max\{0.8, 0.7, 0.9\} = 0.9.$$
$$b_{12} = \max\{0, 0.1, 0\} = 0.1;$$

and so on. Thus $b_{ij} \geq p_{ij}^{t}$; $t = 1, 2, 3$; $1 \leq i \leq 6$ and $1 \leq j \leq 3$. Thus we have $B = (b_{ij})$.

$$B = \begin{array}{c} \\ B_1 \\ B_2 \\ B_3 \\ B_4 \\ B_5 \\ B_6 \end{array} \begin{array}{c} \begin{array}{cccc} O_1 & O_2 & O_3 & O_4 \end{array} \\ \begin{bmatrix} 0.9 & 0.1 & 0 & 0 \\ 0.9 & 0.3 & 0.4 & 0.1 \\ 0.2 & 0.2 & 0.3 & 0.4 \\ 0 & 0.1 & 0.1 & 0.2 \\ 0.9 & 0.2 & 0.2 & 0.4 \\ 0.2 & 0.4 & 0.4 & 0.9 \end{bmatrix} \end{array}$$

is the maximal fuzzy membership matrix. The optimal fuzzy membership matrix O is given by

$$O = \frac{A+B}{2},$$

$$= \left[\frac{(a_{ij})+(b_{ij})}{2}\right],$$

$1 \le i \le 6$, $1 \le j \le 4$.
i.e.,

$$O = \begin{array}{c} \\ B_1 \\ B_2 \\ B_3 \\ B_4 \\ B_5 \\ B_6 \end{array} \begin{array}{c} \begin{array}{cccc} O_1 & O_2 & O_3 & O_4 \end{array} \\ \begin{bmatrix} 0.85 & 0.05 & 0 & 0 \\ 0.7 & 0.25 & 0.35 & 0.05 \\ 0.1 & 0.15 & 0.25 & 0.35 \\ 0 & 0.05 & 0.1 & 0.15 \\ 0.8 & 0.1 & 0.15 & 0.4 \\ 0.15 & 0.3 & 0.35 & 0.8 \end{bmatrix} \end{array}$$

is the matrix such that $a_{ij} \le o_{ij} \le b_{ij}$; $1 \le i \le 6$, $1 \le j \le 3$.
Now we form the average of the three fuzzy membership matrix $P_1$, $P_2$ and $P_3$.

$$\overline{P} = \frac{P_1 + P_2 + P_3}{3}.$$

$$\overline{P} = \begin{array}{c} \\ B_1 \\ B_2 \\ B_3 \\ B_4 \\ B_5 \\ B_6 \end{array} \begin{array}{c} \begin{array}{cccc} O_1 & O_2 & O_3 & O_4 \end{array} \\ \begin{bmatrix} 0.8 & 0.03 & 0 & 0 \\ 0.73 & 0.27 & 0.33 & 0.03 \\ 0.1 & 0.17 & 0.23 & 0.33 \\ 0 & 0.03 & 0.1 & 0.13 \\ 0.8 & 0.1 & 0.17 & 0.4 \\ 0.17 & 0.3 & 0.37 & 0.8 \end{bmatrix} \end{array}$$

$\overline{P} \in [A, B]$ is a fuzzy relational average of membership matrix. Now let X to be views of expert about the industry. We may take all $1 \times 4$ matrices with entries from the fuzzy interval [0.3, 1] i.e., they want atleast to run the factory or the industry with least loss so they do not accept values from the fuzzy interval (0, 0.3).

Thus if [T, U] denotes the set of all $1 \times 4$ matrices with the entries from the interval [0.3, 1] where T = [0.3, 0.3, 0.3, 0.3] and U = [1, 1, 1,1] are the weightages a company expects on the state vector say $R = (r_1, r_2, r_3, r_4)$ the $0.3 \leq r_i \leq 1$; i = 1, 2, 3, 4.

Now we are left with the job of finding the concepts or attributes related with the bonded labourers.

$P \text{ o } Q^T = R$ where 'o' is the max min operation.

Suppose Q = (0.8, 0.6, 0.7, 0.5) is the given vector we shall find the effect of Q on the bonded labourers.

$P_1 \text{ o } Q^T$ = (0.8, 0.8, 0.4, 0.1, 0.8, 0.5),
$P_2 \text{ o } Q^T$ = (0.7, 0.8, 0.3, 0.1, 0.4, 0.5),
$P_3 \text{ o } Q^T$ = (0.8, 0.5, 0.3, 0.2, 0.7, 0.5),
$A \text{ o } Q^T$ = (0.7, 0.5, 0.3, 0.1, 0.7, 0.5),
$B \text{ o } Q^T$ = (0.8, 0.8, 0.4, 0.2, 0.8, 0.5),
$O \text{ o } Q^T$ = (0.8, 0.7, 0.35, 0.15, 0.8, 0.5) and
$\overline{P} \text{ o } Q^T$ = (0.8, 0.73, 0.33, 0.13, 0.8, 0.5).

From the resultant state vector we see when the company / factory chooses to introduce up to 0.8 of modernized textile machines and has a 0.7 degree of free flow or availability of raw goods and if they want their profit to be atleast 0.6 degree and demand for finished goods even though 0.5 degree.

We see the bonded labourers suffer a acute poverty for the main reason they have no knowledge of any other work and hence work for a paltry salary, the power loom and other modern machinery has made them from bad to worse for $B_2$ = 0.8 in majority of the resultant vectors. It is very clear that the government interference has not helped them in any way for it takes a value 0.8, 4 times out of seven times. 0.7; 2 times out of seven times and only once it has taken the value 0.4. Further all entries in $B_6$ is 0.5 which clearly show in all cases majority of the bonded labourers work for more than 10 hours of a day with no holidays with pay. Clearly earning is not mediocre refer [219]. They have no savings and are in permanent debt.

Suppose some experts, wishes to work with a feed back from the bonded labourers and he takes an element, R = [0.2, 0.1, 0.7, 0.8, 0.2, 0.1] we will study the effect of R on the 7 membership matrices from the fuzzy interval of memberships matrices.

Clearly $P^T$ o R = $Q^T$.

$P_1^T$ o R = (0.2, 0.2, 0.3, 0.4),

$P_2^T$ o R = (0.2, 0.2, 0.2, 0.3),

$P_3^T$ o R = (0.2, 0.2, 0.2, 0.3),

$A^T$ o R = (0.2, 0.1, 0.3, 0.3),

$B^T$ o R = (0.2, 0.2, 0.3, 0.4),

$O^T$ o R = (0.2, 0.15, 0.25, 0.35)

and

$\overline{P}^T$ o R = (0.2, 0.17, 0.23, 0.33).

If the bonded labourers have knowledge of other skilled / unskilled labour and if the modernized machinery is not introduced and their earning is mediocre i.e. $\geq$ 0.5 and they have no debt and are free of acute poverty (they can afford a square meal a day) and government after interference help to accommodate them with some employment and they work with 8 hours maximum, holidays with pays; we see from the resultant globalization can not have any impact on the owners of these industries. They have to run their factory only with total loss they will also suffer shortage of raw materials as a craze for the rich and the affordable to buy machine made goods than hand made goods so that the raw materials will not be given to them and they cannot certainly balance the demand with supply. So if the worker has to enjoy as given by the membership matrix R certainly all the companies has to suffer.

It is left for the reader to arrive at an optimum value. As this is only an illustration how the system works we do not go very deep into the solution of the problem. We now try means to find solution for P o Q = R, when Q and R are given i.e. how best we can formulate a feasible membership matrix.

We construct a model for the preference of passengers using the real data got from a transport public corporation form Tamil Nadu [205, 219].

The data supplied to us is;

Hour ending
Q : 6, 7, 8, 9, 10, 11, 12, 13, 14, 15, 16, 17, 18. 19, 20, 21, 22.
Passengers per hour
R : 96, 71, 222, 269, 300, 220, 241, 265, 249, 114, 381, 288, 356, 189, 376, 182, 67.

To convert this data into a Fuzzy Relation Interval Equation (FRIE) model. We partition Q into 3 elements each. Since all concepts are to be realized as fuzzy concepts, we at first stage make the entries of Q ad R to lie between 0 and 1. This is done by multiplying all elements of Q by $10^{-2}$ and the elements of R by $10^{-4}$ respectively. The fuzzy interval of matrix of Q contains

$$Q_1 = [0.06, 0.07, 0.08]^t,$$
$$Q_2 = [0.09, 0.10, 0.11]^t,$$
$$Q_3 = [0.12, 0.13, 0.14]^t,$$
$$Q_4 = [0.15, 0.16, 0.17]^t$$

and

$$Q_5 = [0.18, 0.19, 0.20]^t.$$

Clearly if the interval Q is denoted by [X, Y]; $Q_1$ acts as the minimal element and $Q_5$ as the maximal fuzzy matrix and $Q_3$ = [0.12, 0.13, 0.14] as the optimal fuzzy matrix, the average fuzzy matrix of this interval is different and it is O = [0.10, 0.13, 0.14].

The matrix corresponding to the fuzzy interval matrix of membership is given by $R_1$, $R_2$, $R_3$, $R_4$, $R_5$ where

$R_1$ = $\{(0.0096, 0.0071, 0.0222)\}^t$,
$R_2$ = $\{(0.0269, 0.0300, 0.0220)\}^t$,
$R_3$ = $\{(0.0241, 0.0265, 0.0249)\}^t$,
$R_4$ = $\{(0.0114, 0.0381, 0.0288)\}^t$ and
$R_5$ = $\{(0.0356, 0.0189, 0.0376)\}^t$.

The minimal fuzzy membership matrix is A = $\{(0.0096, 0.0071, 0.0220)\}^t$ and the maximal fuzzy membership matrix B = $\{(0.0356, 0.0381, 0.0376)\}^t$.

The optimal fuzzy membership matrix O = $(0.0226, 0.0226, 0.0298)^t$, the average matrix $\overline{M}^t = (R_1 + R_2 + R_3 + R_4 + R_5)^t$ = $(0.0215, 0.0241, 0.0271)^t$. Now the interval matrix for $R_1^t$, $R_2^t$, …, $\overline{M}^t$ is given by [A, B].

Thus we have from the statistical data calculated the values of R and Q and now our work is to find the preference.

Now we calculate the passengers preference fuzzy matrices $P_1$, $P_2$, $P_3$, $P_4$ and $P_5$ for the pair of fuzzy matrices ($Q_2$, $R_2$), ($Q_3$, $R_3$), ($Q_4$, $R_4$) and ($Q_5$, $R_5$).

We use

$$P_i \text{ o } Q_i = R_i$$

i.e., max $p_{ij}\, q_{jk} = r_{ik}$.

We use the fuzzy relational equation method described in [] and calculate for

$$\begin{matrix} Q_1^t & R_1^t \\ \begin{bmatrix} 0.06 \\ 0.07 \\ 0.08 \end{bmatrix} & \begin{bmatrix} 0.0096 \\ 0.0071 \\ 0.0222 \end{bmatrix} \end{matrix}.$$

We have the fuzzy preference matrix $P_1$ given below:

$$P_1 = \begin{bmatrix} 0.03 & 0.06 & 0.12 \\ 0.0221875 & 0.044375 & 0.08875 \\ 0.069375 & 0.13875 & 0.2775 \end{bmatrix}.$$

For

$$\begin{matrix} Q_2^t & R_2^t \\ \begin{bmatrix} 0.09 \\ 0.10 \\ 0.11 \end{bmatrix} & \begin{bmatrix} 0.0269 \\ 0.0300 \\ 0.0220 \end{bmatrix} \end{matrix}.$$

We have the fuzzy preference matrix $P_2$ given by passengers

$$P_2 = \begin{bmatrix} 0.1345 & 0.269 & 0.06725 \\ 0.15 & 0.3 & 0.075 \\ 0.11 & 0.22 & 0.00605 \end{bmatrix}.$$

Now for the pair $(Q_3, R_3)$

$$\begin{matrix} Q_3^t & R_3^t \\ \begin{bmatrix} 0.12 \\ 0.13 \\ 0.14 \end{bmatrix} & \begin{bmatrix} 0.0241 \\ 0.0265 \\ 0.0249 \end{bmatrix} \end{matrix}.$$

We have given in the following the fuzzy passenger preference matrix $P_3$.

$$P_3 = \begin{bmatrix} 0.2008 & 0.1004 & 0.0502 \\ 0.2208 & 0.1104 & 0.0552 \\ 0.2075 & 0.10375 & 0.051875 \end{bmatrix}.$$

$$\begin{matrix} Q_4^t & R_4^t \\ \begin{bmatrix} 0.15 \\ 0.16 \\ 0.17 \end{bmatrix} & \begin{bmatrix} 0.0114 \\ 0.0381 \\ 0.0288 \end{bmatrix} \end{matrix}$$

the fuzzy passenger preference matrix $P_4$ is as follows.

$$P_4 = \begin{bmatrix} 0.035625 & 0.07125 & 0.0178125 \\ 0.1190625 & 0.23125 & 0.0553125 \\ 0.09 & 0.18 & 0.045 \end{bmatrix}.$$

$$\begin{matrix} Q_5^t & R_5^t \\ \begin{bmatrix} 0.18 \\ 0.19 \\ 0.20 \end{bmatrix} & \begin{bmatrix} 0.0356 \\ 0.0189 \\ 0.0376 \end{bmatrix}, \end{matrix}$$

the fuzzy passenger preference matrix $P_5$ is given below:

$$P_5 = \begin{bmatrix} 0.0445 & 0.089 & 0.178 \\ 0.023625 & 0.04725 & 0.0945 \\ 0.047 & 0.094 & 0.188 \end{bmatrix}.$$

On observing the fuzzy preference matrices $P_1$, $P_2$, …, $P_5$ calculated; we see the preferences correspond to peak hours of the day.

Now having seen and defined the FRIE model we now proceed on to illustrate function in case of fuzzy preference or membership functions are bimatrices and so on.

Suppose there are two sets of experts who want to work on the same problem. But both of them (both the sets) give different order fuzzy preference matrices, then it is not possible to use the fuzzy interval membership matrices model i.e., FRIE model discussed in page 171. So we are forced to build up a new model for this.

We will first explain the new model and give a simple illustrative example of the model.

Suppose we have $n_1$ number of experts giving their $m_1 \times p_1$ fuzzy relational membership matrices as $N_1^1$, $N_2^1$, …, $N_{n_1}^1$ on a certain problem P and for the same problem P we have $n_2$ number of experts who give their fuzzy relational membership matrices are $N_1^2$, $N_2^2$, …, $N_{n_2}^2$; which are $m_2 \times p_2$ fuzzy relational membership matrices. Now how to form the fuzzy interval of membership bimatrices using these sets of matrices $N_1^1$, $N_2^1$, …, $N_{n_1}^1$ and $N_1^2$, $N_2^2$, …, $N_{n_2}^2$. We shall give the method of construction for one set, for the other follows in a very similar and identical way.

Now our main aim is to construct a fuzzy interval of membership matrices using the fuzzy membership matrices $N_1^1$, $N_2^1$, …, $N_{n_1}^1$.

Let $N_t^1 = \{(n_{ij}^t)\}$; $1 \leq t \leq n_1$, $1 \leq i \leq m_1$ and $1 \leq j \leq p_1$. Now we form the minimal and maximal fuzzy membership matrices using these $n_1$ fuzzy membership matrices. Let $A_1$ be the minimal fuzzy membership matrix and $B_1$ be the maximal fuzzy membership matrix. Suppose $A_1 = (a_{ij}^1)$ and $B_1 = (b_{ij}^1)$ how to fill in the $m_1 \times p_1$ elements in $A_1$.

Define

$$a_{ij}^1 = \min \{ n_{ij}^1, n_{ij}^2, \ldots, n_{ij}^{n_1} \};$$

$1 \leq i \leq m_1$ and $1 \leq j \leq p_1$.

i.e.,

$$a_{11}^1 = \min \{ n_{11}^1, n_{11}^2, \ldots, n_{11}^{n_1} \};$$

$$a_{12}^1 = \min \{ n_{12}^1, n_{12}^2, \ldots, n_{12}^{n_1} \};$$

and so on.

$A_1 = (a_{ij}^1)$ will be the minimal fuzzy membership matrix and $a^1_{ij} < n^t_{ij}$; $1 \leq t \leq n_1$ with $1 \leq i \leq m_1$ and $1 \leq j \leq p_1$.

Now let $B_1 = (b_{ij}^1)$, we set the elements of $B_1$ as follows:

$$b_{ij}^1 = \max \{ n_{ij}^1, n_{ij}^2, \ldots, n_{ij}^{n_1} \};$$

with $1 \leq i \leq m_1$ and $1 \leq j \leq p_1$.

i.e.,

$$b_{11}^1 = \max \{ n_{11}^1, n_{11}^2, \ldots, n_{11}^{n_1} \};$$

$$b_{12}^1 = \max \{ n_{12}^1, n_{12}^2, \ldots, n_{12}^{n_1} \};$$

and so on.

Thus we have $b_{ij}^1 \geq n_{ij}^t$; $1 \leq t \leq n_1$ and $1 \leq i \leq m_1$ and $1 \leq j \leq p_1$. Thus $[A_1, B_1]$ will be known as the fuzzy membership of the interval of $m_1 \times p_1$ matrices associated with the fuzzy relational equations i.e., $[A_1, B_1]$ is a FRIE model.

We define for this set of matrices the optimal fuzzy membership matrix of the interval to be

$$O_1 = \frac{A_1 + B_1}{2}.$$

$$= \left[\frac{(a_{ij}^1) + (b_{ij}^1)}{2}\right],$$

$1 \leq i \leq m_1$ and $1 \leq j \leq p_1$.

Clearly $O \in [A_1, B_1]$. We define the mean or the average of the fuzzy membership matrices $N_1^1, N_2^1, ..., N_{n_1}^1$ to be

$$\overline{N^1} = \frac{N_1^1 + N_2^1 + ... + N_{n_1}^1}{n_1}$$

we see $\overline{N^1} \in [A_1, B_1]$. We call $[A_1, B_1]$ to be the associated fuzzy membership of interval of $m_1 \times p_1$ matrices of the FRIE model.

Now the same procedure is adopted for the set of $N_1^2, N_2^2, ..., N_{n_2}^2$ fuzzy membership $m_2 \times p_2$ matrices of the $n_2$ experts associated with the same problem.

We denote this fuzzy interval of fuzzy membership matrices of $m_2 \times p_2$ matrices of the same problem P by $[A_2, B_2]$.

That is $[A_2, B_2]$ is the fuzzy interval of membership $m_2 \times p_2$ matrices of the FRIE, where $A_2$ is the minimal fuzzy membership matrix of the FRIE. $B_2$ is the maximal fuzzy membership matrix of the FRIE. We have if $A_2 = (a_{ij}^2)$ and $B_2 = (b_{ij}^2)$; $1 \leq i \leq m_2$ and $1 \leq j \leq p_2$ if $N_2^t = (n_{ij}^t)$; $1 < t < n_2$ then $a_{ij}^2 \leq n_{ij}^t \leq b_{ij}^2$; $1 \leq t \leq n_2$. $1 \leq i \leq m_2$; $1 \leq j \leq p_2$.

Now let $[A, B] = [A_1, B_1] \cup [A_2, B_2]$; $[A, B]$ contains element N which are fuzzy bimatrices $N = N_1 \cup N_2$ where $N_1 \in [A_1, B_1]$ and $N_2 \in [A_2, B_2]$; we define $A = A_1 \cup A_2$ a fuzzy bimatrix known as the minimal fuzzy membership bimatrix. Similarly $B = B_1 \cup B_2$ is the fuzzy membership bimatrix known as the maximal fuzzy membership bimatrix.

[A, B] is the fuzzy interval of bimatrices known as the fuzzy interval of membership bimatrices associated with the FRIE i.e., FRIBE model. Clearly the optimal O, the fuzzy membership bimatrix is given by $O = O_1 \cup O_2$ where

$$O_2 = \frac{A_2 + B_2}{2}$$

$$= \left[\frac{(a_{ij}^2) + (b_{ij}^2)}{2}\right];$$

$1 \leq i \leq m_2$ and $1 \leq j \leq p_2$.

The average fuzzy membership of matrix $\overline{N} = \overline{N^1} \cup \overline{N^2}$ where

$$\overline{N^2} = \frac{N_1^2 + N_2^2 + \ldots + N_{n_2}^2}{n_2}.$$

We call this fuzzy interval of membership bimatrices as the membership fuzzy relational interval biequation model (FRIBE model). We call fuzzy relational biequation instead of fuzzy relational equation because we get two sets of fuzzy membership matrices of different orders.

We in the same manner define a set of fuzzy bivectors or fuzzy bimatrices associated with the same problem P. Suppose we denote this fuzzy interval of bimatrices by $[X, Y] = [X_1, Y_1] \cup [X_2, Y_2]$ then Q is a fuzzy bimatrix such that $Q = Q_1 \cup Q_2$ where $Q_1 \in [X_1, Y_1]$ and $Q_2 \in [X_2, Y_2]$. We can; given N and Q find the resultant R by the biequation;

$$N \text{ o } Q = R; \ (N \in [A, B] \text{ and } Q \in [X, Y].$$

i.e., $(N_1 \cup N_2) \text{ o } (Q_1 \cup Q_2)$

$$= (N_1 \text{ o } Q_1) \cup (N_2 \text{ o } Q_2)$$
$$= R_1 \cup R_2 = R.$$

It may be possible R is in the fuzzy biinterval of associated membership bimatrices; if the problem P under study had some expected values to form such interval. If not, we can form all the possible equations of the form;

N o Q = R, and the collection of all such R's can be made to form the fuzzy interval of bimatrices related with the FRIBE of the problem P. Thus R $\in$ [T, U] is the fuzzy interval of bimatrices associated with the FRIBE then we have;

[A B] o [X Y] = [T U] is the FRIBE's related dynamical system of fuzzy relational biequational interval of membership bimatrices; bound or satisfied by the operation P o Q = R where P $\in$ [A, B], Q $\in$ [X, Y] and R $\in$ [T, U] and fuzzy bimatrices i.e., $P = P_1 \cup P_2$, $R = R_1 \cup R_2$ and $Q = Q_1 \cup Q_2$.

i.e., $(P_1 \cup P_2)$ o $(Q_1 \cup Q_2) = (R_1 \cup R_2)$
i.e., $(P_1$ o $Q_1) \cup (P_2$ o $Q_2) = (R_1 \cup R_2)$

'o' may be the max min $(p^1_{ij}, q^1_{jk}) \cup$ max min $((p^2_{ij}, q^2_{jk}) = (r^1_{ik}) \cup (r^2_{ik})$. Now it may so happen that in the equation P o Q = R we may be given Q and R and we would be expected to find the fuzzy membership bimatrix of the FRIBE.

Now we will illustrate this situation for P o Q = R given P and Q it is always easy to solve. To give an example of given R and Q to find P is little different and we proceed on to find it.

For this we need two sets of experts to give their views on the same problem we take only the passengers preference problem which we have worked with 5 intervals of times and 5 intervals of passengers per hour. This we call as the first set of experts opinion. Now we divide the time intervals into 3 intervals and the passengers per hour also into 3 intervals.

Already we have calculated the fuzzy interval of preference for this problem using FRIE. Now we shall call them as the first set of solution to the problem, given by the first set of experts.

Now we formulate the second set of solutions to form the fuzzy interval of bimatrices of preference for the FRIBE. Now $Q_1^1$, $Q_2^1$ and $Q_3^1$ be the set of division of the time intervals which are given by

$$\begin{matrix} Q_2^1 \\ \begin{bmatrix} 0.12 \\ 0.13 \\ 0.14 \\ 0.15 \\ 0.16 \end{bmatrix} \end{matrix}, \quad \begin{matrix} Q_1^1 \\ \begin{bmatrix} 0.07 \\ 0.08 \\ 0.09 \\ 0.10 \\ 0.11 \end{bmatrix} \end{matrix} \quad \text{and} \quad \begin{matrix} Q_3^1 \\ \begin{bmatrix} 0.17 \\ 0.18 \\ 0.19 \\ 0.20 \\ 0.21 \end{bmatrix} \end{matrix}.$$

Using these $Q_1^1$, $Q_2^1$ and $Q_3^1$, we see the minimal fuzzy matrix for time interval is $Q_1^1$ and the maximal fuzzy matrix for time interval is $Q_3^1$ the optimal fuzzy matrix is given by

$$Q_1 = \frac{Q_1^1 + Q_3^1}{2} = \begin{bmatrix} 0.12 \\ 0.13 \\ 0.14 \\ 0.15 \\ 0.16 \end{bmatrix} = Q_2^1.$$

Now the average fuzzy matrix is given by

$$\frac{Q_1^1 + Q_2^1 + Q_3^1}{3} = \begin{bmatrix} 0.12 \\ 0.13 \\ 0.14 \\ 0.15 \\ 0.16 \end{bmatrix} = Q^1 = Q_2^1.$$

Thus we see in this set of fuzzy matrices of time interval $Q_1^1$, $Q_2^1$, $Q_3^1$ we see the minimal, maximal and the optimal and their mean / average fall with in the three fuzzy matrices. Thus the fuzzy interval of time matrices falls with in the fuzzy interval $[Q_1^1, Q_3^1]$. Now we have to study the number of fuzzy matrix of passengers $R_1^1$, $R_2^1$, $R_3^1$ where

$$R_1^1 = \begin{bmatrix} 0.0071 \\ 0.0222 \\ 0.0269 \\ 0.0300 \\ 0.0220 \end{bmatrix}, \quad R_2^1 = \begin{bmatrix} 0.0241 \\ 0.0265 \\ 0.0249 \\ 0.0114 \\ 0.0381 \end{bmatrix} \quad \text{and} \quad R_3^1 = \begin{bmatrix} 0.0288 \\ 0.0356 \\ 0.0189 \\ 0.0376 \\ 0.0182 \end{bmatrix}.$$

We find the maximal, minimal, optimal and the average for these three fuzzy matrix $R_1^1$, $R_2^1$, $R_3^1$ to form the fuzzy interval of number of passengers preferring those time intervals. Let $A^1$ be the minimal of the fuzzy matrices $R_1^1$, $R_2^1$, $R_3^1$.

$$A^1 = [(0.0071, 0.0222, 0.0189, 0.0114, 0.0182)]^t.$$

Let $B^1$ be the maximal of the fuzzy matrices $R_1^1$, $R_2^1$, $R_3^1$.

$$B^1 = [(0.0288, 0.0356, 0.0269, 0.0376, 0.0381)]^t.$$

The optimal matrix $O^1$ for this fuzzy interval $[A^1, B^1]$ is defined by

$$O^1 = \frac{A^1 + B^1}{2} = [(0.0180, 0.0289, 0.0229, 0.0245, 0.0282)]^t.$$

The average fuzzy matrix $R^1$ of the three fuzzy matrices $R_1^1$, $R_2^1$, $R_3^1$ is

$$R^1 = [(0.0200, 0.0281, 0.0236, 0.0263, 0.0261)]^t.$$

Clearly $[A^1, B^1]$ contains $R_1^1$, $R_2^1$, $R_3^1$, $R^1$ and $Q^1$. Now we calculate P using the fuzzy relational equation P o Q = R.

Using $R_1^1$ and $Q_1^1$ we get the passengers preference fuzzy matrix $P_1^1$ given below:

$$P_1^1 = \begin{bmatrix} 0.0142 & 0.01775 & 0.02366 & 0.071 & 0.03555 \\ 0.0444 & 0.0555 & 0.074 & 0.222 & 0.111 \\ 0.0538 & 0.06725 & 0.08966 & 0.269 & 0.1345 \\ 0.06 & 0.075 & 0.1 & 0.3 & 0.15 \\ 0.044 & 0.055 & 0.07333 & 0.220 & 0.11 \end{bmatrix}.$$

Now using $R_2^1$ and $Q_2^1$ we get the corresponding passengers preference fuzzy matrices.

$$P_2^1 = \begin{bmatrix} 0.030125 & 0.03765625 & 0.05020833 & 0.0753125 & 0.150625 \\ 0.033125 & 0.04140625 & 0.0552083 & 0.0828125 & 0.165625 \\ 0.031125 & 0.03890625 & 0.051875 & 0.0778125 & 0.155625 \\ 0.01425 & 0.0178125 & 0.02375 & 0.035625 & 0.07125 \\ 0.047625 & 0.05953 & 0.079375 & 0.1190625 & 0.238125 \end{bmatrix}$$

Now using $R_3^1$ and $Q_3^1$ we get the fuzzy matrix of preference $P_3$ is given below;

$$P_3^1 = \begin{bmatrix} 0.0288 & 0.036 & 0.048 & 0.149 & 0.072 \\ 0.0356 & 0.0445 & 0.05933 & 0.178 & 0.089 \\ 0.0189 & 0.023625 & 0.0315 & 0.0945 & 0.04725 \\ 0.0376 & 0.047 & 0.06266 & 0.0188 & 0.094 \\ 0.0182 & 0.02275 & 0.03033 & 0.091 & 0.0455 \end{bmatrix}.$$

Thus the fuzzy interval of membership of bimatrix will give $P_1$ $P_1^1$, $P_2$ $P_2^1$ and $P_3$ $P_3^1$ to be the resultant fuzzy bimatrices; clearly they are mixed square fuzzy bimatrices.

We have shown [X, Y] corresponds to fuzzy time interval matrices. [A, B] which corresponds to fuzzy passenger interval matrices.

[$P_1$, $P_2$, $P_3$, $P_4$, $P_5$] the fuzzy passengers preference $3 \times 3$ matrices.

[$Q_1^1$, $Q_3^1$] fuzzy time interval matrices.

[$A^1$, $B^1$] fuzzy passengers interval matrices.

{$P_1^1$, $P_2^1$ and $P_3^1$} the fuzzy passengers preference matrix.

So the FRIBE of the fuzzy interval of bimatrices with composition rule will be.

$$[X, Y] \cup [Q_1^1, Q_3^1] = Q,$$

$$[A, B] \cup [A^1, B^1] = R.$$

$\{P_1, P_2, P_3, P_4, P_5\} \cup \{P_1^1, P_2^1, P_3^1\}$ = P the resultant fuzzy passenger preference bimatrices calculated using the fuzzy relational biequation;

$$P \text{ o } Q = R.$$

Now it may so happen; that on the same problem P we have more than two sets of experts. Suppose we have say some $t_1$, $t_2$, …, $t_n$, n sets of experts who give their fuzzy membership matrix to be a $m_i \times n_i$ matrix ($m_i \neq m_j$ if $i \neq j$) and i = 1, 2, …, n. Suppose we denote the fuzzy membership matrices of the n sets of experts with their fuzzy membership matrices as

$$\{(P_1^1, P_2^1, ..., P_{t_1}^1), (P_1^2, P_2^2, ..., P_{t_2}^2), \ldots, (P_1^n, P_{2_n}^n, ..., P_{t_n}^n)\}.$$

We will form the fuzzy interval of membership matrices associated with the FRIE as follows. We will explain for the set $\{(P_1^i, P_2^i, ..., P_{t_i}^i)\}$ for i = 1, 2, …, n, how the fuzzy membership interval of $m_i \times n_i$ matrices is constructed. Let $P_i^k = (p_j^{ki})$, k = 1, 2, ..., n, $1 \leq q \leq m_i$ and $1 \leq j \leq n_i$. Let $A^i$ be the minimal fuzzy membership matrix for the fuzzy membership matrices $\{P_1^i, P_2^i, ..., P_{t_i}^i\}$. Suppose $A^i = (a_{qj}^i)$; $1 \leq q \leq m_i$ and $1 \leq j \leq n_i$, we construct the values $(a_{qj}^i)$ as follows:

$$a_{qi}^i = \min\{p_{qj}^{k1}, p_{qj}^{k2}, ..., p_{qj}^{ki}\}.$$
$$a_{11}^i = \min\{p_{11}^{k1}, p_{11}^{k2}, ..., p_{11}^{ki}\},$$
$$a_{12}^i = \min\{p_{12}^{k1}, p_{12}^{k2}, ..., P_{12}^{ki}\}$$

and so on.

Thus $A^i = (a_{qj}^i)$ is the minimal element and is such that

$a_{qj}^i \leq P_{qj}^{ki}$; $1 \leq k \leq i$, $1 \leq q \leq m_i$ and $1 \leq j \leq n_i$ and i = 1, 2, …, n.

On similar lines we construct the maximal fuzzy membership matrix $B^i$ using $\{(P_1^i, ..., P_{t_i}^i)\}$; set of fuzzy membership matrices as follows. Suppose $B^i = (b_{qj}^k)$ define

$$(b^k_{qj}) = \max\{p^1_{qj}, p^2_{qj}, ..., p^{t_i}_{qj}\}.$$

i.e.,
$$(b^k_{11}) = \max\{p^1_{11}, p^2_{11}, ..., p^{t_i}_{11}\},$$
$$(b^k_{12}) = \max\{p^1_{12}, p^2_{12}, ..., p^{t_i}_{12}\},$$

and so on; clearly

$$b^k_{qj} \geq (p^i_{qj}); \ i = 1, 2, \ldots, t_i, \ 1 \leq q \leq m_i \text{ and } 1 \leq j \leq n_i.$$

Now $[A^i, B^i]$ is the fuzzy interval of membership matrices for the fuzzy membership matrices $\{P^i_1, P^i_2, ..., P^i_{t1}\}$. Now we calculate the optimal fuzzy interval of membership matrices

$$O^i = \frac{a^i_{qj} + b^i_{qj}}{2} = \left(\frac{A^i + B^i}{2}\right).$$

We calculate the average fuzzy membership matrix $\overline{N^i}$ as

$$\overline{N^i} = \frac{N^i_1 + N^i_2 + ... + N^i_{t_i}}{t_i}.$$

We see $\overline{N^i} \in [A_i, B_i]$.

Thus we call the fuzzy interval $[A_i, B_i]$ of $m_i \times n_i$ membership matrix to be the FRIE matrix associated with the problem related with the $t_i$ experts. Now for the other sets of experts we find the fuzzy interval of membership matrices $[A^1, B^1]$, $[A^2, B^2]$ and so on.

Thus for these n sets of experts $t_1, t_2, \ldots, t_n$ we have fuzzy interval of membership matrices $[A^1, B^1], [A^2, B^2], \ldots, [A^n, B^n]$ related with the n sets of experts $t_1, t_2, \ldots, t_n$ respectively.

Set $[A, B] = [A^1, B^1] \cup [A^2, B^2] \cup \ldots \cup [A^n, B^n]$ such that $[A, B]$ contains all collection of mixed rectangular n-matrices. If $N \in [A, B]$ then $N = N_1 \cup N_2 \cup \ldots \cup N_n$ where $N_i$ is a $m_i \times n_i$ fuzzy membership matrix from the fuzzy interval of $[A^i, B^i]$ matrices. We set $A = A^1 \cup A^2 \cup \ldots \cup A^n$ and $B = B^1 \cup B^2 \cup \ldots \cup B^n$. The optimal matrix $O = O^1 \cup O^2 \cup \ldots \cup O^n$ belongs to $[A, B]$ and the average matrix $\overline{N} = \overline{N}^1 \cup \overline{N}^2 \cup \ldots \cup \overline{N}^n$.

We define [A, B] to be the fuzzy interval of membership of mixed rectangular n-matrices associated with the FRInEs. This model gives the membership functions collectively.

Now on similar lines we have built for these n classes of fuzzy membership matrices the expected or effect of a collection of matrices, say $Q_1 \cup Q_2 \cup \ldots \cup Q_n = Q \in [X, Y]$, where each Q belongs to the fuzzy interval of expected fuzzy matrix for which, we have to find out the resultant fuzzy interval of n matrix [W, V].

We see [X Y] is a fuzzy interval of n-matrices for which we need to find the resultant. Likewise [W, V] is also a fuzzy interval of n-matrices which is the collection of resultants using the composition of fuzzy intervals;

$$[A, B] \text{ o } [X, Y] = [W, V] \text{ i.e., } P \text{ o } Q = R,$$

P, Q and R are fuzzy n matrices.

We just illustrate how the following the operations are performed.

Suppose

$$P = P_1 \cup P_2 \cup \ldots \cup P_n,$$

$$Q = Q_1 \cup Q_2 \cup \ldots \cup Q_n$$

and

$$R = R_1 \cup R_2 \cup \ldots \cup R_n.$$

$$P \text{ o } Q = R$$

i.e.,

$$(P_1 \cup P_2 \cup \ldots \cup P_n) \text{ o } (Q_1 \cup Q_2 \cup \ldots \cup Q_n) = R_1 \cup R_2 \cup \ldots \cup R_n$$

$$(P_1 \text{ o } Q_1) \cup (P_2 \text{ o } Q_2) \cup \ldots \cup (P_n \text{ o } Q_n) = R_1 \cup R_2 \cup \ldots \cup R_n,$$

i.e., each $P_i$ o $Q_i = R_i$ for i = 1, 2, …, n.

Thus we have seen how the system of fuzzy interval of membership n matrices associated with a FRInE functions.

Now it has become pertinent to mention that FRIE can be worked in different methods for solutions; for if in the equation P o Q = R any two of them out of P, Q, R is known we are not always guaranteed of the solution so we many use neural net works method, neural net work weighted method, genetic algorithm or any other suitable method to arrive at the solution.

We have in [207] used the neutral networks method when the flow rate cannot be determined.

## 3.7 IBAM model and its Generalizations

Now we proceed on to show how the BAM models can exploit the method of interval matrices to obtain a better solution or a solution when several of the experts give their opinion. We have already given a brief description of how a BAM model works in general in chapter one.

Now we use the interval matrices when several of the experts give opinion about the problem using a BAM-model which we call as interval BAM model i.e., IBAM model.

We know a BAM model is always defined using an interval, that the real interval say [–a, a]; a an integer. Since the BAM predicts the time period with in which an event can occur or repeat the occurrence or non occurrence and so on.

Thus if we have the synaptic connection matrix M where entries of M would be from a predetermined fixed interval. It is pertinent to mention that the BAM model uses only matrices and not fuzzy matrices.

We have been using in the models FCM, FCBM, FRM etc only fuzzy interval matrices. Now this is the first instant where we give the application of interval matrices which are not fuzzy.

Further we wish to state that the Bidirectional Associative Memories model is defined by an expert and the expert has the right to choose the scale of interval and the number of attributes connected with the problem. We call the matrix given by the expert as the synaptic connection matrix. Now any synaptic connection matrix associated with the problem will always find its value with in a real interval [–a, a], $a < \infty$. Also given a problem each expert can give his opinion, we do not have any means to compare their views or even the connection matrices given by them. The interval matrices can be used very conveniently in this situation without any bias or modification made to the problem or to the experts opinion.

We first give a IBAM-model before we give the method how it is going to be used in the interval of matrices.

Just the study; cause of vulnerability to HIV/AIDS and factors for migration. The expert said that he wants to work with the scale [–5, 5].

They give 6 attributes relating to the causes for migrant labourers vulnerability to HIV/AIDS.

$A_1$ – No awareness / no education
$A_2$ – Social status
$A_3$ – No social responsibility and social freedom
$A_4$ – Bad company and addictive habits
$A_5$ – Types of profession
$A_6$ – Cheap availability of CSWs.

Factors forcing people for migration are

$C_1$ – Lack of labour opportunities in the home town
$C_2$ – Poverty / seeking better status of life.
$C_3$ – Mobilization of labour contracts
$C_4$ – Infertility of lands due to implementation of wrong research methodologies / failure of monsoon.

The synaptic connection matrix M given by the expert using the scale [–5, 5] is as follows :

$$M = \begin{array}{c} \\ A_1 \\ A_2 \\ A_3 \\ A_4 \\ A_5 \\ A_6 \end{array} \begin{array}{c} \begin{array}{cccc} C_1 & C_2 & C_3 & C_4 \end{array} \\ \begin{bmatrix} 5 & 2 & 4 & 4 \\ 4 & 3 & 5 & 3 \\ -1 & -2 & 4 & 0 \\ 0 & 4 & 2 & 0 \\ 2 & 4 & 3 & 3 \\ 0 & 1 & 2 & 0 \end{bmatrix} \end{array}$$

The expert can spell out any input vector

$X_k$ = {(3, 4, –1, –3, –2, 1)}

at the $k^{th}$ time period,

$S(X_k)$ = (1 1 0 0 0 1)

is the binary state vector of $X_k$.

$$S(X_k)M = (9, 6, 11, 7)$$
$$= Y_{k+1};$$
$$S(Y_{k+1}) = (1\ 1\ 1\ 1)$$

Like wise some one can give the input vector

$$Y_k = \{(-1, 3, 5, 2)\}$$

and find the resultant.

Thus we see the matrix input vector

$$X_k = (x_1, x_2, x_3, x_4, x_5, x_6)$$

are such that $x_i \in [-5, 5]$, for i = 1, 2, 3, 4, 5, 6. Similarly if $Y_k = (y_1, y_2, y_3, y_4)$ is any input vector then $y_j \in [-5, 5]$, j = 1, 2, 3, 4. Also if $M = (m_{ij})$ we see all $m_{ij} \in [-5, 5]$, $1 \leq i \leq 6$ and $1 \leq j \leq 4$.

Thus a BAM model works with three sets of matrices all defined on the same interval [–a, a], $a < \infty$. If M the connection synaptic matrix is given by $M = (m_{ij})$; $m_{ij} \in [-a, a]$; $1 \leq i \leq m$ and $1 \leq j \leq n$.

X is a $1 \times m$ row vector or a row matrix i.e., $X = (x_1, \ldots, x_m)$ where $x_i \in [-a, a]$, i = 1, 2, …, m, Y is a $1 \times n$ row vector or a now matrix i.e., $Y = (y_1, \ldots, y_n)$ where $y_j \in [-a, a]$, j = 1, 2, 3, …, n.

One faces with the problem once a large number of experts give opinion on the same problem using the same set of attributes and the use same interval. Do we have a means to compare them or work them as a collection and not as single entities. To this we use the three sets of interval matrices defined on the same interval [–a, a], $a < \infty$.

Let us now proceed onto give the description of the method for a very general case.

Let P be the problem at hand for which say some t experts are interested in giving their views. Suppose all them agree to work with $A_1, A_2, \ldots, A_m$ attributes and $B_1, B_2, \ldots, B_n$ attributes; where they agree to take $A_1, \ldots, A_m$ along the rows and $B_1, B_2, \ldots, B_n$ along the columns of the synaptic connection matrices $M_i$. They also wish to work with the fixed interval from [–a, a]; $a < \infty$. Let $M_1, M_2, \ldots, M_t$ be the synaptic connection matrix given by the t experts; each $M_i$ is a m × n matrix and they have their entries only from the real interval [–a, a].

Thus if $M_i = \left(a^i_{pq}\right)$, then we have $\left(a^i_{pq}\right) \in [-a, a]$ i.e., $-a \leq \left(a^i_{pq}\right) \leq a$ for i = 1, 2, …, t, and $1 \leq p \leq m$ and $1 \leq q \leq n$. Thus we have a collection of t number of m × n matrices all of them take values from the interval [–a, a].

We making use of these t number of m × n matrices form the interval of matrices [A, B] lying in the interval [–a, a].

For the interval of matrices [A, B] we need to find the minimal and maximal m × n matrix so that all the t matrices $M_1$, …, $M_t$ are in the interval of matrices [A, B] and A is the minimal matrix of this interval and B is the maximal matrix of this interval of matrices [A, B].

Let A = $(a_{pq})$, what should be the entries $(a_{pq})$ so that all $m^i_{pq} \geq a_{pq}$, $1 \leq p \leq m$ and $1 \leq q \leq n$ and $a_{pq} \in [-a, a]$.

Define

$$a_{pq} = \min\{m^1_{pq}, m^2_{pq}, \ldots, m^t_{pq}\}.$$

i.e.,

$$a_{11} = \min\{m^1_{11}, m^2_{11}, \ldots, m^t_{11}\},$$

$$a_{12} = \min\{m^1_{12}, m^2_{12}, \ldots, m^t_{12}\},$$

and so on.

Then A = $(a_{pq})$ is the minimal matrix of the interval [A, B] such that $a_{pq} \leq m^i_{pq}$, i = 1, 2, …, t, $1 \leq p \leq m$ and $1 \leq q \leq n$.

The matrix B = $(b_{pq})$ is the maximal matrix of the interval of matrices and is constructed as follows.

Define

$$(b_{pq}) = \max\{m^1_{pq}, m^2_{pq}, \ldots, m^t_{pq}\},$$

for $1 \leq p \leq m$ and $1 \leq q \leq n$.

Then

$$b_{11} = \max\{m^1_{11}, m^2_{11}, \ldots, m^t_{11}\};$$

$$b_{12} = \max\{m^1_{12}, m^2_{12}, \ldots, m^t_{12}\}$$

and so on.

We have $b_{pq} \geq m^i_{pq}$, $1 \leq i \leq t$, $1 \leq p \leq m$ and $1 \leq q \leq n$.

Now having defined the minimal and maximal matrices with entries from [–a, a] we see [A, B] is the interval of matrices such that all the matrices $M_1, M_2, \ldots, M_t$ are in [A, B]. i.e., all the synaptic connection $m \times n$ matrices given by the t experts with entries from the interval [–a, a] belong to the interval of matrices [A, B].

Now we define optimal synaptic connection $m \times n$ matrix O as

$$O = \frac{A+B}{2} = \frac{a_{pq} + b_{pq}}{2}.$$

Clearly $O \in [A, B]$. We define the average or mean of the t-synaptic connection matrices $M_1, M_2, \ldots, M_t$ to be $\overline{M}$ where

$$\overline{M} = \frac{M_1 + M_2 + \ldots + M_t}{t}.$$

$$= \frac{m^1_{pq} + m^2_{pq} + \ldots + m^t_{pq}}{t};$$

$\overline{M} \in [A, B]$. Thus our interval of synaptic connection matrices can maximum contain t + 4 number of $m \times n$ matrices.

Now what will be the $1 \times m$ matrix of input vector. We define the interval of input matrices as follows :

All $1 \times m$ matrices with entries from [–a, a] with $X_m = (-a, -a, \ldots, -a)$ to be the minimal element and $X_M = (a, a, \ldots, a)$ the $1 \times m$ matrix to be the maximal matrix of the interval of $1 \times m$ matrices defined on the interval / scale [–a, a]. Thus $[X_m, X_M] = \{(x_1, \ldots, x_m) / x_i \in [-a, a], i = 1, 2, \ldots, m, X_m = (-a, -a, \ldots, -a)$ and $X_M = (a, a, \ldots, a)\}$

We call $[X_m, X_M]$ to be the interval of input / resultant vectors of $1 \times m$ matrices.

On similar lines we define the input / resultant vector of $1 \times n$ matrices as $[Y_m, Y_M] = \{(Y_1, Y_2, \ldots, Y_n), Y_i \in [-a, a], i = 1, 2, \ldots, n.\ Y_m = (-a, -a, \ldots, -a)$ and $Y_M = (a, a, a, \ldots, a)\}$. Clearly $Y_m$ is the minimal matrix of the interval of $1 \times n$ matrices. $[Y_m,$

$Y_M$] defined on [–a, a] and $Y_M$ is the maximal matrix of the interval of 1 × n matrices.

Now for a IBAM of t experts to function we have 3 interval matrices defined on the scale [–a, a] viz. [A, B], [$Y_m$, $Y_M$] and [$X_m$, $X_M$]. The following operators on them are defined. Here [A, B] is the interval of m × n synaptic connection matrices containing the t experts opinion, with entries from [–a, a]. [$X_m$, $X_M$] and [$Y_m$, $Y_M$] is the interval of 1 × m and 1 × n resultant / input vector matrices with entries from [–a, a].

Given any input vector $X_k = (x_1, \ldots, x_m) \in [X_m, X_M]$, $S(X_k)$ denotes the binary signal vector formed using $X_k$ i.e., $S(X_k) = (t_1, \ldots, t_m)$ where $t_i$ is 0 or 1 if $t_i = 0$ it implies $x_i$ is 0 or negative i.e., $x_i \in [-a, 0]$ if $t_i = 1$ it implies $x_i \in (0, a]$.

On similar lines we work with $Y_k = (y_1, \ldots, y_n)$.

Thus S ([$X_m$, $X_M$]) o M = S ([$Y_m$, $Y_M$]) and S ([$Y_m$, $Y_M$]) o $M^T$ = S [$X_m$, $X_M$] where just above we have defined how the function S functions.

We will illustrate this by the following example:

We analyze the problem using experts opinion on the factors of migration and the role of government.

Taking the neuronal field $F_x$ as the attributes connected with the factors of migration $C_1$, $C_2$, $C_3$ and $C_4$; the cause for vulnerability of HIV/AIDS and factors for migration $A_1, A_2, \ldots, A_6$, where

| | | |
|---|---|---|
| $A_1$ | – | No awareness / no education |
| $A_2$ | – | Social status |
| $A_3$ | – | No social responsibility and social freedom |
| $A_4$ | – | Bad company and addictive habits |
| $A_5$ | – | Types of profession |
| $A_6$ | – | Cheap availability of CSWs. |

The factors forcing people for migration:

| | | |
|---|---|---|
| $C_1$ | – | Lack of labour opportunities in their home town. |
| $C_2$ | – | Poverty / seeking better status of life |
| $C_3$ | – | Mobilization of labour contractors. |

$C_4$ – Infertility of lands due to implementation of wrong research agricultural methodologies / failure of monsoon

Taking the neuronal field $F_X$ as the attributes connected with the causes of vulnerability resulting in HIV/AIDS and the neuronal field $F_Y$ is taken as factors forcing migration. The 6 × 4 matrix $M_1$ represents the forward synaptic projections from the neuronal field $F_X$ to the neuronal field $F_Y$. The 4 × 6 matrix $M_1^T$ represents the backward synaptic projections $F_X$ to $F_Y$. Now taking $A_1, A_2, ..., A_6$ along the rows and $C_1$, $C_2$, $C_3$ and $C_4$ along the columns we get the synaptic connection matrix $M_1$ given by the first expert and is modeled on the scale [–5, 5].

$$M_1 = \begin{array}{c} \\ A_1 \\ A_2 \\ A_3 \\ A_4 \\ A_5 \\ A_6 \end{array} \begin{array}{c} \begin{array}{cccc} C_1 & C_2 & C_3 & C_4 \end{array} \\ \left[\begin{array}{cccc} 5 & 2 & 4 & 4 \\ 4 & 3 & 5 & 3 \\ -1 & -2 & 4 & 0 \\ 0 & 4 & 2 & 0 \\ 2 & 4 & 3 & 3 \\ 0 & 1 & 2 & 0 \end{array}\right] \end{array}.$$

The synaptic connection matrix $M_2$ given by the second expert and modeled on the same scale [–5, 5] with the same set of attributes along the rows and columns.

$$M_2 = \begin{array}{c} \\ A_1 \\ A_2 \\ A_3 \\ A_4 \\ A_5 \\ A_6 \end{array} \begin{array}{c} \begin{array}{cccc} C_1 & C_2 & C_3 & C_4 \end{array} \\ \left[\begin{array}{cccc} 4 & 3 & 4 & 3 \\ 3 & 4 & 4 & 2 \\ -1 & -3 & 4 & -1 \\ 0 & 3 & 1 & 0 \\ 0 & 3 & 4 & 0 \\ 0 & 1 & 3 & -1 \end{array}\right] \end{array}.$$

The synaptic connection matrix $M_3$ given by the third expert and modeled on the same scale [–5, 5] with the same set of attributes along the same set of rows and column.

$$M_3 = \begin{array}{c} \\ A_1 \\ A_2 \\ A_3 \\ A_4 \\ A_5 \\ A_6 \end{array} \begin{array}{c} \begin{array}{cccc} C_1 & C_2 & C_3 & C_4 \end{array} \\ \left[\begin{array}{cccc} 3 & 2 & 3 & 2 \\ 4 & 3 & 3 & 1 \\ -1 & -2 & 3 & -1 \\ 0 & 2 & 3 & 0 \\ 2 & 2 & 3 & 1 \\ -1 & 0 & 3 & 0 \end{array}\right] \end{array}.$$

We also take the views of the forth expert modeled on the same scale [–5, 5] with the same set of attributes taken along the rows and columns.

Let $M_4$ be the matrix of the fourth expert.

$$M_4 = \begin{array}{c} \\ A_1 \\ A_2 \\ A_3 \\ A_4 \\ A_5 \\ A_6 \end{array} \begin{array}{c} \begin{array}{cccc} C_1 & C_2 & C_3 & C_4 \end{array} \\ \left[\begin{array}{cccc} 4 & 3 & 2 & 1 \\ 3 & 4 & 2 & 0 \\ -1 & -2 & 2 & 0 \\ -1 & 3 & 2 & 1 \\ 1 & 1 & 2 & 0 \\ 2 & 0 & 2 & 0 \end{array}\right] \end{array}.$$

Now using the four experts opinion given by the synaptic connection matrices $M_1$, $M_2$, $M_3$ and $M_4$, we construct the interval of matrices [A, B] such that A is the minimal matrix and B is the maximal matrix and such that the matrices $M_1$, $M_2$, $M_3$ and $M_4$ are in [A, B]. Define $A = (a_{ij})$ as follows:

$$a_{ij} = \min\{m_{ij}^1,\ m_{ij}^2,\ m_{ij}^3, m_{ij}^4\},$$

where $M^t = (m_{ij}^t)$, $1 \leq t \leq 4$, $1 \leq i \leq 6$ and $1 \leq j \leq 4$.

Thus

$$a_{11} = \min \{m^1_{11}, m^2_{11}, m^3_{11}, m^4_{11}\},$$
$$a_{12} = \min \{m^1_{12}, m^2_{12}, m^3_{12}, m^4_{12}\}$$

and so on.

i.e $a_{ij} \le m^t_{ij}$, $1 \le t \le 4$, $1 \le i \le 6$ and $1 \le j \le 4$.

Define B = $(b_{ij})$ as follows, $(b_{ij})$ is the maximal matrix of the interval of matrices [A, B]. Given B = $(b_{ij})$, $b_{ij}$ are defined in the following way $b_{ij} \ge m^t_{ij}$, $1 \le t \le 4$, $1 \le i \le 6$ and $1 \le j \le 4$.

$$b_{11} = \max\{m^1_{11}, m^2_{11}, m^3_{11}, m^4_{11}\}$$
$$b_{12} = \max\{m^1_{12}, m^2_{12}, m^3_{12}, m^4_{12}\}$$

and so on.

Clearly $a_{ij} \le m^t_{ij} \le b_{ij}$ for t = 1, 2, 3, 4; $1 \le i \le 6$ and $1 \le j \le 4$.

Thus [A, B] is an interval of matrices which contains $M_1$, $M_2$, $M_3$, and $M_4$. Define the optimal synaptic connection matrix as

$$O = \frac{(a_{ij}) + (b_{ij})}{2} = (o_{ij}) = \frac{A+B}{2};$$

$1 \le i \le 6$ and $1 \le j \le 4$. O is the optimal synaptic connection matrix associated with the synaptic connection matrices of the interval of matrices [A, B].

The average of synaptic connection matrices $M_1$, $M_2$, $M_3$ and $M_4$ is defined as

$$\overline{M} = \frac{M_1 + M_2 + M_3 + M_4}{4}$$
$$= \frac{(m^1_{ij}) + (m^2_{ij}) + (m^3_{ij}) + (m^4_{ij})}{4};$$

$1 \le i \le 6$ and $1 \le j \le 4$.

Thus $\overline{M} \in$ [A, B]. Now we give explicitly all the synaptic connection matrices A, B, $M_1$, $M_2$, $M_3$, $M_4$, O and $\overline{M}$. We have already given the synaptic connection matrices of the four experts $M_1$, $M_2$, $M_3$ and $M_4$. Now

$$A = \begin{array}{c} \\ A_1 \\ A_2 \\ A_3 \\ A_4 \\ A_5 \\ A_6 \end{array} \begin{array}{c} \begin{array}{cccc} C_1 & C_2 & C_3 & C_4 \end{array} \\ \begin{bmatrix} 3 & 2 & 2 & 1 \\ 3 & 3 & 2 & 0 \\ -1 & -3 & 2 & -1 \\ -1 & 2 & 2 & 0 \\ 1 & 1 & 2 & 0 \\ 0 & 0 & 2 & 0 \end{bmatrix} \end{array}.$$

Having defined the minimal synaptic connection matrix A now we proceed on to define the maximal synaptic connection matrix B i.e.,

$$B = \begin{array}{c} \\ A_1 \\ A_2 \\ A_3 \\ A_4 \\ A_5 \\ A_6 \end{array} \begin{array}{c} \begin{array}{cccc} C_1 & C_2 & C_3 & C_4 \end{array} \\ \begin{bmatrix} 5 & 3 & 4 & 4 \\ 4 & 4 & 5 & 3 \\ -1 & -2 & 4 & 0 \\ 0 & 4 & 3 & 1 \\ 2 & 4 & 4 & 3 \\ 2 & 1 & 3 & 0 \end{bmatrix} \end{array}.$$

Thus B is the maximal synaptic connection matrix of the interval of synaptic connection matrices.

Now we proceed on to give the exact value of the optimal synaptic connection matrix O

$$O = \begin{array}{c} \\ A_1 \\ A_2 \\ A_3 \\ A_4 \\ A_5 \\ A_6 \end{array} \begin{array}{c} \begin{array}{cccc} C_1 & C_2 & C_3 & C_4 \end{array} \\ \begin{bmatrix} 4 & 2.5 & 3 & 2.5 \\ 3.5 & 3.5 & 3.5 & 1.5 \\ -1 & -2.5 & 3 & -0.5 \\ -0.5 & 3 & 2.5 & 0.5 \\ 1.5 & 2.5 & 3 & 1.5 \\ 1 & 0.5 & 2.5 & 0 \end{bmatrix} \end{array}$$

we see O ∈ [A, B]. Now we calculate the average synaptic connection matrix.

$$\overline{M} = \begin{array}{c} A_1 \\ A_2 \\ A_3 \\ A_4 \\ A_5 \\ A_6 \end{array}\begin{bmatrix} 4 & 2.5 & 3.25 & 2.5 \\ 3.5 & 3.5 & 3.5 & 1.5 \\ -1 & -2.5 & 3.5 & -0.5 \\ -0.5 & 2 & 2 & 0 \\ 2 & 2.5 & 3 & 1.5 \\ 0.25 & 0.5 & 2.5 & -0.25 \end{bmatrix},$$

$\overline{M}$ ∈ [A, B]. Clearly $[X_1, X_2]$ = {all 1 × 6 row matrices with entries from the interval [–5, 5] with $X_1$ = [–5, –5, –5, –5, –5, –5] and $X_2$ = [5, 5, 5, 5, 5, 5] as the minimal and maximal input vectors respectively of the interval matrix $[X_1, X_2]$ of the interval of input matrix vectors}.

Similarly $[Y_1, Y_2]$ = {all 1 × 4 row matrices with entires from [–5, 5] with $Y_1$ = [–5, –5, –5, –5] and $Y_2$ = [+5, +5, +5, +5] as the minimal and maximal input vectors of the interval of matrices}. Any input vector $X_k \in [X_1, X_2]$ will be made or transformed into the binary signal vector by the function S as

$$S(X_k) = S(x_1, x_2, \ldots, x_6) = (x'_1, x'_2, \ldots, x'_6).$$

Here $X_k = (x_1, x_2, x_3, \ldots, x_6)$; $x_i$ in $X_k$ will be 0 in $S(x_1, x_2, \ldots, x_6)$ if $x_i \leq 0$ and for $x_i$ in $X_k$ will be 1 in $S(x_1, \ldots, x_6) = S(X_k)$, if $x_i > 0$ thus $(x'_1, \ldots, x'_6)$ is either 0 or 1 i.e., the components in the input vectors are 0 or 1 i.e., OFF or ON.

Like wise for the state vectors in $[Y_1, Y_2]$. Now we will just illustrate the working.

Given $X_k$ an input vector from $[X_1, X_2]$ and M ∈ [A, B], then $S(X_k)$ is the binary state vector of $X_k$.

$$S(X_k)M = Y_{k+1}$$

and

$$S(Y_{k+1}) \in [Y_1, Y_2] \text{ i.e.,}$$
$$S(Y_{k+1})M^T = X_{k+2} \in [X_1, X_2]$$

$S(X_{k+2})M \quad = \quad Y_{k+3}$

and so on till the dynamical system settles in a fixed binary pair. This we show by taking the vector $X_k = (3, 4, –1, –3, –2, 1) \in [X_1, X_2]$ and the resultant from all synaptic connection matrices from [A, B].

Let $X_k$ be the input vector given as $X_k = (3, 4, –1, –3, –2, 1)$ at the $k^{th}$ time period.

The initial vector is given such that illiteracy, lack of awareness, social status and cheap availability of CSWs have stronger impact over vulnerability of migrant workers to become infected by HIV/AIDS.

We suppose that all neuronal state change decisions are synchronous.

The binary signal vector $S(X_k) = (1\ 1\ 0\ 0\ 0\ 1)$. From the activation equation

$$S(X_k)M_1 = (9, 6, 1, 1, 7) = Y_{k+1}.$$

From the activation equation

$$S(Y_{k+1}) = (1\ 1\ 1\ 1).$$

Now

$$S(Y_{k+1})\, M_1^T = (15, 15, 16, 1, 2, 3) = X_{k+2}.$$

From the activation equation

$$S(X_{k+2}) = (1\ 1\ 1\ 1\ 1\ 1).$$

$$S(X_{k+2})M_1 = (10, 12, 20, 10) = Y_{k+3}$$

$$S(Y_{k+3}) = (1\ 1\ 1\ 1).$$

Thus the binary pair {(1 1 1 1 1 1), (1 1 1 1)} represents a fixed point of the dynamical system. Equilibrium of the state has occurred at the time k + 2 when the starting time was k. Thus the fixed point suggest that illiteracy with unawareness, poor social status and cheap availability of CSW lead to the patients remaining or becoming socially free with no social responsibility have all addictive habits added to bad company which directly depends or the types of profession they choose.

Now we study the effect of the same input vector given as $X_k$ = (3, 4, –1, –3, –2, 1) at the $k^{th}$ time period using the second experts dynamical system $M_2$;

| | | |
|---|---|---|
| $S(X_k)$ | = | (1, 1, 0, 0, 0, 1) |
| $S(X_k)M_2$ | = | (7, 8, 11, 4) |
| | = | $Y_{k+1}$ |
| $S(Y_{k+1})$ | = | (1 1 1 1) |
| $S(Y_{k+1}) M_2^T$ | = | (14, 13, –1, 4, 12, 3) |
| | = | $X_{k+2}$ |
| $S(X_{k+2})$ | = | (1 1 0 1 1 1) |
| $S(X_{k+2})M_2$ | = | (10, 14, 20, 6) |
| | = | $Y_{k+3}$ |
| $S(Y_{k+3})$ | = | (1 1 1 1) |
| $S(Y_{k+3}) M_2^T$ | = | $Y_{k+4}$. |

Thus the binary pair {(1 1 0 1 1 1), (1 1 1 1)} represents a fixed point and we see only the attribute $A_3$ remains in the off state. No social responsibility and social freedom has no impact according to this expert, but however has impact with all other attributes as they come to on state.

Next we try the effect of the same input vector $X_k$ = (3, 4, –1, –3, –2, 1) at the $k^{th}$ period using the views of the $3^{rd}$ expert. We use the synaptic connection matrix given by the $3^{rd}$ expert.

| | | |
|---|---|---|
| $S(X_k)$ | = | (1 1 0 0 0 1) |
| $S(X_k)M_3$ | = | (6, 5, 9, 3) |
| | = | $Y_{k+1}$ |
| $S(Y_{k+1})$ | = | (1 1 1 1) |
| $S(Y_{k+1}) M_3^T$ | = | (10, 11, –1, 5, 8, 2) |
| | = | $X_{k+2}$ |
| $S(X_{k+2})$ | = | (1 1 0 1 1 1) |
| $S(X_{k+2})M_3$ | = | (8, 9, 15, 3) |
| | = | $Y_{k+3}$ |
| $S(Y_{k+3})$ | = | (1 1 1 1). |

Thus the resultant is a fixed binary pair given by {(1 1 0 1 1 1), (1 1 1 1)}. We see the opinions of the second and the third expert are identical.

Now we work with the same input vector $X_k$ = (3, 4, –1, –3, –2, 1) and use the fourth expert's synaptic connection matrix $M_4$ of the interval matrix [A, B];

| | | |
|---|---|---|
| $S(X_k)$ | = | (1 1 0 0 0 1) |
| $S(X_k)M_4$ | = | (9, 7, 6, 1) |
| | = | $Y_{k+1}$ |
| $S(Y_{k+1})$ | = | (1 1 1 1) |
| $S(Y_{k+1})\,M_4^T$ | = | (10, 9, –1, 5, 4, 4) |
| | = | $X_{k+2}$ |
| $S(X_{k+2})$ | = | (1 1 0 1 1 1) |
| $S(X_{k+2})M_4$ | = | (9, 11, 10, 2) |
| | = | $Y_{k+3}$ |
| $S(Y_{k+3})$ | = | (1 1 1 1) |
| $S(Y_{k+3})\,M_4^T$ | = | $X_{k+4}$ |
| and $S(X_{k+4})$ | = | (1 1 0 1 1 1). |

Thus we see the resultant or the equilibrium of the dynamical system is given by a binary fixed point {(1 1 0 1 1 1), (1 1 1 1)}. We see the three experts whose synaptic connection matrix is given by $M_2$, $M_3$ and $M_4$ have the same fixed binary pair with only the node $A_3$ in the off state.

Now we proceed on to work with the same input vector $X_k$ = (3, 4, –1, –3, –2, 1); now using the minimal synaptic connection matrix A of the interval matrix [A, B].

| | | |
|---|---|---|
| $S(X_k)$ | = | (1 1 0 0 0 1) |
| $S(X_k)A$ | = | (6, 5, 4, 1) |
| | = | $Y_{k+1}$ |
| $S(Y_{k+1})$ | = | (1 1 1 1) |
| $S(Y_{k+1})\,A^T$ | = | (8, 8, –3, 3, 4, 2) |
| | = | $X_{k+2}$ |
| $S(X_{k+2})$ | = | (1 1 0 1 1 1) |
| $S(X_{k+2})A$ | = | (6, 8, 10, 1) |
| | = | $Y_{k+3}$ |
| $S(Y_{k+3})$ | = | (1 1 1 1). |

Thus the minimal synaptic connection matrix A also gives the same resultant as that of the synaptic connection matrices $M_2$, $M_3$ and $M_4$ given by the second, third and the forth expert respectively.

The resultant given by A is a fixed binary pair {(1 1 0 1 1 1), (1 1 1 1)}.

Next we study the effect of the same input vector $X_k$ = (3, 4, –1, –3, –2, 1) on the synaptic connection matrix B, where B is the maximal synaptic connection matrix of the interval matrix [A, B].

| | | |
|---|---|---|
| $S(X_k)$ | = | (1 1 0 0 0 1) |
| $S(X_k)B$ | = | (11, 8, 12, 7) |
| | = | $Y_{k+1}$ |
| $S(Y_{k+1})$ | = | (1 1 1 1) |
| $S(Y_{k+1})B^T$ | = | (16, 16, 18, 13, 6) |
| | = | $X_{k+2}$ |
| $S(X_{k+2})$ | = | (1 1 1 1 1 1) |
| $S(X_{k+2})B$ | = | (13, 16, 19, 11) |
| | = | $Y_{k+3}$ |
| $S(Y_{k+3})$ | = | (1 1 1 1). |

Thus the resultant given by the maximal connection matrix of the interval matrix [A, B] is a fixed binary pair [(1 1 1 1 1 1), (1 1 1 1)].

The views of the first expert given by the synaptic connection matrix $M_1$ and that of the resultant given by the maximal connection synaptic matrix B coincide. Now we work with the optimal synaptic connection matrix O of the interval matrix [A, B] using the same input vector $X_k$ = (3, 4, –1, –3, –2, 1) on $O \in [A, B]$.

| | | |
|---|---|---|
| $S(X_k)$ | = | (1 1 0 0 0 1) |
| $S(X_k)O$ | = | (8.5, 6.5, 9, 4) |
| | = | $Y_{k+1}$ |
| $S(Y_{k+1})$ | = | (1 1 1 1) |
| $S(Y_{k+1})O^T$ | = | (12, 11, –1, 5.5, 8.5, 4) |
| | = | $X_{k+2}$ |
| $S(X_{k+2})$ | = | (1 1 0 1 1 1) |
| $S(X_{k+2})O$ | = | (10, 12, 14.5, 6) |
| | = | $Y_{k+3}$ |
| $S(Y_{k+3})$ | = | (1 1 1 1). |

Thus the resultant is a fixed binary pair given by {(1 1 0 1 1 1), (1 1 1 1)}. Now we finally proceed on to work with the average synaptic connection matrix $\overline{M}$ for the same input vector

$$X_k = (3, 4, 1, –3, –2, 1).$$

The binary signal vector

$$S(X_k) = (1\ 1\ 0\ 0\ 0\ 1).$$

From the activation equation

| | | |
|---|---|---|
| $S(X_k)\overline{M}$ | = | (7.75, 6.5, 9.25, 3.75) |
| | = | $Y_{k+1}$. |

The binary signal vector

| | | |
|---|---|---|
| $S(Y_{k+1})$ | = | (1 1 1 1) |
| $S(Y_{k+1})\overline{M}^T$ | = | (12.25, 12, -0.5, 3.5 9, 3) |
| | = | $X_{k+2}$. |

The binary signal vector of $X_{k+2}$ is

$$S(X_{k+2}) = (1\ 1\ 0\ 1\ 1\ 1).$$

From the activation equation

| | | |
|---|---|---|
| $S(X_{k+2})\overline{M}$ | = | (9.25, 11, 14.25, 5.25) |
| | = | $Y_{k+3}$ |
| $S(Y_{k+3})$ | = | (1 1 1 1). |

The resultant of the input vector $X_k$ is given by the fixed binary pair {(1 1 0 1 1 1), (1 1 1 1)}. Now we can compare all the fixed binary pairs for the same input vector $X_k$ = (3, 4, –1, –3, –2, 1). Using the first experts opinion i.e., the synaptic connection matrix $M_1$ we get the binary pair for the input vector $X_k$ as {(1 1 1 1 1 1), (1 1 1 1)}. For the dynamical system $M_2$ given by the second expert we get the fixed point to be the equilibrium given by {(1 1 0 1 1 1), (1 1 1 1)}.

For the synaptic connection matrix $M_3$ given by the third expert the resultant is a fixed point given by {(1 1 0 1 1 1), (1 1 1 1)}. The resultant for the input vector same $X_k$ given by the forth expert is also {(1 1 0 1 1 1), (1 1 1 1)}. Now using the minimal synaptic connection matrix A we get the equilibrium to be a fixed binary pair given by {(1 1 0 1 1 1), (1 1 1 1)}.

From the maximal synaptic connection matrix B of the interval of synaptic connection matrices [A, B], we get resultant for the input vector $X_k$ = (3, 4, –1, –3, –2, 1) to be {(1 1 1 1 1 1), (1 1 1 1)}. On similar lines we see the resultant for the input vector $X_k$ using the optimal synaptic connection matrix O is a fixed binary pair given by {(1 1 0 1 1 1), (1 1 1 1)}.

The average synaptic connection matrix $\overline{M}$ of the interval of matrix [A, B] gives the resultant for the input vector $X_k$ to be {(1 1 0 1 1 1), (1 1 1 1)}.

Thus we see the input vector $X_k$ is so powerful that is why on all the synaptic connection matrices of the interval matrix [A, B], we get the same resultant i.e., {(1 1 0 1 1 1), (1 1 1 1)} or {(1 1 1 1 1 1), (1 1 1 1)}.

But it is easily verified that all the eight synaptic connection matrices are very different. Thus we see when more than one expert gives his opinion the interval of synaptic connection matrices helps one to compare results and get the consolidated result as well as give the minimal, optimal and maximal resultants for a given input vector. Thus we see the interval matrices is greatly helpful when the number of experts is more than one and above all equal importance is given to every expert and no bias or favouratism is shown while deriving at the results.

We cannot solve using this model even if two experts give their opinion with different interval period and with different sets of attributes.

Now we understand the interval matrices cannot do any help so we for the first time define the notion of bi-directional biassociative memories (BBAM) when only two experts solve the same problem using different intervals of time and different sets of attributes which may be overlapping; we just describe in general how such problem is solved and then go for an illustration.

Let P be a problem at hand where two experts agree upon to work with it using the BAM model. Both work independently, further they do not agree on all aspects. Thus $E_1$ is the first expert who works with $m_1$ nodes along the rows and $n_1$ nodes along the columns and wishes to take the scale $[-a_1, a_1]$, $a_1$ an integer; $a_1 < \propto$.

The second expert $E_2$ works on the same problem independently taking some $m_2$ nodes along the rows and $n_2$ nodes along the columns. We do not prohibit the expert $E_1$ and $E_2$ having any number of common attributes between $m_1$ and $m_2$ (and / or) $n_1$ and $n_2$. The second expert further wishes to work on the scale from $[-a_2, a_2]$, $a_2 < \propto$ and $a_2$ an integer. Clearly the synaptic connection matrix given by the first expert is a $m_1 \times n_1$ matrix with entries from the interval $[-a_1, a_1]$ and the synaptic connection matrix given by the second expert is a $m_2 \times n_2$ matrix with entries from the interval $[-a_2, a_2]$.

Now let $M_1 = (m^1_{ij})$, $1 \leq i \leq m_1$, $1 \leq j \leq n_1$ be the synaptic connection matrix given by the first expert; $-a_1 \leq (m^1_{ij}) \leq a_1$.

Let $M_2 = (m^2_{ij})$, $1 \leq i \leq m_2$, $1 \leq j \leq n_2$ be the synaptic connection matrix given by the second expert; $-a_2 \leq m^2_{ij} \leq a_2$.

Let $M = M_1 \cup M_2$, M is a bimatrix which is a mixed rectangular bimatrix. We call M to be the synaptic connection bimatrix of the Bidirectional Biassociative memories (BBAM) model. Now we just show how the dynamical system M functions.

Let $X = X_1 \cup X_2$ where $X_1 = (x^1_1, ..., x^1_{m_1})$ and $X_2 = (x^2_1, x^2_2, ..., x^2_{m_2})$ be the input bivector. To find the effect of X on the dynamical bisystem M which is a bimatrix. It is important to note that each $x^1_i$ lies in the interval $-a_1 \leq x^1_i \leq a_1$, $1 \leq i \leq m_1$. Similarly $x^2_j \in X_2$ lies in the interval $-a_2 \leq x^2_j \leq a_2$ with $1 \leq j \leq m_2$.

We first convert the given input bivector into a binary bivector as follows. Given $X = X_1 \cup X_2$ where

$$X_1 = (x^1_1, ..., x^1_{m_1}) \text{ and } X_2 = (x^2_1, x^2_2, ..., x^2_{m_2}) .$$

$$S(X) = (a_1^1, a_2^2, ..., a_{m_1}^1) \cup (a_1^2, a_2^2, ..., a_{m_2}^2).$$

where $a_j^1$ and $a_j^2$ are either 0 or 1. We put $a_j^1 \in X_1$ to be 0 if $x_j^1 \leq 0$ and $a_j^1 = 1$ if $x_j^1 > 0$. This is true for i = 1, 2, …, $m_1$.

Similarly we put for $a_j^2 \in X_2$ to be 0 if $x_j^2 \leq 0$ and $a_j^2 = 1$ if $x_j^2 > 1$; this is true for j = 1, 2, …, $m_2$.

Further when we say $a_i^t = 0$, (t = 1, 2) it means the $i^{th}$ node corresponding to the $i^{th}$ attribute taken along the row is in the off state.

Similarly when we say $a_j^t = 1$, (t = 1, 2) we mean the $j^{th}$ node corresponding to the $j^{th}$ attribute taken along the row of the synaptic matrix $M_t$ is on, (t = 1, 2).

Thus we call S(X) to be the binary state bivector corresponding to the input bivector $X = X_1 \cup X_2$. Thus S is a function which converts any resultant or input vector into a binary state vector.

Now given $M = M_1 \cup M_2$ is the synaptic connection bimatrix and $X = X_1 \cup X_2$ is a input vector at the $k^{th}$ time period $S(X_k)$ is the binary state bivector, clearly $S(X_k) = S(X_1 \cup X_2) = S(X_1) \cup S(X_2)$.

Now how does the dynamical bisystem M function;

$$\begin{aligned} S(X_k)M &= [S(X_1) \cup S(X_2)]\,[M_1 \cup M_2] \\ &= S(X_1)M_1 \cup S(X_2)M_2 \\ &= Y_1 \cup Y_2 \\ &= Y_{k+1} \text{ (say).} \end{aligned}$$

$S(Y_{k+1})$ is a binary state bivector of $Y_{k+1}$; $S(Y_1)$ is a $1 \times n_1$ binary state vector and $S(Y_2)$ is a $1 \times n_2$ binary state vector. $S(Y_{k+1})$ is the $(1 \times n_1 \cup 1 \times n_2)$ bivector at the $(k + 1)^{th}$ time period. Now

$$\begin{aligned} S(Y_{k+1})M^T &= S(Y_1) \cup S(Y_2)\,[M_1^T \cup M_2^T] \\ &= S(Y_1)\,M_1^T \cup S(Y_2)M_2^T \end{aligned}$$

$$= \quad X_{k+2}$$
$$= \quad X_1^1 \cup X_2^1 .$$

Now the binary state bivector of $X_{k+2}$ is given by

$$S(X_{k+2}) \quad = \quad S(X_1^1) \cup S(X_2^1) .$$

Now

$$S(X_{k+2})M \quad = \quad [S(X_1^1) \cup S(X_2^1)][M_1 \cup M_2]$$
$$= \quad S(x_1^1) M_1 \cup S(x_2^1) M_2 ,$$

and so on until we arrive at the equilibrium of the dynamical bisystem which is a fixed binary pair which may occur after the $t^{th}$ period of time say $X_{k+t}$ or $Y_{k+t+1}$. We shall illustrate this with a real model before we proceed on to define the notion of BnAM model ($n \geq 2$).

Let us presume that the problem under investigation is to analyse the causes of migrant labourers vulnerability to HIV/AIDS and the role of the government. Let the attributes associated with the causes of migrant labourers vulnerability to HIV/AIDS:

$A_1$ – No awareness / no education
$A_2$ – Social status
$A_3$ – No social responsibility and social freedom
$A_4$ – Bad company and addictive habits
$A_5$ – Types of profession they do like lorry driver, bore pipe workers; etc
$A_6$ – Cheap availability of CSWs.

The attributes related with the role of government are:

$G_1$ – Alternate job if agriculture fails there by stopping migration
$G_2$ – Awareness clubs in rural areas about HIV/AIDS

$G_3$ – Construction of hospitals in rural areas with HIV/AIDS Counseling cell / compulsory HIV/AIDS test before marriage

$G_4$ – Failed to stop the misled agricultural techniques followed recently by farmers

$G_5$ – No foresight for the government and no precautionary actions taken from the past occurrences

The first expert takes these attributes taking the neuronal field $F_X$ as the attributes connected with the causes of vulnerability resulting in HIV/AIDS and the neuronal field $F_Y$ is taken as the role of government. The 6 × 5 matrix $M_1$ represents the forward synaptic projections from the neuronal field $F_X$ to the neuronal field $F_Y$.

The 5 × 6 matrix $M_1^T$ represents the backward synaptic projections of $F_X$ to $F_Y$. Now taking $A_1, A_2, \ldots, A_6$ along the rows and $G_1, G_2, \ldots, G_5$ as the columns we get the synaptic connection matrix $M_1$ which is modeled in the scale [–5, 5].

$$M_1 = \begin{array}{c} \\ A_1 \\ A_2 \\ A_3 \\ A_4 \\ A_5 \\ A_6 \end{array} \begin{array}{c} \begin{array}{ccccc} G_1 & G_2 & G_3 & G_4 & G_5 \end{array} \\ \begin{bmatrix} 3 & 5 & 1 & 2 & 3 \\ 4 & 4 & 3 & 3 & 2 \\ -2 & 3 & 0 & -2 & 0 \\ 0 & -1 & 1 & -3 & 3 \\ -1 & 0 & 4 & 0 & 1 \\ 5 & 4 & 2 & 3 & 4 \end{bmatrix} \end{array}.$$

However the second expert did not agree to work with the same number of attributes as the first expert. He however choose to take 5 attributes $A_1, \ldots, A_5$, he merged social status with social responsibility and social freedom i.e., $A_3$ and $A_2$ is merged. Thus he worked with $A_1$, $A_2$, $A_4$, $A_5$ and $A_6$ as the attributes related with the vulnerability of the migrant labourers to HIV/AIDS.

Now regarding the role of government he has coupled $G_2$ and $G_3$ thus he workers only with the four attributes $G_1$, $G_2$, $G_3$ and $G_4$ combined as $G_2$, $G_4$ and $G_5$.

Let $M_2$ be the synaptic connection matrix given by the second expert.

$$M_2 = \begin{array}{c} \\ A_1 \\ A_2 \\ A_4 \\ A_5 \\ A_6 \end{array} \begin{array}{c} \begin{matrix} G_1 & G_2 & G_3 & G_4 \end{matrix} \\ \begin{bmatrix} 2 & 4 & -1 & 3 \\ 4 & 3 & 2 & 2 \\ 0 & 2 & -1 & 0 \\ 3 & -2 & 4 & 1 \\ 4 & -2 & 3 & 4 \end{bmatrix} \end{array}$$

modeled in the scale [–4, 4]. $M_2$ is a 5 × 4 synaptic connection matrix. Now the input bivectors $X = X_1 \cup X_2$ are from $F_{X_1} \cup F_{X_2}$, where

$$X = (x_1^1, x_2^1, ..., x_6^1) \cup (x_1^2, x_2^2, x_3^2, x_4^2, x_5^2).$$

The input bivectors from $F_{Y_1} \cup F_{Y_2}$ are of the form

$$\begin{aligned} Y &= Y_1 \cup Y_2 \\ &= (y_1^1, y_2^1, ..., y_5^1) \cup (y_1^2, y_2^2, y_3^2, y_4^2). \end{aligned}$$

The dynamical bisystem associated with the bidirectional biassociative memories model is given by $M = M_1 \cup M_2$, clearly M is a mixed rectangular bimatrix. Just we indicate how a dynamical bisystem functions. Let $X = X_1 \cup X_2$ be a given input bivector; transform the input bivector into a binary state bivector using the S function i.e., $S(X) = S(X_1) \cup S(X_2)$ clearly the entries in $S(X_1)$ and $S(X_2)$ are either 0 or 1 i.e., the attribute corresponding to the $i^{th}$ place is 0 if the node is in the OFF state and 1 if the node is in the ON state. Now we see the effect of S(X) on M i.e., by the activation function equation, we have;

$$\begin{aligned} S(X)M &= [S(X_1) \cup S(X_2)] [M_1 \cup M_2] \\ &= S(X_1) M_1 \cup S(X_2) M_2. \end{aligned}$$

For instance let us consider $X_k = X_1 \cup X_2$ at the $k^{th}$ interval of time $X_k = (4, -2, 0, 1, -3, -2) \cup (3, -1, -2, 1, 0) \in F_{X_1} \cup F_{X_2}$.

$$
\begin{aligned}
S(X_k) &= (1\ 0\ 0\ 1\ 0\ 0) \cup (1\ 0\ 0\ 1\ 0). \\
S(X_k)M_1 &= [S(X_1) \cup S(X_2)]\ [M_1 \cup M_2] \\
&= (1\ 0\ 0\ 1\ 0\ 0)\ M_1 \cup (1\ 0\ 0\ 1\ 0)\ M_2 \\
&= (3, 4, 2, -1, 6) \cup (5, 2, 3, 4) \\
&= Y_1 \cup Y_2 \\
&= Y_{k+1}. \\
S(Y_{k+1}) &= S(Y_1 \cup Y_2) \\
&= S(Y_1) \cup S(Y_2) \\
&= (1\ 1\ 1\ 0\ 1) \cup (1\ 1\ 1\ 1). \\
S(Y_{k+1})M^T &= [S(Y_1) \cup S(Y_2)]\ [M_1^T \cup M_2^T] \\
&= (12, 13, 1, 3, 4, 15) \cup (8, 11, 1, 6, 9) \\
&= X_1^1 \cup X_2^1 \\
&= X_{k+2}. \\
S(X_{k+2}) &= (1\ 1\ 1\ 1\ 1\ 1) \cup (1\ 1\ 1\ 1\ 1). \\
S(X_{k+2})M &= S(Y_1^1)\, M_1 \cup S(Y_2^1)\, M_2 \\
&= (9, 15, 11, 3, 13) \cup (13, 5, 11, 9) \\
&= Y_1^1 \cup Y_2^1 \\
&= Y_{k+3}. \\
S(Y_{k+3}) &= S(Y_1^1) \cup S(Y_2^1) \\
&= (1\ 1\ 1\ 1\ 1) \cup (1\ 1\ 1\ 1). \\
S(Y_{k+3})M^T &= S(Y_1^1)M_1^T \cup S(Y_2^1)M_2^T \\
&= (14, 15, -1, 0, 4, 18) \cup (8, 11, 1, 6, 9) \\
&= X_1^{11} \cup X_2^{11} \\
&= S(X_{k+4}) \\
&= (1\ 1\ 0\ 0\ 1\ 1) \cup (1\ 1\ 1\ 1\ 1). \\
S(X_{k+4}) &= S(X_1^{11})M_1 \cup S(X_2^{11})M_2 \\
&= (11, 13, 10, 8, 10) \cup (13, 5, 11, 10) \\
&= Y_1^{111} \cup Y_2^{111} \\
&= Y_{k+5} \\
S(Y_{k+5}) &= (1\ 1\ 1\ 1\ 1) \cup (1\ 1\ 1\ 1).
\end{aligned}
$$

Thus we see the equilibrium is achieved by the dynamical bisystem.

Now having worked with a input bivector from $F_{X_1} \cup F_{X_2}$ it is left for the reader to work with an input vector from the neuronal bifield $F_{Y_1} \cup F_{Y_2} = F_Y$.

Now having defined a model for a BBAM when two different experts give their opinion on the same problem, now we wish to study or build up a model when n experts give their opinion on the same problem but with varying attributes and scales of intervals ($n \geq 2$). Further we see when we solve any social problem which is dependent on the views of experts more sensitivity and accuracy is arrived only when we work with more number of experts, keeping this in mind we proceed on to work with a problem P on which n experts give their opinion.

Let the first expert wish to take the neuronal field $F_{X_1}$ as the attributes connected with the causes taken along the row of the synaptic matrix and the neuronal field $F_{Y_1}$ is taken as factors connected with the causes of the column of the synaptic connection matrix. Let us assume the first expert works with $m_1$ attributes along the rows and $n_1$ attributes along the columns in the $m_1 \times n_1$ matrix on the scale $[–a_1, a_1]$.

Similarly let the second expert with the attributes connected with neuronal field $F_{X_2}$ which form the rows of the synaptic connection matrix $m_2 \times n_2$ and $F_{Y_2}$ which are associated with the attributes of the neuronal field along the column. Let them give it on the interval $[–a_2, a_2]$ and so on.

Thus the $i^{th}$ expert works on the neuronal field $F_{X_i}$, the attributes of which form the row of the $m_i \times n_i$ matrix and for the column they are taken from the neuronal field $F_{Y_i}$. Thus the synaptic connection $m_i \times n_i$ matrix of the $i^{th}$ expert take their entries from the scale $[–a_i, a_i]$. This is true for i = 1, 2, …, n. Thus let us denote the neuronal n-field $F_X$ by $F_x = F_{X_1} \cup F_{X_2} \cup \ldots \cup F_{X_n}$ and the neuronal n-field $F_Y$ by $F_Y = F_{Y_1} \cup F_{Y_2} \cup \ldots \cup F_{Y_n}$. Let the n-matrix M be denoted by

$$M \quad = \quad M_1 \cup M_2 \cup \ldots \cup M_n$$

where $M_i$ is a $m_i \times n_i$ synaptic connection matrix of the $i^{th}$ expert with entries from $[-a_i, a_i]$, i = 1, 2, …, n. Any input n-vector X from $F_x$ will be of the form $X = X_1 \cup X_2 \cup \ldots \cup X_n$ where $X_i = (x_1^i, ..., x_{m_i}^i)$, i = 1, 2, …, n.

Any input n-vector Y from $F_Y$ will be of the form $Y = Y_1 \cup Y_2 \cup \ldots \cup Y_n$ where $Y_j = (y_1^j, y_2^j, ..., y_{n_j}^j)$, j = 1, 2, …, n. The binary state vector at $k^{th}$ time period, $S(X_k) = S(X_1) \cup S(X_2) \cup \ldots \cup S(X_n)$ where $S(X_j) = (a_1^j, ..., a_{m_i}^j)$ with $a_1^j$ 's zero or 1 for j = 1, 2, …, n.

$$
\begin{aligned}
S(X_k)M &= [S(X_1) \cup S(X_2) \cup \ldots \cup S(X_n)][M_1 \cup M_2 \cup \ldots \cup M_n] \\
&= S(X_1)M_1 \cup S(X_2)M_2 \cup \ldots \cup S(X_n)M_n. \\
&= Y_1 \cup Y_2 \cup Y_3 \cup \ldots \cup Y_n \\
&= Y_{k+1}, \\
S(Y_{k+1}) &= S(Y_1) \cup \ldots \cup S(Y_n).
\end{aligned}
$$

$$
\begin{aligned}
S(Y_{k+1})M^T &= [S(Y_1) \cup \ldots \cup S(Y_n)] [M_1^T \cup M_2^T \cup \ldots \cup M_n^T] \\
&= S(Y_1) M_1^T \cup S(Y_2) M_2^T \cup \ldots \cup S(Y_n) M_n^T \\
&= X_1^1 \cup X_2^1 \cup \ldots \cup X_n^1 \\
&= X_{k+2}^1
\end{aligned}
$$

and so on.

We call $M = M_1 \cup M_2 \cup \ldots \cup M_n$ as the synaptic connection n-matrix of the problem P given by n-experts. Clearly M is a mixed rectangular n-matrix. We denote this by BAnM model. Having seen the model for n-experts we now proceed on to find a suitable model when we have two sets of experts $p_1$, $p_2$ giving their views on the same problem P.

Let us now proceed on to build a new model. Let P be the problem under investigation. On this problem P we have $p_1$ experts using same number of $m_1$ attributes along the rows of the synaptic connection matrix and $n_1$ attributes along the columns of the synaptic connection matrix using the scale $[-a_1,$ $a_1$]. For the same problem P we have $p_2$ experts using the same number say $m_2$ ($m_2 \neq m_1$) along the rows of the synaptic connection matrix and $n_2$ attributes along the columns of the synaptic connection matrix on the scale [–$a_2$, $a_2$].

Thus we get $p_1$ number of $m_1 \times n_1$ synaptic connection matrices defined on the scale [–$a_1$, $a_1$] and $p_2$ number of $m_2 \times n_2$ synaptic connection matrices defined on the scale [–$a_2$, $a_2$].

Now for these two sets of synaptic connection matrices we want to find the interval of synaptic connection matrices we will work for the $p_1$ experts the results in case of $p_2$ experts can be adopted in an identical way.

Given $p_1$ number of synaptic $m_1 \times n_1$ connection matrices in the interval [–$a_1$, $a_1$], to find an interval of synaptic $m_1 \times n_1$ connection matrices on the same interval [–$a_1$, $a_1$]. Let $M_1^1, M_2^1, ..., M_{p_1}^1$ be the $p_1$ number of synaptic $m_1 \times n_1$ connection matrices for the problem P, defined on the interval [–$a_1$, $a_1$]. $M_1^1 = \left(m_{ij}^1\right)$, $M_2^1 = \left(m_{ij}^2\right)$, …, $M_{p_1}^1 = \left(m_{ij}^{p_1}\right)$ where $1 \leq i \leq m_1$, $1 \leq j \leq n_1$. Now to construct an interval of $m_1 \times n_1$ matrices, we have to have a minimal $m_1 \times n_1$ matrix $A_1$ using the $p_1$ number of synaptic connection $m_1 \times n_1$ matrices and a maximal $m_1 \times n_1$ matrix. $B_1$ using the same $p_1$, number of $m_1 \times n_1$ synaptic connection matrices. Let us set $A_1 = (a_{ij}^1)$ we have to find the values $(a_{ij}^1)$, $1 \leq i \leq m_1$ and $1 \leq j \leq n_1$. Set

$$a_{ij}^1 = \min\{m_{ij}^1, m_{ij}^2, ..., m_{ij}^{p_1}\}$$

$1 \leq i \leq m_1$ and $1 \leq j \leq n$.

$$a_{11}^1 = \min\{m_{11}^1, m_{11}^2, ..., m_{11}^{p_1}\}$$

$$a_{12}^1 = \min\{m_{12}^1, m_{12}^2, ..., m_{11}^{p_1}\}$$

and so on.

Thus $A_1$ is the minimal of all the synaptic connection matrices using the $p_1$ number of $m_1 \times n_1$ matrices, built using the interval [–$a_1$, $a_1$] i.e., $a_{ij}^1 \leq m_{ij}^t$ for $1 \leq t \leq p_1$, $1 \leq i \leq m_1$ and $1 \leq j \leq n_1$.

Now on similar lines we construct the maximal synaptic connection $m_1 \times n_1$ matrices $B_1$ using the $p_1$ number $m_1 \times n_1$

matrices $M_1^1, M_2^1, ..., M_{p_1}^1$. Let us put $B = (b_{ij}^1)$ now it is our task to find out $(b_{ij}^1)$. Define

$$b_{ij}^1 = \max \{m_{ij}^1, m_{ij}^2, ..., m_{ij}^{p_1}\}$$

for $1 \leq i \leq m_1$ and $1 \leq j \leq n_1$. That is

$$b_{11}^1 = \max\{m_{11}^1, m_{11}^2, ..., m_{11}^{p_1}\}$$

$$b_{21}^1 = \max\{m_{12}^1, m_{12}^2, ..., m_{12}^{p_1}\}$$

and so on.

Now $B_1 = (b_{ij}^1)$ is the maximal synaptic connection matrix for the matrices $M_t^1$, $1 \leq t \leq p_1$, built using the interval $[-a_1, a_1]$, $1 \leq i \leq m_1$ and $1 \leq j \leq n_1$. Thus we have $m_{ij}^t \leq b_{ij}^1$ for $1 \leq t \leq p_1$, $1 \leq i \leq m_1$ and $1 \leq j \leq n_1$. So we see any entry $m_{ij}^t$ in any of the synaptic connection matrices $(m_{ij}^t)$ given by the experts is such that $a_{ij}^1 \leq m_{ij}^t \leq b_{ij}^1$, $1 \leq t \leq p_1$, $1 \leq i \leq m_1$ and $1 \leq j \leq n_1$.

Thus the interval matrices $[A_1, B_1]$ will contain all the $p_1$ synaptic connection matrices.

Now having found the minimal and maximal elements of the $p_1$ synaptic connection matrices defined on the interval $[-a_1, a_1]$, we now proceed on to find the optimal synaptic connection matrix $O^1$. Define

$$O^1 = \frac{A_1 + B_1}{2} = \frac{(a_{ij}^1) + (b_{ij}^1)}{2},$$

$1 \leq i \leq m_1$ and $1 \leq j \leq n_1$.

Clearly $O^1 \in [A_1, B_1]$. Now we build the average of the synaptic connection matrices $M_1^1, M_2^1, ..., M_{p_1}^1$. Let us denote the average synaptic connection matrix of the $p_1$ experts by $\overline{M}^1$

$$\overline{M}^1 = \frac{(m_1^1 + m_2^1 + ... + m_{p_1}^1)}{p_1}$$

$$= \frac{(m_{ij}^1) + (m_{ij}^2) + ... + (m_{ij}^{p_1})}{p_1}$$

for $1 \leq i \leq m_1$ and $1 \leq j \leq n_1$. Now $\overline{M}^1 \in [A_1, B_1]$.

We define or call $[A_1, B_1]$ to be the interval of synaptic connection matrices i.e., IBAM model defined on the scale $[–a_1, a_1]$ of the $p_1$ experts on the problem P.

On similar lines we build the interval of synaptic connection matrices IBAM model; $[A_2, B_2]$ on the interval $[–a_2, a_2]$ using the views of the $p_2$ experts on the same problem P.

Here $A_2$ is the minimal synaptic connection matrix of the synaptic connection matrices $M_1^2, M_2^2, ..., M_{p_2}^2$ and $B_2$ is the maximal synaptic connection matrices of the $p_2$ experts. Here also

$$O^2 = \frac{A_2 + B_2}{2}$$

will be the optimal synaptic connection matrix of the interval of synaptic connection matrices $[A_2, B_2]$. $\overline{M}^2$ will be the average of the synaptic connection matrices of the $p_2$ experts $M_1^2, M_2^2, ..., M_{p_2}^2$.

Set $[A, B] = [A_1, B_1] \cup [A_2, B_2]$, we see elements, on the right are two interval matrices. Thus any element M in [A, B] will be of the form $M = M_1 \cup M_2$ where $M_1 \in [A_1, B_1]$ and $M_2 \in [A_2, B_2]$ i.e., M is a bimatrix. Define $A = A_1 \cup A_2$ and $B = B_1 \cup B_2$ which we choose to call as the minimal and maximal bimatrices of the bi-interval [A, B]. We define the optimal element O of [A, B] to be $O_1 \cup O_2$ and the average element to be $\overline{M} = \overline{M}^1 \cup \overline{M}^2$.

Thus we call [A, B] to be the bi-interval of the synaptic connection of bimatrices defined on the bi-interval $[–a_1, a_1] \cup [–a_2 \cup a_2]$ related with the problem P which is denoted as the BIBAM model.

Now having defined the notion of bi-interval of synaptic connection matrices we have to find a method to show how they function.

Now we know for the interval $[A_1, B_1]$ of synaptic connection of $m_1 \times n_1$ matrices of $p_1$ experts for the problem P defined on the interval $[–a_1, a_1]$ any input vector will be from the neuronal field $F_{X_1}$ and every vector in $F_{X_1}$ will be a $1 \times m_1$ matrix with entries from the interval [–$a_1$, $a_1$]. Like wise the input vectors related with the attributes taken along the column, i.e., any input vector will be from the neuronal field $F_{Y_1}$ and will be a 1 × $n_1$ matrices with entries from [–$a_1$, $a_1$]. We shall denote the interval of input vector from the neuronal field $F_{X_1}$, we denote the interval of 1 × $m_1$ matrices from the neuronal field $F_{X_1}$ with entries from the interval [–$a_1$, $a_1$] by $\left[X_1^1, X_2^1\right]$ where $X_1^1 = (-a_1, -a_1, ..., -a_1)$ the minimal 1 × $m_1$ input matrix of the interval of input matrices $\left[X_1^1, X_2^1\right]$ and $X_2^1$ = (+$a_1$, + $a_1$, …, + $a_1$) the maximal 1 × $m_1$ input matrix of the interval of input matrices $\left[X_1^1, X_2^1\right]$.

We see the synaptic connection matrix of the problem in general cannot work with any arbitrary input vector from the interval of input matrices $\left[X_1^1, X_2^1\right]$. For the dynamical system in general can recognize only binary vectors we say a 1 × $m_1$ matrix or a vector is a binary vector if the entries in that binary matrix is either 0 or 1.

But the collection of input vectors from the interval of 1 × $m_1$ matrices are in general not binary vectors we use the synchronizing function S. Suppose $X = (x_1^1, x_2^1, ..., x_{m_1}^1)$ input vector from the interval of matrices $(X_1^1, X_2^1)$.

Then $S(X) = S(x_1^1, x_2^1, ..., x_{m_1}^1)$ = $(t_1^1, t_2^1, ..., t_{m_1}^1)$. We know $(x_1^1, x_2^1, ..., x_{m_1}^1)$ is a 1 × $m_1$ matrix with entries from [–$a_1$, $a_1$]. We define $t_i^1$, $1 \leq i \leq m_1$ as follows. $t_i^1 = 0$ if $x_i^1 \leq 0$ and $t_j^1 = 1$ if $x_i^1 > 0$. Thus each $t_i^1$ for $1 \leq i \leq m_1$ will either be 1 or 0. Thus the value of the entries take in S(X) will either be 0 or 1.

This 1 × $m_1$ matrix $S(X) = (t_1^1, t_2^1, ..., t_{m_1}^1)$ is called the binary state vector for it denotes the ON of OFF state of the nodes, taken along the rows of the synaptic connection $m_1$ × $n_1$ matrix.

Now on similar lines if the input vectors are from the nodes of the column entries then we denote by $\left[Y_1^1, Y_2^1\right]$ the collection of all input 1 × $n_1$ matrices with entries from the interval [$-a_1$, $a_1$]. Thus $\left[Y_1^1, Y_2^1\right]$ denotes the interval of input 1 × $n_1$ matrices with entries from the scale [$-a_1$, $a_1$]. Here $Y_1^1$ = ($-a_1$, …, $-a_1$) is called the minimal element of the input vectors from the interval of 1 × $n_1$ matrices $\left[Y_1^1, Y_2^1\right]$. Similarly $Y_2^1 = \left(a_1, a_1, ..., a_1\right)$ is called as maximal input vector from the interval of 1 × $n_1$ matrices $\left[Y_1^1, Y_2^1\right]$. Thus any 1 × $n_1$ matrix P with entries from the scale [$-a_1$, $a_2$] will be in the interval of $\left[Y_1^1, Y_2^1\right]$, that is any 1 × $n_1$ matrices P ∈ $\left[Y_1^1, Y_2^1\right]$ if and only if $P = \left(p_1, ..., p_{n_1}\right)$, then $-a_1 \leq p_i \leq a_1$, i = 1, 2, …, $p_1$.

Thus any 1 × $n_1$ matrix from the interval of input vectors $\left[Y_1^1, Y_2^1\right]$ can be made into a binary state vector using the function S. i.e., if $Y = \left(y_1^1, y_2^1, ..., y_{n_1}^1\right)$ and Y ∈ $\left[Y_1^1, Y_2^1\right]$ then S(Y) = $\left(p_1^1, ..., p_{n_1}^1\right)$ where $p_i$ is either zero or one i.e., $p_i^1 = 0$ if $y_i^1 \leq 0$ and $p_j^1 = 1$ if $y_j^1 > 0$.

Thus the binary state vector gives the ON or OFF state of the nodes of attributes taken along the column of the synaptic connection matrix at the $k^{th}$ time period. Further it is to be noted that any dynamical system of synaptic connection matrices can only recognize the binary state vector of the input vector at the $k^{th}$ time period.

Now on similar lines for the synaptic connection $m_2 \times n_2$ matrices given by the second set of $p_2$ experts given by the interval of connection matrices [$A_2$, $B_2$]. Any input vector from the row attributes will be a 1 × $m_2$ row matrix / vector the entries are from [$-a_2$, $a_2$]. The input vector are from the neuronal field . We denote the interval of input 1 × $m_2$ vector from the interval field $F_{X_1}$ by $\left[X_1^2, X_2^2\right]$. Here the 1 × $m_2$ matrices from the interval of matrices $\left[X_1^2, X_2^2\right]$ take the values from the scale $[-a_2, a_2]$. Further $X_1^2$ is the minimal input vector of $\left[X_1^2, X_2^2\right]$ given by $[(-a_2, -a_2, \ldots, -a_2)]$ and the maximal input vector of $\left[X_1^2, X_2^2\right]$ is $X_2^2$ given by $X_2^2$ = $[(a_2, a_2, \ldots, a_2)]$. Any vector $X^1 \in \left[X_1^2, X_2^2\right]$ is a $1 \times m_2$ matrix taking its values from $[-a_2, a_2]$.

Similarly using the neuronal field $F_{Y_2}$ which contains the set of all input $1 \times n_2$ vectors / matrices from the scale $[-a_2, a_2]$.

Let $\left[Y_1^2, Y_2^2\right]$ denote the interval of input $1 \times n_2$ vectors / matrices with entries from $[-a_2, a_2]$. Any $1 \times n_2$ vector $Y^1 \in \left[Y_1^2, Y_2^2\right]$ will take its entries from $[-a_2, a_2]$. Here for this interval of $1 \times n_2$ matrices $\left[Y_1^2, Y_2^2\right]$ we see $Y_1^2$ = $(-a_2, -a_2, \ldots, -a_2)$ is the minimal $1 \times n_2$ matrix. Thus the triple $\left[X_1^1, X_2^1\right]$, $[A_1, B_1]$ $\left[Y_1^2, Y_2^2\right]$ will be the multi expert i.e., $p_1$ expert IBAM model. Similarly, $[X_1^2, X_2^2], [A_2, B_2], [Y_1^2, Y_2^2]$ will be the multi $p_2$-expert IBAM model for the problem P.

Now set $[A, B] = [A_1, B_1] \cup [A_2, B_2]$,

$$[X_1, X_2] = [X_1^1\ X_2^1] \cup [X_1^2, X_2^2]$$

and

$$[Y_1, Y_2] = [Y_1^1, Y_2^1] \cup [Y_1^2, Y_2^2],$$

then the triple $\{[A, B], [X_1, X_2], [Y_1, Y_2]\}$ is called as the synaptic connection interval of bimatrices together with the input vectors defined on the bintervals $[-a_1, a_1] \cup [-a_2, a_2]$ for the problem P given by the two sets of experts $p_1$, $p_2$ of the Interval of Bidirectoral Associative Memories Bimodel [IBBAM].

$$S([X_1, X_2])\ [A, B] = S\ [Y_1, Y_2]$$

or

$$S\ [Y_1, Y_2]\ [A, B]^T = S([X_1, X_2])$$

where S is the function which makes any input vector into a state binary vector.

i.e., $[S[X_1^1, X_2^1] \cup S[X_1^2, X_2^2][A_1, B_1] \cup [A_1, B_2]$

$$= \quad S[Y_1^1, Y_2^1] \cup S[Y_1^2, Y_2^2]$$

i.e., $S[X_1^1 X_2^1][A_1, B_1] \cup S[X_1^2, X_2^2][A_2, B_2]$

$$= \quad S[Y_1^1, Y_2^1] \cup [Y_1^2, Y_2^2]$$

We will show how the IBBAM bimodel functions by an illustration from the real data model.

Let us consider the problem of analyzing the factors forcing people for migration and the role of government in the context of HIV/AIDS infected poor migrant labourers.

The attributes taken by the first set of experts is as follows.

F – Factors forcing people for migration.

$F_1$ – Lack of labour opportunities in their hometown

$F_2$ – Poverty seeking better status of life

$F_3$ – Mobilization of labour contractors

$F_4$ – Infertility of lands due to implementation of wrong research methodologies /failure of monsoon.

The attributes given by this first set of experts regarding the role of government.

$G_1$ – Alternate employment if the harvest fails there by stopping migration

$G_2$ – Awareness clubs in rural areas about HIV/AIDS

$G_3$ – Construction of hospitals in rural areas with HIV/AIDS counseling cell/ compulsory HIV/AIDS test before marriage

$G_4$ – Failed to stop the misled agricultural techniques followed recently by farmers

$G_5$ – No foresight for the government and no precautionary actions taken from the past occurrences.

Taking the neuronal field $F_{X_1}$ as the role of government and the neuronal field $F_{Y_1}$ as the attributes connected with, forcing people for migration we form the synaptic projections from $F_{X_1}$ to $F_{Y_1}$.

Taking along the columns of the synaptic matrix the factors forcing people for migration and along the rows of the synaptic matrix the factors related with the role of government. In this problem we have just 3 experts giving their opinion.

Let $M_1^1$ denote the synaptic connection matrix related with the first expert on the scale [–4, 4]

$$M_1^1 = \begin{array}{c} \\ G_1 \\ G_2 \\ G_3 \\ G_4 \\ G_5 \end{array}\begin{array}{c} \begin{array}{cccc} F_1 & F_2 & F_3 & F_4 \end{array} \\ \begin{bmatrix} 3 & 2 & 3 & 2 \\ 0 & -2 & 0 & 1 \\ 4 & -4 & 4 & 0 \\ 2 & 4 & -1 & 4 \\ 3 & 3 & -2 & 3 \end{bmatrix} \end{array}.$$

Let $M_2^1$ denote the synaptic connection matrix given by the second expert on the same set of attributes and on the same scale [–4, 4]

$$M_2^1 = \begin{array}{c} \\ G_1 \\ G_2 \\ G_3 \\ G_4 \\ G_5 \end{array}\begin{array}{c} \begin{array}{cccc} F_1 & F_2 & F_3 & F_4 \end{array} \\ \begin{bmatrix} 4 & 3 & 2 & 2 \\ 0 & -3 & 0 & 0 \\ 4 & -4 & 3 & 0 \\ 3 & 3 & -2 & 4 \\ 4 & 2 & -2 & 3 \end{bmatrix} \end{array}.$$

Now the 3$^{rd}$ expert also choose to work on the same collection of attributes and on the same scale [–4, 4]. Let $M_3^1$ denote the synaptic connection matrix of the third expert.

$$M_3^1 = \begin{array}{c} \\ G_1 \\ G_2 \\ G_3 \\ G_4 \\ G_5 \end{array} \begin{array}{c} \begin{array}{cccc} F_1 & F_2 & F_3 & F_4 \end{array} \\ \begin{bmatrix} 4 & 3 & 3 & 2 \\ 0 & -3 & -1 & 1 \\ 3 & -4 & 4 & 0 \\ 3 & 4 & -2 & 4 \\ 4 & 3 & -3 & 4 \end{bmatrix} \end{array}.$$

Now on the same problem four experts from the second set of experts wish to give their opinion and they choose to change the number of attributes as well the scale from [–4, 4] to [–5, 5].

They work with the same number of column attributes viz. $F_{Y_2}$ but this second set of experts have coupled the attributes $G_2$ and $G_3$ as a single node. Thus they work only with four attributes related with the government role in the problems of migrant labourers and their vulnerability to HIV/AIDS.

Now let $M_i^2$ represent the opinion of the four experts from the second set i.e., i = 1, 2, 3, 4. Let the synaptic projection matrices from $F_{X_2}$ to $F_{Y_2}$ be denoted by $M_1^2$, $M_2^2$, $M_3^2$ and $M_4^2$. Thus as in case of the first set of experts the second set of experts also take the role of government along the rows of the synaptic connection matrix and the factors forcing migration along the columns of the synaptic connection matrix. This second set of experts work with the scale [–5, 5].

Let the synaptic connection matrix given by the first expert of the second set be given by $M_1^2$

$$M_1^2 = \begin{array}{c} \\ G_1^1 \\ G_2^1 \\ G_3^1 \\ G_4^1 \end{array} \begin{array}{c} \begin{array}{cccc} F_1 & F_2 & F_3 & F_4 \end{array} \\ \begin{bmatrix} 4 & 3 & 3 & 2 \\ 3 & -4 & 2 & -1 \\ 3 & 5 & -3 & 5 \\ 5 & 3 & -3 & 4 \end{bmatrix} \end{array}.$$

Now we give the synaptic connection matrix $M_2^2$, the views of the second expert from the second set of experts on the same scale [–5, 5].

$$M_2^2 = \begin{array}{c} \\ G_1^1 \\ G_2^1 \\ G_3^1 \\ G_4^1 \end{array} \begin{array}{c} \begin{array}{cccc} F_1 & F_2 & F_3 & F_4 \end{array} \\ \begin{bmatrix} 4 & 3 & 4 & 3 \\ 4 & -4 & 3 & 0 \\ 3 & 5 & -3 & 5 \\ 5 & 3 & -4 & 5 \end{bmatrix} \end{array}.$$

Let $M_3^2$ correspond to the opinion of the 3rd expert from the second set experts on the same scale [–5, 5] with the attributes $G_1^1, G_2^1, G_3^1$ and $G_4^1$ along the row and $F_1$, $F_2$, $F_3$ and $F_4$ along the column.

$$M_3^2 = \begin{array}{c} \\ G_1^1 \\ G_2^1 \\ G_3^1 \\ G_4^1 \end{array} \begin{array}{c} \begin{array}{cccc} F_1 & F_2 & F_3 & F_4 \end{array} \\ \begin{bmatrix} 5 & 3 & 3 & 3 \\ 0 & -4 & 2 & 1 \\ 4 & 4 & -3 & 5 \\ 5 & 3 & -4 & 4 \end{bmatrix} \end{array}.$$

Let $M_4^2$ denote the synaptic connection matrix given by the forth expert on the same set of attributes and on the same scale [–5, 5].

$$M_4^2 = \begin{array}{c} \\ G_1^1 \\ G_2^1 \\ G_3^1 \\ G_4^1 \end{array} \begin{array}{c} \begin{array}{cccc} F_1 & F_2 & F_3 & F_4 \end{array} \\ \begin{bmatrix} 4 & 4 & 3 & 2 \\ 2 & -4 & 2 & 0 \\ 4 & 4 & -4 & 4 \\ 5 & 4 & -3 & 3 \end{bmatrix} \end{array}.$$

Now we have to give a interval of bimatrix model to this problem.

To achieve this first we denote how we work to obtain the interval of synaptic connection matrices associated with the first set of experts. For the second set of experts we just produce the values of the interval of synaptic connection of matrices.

We are given the synaptic connection matrices $M_1^1, M_2^1$ and $M_3^1$ of the first set of experts. Now we will construct the interval of synaptic connection of matrices $[A_1, B_1]$ on the scale [–4, 4].

Clearly if $A_1 = (a_{ij}^1)$ is the minimal element of $[A_1, B_1]$ then we fill in the elements $(a_{ij}^1)$ as follows.

$$(a_{ij}^1) = \min(m_{ij}^1, m_{ij}^2, m_{ij}^3)$$

where

$$M_1^1 = (m_{ij}^1),\ M_2^1 = (m_{ij}^2) \text{ and } M_3^1 = (m_{ij}^3);$$

$1 \le i \le 5$ and $1 \le j \le 4$.

Thus

$$A_1 = \begin{array}{c} \\ G_1 \\ G_2 \\ G_3 \\ G_4 \\ G_5 \end{array} \begin{array}{c} \begin{array}{cccc} F_1 & F_2 & F_3 & F_4 \end{array} \\ \begin{bmatrix} 3 & 2 & 2 & 2 \\ 0 & -3 & -1 & 0 \\ 3 & -4 & 4 & 0 \\ 2 & 3 & -2 & 4 \\ 3 & 2 & -3 & 3 \end{bmatrix} \end{array}$$

is defined to be the minimal synaptic connection matrix of the interval of matrices related to the synaptic connection matrices $M_1^1, M_2^1$ and $M_3^1$ we see $a_{ij}^1 \le m_{ij}^t$, t = 1, 2, 3, $1 \le i \le 5$ and $1 \le j \le 4$.

Let $B_1$ denote the maximal element of the synaptic connection matrices $M_1^1, M_2^1$ and $M_3^1$.

Let $B_1 = (b_{ij}^1)$, $1 \le i \le 5$ and $1 \le j \le 4$. Now how do the elements of the maximal matrix $B_1$ defined.

$$b_{ij}^1 = \max\{m_{ij}^1, m_{ij}^2, m_{ij}^3\},$$

$1 \le i \le 5$ and $1 \le j \le 4$.

Thus

$$B_1 = \begin{array}{c} \\ G_1 \\ G_2 \\ G_3 \\ G_4 \\ G_5 \end{array} \begin{array}{c} \begin{array}{cccc} F_1 & F_2 & F_3 & F_4 \end{array} \\ \begin{bmatrix} 4 & 3 & 3 & 2 \\ 0 & -2 & 0 & 1 \\ 4 & -4 & 4 & 0 \\ 3 & 4 & -1 & 4 \\ 4 & 3 & -2 & 4 \end{bmatrix} \end{array}$$

Clearly $b^1_{ij} \geq m^t_{ij}$ , t = 1, 2, 3. Thus we call $[A_1, B_1]$ to be the interval of synaptic connection of matrices $M^1_1, M^1_2$ and $M^1_3$ on the scale [–4, 4].

Now we define the optimal synaptic connection matrix $O^1$ for the interval of synaptic connection of matrices by.

$$O_1 = \frac{A_1 + B_1}{2}$$

i.e.,

$$= \frac{(a^1_{ij}) + (b^1_{ij})}{2}.$$

$$O^1 = \begin{array}{c} \\ G_1 \\ G_2 \\ G_3 \\ G_4 \\ G_5 \end{array} \begin{array}{c} \begin{array}{cccc} F_1 & F_2 & F_3 & F_4 \end{array} \\ \begin{bmatrix} 3.5 & 2.5 & 2.5 & 2 \\ 0 & -2.5 & -0.5 & -0.5 \\ 3.5 & -4 & 4 & 0 \\ 2.5 & 3.5 & -1.5 & 4 \\ 3.5 & 2.5 & -2.5 & 3.5 \end{bmatrix} \end{array}.$$

Clearly $O^1 \in [A_1, B_1]$, it is the optimal synaptic connection matrix of the interval of synaptic connection matrices $[A_1, B_1]$. Now we calculate the average/mean of the synaptic connection matrices $\overline{M}^1$ in $[A_1, B_1]$.

$$\overline{M}^1 = \frac{M^1_1 + M^1_2 + M^1_3}{3}$$

$$= \frac{(m_{ij}^1) + (m_{ij}^2) + (m_{ij}^2)}{3}$$

with $1 \le i \le 5$ and $1 \le j \le 4$ on the same interval [–4, 4]. Clearly $\overline{M}^1 \in [A_1, B_1]$.

$$\overline{M}^1 = \begin{array}{c} \\ G_1 \\ G_2 \\ G_3 \\ G_4 \\ G_5 \end{array}\begin{array}{c} \begin{array}{cccc} F_1 & F_2 & F_3 & F_4 \end{array} \\ \begin{bmatrix} 3.7 & 2.7 & 2.7 & 2 \\ 0 & -2.7 & -0.33 & 0.66 \\ 3.7 & -4 & 3.7 & 0 \\ 2.7 & 3.7 & -1.7 & 4 \\ 3.7 & 2.7 & -2.3 & 3.3 \end{bmatrix} \end{array}.$$

Thus the interval of synaptic connection matrices $[A_1, B_1]$ contains the matrices $A_1$, $B_1$, $O^1$, $\overline{M}^1, M_1^1, M_2^1$ and $M_3^1$ defined on the scale [–4, 4], related to the first set of experts.

Now we just give the interval of synaptic connection of matrices $[A_2, B_2]$ using the synaptic connection matrices of 4 experts on the same problem defined on the interval [–5, 5].

Let $M_t^2 = [m_{ij}^t]$, t = 1, 2, 3, 4, $1 \le i \le 4$ and $1 \le j \le 4$.

$$A_2 = \begin{array}{c} \\ G_1^1 \\ G_2^1 \\ G_3^1 \\ G_4^1 \end{array}\begin{array}{c} \begin{array}{cccc} F_1 & F_2 & F_3 & F_4 \end{array} \\ \begin{bmatrix} 4 & 3 & 3 & 2 \\ 0 & -4 & 2 & -1 \\ 3 & 4 & -4 & 4 \\ 5 & 3 & -4 & 3 \end{bmatrix} \end{array}.$$

$A_2 = (a_{ij}^2)$ is defined as

$$a_{ij}^2 = \min\{m_{ij}^1, m_{ij}^2, m_{ij}^3, m_{ij}^4\};$$

$1 \le i \le 4$ and $1 \le j \le 4$.

Further $a_{ij}^2 \le m_{ij}^t$, t = 1, 2, 3, 4; $1 \le i \le 4$ and $1 \le j \le 4$.

Thus $A_2$ is the minimal synaptic connection matrix of the interval containing the synaptic connection matrices $M_1^2, M_2^2, M_3^2$ and $M_4^2$.

Let $B_2$ be the maximal synaptic connection matrix of the interval containing the synaptic connection matrices. $M_1^2, M_2^2, M_3^2$ and $M_4^2$.

Let us define $B_2 = (b_{ij}^2)$; $1 \leq i \leq 4$, $1 \leq j \leq 4$, where

$$b_{ij}^2 = \max\{m_{ij}^1, m_{ij}^2, m_{ij}^3, m_{ij}^4\}$$

where $1 \leq i \leq 4$ and $1 \leq j \leq 4$.

$$B_2 = \begin{array}{c} \\ G_1^1 \\ G_2^1 \\ G_3^1 \\ G_4^1 \end{array} \begin{array}{c} \begin{array}{cccc} F_1 & F_2 & F_3 & F_4 \end{array} \\ \begin{bmatrix} 5 & 4 & 4 & 3 \\ 4 & -4 & 3 & 1 \\ 4 & 5 & -3 & 5 \\ 5 & 4 & -3 & 5 \end{bmatrix} \end{array}.$$

$B_2$ is the maximal synaptic connection matrix of the interval of synaptic matrices $[A_2, B_2]$. Further $a_{ij}^2 \leq m_{ij}^t \leq b_{ij}^2$, $t = 1, 2, 3, 4$; $1 \leq i \leq 4$ and $1 \leq j \leq 4$, where $(m_{ij}^t) = M_t^2$; $t = 1, 2, 3, 4$.

Now as in case of the first set of experts we define the optimal synaptic connection matrix

$$O^2 = \frac{(A_2 + B_2)}{2} = \left(\frac{a_{ij}^2 + b_{ij}^2}{2}\right);$$

$1 \leq i \leq 4$; $1 \leq j \leq 4$.

$$O^2 = \begin{array}{c} \\ G_1^1 \\ G_2^1 \\ G_3^1 \\ G_4^1 \end{array} \begin{array}{c} \begin{array}{cccc} F_1 & F_2 & F_3 & F_4 \end{array} \\ \begin{bmatrix} 4.5 & 3.5 & 3.5 & 2.5 \\ 2 & -4 & 2.5 & 0 \\ 3.5 & 4.5 & -3.5 & 4.5 \\ 5 & 3.5 & -3.5 & 4 \end{bmatrix} \end{array}$$

$O^2 \in [A_2, B_2]$ and is the optimal synaptic connection matrix.

Now let $\overline{M}^2$ denotes the mean of the synaptic connection matrices given by the four experts from the second set ;

$$\overline{M}^2 = \frac{M_1^2 + M_2^2 + M_3^2 + M_4^2}{4}$$

$$= \left( \frac{(m_{ij}^1) + (m_{ij}^2) + (m_{ij}^3) + (m_{ij}^4)}{4} \right).$$

$$\overline{M}^2 = \begin{array}{c} \\ G_1^1 \\ G_2^1 \\ G_3^1 \\ G_4^1 \end{array} \begin{array}{c} \begin{array}{cccc} F_1 & F_2 & F_3 & F_4 \end{array} \\ \begin{bmatrix} 4.25 & 3.25 & 3.25 & 2.50 \\ 4.25 & -4 & 2.25 & -0.5 \\ 4.25 & 4.5 & -3.25 & 4.25 \\ 4.75 & 4.25 & -3.25 & 4 \end{bmatrix} \end{array}.$$

$\overline{M}^2 \in [A_2, B_2]$. Further we see the interval $[A_2, B_2]$ has 8 matrices 4 of them are synaptic connection matrices given by the four experts and the remaining, the related synaptic connection matrices.

Further we see that $F_{X_2}$ and $F_{Y_2}$ are the neuronal fields related with the interval of synaptic connection matrices defined on the interval [–5, 5]. $F_{X_2}$ and $F_{Y_2}$ are also defined on the same scale [–5, 5].

Further we see that the interval $[A_1, B_1]$ corresponds to the same problem but with a very different set of experts and also the different set of attributes and also with a different scale. This interval of synaptic connection matrices works with 3 experts with different of set experts and on the scale [–4, 4]. As already said the input vectors belong to neuronal fields $F_{X_1}$ and $F_{Y_1}$ which are related to the interval of synaptic connection matrices defined on the interval [–4, 4].

Let $[A, B] = [A_1, B_1] \cup [A_2, B_2]$ be the collection of all synaptic connection bimatrices M such that $M = M_1 \cup M_2$ where $M_i \in [A_i, B_i]$, i = 1, 2. [A, B] is know as the interval of synaptic connection of bimatrices related with the problem on the binterval [–4, 4] ∪ [–5, 5]. Clearly the input bivector X = $X_1$ ∪ $X_2$ is such that $X_1 \in F_{X_1}$ and $X_2 \in F_{X_1}$ and X ∈ $F_{X_1} \cup F_{X_2}$. Likewise $F_Y = F_{Y_1} \cup F_{Y_2}$ is the collection of all input bivectors defined on the bi-interval [–4, 4] ∪ [–5, 5]. This model is defined as the interval Bimatrix of bi-directional associative memories [IBBAM].

Now we will just illustrate how a single input bivector works on the interval [A, B] of synaptic connection bimatrices; defined on the bi-interval [–4, 4] ∪ [–5, 5]. Here A = $A_1 \cup A_2$, B = $B_1 \cup B_2$, O = $O^1 \cup O^2$ and $\overline{M} = \overline{M}^1 \cup \overline{M}^2$ denotes the minimal, maximal optimal and average/mean synaptic connection bimatrix of the interval of bimatrices [A, B] defined on the bi-interval [–4, 4] ∪ [–5, 5] respectively.

Suppose the experts uniformly agree to work on the input bivector at the k^th time period $X_k = X_1 \cup X_2$ ($X_1 \in F_{X_1}$ and $X_2 \in F_{X_2}$) i.e., X = (3, –0.2, 0, –1, 2) ∪ (1, –1, –4, 0).

Using the function S we first transform the input bivector into the binary bivector so that the bisystem IBBAM can recognize the input vector.

Let us take A ∈ [A, B], i.e., A = $A_1 \cup A_2$ ($A_1 \in [A_1, B_1]$ and $A_2 \in [A_2, B_2]$ the minimal matrices of the interval of matrices i.e.

$$A = \begin{array}{c} \\ G_1 \\ G_2 \\ G_3 \\ G_4 \\ G_5 \end{array}\!\!\begin{array}{c} \begin{array}{cccc} F_1 & F_2 & F_3 & F_4 \end{array} \\ \begin{bmatrix} 3 & 2 & 2 & 2 \\ 0 & -3 & -1 & 0 \\ 3 & -4 & 4 & 0 \\ 2 & 3 & -2 & 4 \\ 3 & 2 & -3 & 3 \end{bmatrix} \end{array} \cup \begin{array}{c} \\ G_1^1 \\ G_2^1 \\ G_3^1 \\ G_4^1 \end{array}\!\!\begin{array}{c} \begin{array}{cccc} F_1 & F_2 & F_3 & F_4 \end{array} \\ \begin{bmatrix} 4 & 3 & 3 & 2 \\ 0 & -4 & 2 & -1 \\ 3 & 4 & -4 & 4 \\ 5 & 3 & -4 & 3 \end{bmatrix} \end{array}.$$

| | | |
|---|---|---|
| $S(X_k)$ | = | $S(X_1) \cup S(X_2)$ |
| | = | $S(X_1 \cup X_2)$ |
| | = | S[(3, –2, 0, –1, 2)] ∪ S[(1, –1, –4, 0)] |
| | = | (1 0 0 0 1) ∪ (1 0 0 0). |
| $S(X_k)A$ | = | (6, 4, –1, 5) ∪ (4, 3, 3, 2) |

$$
\begin{aligned}
& = && Y_1 \cup Y_2 \\
& = && Y_{k+1} \\
S(Y_{k+1}) & = && (1\ 1\ 0\ 1) \cup (1\ 1\ 1\ 1). \\
S(Y_{k+1})A^T & = && (7, –3, –4, 0, 8) \cup (12, –3, 7, 7) \\
& = && X_1^1 \cup X_2^1 \\
& = && X_{k+2} \\
S(X_{k+2}) & = && (1\ 0\ 0\ 0\ 1) \cup (1\ 0\ 1\ 1). \\
S(X_{k+2})A & = && (6, 4, –1, 5) \cup (12, 10, –5, 9) \\
& = && Y_1^1 \cup Y_2^1 \\
& = && Y_{k+3}.
\end{aligned}
$$

Now

$$
\begin{aligned}
S(Y_{k+3}) & = && (1\ 1\ 0\ 1) \cup (1\ 1\ 0\ 1). \\
S(Y_{k+3})A^T & = && (7, –3, –4, 0, 8) \cup (9, –5, 11, 11) \\
& = && X_{k+4} \\
& = && X_1^{11} \cup X_2^{11} \\
S(X_{k+4}) & = && (1\ 0\ 0\ 0\ 1) \cup (1\ 0\ 1\ 1).
\end{aligned}
$$

Thus we see the resultant is a fixed binary bipair. One can work with any input vector at any desired time period k, from the neuronal bifield $F_Y = F_{Y_1} \cup F_{Y_2}$ .

Further as we have only just venutured to illustrate with examples the new IBBAM model, we do not proceed to work with more and more input vectors or with various synaptic connection bimatrices of the interval of bimatrices [A, B].

Next we proceed on to work with the interval of n-matrices n ≥ 3; [n = 2 will be the interval of bimatrices].

Let us first define the notion of interval of n-matrices

**DEFINITION 3.7.1:** *Suppose we have n-intervals of matrices say $[A_1, B_1]$, $[A_2, B_2]$, …, $[A_n, B_n]$ where each $[A_i, B_i]$ is an interval of $n_i \times m_i$ matrix $n_i \neq n_j$ if $(i \neq j)$; for i = 1, 2, …, n. We define the interval of n-matrices $[A, B] = [A_1, B_1] \cup [A_2, B_2] \cup \ldots \cup [A_n, B_n]$ where [A, B] is the collection of all n-matrices i.e., if $M \in [A, B]$ then $M = M_1 \cup M_2 \cup \ldots \cup M_n$, where $M_i \in [A_i, B_i]$; i = 1, 2, …, n, $M_i \neq M_j$ if $i \neq j$ with $A = A_1 \cup A_2 \cup \ldots \cup A_n$ and $B = B_1 \cup B_2 \cup \ldots \cup B_n$.*

So when one speaks of interval of n-matrices [A, B] then one has the above conditions to be true; we just illustrate this by a very simple example.

***Example 3.7.1:*** Let [A, B] = $[A_1, B_1] \cup [A_2, B_2] \cup [A_3, A_3] \cup [A_4, B_4]$ be the interval of 4-matrices.

$[A_1, B_1]$ = $\{[a_1, a_2, a_3]$ collection of all $1 \times 3$ matrices with entries from $Z^+\}$;

$$[A_2, B_2] = \left\{\begin{pmatrix} a & b \\ c & d \end{pmatrix} / a, b, c, d \in Z\right\};$$

$$[A_3, B_3] = \left\{\begin{bmatrix} a_1 \\ a_2 \\ a_3 \\ a_4 \end{bmatrix} / a_i \in Q; i = 1, 2, 3, 4\right\} \text{ and}$$

$$[A_4, B_4] = \left\{\begin{pmatrix} a_{11} & a_{12} & a_{13} & a_{14} \\ a_{21} & a_{22} & a_{23} & a_{24} \\ a_{31} & a_{32} & a_{33} & a_{34} \end{pmatrix} / a_{ij} \in Q; 1 \leq i \leq 3 \text{ and } 1 \leq j \leq 4\right\}.$$

Now any typical element in [A, B] will be of the form;

$M = M_1 \cup M_2 \cup M_3 \cup M_4$, where $M_i \in [A_i, B_i]$, i = 1, 2, 3, 4.

$$M = [3, 5, 1] \cup \begin{bmatrix} 0 & -3 \\ 2 & 9 \end{bmatrix} \cup \begin{bmatrix} 0.5 \\ 3 \\ -2 \\ 17/5 \end{bmatrix} \cup \begin{bmatrix} -3 & 2 & 1 & 0 \\ 4 & 6 & -5 & 2 \\ 3/7 & 0 & -6 & 7 \end{bmatrix}.$$

Now we see each matrix in the interval of matrices are different form one interval of matrices is a row matrix, another interval is a square matrix and one interval of matrix is a column matrix and another a rectangular matrix. Thus we have several types of matrices to be present in an interval of n-matrices unlike the interval of matrices.

How to distinguish them? We define for this these types of interval of matrices.

We say an interval of n-matrices [A, B] = [$A_1$, $B_2$] ∪ [$A_2$, $B_2$] ∪ … [$A_n$, $B_n$] to be a mixed square interval of n-matrices if each of the interval of matrices [$A_i$, $B_i$]; i = 1, 2, …, n is a $m_i \times m_i$ square matrix. It is to be noted if i ≠ j then $M_i \neq M_j$. A very simple illustration would be.

***Example 3.7.2:*** Let [A, B] = [$A_1$, $B_1$] ∪ [$A_2$, $B_2$] ∪ [$A_3$, $B_3$] ∪ [$A_4$, $B_4$] where

[$A_1$, $B_1$] is the interval of 2 × 2 matrices with entries from Q.
[$A_2$, $B_2$] is the set of all 5 × 5 matrices with entries from $Z^+$,
[$A_3$, $B_3$] is the set of all 3 × 3 matrices with entries from $Z_{10}$ (ring of integers modulo 10) and
[$A_4$, $B_4$] is the collection of all 6 × 6 matrices with entries from $Z_2$ the prime field of characteristic 2.

Thus any element M in [A, B] would be of the form, M = $M_1 \cup M_2 \cup M_3 \cup M_4$.

$$M = \begin{bmatrix} 2 & -5 \\ 7/3 & 4 \end{bmatrix} \cup \begin{bmatrix} 9 & 1 & 2 & 4 & 3 \\ 5 & 6 & 1 & 1 & 2 \\ 1 & 7 & 2 & 2 & 1 \\ 2 & 8 & 3 & 9 & 20 \\ 5 & 12 & 4 & 5 & 12 \end{bmatrix} \cup$$

$$\begin{bmatrix} 9 & 0 & 1 \\ 2 & 1 & 3 \\ 1 & 7 & 5 \end{bmatrix} \cup \begin{bmatrix} 1 & 0 & 1 & 0 & 0 & 1 \\ 1 & 1 & 1 & 1 & 0 & 0 \\ 0 & 1 & 1 & 1 & 0 & 1 \\ 0 & 1 & 1 & 0 & 1 & 1 \\ 1 & 1 & 1 & 0 & 0 & 0 \\ 0 & 0 & 0 & 1 & 1 & 1 \end{bmatrix}$$

We define an interval of n-matrices to be a mixed rectangular n-matrices if each of the matrices in each of the interval $[A_i, B_i]$ where $[A, B] = [A_1, B_1] \cup [A_2, B_2] \cup \ldots \cup [A_n, B_n]$, for i = 1, 2, …, n is a $m_i \times n_i$ rectangular matrix; clearly if $i \neq j$ then $m_i \neq m_j$. As in case of mixed square interval of n-matrices one can construct examples.

We call an interval of n-matrices $[A, B] = [A_1, B_1] \cup [A_2, B_2] \cup \ldots \cup [A_n, B_n]$ to be an interval of mixed n-matrices if the interval of matrices $[A_i, B_i]$ is either a square $m_i \times m_i$ matrix and $[A_j, B_j]$ is a $t_j \times s_j$ matrix $t_j \neq s_j$. $(i \neq j)$ $1 \leq i, j \leq n$.

We call $[A\ B] = [A_1, B_1] \cup [A_2, B_2] \cup \ldots \cup [A_n, B_n]$ to be an interval of mixed row matrices if each of the interval of matrices in $[A_i, B_i]$ is a row matrix.

Like wise we call $[A, B] = [A_1, B_1] \cup \ldots \cup [A_n, B_n]$ is an interval of mixed column matrix if each of the matrices in each of the interval of matrices $[A_i, B_i]$ is a column matrix for i = 1, 2, …, n.

***Note:*** Some times all the matrices in these types of matrices can be matrices of same order in $[A, B] = [A_1, B_1] \cup \ldots \cup [A_n, B_n]$ but can be defined on different domains then how to distinguish them.

Just we illustrate it with an example before we give the essential definition.

***Example 3.7.3:*** Let $[A, B] = [A_1, B_1] \cup [A_2, B_2] \cup \ldots \cup [A_5, B_5]$ where $[A_1, B_1]$ = set of all $2 \times 2$ matrices with entries from the interval [0, 5]. $[A_2, B_2]$ is the set of $2 \times 2$ matrices with entries from $Z_2$ ring of integers modulo 2. $[A_3, B_3]$ is the set of all $2 \times 2$ matrices with entries from the interval [–1, 1], $[A_4, B_4]$ is the set of all $2 \times 2$ matrices with entries from $Z_5$ ring of integers modulo 5. $[A_5, B_5]$ is the set of all $2 \times 2$ matrices with entries from $Q \setminus \{0\ 2\}$.

We see all matrices in each of the interval of matrices $[A_i, B_i]$; i = 1, 2, 3, 4, 5 are only $2 \times 2$ matrices with entries from very different intervals.

We define them as follows:

**DEFINITION 3.7.2:** *Let [A, B] = [$A_1$, $B_1$] $\cup$ [$A_2$, $B_2$] $\cup \ldots \cup$ [$A_n$, $B_n$]. If all the n interval of matrices [$A_i$, $B_i$] are $m \times m$ matrices defined over different intervals then we call [A, B] to be interval of $m \times m$ square n-matrices. Like wise we can define interval of $m \times p$ rectangular n-matrices and so on.*

Now having defined the interval of n-matrices we proceed on to define the nIBAM model $n \geq 3$.

Let us assume we have a problem P on hand, which mainly deals with an unsupervised data. Several experts are interested to give their solutions as they are involved with it in one way or the other. Further they have grouped them into say t-groups to work on it. Thus we have sets of t-groups of experts working on the problem P. Each of the experts in each of the groups $g_1, \ldots, g_t$ agree upon to work with some $m_i$ attributes from the fuzzy neuronal field $F_{X_i}$ and $n_i$ attributes from the fuzzy neuronal field $F_{Y_i}$ on the scale $[-a_i, a_i]$; i = 1, 2, …, n.

Suppose the group $g_i$ has say $r_i$ experts. This is true for i = 1, 2, …, t. So the group $g_1$ has $r_1$ experts and they will work with $m_1$ attributes from the neutronal field $F_{X_1}$ and $n_1$ attributes from the neuronal field $F_{Y_1}$ on the scale $[-a_1, a_1]$. Let the association of the attributes given by any expert from group $g_1$ from the neutronal field $F_{X_1}$ to $F_{Y_1}$ be denoted by the synaptic connection matrices $M_s^1$, where $1 \leq s \leq r_1$, on the scale $[-a_1, a_1]$.

Thus we have the $r_1$ experts giving the $r_1$ number of synaptic connection matrices denoted by $M_1^1, M_2^1, \ldots, M_{r_1}^1$. Now we will construct the interval of these $r_1$ number of $m_1 \times n_1$ synaptic connection matrices given by these $r_1$ experts. To find this we have to first find an interval say [$A_1$, $B_1$] where $A_1$ corresponds to the constructed minimal synaptic connection matrix using these $r_1$ experts and $B_1$ corresponds to the constructed maximal synaptic connection matrices related with the views of the $r_1$ experts.

Define $A_1 = (a^1_{ij})$ now we have to find the values of $(a^1_{ij})$ in the interval $[–a_1, a_1]$ using the $r_1$ synaptic connection matrices. Define

$$a^1_{ij} = \min\{m^1_{ij}, m^2_{ij}, ..., m^{r_1}_{ij}\};$$

$1 \leq i \leq m_1$ and $1 \leq j \leq n_1$; with $m^1_t = (m^t_{ij})$; $1 \leq t \leq r_1$ i.e.,

$$a^1_{11} = \min\{m^1_{11}, m^2_{11}, ..., m^{r_1}_{11}\}$$

and

$$a^1_{12} = \min\{m^1_{12}, m^2_{12}, ..., m^{r_1}_{12}\}$$

and so on.

Thus $A_1 = (a^1_{ij})$ is such that $a^1_{ij} \leq m^t_{ij}$, $1 \leq t \leq r$, $1 \leq i \leq m_1$ and $1 \leq j \leq n_1$.

Now we build $B_1 = (b^1_{ij})$ as follows:

Define each

$$b^1_{ij} = \max\{m^1_{ij}, m^2_{ij}, ..., m^{r_1}_{ij}\}$$

for $1 \leq i \leq m_1$ and $1 \leq j \leq n_1$;

$$b^1_{11} = \max\{m^1_{11}, m^2_{11}, ..., m^{r_1}_{11}\},$$

$$b^1_{12} = \max\{m^1_{12}, m^2_{12}, ..., m^{r_1}_{12}\}$$

and so on.

Clearly $b^1_{ij} \geq m^t_{ij}$ for $1 \leq t \leq r_1$ and $1 \leq i \leq m_1$ and $1 \leq j \leq n_1$; Thus $a^1_{ij} \leq m^t_{ij} \leq b^1_{ij}$; $t = 1, 2, ..., r_1$.

We have $[A_1, B_1]$ to be the interval of synaptic $m_1 \times n_1$ matrices containing the $r_1$ synaptic connection matrices given by the $r_1$ experts.

Now we define for this interval of synaptic $m_1 \times n_1$ matrices the related optimal and average synaptic connection matrices as follows:

$$O^1 = \frac{A_1 + B_1}{2}$$

$$= \frac{(a^1_{ij}) + (b^1_{ij})}{2}$$

$$= (o^1_{ij})\ ;$$

$1 \leq i \leq m_1$ and $1 \leq j \leq n_1$ .

$$\overline{M}^1 = \frac{m^1_1 + m^1_2 + \ldots + m^1_{r_1}}{r_1}$$

$$= \frac{(m^1_{ij}) + (m^2_{ij}) + \ldots + (m^{r_1}_{ij})}{r_1}$$

$$= (\overline{m}_{ij})\ .$$

Clearly $O^1$ and $\overline{M}^1 \in [A_1, B_1]$ by the very construction.

It is to be noted that the interval can contain a maximum number of $r_1 + 4$ number of $m_1 \times n_1$ matrices. It may so happen that these matrices $A_1$, $B_1$, $\overline{M}^1$ or $O^1$ may concide with any of the $M^1_1, M^1_2, \ldots, M^1_{r_1}$ matrices also. That is why, we say it can contain a maximum number of $r_1 + 4$ matrices.

Now the same procedure is adopted for the $r_2$ synaptic connection $m_2 \times n_2$ matrices given by the second set of $r_2$ experts on the same problem on the scale $[–a_2, a_2]$. We denote this interval of synaptic connection matrices by $[A_2, B_2]$ where $A_2$ is the minimal synaptic connection matrix of the interval of matrices $[A_2, B_2]$ related with the synaptic connection matrices $M^2_1, M^2_2, \ldots, M^2_{r_2}$ .

Similarly $B_2$ is the maximal synaptic connection matrix and

$$O^2 = \frac{A_2 + B_2}{2}$$

$$= \frac{(a^2_{ij}) + (b^2_{ij})}{2}$$

will be the associated optimal synaptic connection matrix. The average / mean synaptic connection matrix

$$\overline{M}^2 = \frac{m_1^2 + m_2^2 + ... + m_{r_2}^2}{r_2}$$

$$= \frac{(m_{ij}^1) + (m_{ij}^2) + ... + (m_{ij}^{r_2})}{r_2}.$$

The same procedure is repeated for all t-sets of experts $r_1$, $r_2$, … $r_t$.

So we have t-intervals of $m_i \times n_i$ ($1 \le i \le t$) synaptic connection matrices associated with the problem P, given by $[A_1, B_1]$, $[A_2, B_2]$, …, $[A_t, B_t]$.

Now set

$$[A, B] = [A_1, B_1] \cup [A_2, B_2] \cup \ldots \cup [A_t, B_t]$$

where [A, B] is the collection of mixed rectangular t-matrices such that if P is any typical element of the set [A, B] then $P = P_1 \cup P_2 \cup \ldots \cup P_t$ where $P_i \in [A_i, B_i]$, i = 1, 2, …, t. Further $A = A_1 \cup A_2 \cup \ldots \cup A_t$; $B = B_1 \cup B_2 \cup \ldots \cup B_t$, $O = O_1 \cup O_2 \cup \ldots \cup O_t$ and $\overline{M} = \overline{M}^1 \cup \overline{M}^2 \cup ... \cup \overline{M}^t$ where A is the minimal synaptic connection t-matrix of the interval of synaptic connection of t-matrices [A, B]. Like wise B is the maximal synaptic connection t-matrix. O the optimal synaptic connection t-matrix and $\overline{M}$ is the average/ mean of the synaptic connection t-matrix of the t sets of experts associated with the problem P.

We call [A, B] to be the interval of synaptic connection of t-matrices of the bidirectional t-associative memories model (t-IBAM model).

This model also works as the Bidirectional 2-associative memories models (BIBAM model) when t =2 which has been very clearly explained.

The Bidirectional t-associative memories model is a Bi directional bi associative model when t = 2. Thus we have seen the notion of interval of n-matrices finds its applications in the BAM models and their generalizations when the experts choose to work with any random number of variables i.e., attributes or concepts and also they choose different scales of intervals. This new model will be shortly represented by t-IBAM model.

## 3.8 Introduction to Interval Neutrosophic Matrices and use of these Matrices in Neutrosophic Models

Next we proceed on to introduce the notion of interval of neutrosophic matrices for the first time and proceed onto give their applications to the real world problems.

For the brief notion about neutrosophy and neutrosophic matrices please refer [220], also a brief introduction is given in chapter I. Throughout the book 'I' will denote the indeterminate and $I^2 = I$. Further I + I = 2I, I – I = 0, I + I + … + I, n times is n I and so on.$\langle Q \cup I\rangle$ is the neutrosophic field of rationals, $\langle Z \cup I\rangle$ the neutrosophic integral domain or the neutrosophic ring of integers, $\langle R \cup I\rangle$ and $\langle C \cup I\rangle$ are the neutrosophic field of reals and complex numbers respectively.

All these neutrosophic fields are of characteristic 0. $\langle Z_n \cup I\rangle$ is the neutrosophic ring of integers modulo n and of characteristic n.

The set of matrices $M_{m\times n} = \{(a_{ij}) / a_{ij} \in \langle Z \cup I\rangle\}$ is the collection of neutrosophic matrices with entries from the neutrosophic ring of integers. $T_{q\times q} = \{(t_{ij}) / t_{ij} \in \langle Q \cup I\rangle\}$ is neutrosophic ring of matrices under usual matrix addition and matrix multiplication.
Thus for more about neutrosophic matrices please refer [220].

***Example 3.8.1:*** Let

$$M = \begin{bmatrix} 3+I & 2I & -8I+1 \\ 7 & 8+5I & 3I-1 \\ 2 & -30 & 11 \end{bmatrix}$$

M is a 3 × 3 neutrosophic square matrix with entries from the ring of neutrosophic integers $\langle Z \cup I\rangle$.

***Example 3.8.2:*** Let $A = A_1 \cup A_2$ where

$$A_1 = \begin{bmatrix} 3I-2 & 7I & 5I+7 \\ -3 & 8-2I & 11I \\ 0 & 3+17I & 5 \end{bmatrix}$$

and

$$A_2 = \begin{bmatrix} 3I+7 & 2 \\ 13 & 10I \end{bmatrix}$$

is a mixed square neutrosophic bimatrix.

***Example 3.8.3:*** Let B = $B_1 \cup B_2$ where

$$B_1 = \begin{bmatrix} 2I-7 & 0 & 9I \\ 4 & 14I-3 & 1 \end{bmatrix}$$

and

$$B_2 = \begin{bmatrix} 12 & 10I-1 & 17 & 14I & 16 & -10I \end{bmatrix}$$

is a mixed rectangular neutrosophic bimatrix.

***Example 3.8.4:*** Let M = $M_1 \cup M_2 \cup \ldots \cup M_n$ where each $M_i$ is either a square neutrosophic matrix or a rectangular neutrosophic matrix for i = 1, 2, …, n, then we call M the mixed neutrosophic n-matrix.

For more about these structures please refer [220-222].

Now we proceed on to define interval of neutrosophic matrix for this we need to define an order on the neutrosophic ring or field from which entries will be taken. Further throughout this book we will be considering only real neutrosophic fields or rings. So if $\langle R \cup I \rangle$ is the real neutrosophic field and if x and y are in $\langle R \cup I \rangle$ and x and y are real i.e., x, y $\in$ R then we can say either $x \leq y$ or $x \geq y$. i.e., we can compare them or to be more mathematically 'order' them under a linear order or still say min {x, y] or max {x, y}, can be easily obtained. Note if m, n $\in \langle R \cup I \rangle$ and m = tI and n = sI where s, t $\in$ R then also, we can say min {m, n} or max of {m, n}.

***Note:*** Min of {m, n} is the min of [s, t} × I (max of {m, n} is the max of {s, t} × I} i.e. If m = 5I and n = 0.9I clearly m > n, and min {5I, 0.9I} in 0.9I and max {5I, 0.9I} is 5I i.e., min {5I, 0.9I} = min {5, 9} × I = 0.9 × I. Like wise max {5I, 0.9I} = max {5, 0.9} × I = 5 × I.

Now the problem arises when we have a pair (x, y) ∈ ⟨R ∪ I⟩ and x = 7I and y = 3 or when x = 8I – 2 and y = 5I or when x = 5 + 7I and y = 2 + 11I but when we have x = 8I + 11 and y = 12I + 17 then we say max {x, y} = 12I + 17 and min {x, y} = 8I + 11, so we define in case of neutrosophic matrices two concepts called interval of neutrosophic matrices and pseudo interval of neutrosophic matrices.

We shall see how they differ and how they would be used in the study of the problem. We will also define notions like pseudo neutrosophic maximum and pseudo neutrosophic minimum, pseudo real maximum and pseudo real minimum to cater to the needs of the problem.

Suppose we have {7I, 25} are elements from ⟨Z ∪ I⟩ the pseudo neutrosophic minimum is 7I. The pseudo neutrosophic maximum is 25. Suppose we have {20I + 17, 15I + 42}.

The pseudo neutrosophic minimum is {15I + 42}. The pseudo real minimum is 20I +17. The pseudo real maximum 15I + 42. The pseudo neutrosophic maximum is 20I +17.

Thus according to the need of the problem we will formulate either a pseudo neutrosophic order or the pseudo real order on the interval of neutrosophic matrices depending on the situation of the problem.

So we proceed to define these types of interval of neutrosophic matrices.

**DEFINITION 3.8.1:** *Let [A, B] be a collection of all neutrosophic m × n matrices, we say C a neutrosophic matrix in this collection if and only if C = ($c_{ij}$), (A = ($a_{ij}$) and B = ($b_{ij}$)) then $a_{ij} \leq c_{ij} \leq b_{ij}$. We call [A, B] the interval of neutrosophic matrices. With A = ($a_{ij}$) the minimal matrix of the interval and B = ($b_{ij}$) to be the maximum matrix of the interval.*

Thus '≤' is the usual order and any two elements in their $ij^{th}$ places are comparable. First we illustrate this by the following example.

***Example 3.8.5:*** Let us consider the neutrosophic interval of matrices [A, B] where

$$A = \begin{bmatrix} 5I & 3+2I \\ 17 & 10-3I \end{bmatrix}$$

and

$$B = \begin{bmatrix} 22+15I & 10+12I \\ 40+5I & 50+17I \end{bmatrix}.$$

Now take

$$C = \begin{bmatrix} 6I+1 & 4+3I \\ 18 & 20+I \end{bmatrix}.$$

Clearly C is in the interval of neutrosophic matrices [A, B]. All matrices $m = (m_{ij})$

$$m = \begin{bmatrix} m_{11} & m_{12} \\ m_{21} & m_{22} \end{bmatrix}$$

with $5I \leq m_{11} \leq 22 + 15I$, $3 + 2I \leq m_{12} \leq 10 + 12I$, $17 \leq m_{21} \leq 40 + 5I$ and $10 - 3I \leq m_{22} \leq 50 + 17I$, are in the interval of neutrosophic matrices.

Take

$$N = \begin{bmatrix} 2I & 16+18I \\ 16 & 10+2I \end{bmatrix}$$

a neutrosophic matrix. Clearly N is not an element of the interval [A, B] of neutrosophic matrices.

Now we can have the interval of neutrosophic matrices to be an interval of square $n \times n$ neutrosophic matrix or an interval of rectangular $m \times n$ ($m \neq n$) neutrosophic matrix or an interval of column neutrosophic vector or an interval of row neutrosophic vector.

Now having defined the notion of interval of neutrosophic matrices we proceed on to define interval of pseudo neutrosophic matrices.

**DEFINITION 3.8.2:** *Let [A, B] be a collection of neutrosophic m × n matrices, where A = ($a_{ij}$) and B = ($b_{ij}$). We call [A, B] the pseudo neutrosophic interval of neutrosophic m × n matrices if and only if C = ($c_{ij}$) is in the interval [A, B] then pseudo neutrosophic min of {$c_{ij}$, $a_{ij}$} = $a_{ij}$ and pseudo neutrosophic max of {$c_{ij}$, $b_{ij}$} = $b_{ij}$.*

We just illustrate this by a very simple example.

***Example 3.8.6:*** Let [A, B] be the pseudo neutrosophic interval of neutrosophic 3 × 2 matrices where

$$A = \begin{bmatrix} 2I-4 & 7 \\ 8I & 4+I \\ I+7 & 20 \end{bmatrix}$$

and

$$B = \begin{bmatrix} 20I-2 & 27 \\ 12I+4 & 2+8I \\ 6+3I & 41 \end{bmatrix};$$

consider

$$C = \begin{bmatrix} 8I-4 & 20 \\ 8I-7 & 2+2I \\ I+4 & 22 \end{bmatrix},$$

C is in the pseudo neutrosophic interval of neutrosophic 3 × 2 matrices [A, B].

Now we proceed on to define the notion of pseudo real interval of neutrosophic matrices.

**DEFINITION 3.8.3:** *Let [A, B] be the collection of all m × n neutrosophic matrices where A = ($a_{ij}$) and B = ($b_{ij}$). We all [A,*

*B] the pseudo real interval of neutrosophic matrices if $C = (c_{ij})$ is any matrix of the interval [A, B], then;*

$$\text{pseudo real min of } \{c_{ij}, a_{ij}\} = a_{ij}$$

*and*

$$\text{pseudo real max of } \{c_{ij}, b_{ij}\} = b_{ij}.$$

The reader is given the task of finding an example of pseudo real interval of neutrosophic matrices. Having defined 3 types of interval of neutrosophic matrices we now proceed on to define the notion of interval of neutrosophic bimatrices.

**DEFINITION 3.8.4:** *Let $[A_1, B_1]$ and $[A_2, B_2]$ be two intervals of neutrosophic matrices. Define $[A, B] = \{M_1 \cup M_2 / M_1 \in [A_1, B_1]$ and $M_2 \in [A_2, B_2]\}$ with $A = A_1 \cup A_2$, the minimal neutrosophic bimatrix and $B = B_1 \cup B_2$ the maximal neutrosophic bimatrix to be the interval of neutrosophic bimatrix.*

We call [A, B] the interval of neutrosophic bimatrices. If $[A_i, B_i]$ is a interval of neutrosophic $m_i \times m_i$ square matrix i = 1, 2 then we call [A, B] the mixed square bimatrix. If $[A_i, B_i]$ is an interval of neutrosophic $m_i \times n_i$ rectangular matrix i = 1, 2; then we call [A, B] to be the interval of neutrosophic mixed rectangular bimatrix. If $[A_1, B_1]$ is an interval of neutrosophic square bimatrix and $[A_2, B_2]$ is an interval of neutrosophic rectangular matrix then we define [A, B] the interval of mixed neutrosophic bimatrices.

We illustrate them by the following examples:

***Example 3.8.7:*** Let $[A, B] = [A_1, B_1] \cup [A_2, B_2]$ be an interval of neutrosophic bimatrices. In the interval of neutrosophic matrices $[A_1, B_1]$ where

$$A_1 = \begin{bmatrix} a_{11} & a_{12} \\ a_{21} & a_{22} \end{bmatrix}$$

with $a_{11} = 3I + 2$, $a_{12} = 12I$, $a_{21} = 3$, $a_{22} = 7I – 1$ and

$$B_1 = \begin{bmatrix} b_{11} & b_{12} \\ b_{21} & b_{22} \end{bmatrix}$$

where $b_{11} = 5I + 3$, $b_{12} = 18I$, $b_{21} = 18$ and $b_{22} = 8I + 11$.

For the interval of neutrosophic matrices $[A_2, B_2]$ where

$$A_2 = \begin{bmatrix} a_{11} & a_{12} \\ a_{21} & a_{22} \\ a_{31} & a_{32} \end{bmatrix}$$

and

$$B_2 = \begin{bmatrix} b_{11} & b_{12} \\ b_{21} & b_{22} \\ b_{31} & b_{32} \end{bmatrix};$$

with
$a_{11} = 3I – 2$, $a_{12} = 12$, $a_{21} = 4I + 6$, $a_{22} = 14I$, $a_{31} = –11 + I$ and $a_{32} = 6 – 5I$.
$b_{11} = 2 – 4I$, $b_{12} = 12I +13$, $b_{21} = 7I + 8$, $b_{22} = 40I$, $b_{31} = 5 + 12I$ and $b_{32} = 17 + 3I$.

Thus $[A, B] = [A_1, B_1] \cup [A_2, B_2]$ is an interval of neutrosophic bimatrices.

We give the definition of mixed interval of neutrosophic bimatrices.

**DEFINITION 3.8.5:** *Let $[A, B] = [A_1, B_1] \cup [A_2, B_2]$ where $[A_1, B_1]$ is the interval of neutrosophic matrices and $[A_2, B_2]$ is the pseudo neutrosophic interval of neutrosophic matrices or if $[A_2, B_2]$ is the pseudo real interval of neutrosophic matrices or $[A_1, B_1]$ is itself a pseudo real interval of neutrosophic matrices and so on, then we call $[A, B]$ the mixed interval of neutrosophic bimatrices.*

Now we proceed onto define the notion of interval of neutrosophic n-matrices and mixed interval of neutrosophic n-matrices.

**DEFINITION 3.8.6:** *Let $[A_1, B_1]$, $[A_2, B_2]$, …, $[A_n, B_n]$ be the n intervals of neutrosophic matrices set, $[A, B] = [A_1, B_1] \cup [A_2, B_2] \cup \ldots \cup [A_n, B_n] = \{M = M_1 \cup M_2 \cup \ldots \cup M_n / M_i \in [A_i, B_i]$ with $A = A_1 \cup A_2 \cup \ldots \cup A_n$, $B = B_1 \cup B_2 \cup \ldots \cup B_n$ the set of all neutrosophic n matrices\}.*

*We call [A, B] to be the interval of neutrosophic n-matrices. If all the interval of neutrosophic matrices $[A_i, B_i]$ are square matrices we call [A, B] to be the interval of neutrosophic mixed square n-matrices.*

*If the intervals of neutrosophic matrices $[A_i, B_i]$ are rectangular matrices then we call [A, B] to be the interval of neutrosophic mixed rectangular matrices. Like wise when the intervals of matrices $[A_i, B_i]$ are either neutrosophic square matrices or neutrosophic rectangular matrices then we call [A, B] an interval of neutrosophic mixed matrices.*

We can imagine an interval of neutrosophic n-matrices similar to the notion of the interval of n-matrices where every general matrix in the n-matrices is replaced by the neutrosophic matrix.
Now we proceed on to define the notion of interval of fuzzy neutrosophic matrices. We know a fuzzy neutrosophic matrix is a matrix in which every entry is from the set $\langle [0, 1] \cup [0, I] \rangle = \{x + yI / y, x \in [0, 1]\}$.

Here the three types of order can be defined on the set $\langle [0, 1] \cup [0, I] \rangle$. If $x, y \in \langle [0, 1] \cup [0, I] \rangle$, min\{x, y\} and max \{x, y\} provided $x \leq y$ or $x \geq y$. If $x \geq y$ or $x \leq y$ cannot be said then we define pseudo real maximum or real minimum as in case of real neutrosophic numbers, i.e., if x = 0.7 + 0.8I and y = 0.9 + 0.6 I then pseudo real of maximum \{x, y\} is 0.9 + 6I, pseudo real minimum of \{x, y\} is 0.7 + 8I, pseudo neutrosophic maximum of \{x, y\} is 0.7 + 8I and pseudo neutrosophic minimum of \{x, y\} is 0.9 + 0.6I. Clearly if x = 0.8I and y = 0.8I + 0.2 or $y_1$ = 0.8I – 0.2 then $x \leq y$ and $y_1 \leq x$.

Keeping these three types of orders in mind we proceed on to define interval of fuzzy neutrosophic matrices.

**DEFINITION 3.8.7:** *Let [A, B] = {set of all n × n fuzzy neutrosophic matrices with entries from ⟨[0, I] ∪ [0, 1]⟩ where A = ($a_{ij}$), and B = ($b_{ij}$)} are n × n square matrices. We have for any C = ($c_{ij}$) ∈ [A, B]; if and only if $a_{ij} \leq c_{ij} \leq b_{ij}$. Then we call [A, B] to be the interval of fuzzy neutrosophic square matrices.*

If instead of taking fuzzy neutrosophic square matrices we have [A, B] to be the collection of all m × n (m ≠ n) fuzzy neutrosophic rectangular matrices then we call [A, B] to be the interval of fuzzy neutrosophic rectangular matrices.

***Example 3.8.8:*** Let [A, B] be the interval of fuzzy neutrosophic 2 × 2 matrices where

$$A = \begin{bmatrix} 0.I & 0.2I + 0.3I \\ 0.5 + 0.3I & 0.7 + 0.5I \end{bmatrix}$$

and

$$B = \begin{bmatrix} I & 0.3 + I \\ 1 + I & 0.3I + 1 \end{bmatrix}.$$

Take

$$C = \begin{bmatrix} 0.9I & 0.7 + 0.5I \\ 0.9 + 0.4I & 0.8 + 0.6I \end{bmatrix}$$

clearly C is in [A, B]. If we take

$$D = \begin{bmatrix} 0.3 + I & 0.3I + I \\ 1 + I & 0.9I \end{bmatrix}$$

D is a 2 × 2 fuzzy neutrosophic matrix but D ∉ [A, B].

Thus [A, B] is an interval of fuzzy neutrosophic 2 × 2 matrices.

Now having defined the interval of fuzzy neutrosophic matrices and give an example, the reader is expected to give more examples of interval of fuzzy neutrosophic matrices.

**DEFINITION 3.8.8:** *Let [A, B] be the collection of all n × n fuzzy neutrosophic matrices where A = ($a_{ij}$) and B = ($b_{ij}$). For any C = ($c_{ij}$), a n × n fuzzy neutrosophic matrix is said to belong to the collection [A, B] if and only if pseudo neutrosophic max {$a_{ij}$, $c_{ij}$} = $c_{ij}$, pseudo neutrosophic min {$a_{ij}$, $c_{ij}$} = $a_{ij}$ pseudo neutrosophic min {$c_{ij}$, $b_{ij}$} = $c_{ij}$ and pseudo neutrosophic max {$c_{ij}$, $b_{ij}$} = $b_{ij}$.*

We call [A, B] the pseudo neutrosophic interval of n × n fuzzy neutrosophic matrices. If [A, B] is the collection of m × n rectangular fuzzy neutrosophic matrices then we call [A, B] to be the pseudo neutrosphic interval of rectangular m × n fuzzy neutrosophic matrices.

**DEFINITION 3.8.9:** *Let [A, B] be the collection of all n × n fuzzy neutrosophic matrices, where A = ($a_{ij}$) and B = ($b_{ij}$). We say a n × n fuzzy neutrosophic matrix C = ($c_{ij}$) is in [A, B] if and only if pseudo real max {$a_{ij}$, $c_{ij}$} = $c_{ij}$ and pseudo real min {$a_{ij}$, $c_{ij}$} = $a_{ij}$; pseudo real min {$c_{ij}$, $b_{ij}$} = $c_{ij}$ and pseudo real max {$c_{ij}$, $b_{ij}$} = $b_{ij}$.*

*Let [A, B] denote such collection of matrices C. We call [A, B] to be the pseudo real interval of the fuzzy neutrosophic n × n square matrices. Now if we consider instead of n × n square matrix a m × n (m ≠ n) rectangular fuzzy neutrosophic matrices then we call the set [A, B] the pseudo real interval of fuzzy neutrosophic m × n rectangular matrices where A = ($a_{ij}$) and B = ($b_{ij}$) are m × n rectangular matrices.*

The reader is expected to construct more and more examples of pseudo neutrosophic interval of fuzzy neutrosophic matrices and pseudo real interval of fuzzy neutrosophic matrices.

Now we proceed on to define the notion of interval of fuzzy neutrosophic bimatrices and interval of fuzzy neutrosophic n-matrices (n > 2).

**DEFINITION 3.8.10:** *Let $[A_1, B_1]$ and $[A_2, B_2]$ be two intervals of fuzzy neutrosophic matrices. We set $[A, B] = [A_1, B_1] \cup [A_2, B_2]$ = {Collection of all fuzzy neutrosophic bimatrices $M = M_1 \cup M_2$ where $M_1 \in [A_1, B_1]$ and $M_2 \in [A_2, B_2]$ with $A = A_1 \cup A_2$ and $B = B_1 \cup B_2$} and define [A, B] to be the interval of fuzzy neutrosophic bimatrices.*

*If the matrices in both the interval of fuzzy neutrosophic matrices are square then, we call $[A, B] = [A_1, B_1] \cup [A_2, B_2]$ to be the interval of fuzzy neutrosophic mixed square bimatrices; if on the other hand both the interval of fuzzy neutrosophic matrices are rectangular then we call [A, B] to be an interval of fuzzy neutrosophic mixed rectangular bimatrices. If one of the interval of fuzzy neutrosophic matrices is square and the other interval of fuzzy neutrosophic matrices is rectangular then we call $[A, B] = [A_1, B_1] \cup [A_2, B_2]$ to be the interval of mixed fuzzy neutrosophic bimatrices.*

The reader is expected to construct examples of each type, for it is a matter of routine and can be constructed as in case of interval of fuzzy bimatrices.

**DEFINITION 3.8.11:** *Let $[A_1, B_1], [A_2, B_2], \ldots, [A_n, B_n]$ be a collection of n-intervals of fuzzy neutrosophic matrices . Set $[A, B] = [A_1, B_1] \cup [A_2, B_2] \cup \ldots \cup [A_n B_n] = \{M_1 \cup M_2 \cup \ldots \cup M_n / M_i \in [A_i, B_i], i = 1, 2, \ldots, n$ and $A = A_1 \cup A_2 \cup \ldots \cup A_n$ and $B = B_1 \cup B_2 \cup \ldots \cup B_n\}$.*

*We define [A, B], the collection of all n-matrices from the n intervals of fuzzy neutrosophic matrices as the interval of fuzzy neutrosophic n-matrices ($n \geq 2$). When each of the n interval of fuzzy neutrosophic matrices $[A_i, B_i]$ are only square matrices for $i = 1, 2, \ldots, n$; then we call [A, B] the interval of fuzzy neutrosophic mixed square n-matrices.*

*If each of the n-interval of fuzzy neutrosophic matrices $[A_i, B_i]$ happen to be rectangular matrices for $i = 1, 2, \ldots, n$ then we call [A, B] to be the interval of fuzzy neutrosophic mixed rectangular n-matrices.*

*If in the collection of interval of fuzzy neutrosophic n-matrices $[A, B] = [A_1, B_1] \cup [A_2, B_2] \cup \ldots \cup [A_n, B_n]$ some of the intervals $[A_i, B_i]$ happen to be rectangular fuzzy* *neutrosophic matrices and some of the intervals $[A_j, B_j]$ happen to be only square fuzzy neutrosophic matrices then we call [A, B] the interval of fuzzy neutrosophic mixed n-matrices.*

*Now if in the interval of fuzzy neutrosophic n-matrices [A, B] = $[A_1, B_1] \cup [A_2, B_2] \cup \ldots \cup [A_n, B_n]$ we have some interval $[A_i, B_i]$ to be a pseudo real interval of fuzzy neutrosophic matrices and some intervals $[A_j, B_j]$ to be a pseudo neutrosophic matrices and some intervals $[A_k, B_k]$ to be just interval of fuzzy neutrosophic matrices then we call [A, B] to be the mixed interval of mixed fuzzy neutrosophic n-matrices.*

As in case of interval of fuzzy mixed n-matrices the reader can construct examples by replacing a fuzzy matrix by a fuzzy neutrosophic matrix.

Now we proceed on to give their applications.

## 3.9 Application of Interval of Neutrosophic Matrices to Neutrosophic Models

Bidirectional Associative Memories (BAM) model has been just recalled in the chapter one. For more refer [96, 214]. Further several existing results and new results have been introduced in the book to appear.

However for the sake of completeness we just recall the description and the working of the model.

Suppose we have a problem P at hand which comprises only of the unsupervised data i.e., no statistical data is available. Suppose we have $\langle X \cup I \rangle$, n attributes or neurons associated with the problem P in the neuronal neutrosophic field $F_{\langle X \cup I \rangle}$; m other related attributes or concepts which makes the problem to function or highly dependent on the working or in the determination of properties associated with the solution of the problem in the neuronal neutrosophic field $F_{\langle Y \cup I \rangle}$.

Let us consider the neuronal field $F_X$ associated with the n-attributes and the neuronal field $F_Y$ associated with the m-attributes. Now we make the neuronal field to be generated by $\langle F_X \cup I \rangle$ which we denote by $F_{\langle X \cup I \rangle}$ and call $F_{\langle X \cup I \rangle}$ to be the neutrosophic neuronal field associated with the n attributes. On similar lines we define $F_{\langle X \cup I \rangle}$.

Here three things happen, the attributes themselves can be neutrosophic or the concepts may turn to be neutrosophic depending on the external situation or its relation with other concepts may be neutrosophic.

Let M denote the synaptic connection matrix of functions from the neutrosophic neuronal field $F_{\langle X \cup I \rangle}$ to $F_{\langle Y \cup I \rangle}$ taken on the scale say $\langle [-a_1, a_1] \cup [-a_1 I, a_1 I] \rangle$. M = $(a_{ij})$ is clearly a neutrosophic matrix with entries from $\langle [-a_1, a_1] \cup [-a_1 I, a_1 I] \rangle$, $1 \leq i \leq m$ and $1 \leq j \leq n$. Now the entries in the neutrosophic neuronal field $F_{\langle X \cup I \rangle}$ are taken as the input neutrosophic vectors for which the synaptic connection neutrosophic matrix M would be used. For the other set of attributes the input neutrosophic vectors are taken from the neutrosophic neuronal field $F_{\langle Y \cup I \rangle}$. Both these neutrosophic neutronal fields $F_{\langle X \cup I \rangle}$ and $F_{\langle Y \cup I \rangle}$ take their entries from the neutrosophic interval $\langle [-a_1, a_1] \cup [-a_1 I, a_1 I] \rangle$. The input neutrosophic vectors are taken from these neutrosophic fields $F_{\langle X \cup I \rangle}$ and $F_{\langle Y \cup I \rangle}$.

This synaptic connection neutrosophic matrix M forms the dynamical system of the problem. When we input any neutrosophic vector from the neutrosophic field we get the resultant neutrosophic vector which gives the equilibrium of the system.

As in case of BAM even in the case of Neutrosophic Bidirectional Associative memories (NBAM) we see the input neutrosophic vector from the neuronal neutrosophic fields cannot be recognized by the dynamical system NBAM. We define the neutrosophic trivector system, for which the function S is used.

If X = $(x_1, \ldots, x_n)$; $x_i \in \langle [-a_1, a_1] \cup [-a_1 I, a_1 I] \rangle$ then,

$$
\begin{aligned}
S(x_i) &= 0 \text{ if } x_i \leq 0 \\
&= 1 \text{ if } x_i > 0 \\
&= 0 \text{ if } x_i \text{ is negative neutrosophic} \\
&= I \text{ if } x_i \text{ is +ve neutrosophic.}
\end{aligned}
$$

For example if

$$X = (2I, -I, 0, -1, 4, -6, 3I)$$
$$S(X) = (I\ 0\ 0\ 0\ 1\ 0\ I).$$

Now if $X_k$ is any input vector at the $k^{th}$ period of time from $F_{\langle X \cup I\rangle}$ then its effect on the dynamical system M is given by

$$\begin{aligned} S(X_k)M &= (y_1, \ldots, y_n) \\ &= Y_{k+1} \\ S(Y_{k+1})M^T &= \left(x_1^1, ..., x_m^1\right) \\ &= X_{k+2}. \\ S(X_{k+2})M &= \ldots \end{aligned}$$

and so on till the equilibrium of the system is arrived resulting in a fixed trinary pair. For more please refer [96].

As this notion is new we just illustrate it with a simple NBAM model.

The study of the problem of a sample group is considered. This group consists of HIV / AIDS infected migrant labourers in the age group 20 to 58 and they were involved in a variety of deregulated labour such as transport, truck or cab drivers, construction labourers, daily wagers etc.

The two heads under which we wish to study this problem is as follows:

A : Causes for Migrant labourer's vulnerability to HIV/AIDS.

$A_1$ – No awareness / education
$A_2$ – Social status
$A_3$ – No social responsibility enormous freedom
$A_4$ – Bad company and additive habits
$A_5$ – Type of profession
$A_6$ – Cheap availability of CSWs.

B: Factors forcing them for migration.

$B_1$ - Lack of labour opportunities in their hometown
$B_2$ - Poverty / seeking better status in life

$B_3$ - Mobilization of contract labourers by middlemen

$B_4$ - Infertility of land due to wrong means of agricultural research methodologies / failure of monsoon

$B_5$ - Globalization industrialization / modernization.

Taking the neutrosophic neuronal field $F_{\langle X \cup I \rangle}$ as the attributes connected with the causes of vulnerability resulting in HIV/AIDS and the neutrosophic neuronal field $F_{\langle Y \cup I \rangle}$ is taken as the factors forcing people for migration.

The 6 × 5 neutrosophic matrix M represents the forward synaptic projections from the neutrosophic neuronal field $F_{\langle X \cup I \rangle}$ to the neutrosophic neuronal field $F_{\langle Y \cup I \rangle}$. The 5 × 6 neutrosophic matrix $M^T$ represents the backward synaptic projections $F_{\langle X \cup I \rangle}$ to $F_{\langle Y \cup I \rangle}$. Now taking $A_1, A_2, \ldots, A_6$ along the rows and $B_1, B_2 \ldots, B_5$ along the columns we get the synaptic connection neutrosophic matrix M which is modeled on the neutrosophic scale $\langle[-5, 5] \cup [-5I, 5I]\rangle$.

$$M = \begin{array}{c} \\ A_1 \\ A_2 \\ A_3 \\ A_4 \\ A_5 \\ A_6 \end{array} \begin{array}{c} \begin{array}{ccccc} B_1 & B_2 & B_3 & B_4 & B_5 \end{array} \\ \begin{bmatrix} 5 & 2 & 4 & 4 & I \\ 4 & 3 & 5 & 3 & 0 \\ -1 & -2 & 4 & I & 0 \\ 0 & 4 & 2 & 0 & 0 \\ 2 & 4 & 3 & 3 & 4 \\ 0 & 2 & 0 & 0 & 0 \end{bmatrix} \end{array}.$$

Let $X_k$ be the input vector given by (3, 4, –1, –3, –2, 1) at the $k^{th}$ time period. The initial vector is given such that literacy, lack of awareness, social status and cheap availability of CSWs have stronger impact over migration. We suppose that all neuronal state change decisions are synchronous.

The trinary vector

$$S(X_k) = (1\ 1\ 0\ 0\ 0\ 1).$$

From the activation equation

| | | |
|---|---|---|
| $S(X_k)M$ | = | (9, 7, 9, 7, I) |
| | = | $Y_{k+1}$. |
| $S(Y_{k+1})$ | = | (1 1 1 1 I) |
| $S(Y_{k+1})M^T$ | = | (15 + I, 15, 4 + I, 6, 10 + I, 2) |
| | = | $X_{k+2}$. |
| $S(X_{k+2})$ | = | (1 1 1 1 1 1) |
| $S(X_{k+2})M$ | = | (6, 13, 18, 10 + I, 4 + I) |
| | = | $Y_{k+3}$. |
| $S(Y_{k+3})$ | = | (1 1 1 1 1) |
| $S(Y_{k+3})M^T$ | = | (15 + I, 15 + I, 6, 12, 2) |
| | = | $X_{k+4}$. |
| $S(X_{k+4})$ | = | (1 1 I 1 1 1) |
| $S(X_{k+4})M$ | = | (7 – I, 15 – 2I, 14 + I, 10 + I, 4 + I) |
| | = | $Y_{k+5}$. |
| $S(Y_{k+5})$ | = | (1 1 1 1 1). |

But it is pertinent to mention here that as they proceed on with varying time the system predicts the increase in the indeterminacy factor. Thus one can analyze why the indeterminacy increases with time for certain attributes.

Next we proceed onto build the notion of interval of neutrosophic synaptic connection matrices associated with the NBAM model. We have just introduced and analyzed how the NBAM model functions. Also in this chapter we have just introduced the new concept of interval of neutrosophic square (rectangular and mixed) matrices and also the notion of mixed interval of neutrosophic matrices.

Now we will just introduce the notion of neutrosophic bidirectional bi-associative memories (NBBAM) model and use the notion of neutrosophic bimatrices as the dynamical bisystem of this bimodel.

Let P be a problem at hand and let the related data associated with the problem be only an unsupervised data. Suppose two experts wish to work in the data and choose to work with two different sets of attributes and also on two different scales.

We have to formulate a method so that both of them are given equal importance and also their views are comparable at each stage or time interval. To this end let us say the first expert wishes to work with ($m_1$, $n_1$) sets of attributes on the scale $\langle[-a_1, a_1] \cup [-a_1I, a_1I]\rangle$. Let the neutrosophic neuronal field associated in the study be given by $F_{\langle X_1 \cup I\rangle}$ and $F_{\langle Y_1 \cup I\rangle}$; on the same interval $\langle[-a_1, a_1] \cup [-a_1I, a_1I]\rangle$.

Let $M_1$ be the neutrosophic synaptic connection matrix given by the first expert. $M_1 = \left(m_{ij}^1\right)$ and $m_{ij}^1 \in \langle[-a_1, a_1] \cup [-a_1I, a_1I]\rangle$ where $1 \le i \le m_1$ and $1 \le j \le n_1$.

On similar lines let $M_2$ be the neutrosophic synaptic connection matrix given by the second expert on the scale $\langle[-a_2, a_2] \cup [-a_2I, a_2I]\rangle$ i.e., $M_2 = \left(m_{ij}^2\right)$ with $m_{ij}^2 \in \langle[-a_2, a_2] \cup [-a_2I, a_2I]\rangle$, $1 \le i \le m_2$ and $1 \le j \le n_2$. $F_{\langle X_2 \cup I\rangle}$ and $F_{\langle Y_2 \cup I\rangle}$ be the neutrosophic neuronal fields associated with the second experts input vectors and are defined over the scale $\langle[-a_2, a_2] \cup [-a_2I, a_2I]\rangle$.

Now set $M = M_1 \cup M_2$, this represents a neutrosophic synaptic connection bimatrix of both the experts. M will be called as the dynamical bisystem of the Bidirectional biassociative memories (NBBAM) model.

If $X_k = X_1 \cup X_2 \in F_{\langle X_1 \cup I\rangle} \cup F_{\langle X_2 \cup I\rangle}$ is any input bivector at the $k^{th}$ time period from the bi-interval $\langle[-a_1, a_1] \cup [-a_1I, a_1I]\rangle \cup \langle[-a_2, a_2] \cup [-a_2I, a_2I]\rangle$; then $S(X_k)$ is the trivector with elements from the set {+1, I, 0}.

Now

$$\begin{aligned} S(X_k)M &= [S(x_1) \cup S(x_2)]\,[M_1 \cup M_2] \\ &= S(X_1)\,M_1 \cup S(X_2)\,M_2 \\ &= Y_1 \cup Y_2 \\ &= Y_{k+1}. \end{aligned}$$

We work with $S(Y_{k+1})M^T$ and so on. Thus we can arrive at the equilibrium of the dynamical bisystem.

Thus when we have only two experts we can work with the neutrosophic synaptic connection bimatrix.

When we have more than 2 experts say n experts and each of them choose to work with varying number of attributes and with varying scales how to construct a new system which can work with all the n-experts simultaneously.

Let P be the problem at hand and suppose n experts wish to work with the problem with varying number of attributes and also with varying scales. Suppose we have the first expert $E_1$ working with $(m_1, n_1)$ attributes on the neutrosophic neuronal scale $\langle[-a_1, a_1] \cup [-a_1I, a_1I]\rangle$ using a $m_1 \times n_1$ neutrosophic matrix. Let $M_1 = (m_{ij}^1)$ denote the $m_1 \times n_1$ neutrosophic synaptic connection matrix given by the first expert on the scale $I_1 = \langle[-a_1, a_1] \cup [-a_1I, a_1I]\rangle$; $\left(m_{ij}^1\right) \in I_1$.

On similar lines let the second expert $E_2$ work with $(m_2, n_2)$ attributes on the neutrosophic neuronal scale $\langle[-a_2, a_2] \cup [-a_2I, a_2I]\rangle = I_2$, and let $M_2 = \left(m_{ij}^2\right)$ denote the $m_2 \times n_2$ neutrosophic synaptic connection matrix with $m_{ij}^2 \in I_2$; $1 \leq i \leq m_2$ and $1 \leq j \leq n_2$; and so on. Thus the n experts $E_1, E_2, \ldots, E_n$ give their views in the n- synaptic connection neutrosophic matrices $M_1, M_2, \ldots, M_n$ on the neutrosophic scales $I_1, I_2, I_3, \ldots, I_n$ respectively.

Let $M = M_1 \cup M_2 \cup \ldots \cup M_n$ be the neutrosophic n-matrix, we define M to be the neutrosophic synaptic connection n-matrix of the dynamical system. Let $F_{\langle X_i \cup I\rangle}$, $F_{\langle Y_i \cup I\rangle}$ be the neutrosophic neuronal fields given by the expert $E_i$ with the neutrosophic scale $I_i$; $1 \leq i \leq n$.

Thus any input neutrosophic n-vector in this dynamical system will be $X = X_1 \cup X_2 \cup \ldots \cup X_n$ where $X_i \in F_{\langle X_i \cup I\rangle}$ and i = 1, 2, …, n.

As in case of other NBAMs the dynamical system cannot in general recognize the neutrosophic n-vector. So we use the S-function,

$$S(X) = S(X_1) \cup S(X_2) \cup \ldots \cup S(X_n)$$

where each $S(X_i)$ will be a trinary vector i.e., its entries will be from the set {0, 1, I}.

Now we will just denote the functions of the dynamical system. Let $X_k = X_1 \cup X_2 \cup \ldots \cup X_n$ be the input neutrosophic n-vector from $I_1 \cup I_2 \cup \ldots \cup I_n$ at any $k^{th}$ interval of time. $S(X_k)$ be the changed trinary vector using the S function, $S(X_k) = S(X_1) \cup S(X_2) \cup \ldots \cup S(X_n)$.

To study the effect of $X_k$ on the dynamical n-system M, we find

$$
\begin{aligned}
S(X_k)M &= [S(X_1) \cup S(X_2) \cup \ldots S(X_n)]\ [M_1 \cup M_2 \cup \ldots \cup M_n] \\
&= S(X_1)M_1 \cup S(X_2)M_2 \cup \ldots \cup S(X_n)M_n \\
&= Y_1 \cup Y_2 \cup \ldots \cup Y_n \\
&= Y_{k+1}.
\end{aligned}
$$

$$
\begin{aligned}
S(Y_{k+1}) &= S(Y_1) \cup S(Y_2) \cup \ldots \cup S(Y_n). \\
S(Y_{k+1})M^T &= [S(Y_1) \cup S(Y_2) \cup \ldots \cup S(Y_n)]\ [M_1^T \cup M_2^T \cup \ldots \cup M_n^T] \\
&= S(Y_1)\ M_1^T \cup S(Y_2)\ M_2^T \cup \ldots \cup S(Y_n)\ M_n^T \\
&= X_1^1 \cup X_2^1 \cup \ldots \cup X_n^1 \\
&= X_{k+2};
\end{aligned}
$$

and so on; until we arrive at the equilibrium of the dynamical n-system. Now this model will help us to study all the effects of the input neutrosophic n-vectors and compare them at every interval of time. Only this system can show the effect at each moment of time so that each and every experts opinion cannot only be compared but relative opinion can also be formed without any effort. Thus this system has all advantages over the existing systems.

Now we proceed on the define the nINBAM model when more than one expert gives his opinion having the same number of attributes and also the same set of attributes and further agree upon to work with the same scale of interval. In this case we make use of the notion of interval of neutrosophic matrices. Let n number of experts work with the problem P and let all of them choose to work with (m, n) attributes on the same neutrosophic interval $\langle[-a, a] \cup [-aI, aI]\rangle$. Let the n number of $m \times n$ matrices given by the n experts be denoted by $E_1, E_2, \ldots, E_n$ we have to find the interval having a mixed collection of $m \times n$ neutrosophic matrices, for we may not be able to find min or max. It can be also be a pseudo neutrosophic or pseudo real interval. We first decide whether to order using pseudo neutrosophic interval of neutrosophic matrices or pseudo real interval of neutrosophic matrices or is it just an interval of neutrosophic matrices.

Let $E_t = \left(e_{ij}^t\right)$; t = 1, 2, …, n, 1 ≤ I ≤ m and 1 ≤ j ≤ n then if [A, B] is the interval of neutrosophic matrices. Let A = $(a_{ij})$ and B = $(b_{ij})$ set

$$a_{ij} = \min\left\{e_{ij}^1, e_{ij}^2, ..., e_{ij}^n\right\}$$

if it is accepting a usual order or

$$a_{ij} = \text{pseudo neutrosophic min}\left\{e_{ij}^1, e_{ij}^2, ..., e_{ij}^n\right\};$$

1 ≤ i ≤ m and 1 ≤ j ≤ n.

$$a_{ij} = \text{pseudo real min}\left\{e_{ij}^1, e_{ij}^2, ..., e_{ij}^n\right\};$$

1 ≤ i ≤ m and 1 ≤ j ≤ n.
On similar lines we define B = $(b_{ij})$

$$b_{ij} = \max\left\{e_{ij}^1, e_{ij}^2, ..., e_{ij}^n\right\};$$

1 ≤ i ≤ m, 1 ≤ j ≤ n;

$$b_{ij} = \text{pseudo neutrosophic max}\left\{e_{ij}^1, e_{ij}^2, ..., e_{ij}^n\right\};$$

1 ≤ i ≤ m,

$$b_{ij} = \text{pseudo real max}\left\{e_{ij}^1, e_{ij}^2, ..., e_{ij}^n\right\};$$

1 ≤ i ≤ m, 1 ≤ j ≤ n.

Thus [A, B] is the interval of synaptic connection of neutrosophic matrices defined on the scale ⟨[–a a] ∪ [–aI, aI]⟩. The optimal synaptic connection of neutrosophic matrices of the interval of synaptic connection of neutrosophic matrices is defined by

$$O = \frac{A+B}{2},$$

$$= \frac{(a_{ij})+(b_{ij})}{2} \in [A, B].$$

Let us find the mean of the synaptic connection of neutrosophic matrices of the interval [A, B].

$$E = \frac{E_1 + E_2 + ... + E_n}{n}$$

$$= \frac{\left(e_{ij}^{1}\right) + \left(e_{ij}^{2}\right) + ... + \left(e_{ij}^{n}\right)}{n}.$$

Clearly $E \in [A, B]$.

[A, B] is defined as the interval of synaptic connection of neutrosophic matrices of n experts of the nINBAM model. [A, B] can be interval of neutrosophic synaptic connection of matrices or pseudo neutrosophic interval of synaptic connection of neutrosophic matrices or pseudo real interval of synaptic connection of neutrosophic matrices of the nINBAM model.

Now when we have several experts and they do not agree upon the same set of attributes then how to tackle the situation. We describe it as follows:

Let us consider the problem P and suppose we have n sets of experts working on the problem. They choose to work with different sets of neutrosophic intervals and they also choose to work with different sets of attributes; how to form a NBAM model which will deal with this situation simultaneously. To this end we construct the following model.

Let the sets of experts be denoted by $t_1, \ldots, t_n$. So the first set consists of $t_1$ experts, second set consists of $t_2$ experts and so on. Now let the $t_i$ set of experts work with $(m_i, n_i)$ set of attributes on the neutrosophic scale $I_i = \langle[-a_i, a_i] \cup [-a_iI, a_iI]\rangle$ i = 1, 2, …, n. Let $F_{\langle X_i \cup I\rangle}$ and $F_{\langle Y_i \cup I\rangle}$ denote the neutrosophic neuronal field from which the backward synaptic function neutrosophic matrices and forward synaptic function neutrosophic matrices are formed.

Thus these $t_i$ neutrosophic matrices would only be a $m_i \times n_i$ neutrosophic matrices with entries from the neutrosophic interval $I_i$; i = 1, 2, …, n.

Thus we have a set of n, $m_i \times n_i$; $t_i$ number of synaptic connection neutrosophic matrices of the $t_i$ experts defined over the scale $I_i$ (true for i = 1, 2, … , n). Now using these $t_i$ neutrosophic $m_i \times n_i$ matrices we form the interval $[A_i, B_i]$ of synaptic neutrosophic connection $m_i \times n_i$ matrices with entries from $I_i$. This interval of synaptic neutrosophic connection matrices may be pseudo neutrosophic interval of neutrosophic matrices or pseudo real interval of neutrosophic synaptic connection matrices or just the interval of neutrosophic connection matrices. Thus we have n sets of intervals of synaptic connection of neutrosophic matrices viz. $[A_1, B_1]$, $[A_2, B_2]$, …, $[A_n, B_n]$.

Now set $[A, B] = [A_1, B_1] \cup [A_2, B_2] \cup \ldots \cup [A_n, B_n]$ where $A = A_1 \cup A_2 \cup \ldots \cup A_n$ and $B = B_1 \cup B_2 \cup \ldots \cup B_n$ and the O, the optimal synaptic connection neutrosophic matrix is $O = O_1 \cup O_2 \cup \ldots \cup O_n$ and the mean synaptic connection neutrosophic matrix of the system

$$\bar{M} = \bar{M}^1 \cup \bar{M}^2 \cup \ldots \cup \bar{M}^n$$

where

$$\bar{M}^i \in [A_i, B_i];\ i = 1, 2, \ldots, n.$$

Thus the set [A, B] consists of all neutrosophic $m_i \times n_i$ rectangular n-matrices called the mixed interval of synaptic connection of neutrosophic n-matrices associated with the nINBAM model; of n sets of experts.

This model can function as a combined nINBAM model giving opinion of n-sets of experts at any $k^{th}$ period of time. The input neutrosophic n vectors will be taken from $I_1 \cup I_2 \cup \ldots \cup I_n$ or to be more specific from the neutrosophic neuronal n-fields

$$F_{\langle X_1 \cup I \rangle} \cup F_{\langle X_2 \cup I \rangle} \cup \ldots \cup F_{\langle X_n \cup I \rangle}$$

and

$$F_{\langle Y_1 \cup I \rangle} \cup F_{\langle Y_2 \cup I \rangle} \cup \ldots \cup F_{\langle Y_n \cup I \rangle}.$$

Just any synaptic connection neutrosophic n-matrix M in the interval [A, B] would be of the form $M = M_1 \cup M_2 \cup \ldots \cup M_n$ where $M_i \in [A_i, B_i]$ ; i = 1, 2, …, n. Thus if $X_k \in F_{\langle X_1 \cup I \rangle} \cup F_{\langle X_2 \cup I \rangle} \cup \ldots \cup F_{\langle X_n \cup I \rangle}$, then $X_k = X_1 \cup X_2 \cup \ldots \cup X_n$ as in case of other NBAMs, $X_k$ will not be recognized by the dynamical system of the nINBAM, so we use the usual S function; thus

$$S(X_k) = S(X_1) \cup \ldots \cup S(X_n)$$

where $S(X_t)$ takes values only from the set {0, 1, I}.

Thus

$$\begin{aligned} S(X_k)M &= [S(X_1) \cup \ldots \cup S(X_n)][M_1 \cup M_2 \cup \ldots \cup M_n] \\ &= S(X_1)M_1 \cup S(X_2)M_2 \cup \ldots \cup S(X_n)M_n \\ &= Y_1 \cup Y_2 \cup \ldots \cup Y_n \\ &= Y_{k+1}. \end{aligned}$$

$$S(Y_{k+1}) = S(Y_1) \cup S(Y_2) \cup \ldots \cup S(Y_n).$$

Thus

$$\begin{aligned} S(Y_{k+1})\, M^T &= [S(Y_1) \cup S(Y_2) \cup \ldots \cup S(Y_n)] \\ &\quad \left[M_1^T \cup M_2^T \cup \ldots \cup M_n^T\right] \\ &= S(Y_1)\, M_1^T \cup S(Y_2)M_2^T \cup \ldots \cup S(Y_n)M_n^T \\ &= X_1^1 \cup X_2^1 \cup \ldots \cup X_n^1 \\ &= X_{k+3}; \end{aligned}$$

and so on. The same procedure is repeated until we arrive at the equilibrium of the dynamical system $M \in [A, B]$. The equilibrium of the system gives a fixed n-point, using which interpretation can be made as in case of BmAMs. Having seen the use of interval of neutrosophic n matrices now we proceed on to give the application of fuzzy neutrosophic matrices, fuzzy neutrosophic n-matrices and the interval of fuzzy interval of fuzzy neutrosophic n-matrices, $n \geq 2$.

The fuzzy neutrosophic interval matrices can be used in the other neutrosophic models like NCMs and NRMs to perform multifold actions.

The NCM model has been introduced in the book [216]. The notion of bi-NCM model and n-NCM model have been described in [223]. Now we will just describe how the new model makes use of several experts opinion.

Suppose we have a problem P which has only unsupervised data and some n experts work on it using NCMs, how to study this without using the notion of combined NCMs (CNCMs). Let us assume all the n experts want to work with the same number of m attributes associated with the same problem P.

Let us further assume that they wish to work with a non simple NCM. Let $M_1, \ldots, M_n$ be the opinion given in the form of connection matrices related with the problem; i.e., $M_t = (m_{ij}^t)$; $1 \leq t \leq n$ and $1 \leq i, j \leq m$. All the n matrices are $m \times m$ matrices with entries from $\langle [0, 1] \cup [0, I] \rangle$. Form the interval of neutrosophic matrices by defining the max and min for these n matrices as $B = (b_{ij})$ and $A = (a_{ij})$ respectively. The $a_{ij}$'s are defined as

$$a_{ij} = \min \{m_{ij}^1, m_{ij}^2, \ldots, m_{ij}^n\}$$

if $\min \{m_{ij}^1, m_{ij}^2, \ldots, m_{ij}^n\}$; does not exists then define

$$a_{ij} = \text{pseudo neutrosophic } \min \{m_{ij}^1, m_{ij}^2, \ldots, m_{ij}^n\};$$

$1 \leq i, j \leq m$ or

$$a_{ij} = \text{pseudo real } \min \{m_{ij}^1, m_{ij}^2, \ldots, m_{ij}^n\}$$

which ever is preferred by the experts.

Likewise for $B = (b_{ij})$ define

$$b_{ij} = \max \{m_{ij}^1, m_{ij}^2, \ldots, m_{ij}^n\};$$

$1 \leq i, j \leq m$ if such $\max \{m_{ij}^1, m_{ij}^2, \ldots, m_{ij}^n\}$ does not exist define pseudo neutrosophic $\max \{m_{ij}^1, m_{ij}^2, \ldots, m_{ij}^n\}$; $1 \leq i, j \leq m$ or pseudo real $\max \{m_{ij}^1, m_{ij}^2, \ldots, m_{ij}^n\}$.

We define this interval [A, B] to be the interval of neutrosophic $m \times m$ connection matrices of the NCM denoted by INCM.

Define optimal NCM, O of the system as

$$O = \frac{(a_{ij} + b_{ij})}{2} = \frac{A + B}{2}.$$

Define the mean of the NCM as

$$M = \frac{M_1 + ... + M_n}{n}$$
$$= \frac{\left(m_{ij}^1\right) + \left(m_{ij}^2\right) + ... + \left(m_{ij}^n\right)}{n}.$$

Clearly O, M ∈ [A, B]. [A, B] is called the associated interval of fuzzy neutrosophic matrices of the NCM (INCM).

Now if on the other hand we have some t-sets of experts giving their opinion about the problem with varied sets of attributes how to consolidate this study and how to construct a single dynamical system to analyze such a problem. We see the interval of fuzzy neutrosophic t-matrices will serve the purpose.

Let $p_1, \ldots, p_t$ set of experts give their opinion on the same problem P. Let all the experts in the set $p_i$ work with $m_i$ attributes and let them choose to work with non simple NCM on the neutrosophic interval $\langle [0, 1] \cup [0, I] \rangle$ for i = 1, 2, …, t.

Let us as described in earlier section of this chapter the fuzzy model, construct the interval $[A_i, B_i]$ of fuzzy neutrosophic $m_i \times m_i$ matrices with entries from $\langle [0, 1] \cup [0, I] \rangle$, for i = 1, 2, …, t.

Set $[A, B] = [A_1, B_1] \cup [A_2, B_2] \cup \ldots \cup [A_t, B_t]$ with $A = A_1 \cup A_2 \cup \ldots \cup A_t$ and $B = [B_1 \cup B_2 \cup \ldots \cup B_t]$ i.e. A and B may be minimal and maximal n matrices of the interval of fuzzy neutrosophic n-matrices.

Define the optimal of fuzzy neutrosophic n-matrices as

$$O = \frac{A + B}{2}$$
$$= \frac{(A_1 + B_1)}{2} \cup \frac{(A_2 + B_2)}{2} \cup ... \cup \frac{(A_t + B_t)}{2}$$
$$= O_1 \cup O_2 \cup \ldots \cup O_t;$$

$O_i \in [A_i, B_i]$;

O will be known as the fuzzy neutrosophic optimal connection n-matrix of the interval [A, B]. Define the mean fuzzy neutrosophic connection n-matrix to be $M = M_1 \cup M_2 \cup \ldots \cup$ $M_t$; $M_i$ are the mean connection fuzzy neutrosophic n-matrices; $M_i \in [A_i, B_i]$; i = 1, 2, …, t. This interval [A, B] will be known as the interval of fuzzy neutrosophic n-matrices of the INn-CM ($n \geq 2$). When n = 2 we get the interval of fuzzy neutrosophic bimatrices of the INBCM or IN2-CM.

Likewise we can define interval of fuzzy neutrosophic matrices for the INnCM (neutrosophic n-relational maps) model ($n \geq 2$).

The case of FnRMs have been very clearly explained with illustrations. In case of INnRMs we need to follow two different steps.

1. Replace the fuzzy relational matrices by the fuzzy neutrosophic relational matrices.
2. In case of the interval of relational fuzzy neutrosophic matrices related with a IN-nRM we may have three types of order in the interval of n-fuzzy neutrosophic matrices.
3. The usual order min-max of the fuzzy neutrosophic matrices on the interval of n-matrices $[A, B] = [A_1, B_1] \cup \ldots \cup [A_n, B_n]$.
4. The pseudo neutrosophic order min-max of fuzzy neutrosophic matrices on the interval of n-matrices $[A, B] = [A_1, B_1] \cup \ldots \cup [A_n, B_n]$.
5. The pseudo or real min-max of fuzzy neutrosophic matrices on the interval of n-matrices $[A, B] = [A_1, B_1] \cup [A_2, B_2] \cup \ldots \cup [A_n, B_n]$.

All neutrosophic models can be constructed following the above changes. Thus the interested reader is requested to refer [219].

# INDEX

## A

## B

## C

## H

## I

## L

## M

## N

## O

## P

## Q

## R

## S

## T

# ABOUT THE AUTHORS

**Dr.W.B.Vasantha Kandasamy** is an Associate Professor in the Department of Mathematics, Indian Institute of Technology Madras, Chennai. In the past decade she has guided 11 Ph.D. scholars in the different fields of non-associative algebras, algebraic coding theory, transportation theory, fuzzy groups, and applications of fuzzy theory of the problems faced in chemical industries and cement industries. Currently, four Ph.D. scholars are working under her guidance.

She has to her credit 633 research papers. She has guided over 51 M.Sc. and M.Tech. projects. She has worked in collaboration projects with the Indian Space Research Organization and with the Tamil Nadu State AIDS Control Society. This is her 28$^{th}$ book.

On India's 60th Independence Day, Dr.Vasantha was conferred the Kalpana Chawla Award for Courage and Daring Enterprise by the State Government of Tamil Nadu in recognition of her sustained fight for social justice in the Indian Institute of Technology (IIT) Madras and for her contribution to mathematics. (The award, instituted in the memory of Indian-American astronaut Kalpana Chawla who died aboard Space Shuttle Columbia). The award carried a cash prize of five lakh rupees (the highest prize-money for any Indian award) and a gold medal.
She can be contacted at vasanthakandasamy@gmail.com
You can visit her on the web at: http://mat.iitm.ac.in/~wbv or: http://www.vasantha.net

---

**Dr. Florentin Smarandache** is an Associate Professor of Mathematics at the University of New Mexico in USA. He published over 75 books and 100 articles and notes in mathematics, physics, philosophy, psychology, literature, rebus. In mathematics his research is in number theory, non-Euclidean geometry, synthetic geometry, algebraic structures, statistics, neutrosophic logic and set (generalizations of fuzzy logic and set respectively), neutrosophic probability (generalization of classical and imprecise probability). Also, small contributions to nuclear and particle physics, information fusion, neutrosophy (a generalization of dialectics), law of sensations and stimuli, etc.
He can be contacted at smarand@unm.edu